\documentclass[12pt]{article}
\usepackage{defs1}
  
\title{Differential cohomology}
\author{Ulrich Bunke\thanks{NWF I - Mathematik,
Universit{\"a}t Regensburg,
93040 Regensburg,
GERMANY, ulrich.bunke@mathematik.uni-regensburg.de}  
}

\newcommand{\bku}{{\mathbf{ku}}}
\newcommand{\Free}{{\tt Free}}

\newcommand{\fibre}{{\tt kernel}}

\newcommand{\Filt}{{\tt Filt}}
\newcommand{\aut}{{\tt aut}}
\newcommand{\ed}{{\tt end}}
\newcommand{\tot}{{\tt tot}}
\newcommand{\bs}{{\mathbf{s}}}
\newcommand{\desusp}{{\tt desusp}}
\newcommand{\GeomFam}{{\tt GeomFam}}

\newcommand{\harm}{{\tt harm}}
\newcommand{\Cycle}{{\tt Cycle}}
\newcommand{\cN}{{\mathcal{N}}}
\newcommand{\hbMU}{{\widehat{\mathbf{MU}}}}
\newcommand{\bMU}{{\mathbf{MU}}}
\newcommand{\hcyc}{{\tt \widehat{cycl}}}
 \newcommand{\hbKU}{\widehat{\mathbf{KU}}} 
 \newcommand{\bKU}{\mathbf{KU}}
\newcommand{\cw}{{\tt cw}}

\newcommand{\Simp}{{\mathbf{Simp}}}
\newcommand{\cK}{{\mathcal{K}}}
\newcommand{\su}{{\mathbf{su}}}
\newcommand{\Line}{{\mathrm{Line}}}

\newcommand{\bu}{{\mathbf{u}}}
\newcommand{\pt}{{\mathrm{pt}}}
\newcommand{\bH}{\mathbb{H}}

\newcommand{\hol}{\mathrm{hol}}
\newcommand{\taaa}{{\mathfrak{t}}}
\newcommand{\cs}{{\mathbf{cs}}}
\newcommand{\MSO}{{\mathbf{MSO}}}
\newcommand{\CS}{{\mathbf{CS}}}
\newcommand{\so}{{\mathbf{so}}}
\newcommand{\vol}{{\mathrm{vol}}}
\newcommand{\Ad}{{\mathrm{Ad}}}
\newcommand{\gaaa}{{\mathfrak{g}}}

\newcommand{\cI}{{\mathcal{I}}}

\newcommand{\Rings}{{\mathbf{Rings}}}

\newcommand{\cL}{{\mathcal{L}}}
\newcommand{\cW}{{\mathcal{W}}}

\newcommand{\hcycl}{\widehat{\tt cycl}}

\newcommand{\bK}{{\mathbf{K}}}

\newcommand{\cU}{{\mathcal{U}}}

\newcommand{\cD}{{\mathcal{D}}}

 \newcommand{\Cone}{{\tt Cone}}
 
 \newcommand{\Vect}{{\tt Vect}}
 
 \newcommand{\CommMon}{{\mathbf{CommMon}}}
 \newcommand{\Cat}{{\mathbf{Cat}}}

\newcommand{\cE}{{\mathcal{E}}}

\renewcommand{\Im}{\operatorname{Im}}

\begin{document}
\maketitle
\begin{abstract}

These course note first provide an introduction to secondary characteristic classes and differential cohomology. They continue with  a presentation of a stable homotopy theoretic approach to the  theory of differential extensions of generalized cohomology theories including products and Umkehr maps. 
\end{abstract}

\tableofcontents

\section{Introduction}

This is a script based on my introductory course on differential cohomology tought in spring 2012
at the Regensburg University.  The course was designed for students having a profound background in differential geometry, algebraic topology and homotopy theory. 

The course starts with a detailed introduction to characteristic forms for vector bundles. I explain how one can use locality  and integrality of these forms to construct secondary invariants.  
The main emphasis in this part is put on explicit calculations of examples. I think that this experience is important to understand and work with the fine structure of differential cohomology to be introduced later.

 To capture these secondary invariants for complex vector bundles systematically and refine them even more 
I first introduce differential integral cohomology. More specifically, I explain 
   the sheaf theoretic definition of differential cohomology usually called smooth Deligne cohomology. The main theoretical result here is the construction of the  differential refinement of the Chern-Weyl homomorphism due to Cheeger-Simons.  I give a detailed discussion of various structures like integration and products.  Up to this point the course reviews the classical part of the theory. The theory is complemented  by a variety of examples. It is again important to do these exercises to get an idea how the theory works, and to learn some calculational tricks.

The course then continues with a presentation of a full-fledged theory of  differential extensions of generalized cohomology theories including products and integration. This stable homotopy theoretic  approach  is new and complements many special constructions for particular cases. The details of the homotopy theoretic approach are developed in Section \ref{aug2040} which can be read independently of the previous sections. I again tried to add, as much as possible, interesting explicit calculations.  As examples of differential extensions of generalized cohomology theories  I 
discuss the cases of  complex $K$-theory and  complex bordism in some detail.

The script contains a variety of exercises ranging from explicit calculations with numerical outputs to 
verifications of statements made without detailed proofs in the course. While these theoretical 
pieces should be more or less straight-forward,  some of the explicit calculations are more involved and
were major pieces of original works or Phd theses. This case will be indicated by giving appropriate references. For many problems I have added a sketch of a proof.
The exercise sheets 1-13 were discussed in the weekly exercise classes which complemented the course. The Sheet 14 (cf. \ref{sheet14}) collects some problems which were left open during the course. Some of them may be interesting research  projects.

\bigskip

{\em Acknowledgements: I am grateful to the participants of the course (CM. Ruderer, M. Spitzweck, J. Sprang, A. Strack, G. Tamme, M. Voelkl)   for their patience and critical questions which lead to many crucial corrections. Furthermore, I thank Th. Nikolaus and D. Gepner for various helpful discussions, in particular about $\infty$-categorical aspects. 
Finally I thank Ch. Becker and  A. Kuebel for their support to eliminate typos and small mistakes. }

\section{Characteristic forms}
 \subsection{Linear Connections}
 
We consider a smooth manifold $M$. By $\Omega(M)$ we denote the de Rham complex of $M$. For its complexification we use the notation $\Omega(M,\C)$. Let
 $E\to M$ be a complex vector bundle. Then we have the  $\Omega(M,\C)$-module  $\Omega(M,E)$ of differential forms with coefficients in $E$. In particular, $\Omega^{0}(M,E)$ is the space of sections of $E$. Locally, an element of $\Omega(M,E)$ can be written as a finite sum of elements of the form $\omega\otimes \phi$ with $\phi\in \Omega(M,\C)$
 and $\phi\in \Omega^{0}(M,E)$. The product of $\alpha\in \Omega(M,\C)$ with
 $\omega\otimes \phi$ is then given by $\alpha\wedge  (\omega\otimes \phi):=(\alpha\wedge \omega)\otimes \phi$. 
 
 \begin{ddd}
A connection on $E$ is a map
$$\nabla:\Omega^{0}(M,E)\to \Omega^{1}(M,E)$$
satisfying the Leibnitz rule
$$\nabla (f\wedge \phi)=df\wedge \phi +f\wedge \nabla\phi\ , \quad \forall f\in \Omega^{0}(M,\C)\ , \forall \phi\in \Omega^{0}(M,E)\ .$$
\end{ddd}
In the following we discuss some constructions with connections and observe that connections always exist.

A trivialization $\phi$ of $E$ is a collection of sections $\phi=(\phi_{\alpha})_{\alpha=1}^{k}$ of $E$ whose evaluations at each point $m\in M$ form a basis of the fibre $E_{m}$ of $E$ at $m$.
A trivialization $\phi$
 induces a connection $\nabla^{\phi}$ as follows. 
\begin{prob}\label{apr1202} Show that
there exists a unique connection  
$\nabla^{\phi}$   on $E$ characterized by $\nabla^{\phi}(\phi_{\alpha})=0$ for all $\alpha\in \{1,\dots k\}$.
 \end{prob}
Connections can be glued with the help of a  partition of unity.
Let $(\nabla_{i})_{i\in I}$  be a collection of connections on $E$ and
$(\chi_{i})_{i\in I}$ be a partition of unity. For a function $\chi\in \Omega^{0}(M,\C)$ and a connection $\nabla$ we write $\chi\nabla$ for the composition of the operators $\nabla$ and multiplication by $\chi$.
\begin{prob}
Show that
the sum $\sum_{i\in I} \chi_{i} \nabla_{i}$ is a connection on $E$.
\end{prob}

\begin{lem}
Every vector bundle admits a connection
\end{lem}
\proof
Since a vector bundle is locally trivial, it locally admits connections. We get a global connection by glueing the local ones.
 \begin{prob}
Work out the details.
\end{prob}
\hB 
We consider the bundle of algebras $\End(E)\to M$. It induces a graded algebra  
$\Omega(M,\End(E))$ over $\Omega(M,\C)$ with product given by 
$$(\alpha\otimes \Psi)\wedge (\beta\otimes \Phi):=(\alpha\wedge \beta)\otimes \Psi\circ\Phi\ ,$$
This algebra acts on $\Omega(M,E)$ by $$(\alpha\otimes \Psi)\wedge (\omega\otimes \phi):=(\alpha\wedge \omega)\otimes\Psi(\phi)\ .$$

 \begin{prob}\label{may0601} Show that
 the set of connections  on $E$ has the structure of an affine space over $\Omega^{1}(M,\End(E))$ given by  
  $$(\nabla+\omega)(\phi)=\nabla\phi + \omega\wedge\phi\ , \quad \forall \omega\in \Omega^{1}(M,\End(E))\ , \forall \phi\in \Omega^{0}(M,E)\ .$$
  \end{prob}
  
  \begin{prob}\label{apr1221}
  Show that a  connection on $E$ uniquely extends to a connection on $\End(E)$ which is  characterized by 
\begin{equation}\label{mar2802}(\nabla\omega)\wedge \phi=\nabla(\omega\wedge \phi)-\omega\wedge \nabla \phi\end{equation}
for $\omega\in \Omega^{0}(M,\End(E))$, $\phi\in \Omega^{0}(M,E)$.
\end{prob}

\begin{prob}\label{apr1610}
Show that a connection on $E$ uniquely extends to a linear map $\nabla:\Omega(M,E)\to \Omega(M,E)$
satisfying  the Leibnitz rule
\begin{equation}\label{mar2801}\nabla(\omega\wedge \phi)=d\omega\wedge \phi+(-1)^{k}\omega\wedge \nabla\phi\ ,\quad \forall \omega\in \Omega^{k}(M,\C)\ ,\forall \phi\in \Omega(M,E)\ .\end{equation}
\end{prob}

Using this extension we consider $\nabla^{2}:=\nabla\circ \nabla:\Omega(M,E)\to \Omega(M,E)$. 

\begin{prob}\label{apr1201}
Show that $\nabla^{2}$ is given by muliplication by a uniquely determined element
$R^{\nabla}\in \Omega^{2}(M,\End(E))$.\end{prob}
\begin{ddd}\ The element 
 $R^{\nabla}\in \Omega^{2}(M,\End(E))$ characterized in Problem \ref{apr1201}
 is called the curvature of $\nabla$. A connection is called flat if its curvature vanishes.
 \end{ddd}
 The connection associated to a trivialization of $E$ (Problem \ref{apr1202}) is flat.

 Note that $\nabla$ is an odd operator. Hence we can write
 $\nabla^{2}=\frac{1}{2}[\nabla,\nabla]$. Similarly, the application of the extension of $\nabla$
 defined in Problems \ref{apr1221} and \ref{apr1610} to $\Phi\in \Omega^{ev}(M,\End(E))$ can be written as
 $\nabla\Phi=[\nabla,\Phi]$. It follows the Bianchi identity  \begin{equation}
 \label{apr1224}\nabla R^{\nabla}=\frac{1}{2}[\nabla,[\nabla,\nabla]]=0\ .\end{equation}

Let $f:M^{\prime}\to M$ be a smooth map and consider the pull-back diagram
$$\xymatrix{E^{\prime}\ar[d]\ar[r]^{F}&E\ar[d]\\M^{\prime}\ar[r]^{f}&M}$$
of vector bundles. We have a pull-back operation of sections
$(F,f)^{*}:\Omega(M,E)\to \Omega(M^{\prime},E^{\prime})$
such that
$$F(((F,f)^{*}\omega)(m^{\prime}))=\omega(m)\circ \Lambda^{k} df(m^{\prime})$$
as linear maps $\Lambda^{k}T_{m^{\prime}}M^{\prime}\to E_{m}$
for $\omega\in \Omega^{k}(M,E)$, $m^{\prime}\in M^{\prime}$ and $m:=f(m^{\prime})$.

  We consider a connection $\nabla$ on $E$.
\begin{prob}
Show that there is a unique connection $\nabla^{\prime}:=(F,f)^*\nabla$ on $E^{\prime}$  which is  characterized by 
$$\nabla^{\prime} ((F,f)^{*}\phi)=(F,f)^{*}(\nabla\phi)\ ,\quad \forall \phi\in \Omega^{0}(M,E)\ .$$
\end{prob}
If $f=\id$, then we write $F^{*}:=(F,\id)^{*}$. If $E^{\prime}=f^{*}E$ and $F$ is the canonical map, then we write $f^{*}\nabla:=(F,f)^{*}\nabla$.
\begin{prob}
Verify that
if $f=\id_M$ and $F\in \Omega^{0}(M,\End(E))$ is an automorphism of $E$, then  we have
$$F^{*}\nabla=\nabla-F^{-1}\wedge \nabla F\ .$$
\end{prob}
\proof
Note that $F^{*}\phi=F^{-1}\wedge \phi$. \hB 
We let $\iota_{F}:f^{*}\End(E)\stackrel{\sim}{\to}\End(E^{\prime})$ be the isomorphism of endomorphism bundles induced by $F$.
\begin{prob}\label{apr1220} Show that
the curvature satisfies
$$R^{(F^{*},f)\nabla}=\iota_{F}(F,f)^{*}R^{\nabla}\ .$$
\end{prob}

A connection on $E$ induces a parallel transport $\|_\gamma:E_{\gamma(0)}\to E_{\gamma(1)}$ along curves $\gamma:[0,1]\to M$. This material is assumed as a prerequisite.

Later we need the following variants.

\begin{ex}{\em
If $E\to M$ is a real vector bundle, then 
a connection is a map
$$\nabla:\Omega^{0}(M,E)\to \Omega^{1}(M,E)$$
satisfying the Leibnitz rule for real functions $f\in C^{\infty}(M)$.
Here the real vector space $\Omega(M,E)=\Gamma(M;\Lambda^{*}T^{*}M\otimes E)$ is an $\Omega(M)$-module.}
\end{ex}

Later in this course (e.g in \ref{may2004}) we will assume knowledge of the theory of connections on principal bundles.
We refer to \cite{MR1393940} for details. There is a parallel theory of existence, glueing, and curvature.

\subsection{Invariant connections and quotients}

 Let $G$ be a Lie group which acts on $M$. Assume that this action lifts to an action on the bundle $E\to M$. For every $g\in G$ we get a pull-back diagram
 $$\xymatrix{E\ar[d]\ar[r]^{\tilde g}&E\ar[d]\\M\ar[r]^{g}&M}\ ,$$
 and the map $g\mapsto \tilde g$ satifies an associativity rule.
 
 \begin{ddd}
 A connection $\nabla$ on $E$ is called $G$-invariant if 
 $(\tilde g,g)^{*}\nabla =\nabla$ for all $g\in G$.
 \end{ddd}

 \begin{prob}\label{apr1207} 
 Show that  if $G$ acts properly on $M$, then there exists a $G$-invariant connection on $E$.
 \end{prob}
\proof By assumption the map
$(\mu,\pr_M):G\times M \to M\times M$ is proper.
There exists a function $\chi\in C^{\infty}(M)$ such that
$(\mu,\pr_M)^{-1}(\supp(\chi)\times K)$ is compact for all compact $K\subseteq M$ and $\int_{G}\: g^{*}\chi \:dg\equiv 1$.  Let $\tilde \nabla$ be any connection on $E$. 
Then
$$\nabla:=\int_{G}\: (\tilde g,g)^{*} (\chi\tilde  \nabla)\:dg$$
is an invariant connection. Observe the convergence of this integral.

For the construction of the function $\chi$ one could proceed as follows.
Let $(\phi_i)_{i\in I}$ be a collection of non-negative compactly supported functions 
  whose supports cover $M$ and set
$\bar \phi_i:=\int_G g^{*}\phi_i dg$. The supports of these functions cover $M$ and  by paracompactness  we can choose a subset $J\subseteq I$ such that
$(\supp(\bar \phi_i))_{i\in J}$ is a locally finite covering of $M$.
We set $$\chi:=\frac{1}{\sum_{i\in J}\bar \phi_j}\sum_{i\in J} \phi_i\ .$$
 \hB

Assume that $G$ acts freely and properly on $M$. Then we have a $G$-principal bundle
$M\to \bar M:=G\backslash M$. Furthermore, $G$ acts freely on $E$ such that $\bar E:=G\backslash E\to \bar M$ is a vector bundle which fits into  the pull-back diagram 
 $$\xymatrix{E\ar[d]\ar[r]^{P}&\bar E\ar[d]\\M\ar[r]^{p}&\bar M}\ .$$
 By $\Omega^{0}(M,E)^{G}\subseteq \Omega^{0}(M,E)$ we denote the $G$-invariant sections of $E$.
 \begin{prob}
 Show that there is a natural isomorphism
 $$(p,P)^{*}:\Omega^{0}(\bar M,\bar E)\cong \Omega^{0}(M,E)^{G}\ .$$  
 \end{prob}

 \begin{ddd}
 We say that an invariant connection $\nabla$ on $E$ descends  if there exists
 a connection $\bar \nabla$ on $\bar E\to \bar M$ such that
 $\nabla=(P,p)^{*}\bar \nabla$.
 \end{ddd}
 
 \begin{prob}
 If $\nabla$ descends, then $\bar \nabla$ is uniquely determined.
 \end{prob}
 
 \begin{prob} 
 Assume that $G$ is discrete. Then
   every invariant connection    $\nabla$ on $E$ descends. \end{prob}
 
 If $G$ is not discrete, then there is an obstruction against descending an invariant  connection $\nabla$.
 Let $\gaaa$ be the Lie algebra of $G$. By $\cX(M)$ we denote the set of  vector fields on $M$.
 \begin{ddd} For $A\in \gaaa$  the fundamental vector field
 $A^{\sharp}\in \cX(M)$  is defined by
 $$A^{\sharp}(m)=\frac{d}{dt}\exp(tA)m\ .$$
 \end{ddd}
 
 \begin{ddd}
For $A\in \gaaa$ the  Lie derivative  $\cL_A:\Omega(M,E)\to \Omega(M,E)$
 is defined by 
$$\cL_{A}\phi:=\frac{d}{dt}_{|t=0}(\widetilde{\exp(tA)},\exp(tA))^{*}\phi\ ,  \phi\in \Omega(M,E)\ .$$
\end{ddd}

Let $\nabla$ be an invariant connection.
\begin{ddd}\label{apr2401}
The moment map $\mu^{\nabla}\in \Hom(\gaaa, \Omega^{0}(M,\End(E))^{G}$
 is  defined by 
 $$\mu^{\nabla}(A)\wedge \phi:=\nabla_{A^{\sharp}} \phi-\cL_{A}
\phi\ , \quad  \phi\in \Omega^{0}(M,E)\ .$$
\end{ddd} 
\begin{prob}
 Show that $\mu^{\nabla}$ is well-defined. Further show that 
 $\nabla$ descends if and only if $\mu^{\nabla}=0$.
\end{prob}

\begin{ex}{\em 
Assume that $G$ is a compact Lie group.
We consider the  action of $G$ on $TG\to G$ by left multiplication.
 We consider a biinvariant  Riemannian
metric on $G$ and let $\nabla$ be the Levi-Civita connection on $TG$.
Furthermore let $\nabla^{l}$ be the connection defined by the trivialization of $TG$ by 
left-invariant vector fields.
\begin{prob}Calculate the moment maps of $\nabla$ and $\nabla^{l}$.
\end{prob}
We shall see that $\nabla^{l}$ descends, while $\nabla$ does not.

}\end{ex}

\begin{ex}\label{apr1203} 
{\em
 Let $\rho:G\to \End(V)$ be a linear representation of $G$ on  a finite-dimensional  complex vector space
 $V$ and assume that $G$ acts freely and properly on $M$. Then we get an action
 of $G$ on the trivial bundle $E:=M\times V$ by $g(m,v):=(gm,\rho(g)v)$. The trivial connection on $E\to M$ is $G$-invariant. \begin{prob}Calculate the moment map.\end{prob}   
We shall see that $\nabla$ descends if and only if $d\rho=0$. }
 \end{ex}

\begin{ex}{\em
In this example we discuss the relation of the moment map in the sense  of Definition \ref{apr2401}
with the moment map in symplectic geometry.
Let $(M,\omega)$ be a symplectic manifold  such that $\frac{1}{2\pi i}\omega$ is integral (Definition \ref{apr2402}).
For simplicity we assume that $M$ is simply connected. Then there exists a line bundle
$L\to M$ with connection $\nabla$ such that $R^{\nabla}=-\omega$. In fact, $(L,\nabla)$ is uniquely determined up to isomorphism by these conditions. We now assume that $G$ acts on $M$ in a Hamiltonian fashion with symplectic moment map $\mu^{\omega}\in \Hom(\gaaa,C^{\infty}(M))^{G}$.
It is characterized up to an element of $(\gaaa^{*})^{G}$ by
$$d\mu^{\omega}(A)=i_{A^{\sharp}}\omega\ .$$

We now observe that $\End(L)=M\times \C$ so that $\mu^{\nabla}\in \Hom(\gaaa,C^{\infty}(M,\C))^{G}$.
The action of $G$ on $M$ extends to an action of a $\C^{*}$-central extension $\tilde G\to G$
on $L\to M$. For simplicity we assume that this extension splits. Two splits differ by a character $\chi:G\to \C^{*}$.  Note that the corresponding moment maps differ by $d\chi\in (g^{*}\otimes \C)^{G}$.
We calculate for $A\in \gaaa$ and $X\in \cX(M)$
$$d\mu^{\nabla}(A)(X)=[\nabla_X,\nabla_{A^{\sharp}}-\cL_A]=-\omega(X,A^{\sharp})=d\mu^{\omega}(A)(X)\ .$$
We conclude that we can adjust the action of $G$ on $L$ such that
$$\mu^{\nabla}=\mu^{\omega}\ .$$
In the theory of quantization one wants to define a quotient
$\bar L\to \bar M$ with connection $\bar \nabla$ whose curvature is symplectic.
This is obstructed exactly by the moment map. The way out is to restrict to its
zero set $M_0:=(\mu^{\omega})^{-1}(0)$ and to take $\bar L_0\to \bar M_0$ with the induced connection (assuming that $M_0$ is smooth and $G$ acts freely on $M_0$). 
This construction is called symplectic reduction.
}
\end{ex}

\begin{ex}\label{apr1204}  {\em
Let $M$ be a compact connected manifold with base point $m\in M$ and $\tilde M\to M$ be the associated universal covering.  Then $\pi:=\pi_{1}(M,m)$ is finitely generated and $J(M):=\Hom(\pi,U(1))$ is  finite-dimensional abelian Lie group, the Jacobian of $M$.
Its connected component of the identity is a torus.
  We consider the action of $\pi$ on
$E:=J(M)\times \tilde M\times \C$  given by  $\gamma(\rho,\tilde m,z):=(\rho,\tilde m \gamma^{-1},\rho(\gamma) z)$. This action covers the free and proper action of $\pi$ on $J(M)\times \tilde M$ (on the second factor). Hence we have a quotient $P\to J(M)\times M$ which is called the Poincar\'e bundle. 

For $\gamma\in \pi$ we consider the $U(1)$-valued function $\rho\mapsto \rho(\gamma)$ on $J(M)$.
We have $ \rho(\gamma)^{-1}d\rho(\gamma)\in \Omega^{1}_{cl}(J(M),\C)$. This is actually a homomorphism
$ \rho^{-1}d\rho:\pi\to  \Omega^{1}_{cl}(J(M),\C)$, hence a cohomology class
$[\theta]\in H^{1}_{dR}(M;\C)\otimes_{\C}  \Omega_{cl}^{1}(J(M),\C)$
represented by a closed form $\theta\in \Omega^{1}(M;\C)\otimes_{\C} \Omega^{1}_{cl}(J(M),\C)\subseteq \Omega_{cl}^{2}(J(M)\times M;\C)$. 
 Since $H^{1}(\tilde M;\C)=0$ we can choose a function
$\alpha\in \Omega^{0}(\tilde M;\C)\otimes_{\C}  \Omega_{cl}^{1}(J(M),\C)\subset \Omega^{1}(J(M)\times \tilde M,\C)$ such that
$d \alpha$ is $\pi$-invariant and its descent to $M$ equals $\theta$. Then we define the connection $\nabla$ on the trivial bundle $E$ by
$$\nabla=d - \alpha \ .$$  
\begin{prob}
Show that this connection is $\pi$-invariant.
\end{prob}
\proof
Let $[\gamma]\in \pi$ and $\phi$ be the constant section of $E$ with value $1$.
Then $([\gamma]^{*}\nabla)[\gamma]^{*}(\phi)=[\gamma]^{*}(\nabla \phi)$. We shall first compute the r.h.s. We have $\nabla \phi=-\alpha$. Furthermore $[\gamma]^{*}\alpha-\alpha=\int_{[\gamma]} d\alpha=  \rho(\gamma)^{-1}d\rho(\gamma)$. Hence $[\gamma]^{*}(\nabla \phi)=-\rho^{-1}(\gamma)\alpha- \rho(\gamma)^{-2}d\rho(\gamma)$.  On the other hand
$\nabla ([\gamma]^{*}\phi)=  -\rho(\gamma)^{-2}d\rho(\gamma)  -\rho(\gamma)^{-1}\alpha $. 
Hence $\nabla=[\gamma]^{*}\nabla$. \hB

Hence the choice of $\alpha$ determines a connection $\nabla^{P,\alpha}$ on the Poincar\'e bundle.
\begin{prob}\label{apr1205}
Calculate the curvature of $\nabla^{P}$.
\end{prob}
\proof
We have
$R^{\nabla}=-d\alpha=- \pi^{*} \theta$.
It follows that
$R^{\nabla^{P}}=-\theta\in \Omega^{2}(J(M)\times M,\C)$, where we identify
$M\times \C\cong \End(P)$.
\hB 

\begin{prob}\label{apr1206}
Make this construction explicit in the case where $M$ itself is a torus $T^{n}:=\R^{n}/\Z^{n}$. In this case we can take $\alpha$ to be linear. 
\end{prob}
\proof
We identify $\R^{n}$ with its dual with respect to the standard scalar product.
Then $J(T^{n})\cong \R^{n}/\Z^{n}$ with the identification
$[y](u)=\exp(2\pi i\langle y,u\rangle)$ for $u\in \Z^{n}\cong \pi_1(T^{n})$,
$y\in \R^{n}$. We have
$\rho^{-1}d\rho(u)= 2\pi i   \langle  dy,u\rangle$. 
We can take
$\theta:=2\pi i \langle  dy\wedge dx\rangle$. Finally, for $\alpha$ we chose
$\alpha:=2\pi i \langle dy,x\rangle$.
The curvature of $\nabla^{P}$ is given by
$R^{\nabla^{P}}=-2\pi i \langle dy\wedge dx\rangle$.
}
\end{ex}

\begin{ex}{\em
Let $M$ be a compact connected manifold with base point $m_{0}$.
We consider the functor
$F:\Mf^{op}\to \Set$
which associates to a manifold
$T$ the set of isomorphism classes of line bundles $L\to T\times M$
with a flat partial connection $\nabla:\Gamma(T\times M,L)\to \Gamma(T\times M,T^{*}M\otimes L)$ along the foliation given by the fibres of the projection $T\times M\to T$
and a trivialization of these structures restricted to   $T\times \{m_{0}\}$.
For $f:T^{\prime}\to T$ we let $F(f):F(T)\to F(T^{\prime})$ be given by  pull-back.

\begin{prob}
Show that this functor is representable.
More concretely, consider the Poincar\'e
 bundle $P\to J(M)\times M$ (\ref{apr1204}) with the partial connection $\nabla$ induced by $\nabla^{P}$.
 Observe that this determines
 an object of $F(J(M))$ which is independent of the choices involved in the construction. It induces a natural transformation
 $$\Hom_{\Mf}(\dots,J(M))\to F$$
 
 which turns out to be an isomorphism of functors.
 \end{prob}}\end{ex}

\begin{ex}\label{apr1260}{\em 
We have a natural  inclusion
$U(1)\subset SU(2)$ as a maximal torus. We have an identification $SU(2)/U(1)\cong \C\P^{1}$ and
$SU(2)\times_{U(1),\id} \C\cong L$, where $L\to \C\P^{1}$ is the tautological bundle.
The group $SU(2)$ still acts from the left on the quotient $L\to \C\P^{1}$.
\begin{prob}\label{apr1210}
Show that there is a unique $SU(2)$-invariant connection $\nabla^{L}$ on
$L\to \C\P^{1}$. Calculate its curvature form explicitly.
\end{prob}
\proof
Existence is ensured by \ref{apr1207}.
The difference of two invariant connections is an element of
$\Omega^{1}(\C\P^{1},\C)^{SU(2)}$. Show that this space is trivial.
\hB
The result of the  curvature calculation is 
$$R^{\nabla^{L}}=2\pi i \vol_{\C\P^{1}}\ ,$$
where $\vol_{\C\P^{1}}$ is the unique normalized $S^{3}$-invariant volume form on $\C\P^{1}$.
}
\begin{prob}
Generalize \ref{apr1210} to $\C\P^{n}$ with action of $SU(n+1)$.
\end{prob}
  \end{ex}

\subsection{Characteristic forms, classes and transgression - theory}

We consider characteristic forms for complex vector bundles. 
 
\begin{ddd}
A characteristic form $\omega$ of degree $n$ 
associates to each connection $\nabla$ on a complex vector bundle bundle $E\to M$ a closed form
$\omega(\nabla)\in \Omega^{n}_{cl}(M,\C)$
such hat
$\omega((F,f)^{*}\nabla)=f^{*}\omega(\nabla)$
for all pull-back diagrams
$$\xymatrix{E^{\prime}\ar[d]\ar[r]^{F}&E\ar[d]\\M^{\prime}\ar[r]^{f}&M}\ .$$
\end{ddd}

\begin{prob}\label{apr1601}
Let $\omega$ be a characteristic form of degree $\ge 1$.
Show that
$\omega(\nabla)=0$ if
$\nabla$ is flat.
\end{prob}
\proof
Reduce to a local question. Locally a flat bundle is pulled back from a point.
On a point there are no higher-degree forms.
\hB 

\begin{ex}{\em
The form
$\dim$ which associates to $(E,\nabla)$ the function $\dim(E)\in \Omega^{0}_{cl}(M;\C)$
is a characteristic form of degree $0$.
}
\end{ex}

Our main examples are Chern classes and the components of the Chern character. 
We start with Chern classes. The determinant is a fibrewise polynomial function
$\det:\End(E)\to M\times \C$ and therefore extends to
$$\det:\Omega^{ev}(M,\End(E))\to \Omega^{ev}(M,\C)\ .$$
For example,
$$\det\left(\begin{array}{cc}dx\wedge dy&dy\wedge dz\\dx\wedge dz&dz\wedge du\end{array}\right) =dx\wedge dy\wedge dz\wedge du\ .$$

The curvature $R^{\nabla}$ (cf. Definition \ref{apr1201}) of the connection $\nabla$ is an element $$R^{\nabla}\in \Omega^{2}(M,\End(E))\subset \Omega^{ev}(M,\End(E))\ .$$ 
\begin{ddd} The homogeneous pieces of the total Chern form which is defined by $$c(\nabla):=\det(1-\frac{1}{2\pi i} R^{\nabla})=1+c_{1}(\nabla)+\dots + c_{n}(\nabla)\ , \quad  c_{n}(\nabla)\in \Omega^{2n}(M,\C)$$
are called the Chern forms of $\nabla$.
\end{ddd}
For example $$c_{1}(\nabla)=-\frac{1}{2\pi i}\Tr(R^{\nabla})\ .$$
\begin{lem}
The Chern forms are characteristic forms.
\end{lem}\label{apr1230}
\proof Use \ref{apr1220} for compatibility with pull-back. 
In order to see closedness we use that for
  $\Phi\in \Omega^{\ge 2}(M,\End(E))$ we have $$d\det(1+\Phi)=\det(1+\Phi)\Tr \left(
  (1+\Phi)^{-1}
  [\nabla ,\Phi]\right)\ .$$ 
\begin{prob}\label{apr1225}
Prove this formula.
\end{prob}
\proof
Note that an element $\Phi\in \Omega^{\ge 1}(M,\End(E))$ is nilpotent. We understand $\exp (\Phi)$ or $\log(1+\Phi)$ in the sense of formal power series. 
We write
$$\log\det(1+\Phi)=\Tr(\log(1+\Phi))\ .$$
Using that $d$ is a derivation    we get
$$d\log \det(1+\Phi)=d \Tr (\log(1+\Phi))\ .$$
We now use 
$d \Tr(\Psi)= \Tr( [\nabla,\Psi])$ which implies
$d \Tr \Phi^{n}= n\Tr \left( \Phi^{n-1} [\nabla,\Phi]\right)$, and finally $(1+\Phi)^{-1}=\sum_{n=0}^{\infty} (-1)^{n} \Phi^{n}$. \hB 

If we insert $\Phi:=-\frac{1}{2\pi i}R^{\nabla}$ and the Bianchi identity \eqref{apr1224} we get $d c(\nabla)=0$.\hB

If $(E,\nabla^{E})$ and $(F,\nabla^{F})$ are two vector bundles with connection on $M$, then
we obtain an induced connection $\nabla^{E\oplus F}$ on $E\oplus F$.
It is given by 
$$\nabla^{E\oplus F}(\phi\oplus \psi)=\nabla^{E}\phi\oplus \nabla^{F} \psi\ , \quad \phi\in \Omega^{0}(M,E)\ , \: \psi\in \Omega^{0}(M,F)\ .$$

\begin{prob}
Show that
$c(\nabla^{E\oplus F})=c(\nabla^{E})\wedge c(\nabla^{F})$.
\end{prob}

\begin{ddd}\label{may2803}
The  Chern character form is given by 
$$\ch(\nabla):=\Tr \exp(-\frac{1}{2\pi i} R^{\nabla})\in \Omega^{ev}(M)$$
\end{ddd}

We write $$\ch=\ch_{0}\oplus \ch_{2}\oplus \dots\ ,\quad \ch_{2n}\in \Omega^{2n}(M,\C)$$
for its decomposition into homogeneous pieces.
 \begin{lem}
For all $k\ge 0$ the homogeneous piece $\ch_{2k}$ of the Chern character form  is a characteristic form.
\end{lem}
\proof
We argue as in the proof of Lemma \ref{apr1230}.
For closedness we use again the Bianchy identity and
the relation
$$d \Tr(\exp(\Phi))=\Tr(\exp(\Phi)[\nabla,\Phi] )\ .$$ \hB

In the following we show that the Chern character forms are compatible with the tensor operations.
If $(E,\nabla^{E})$ and $(F,\nabla^{F})$ are two vector bundles with connection on $M$, then the tensor product $E\otimes F$  and the bundle $\Hom(E,F)$ have induced connections $\nabla^{E\otimes F}$ and $\nabla^{\Hom(E,F)}$.
\begin{prob}
Show that $\nabla^{E\otimes F}$ and $\nabla^{\Hom(E,F)}$ are  uniquely determined by the conditions that
$$\nabla^{E\otimes F}(\phi\otimes \psi)=\nabla^{E} \phi\otimes \psi+\phi\otimes \nabla^{F} \psi\ , \quad \phi\in \Omega^{0}(M,E)\ ,\: \psi\in \Omega^{0}(M,F)\ ,$$
$$\nabla^{\Hom(E,F)}\Psi=\nabla^{F}\circ \Psi-\Psi\circ \nabla^{E}\ ,\quad \Psi\in \Omega^{0}(M,\Hom(E,F))\ .$$ 
\end{prob}
We equip $M\times \C\to M$ with the trivial connection and define $\nabla^{E^{*}}$ as the connection induced by $\nabla^{E}$ on
$E^{*}:=\Hom(E,M\times \C)$.

\begin{prob}
Show that
$R^{\nabla^{E\otimes F}}=R^{\nabla^{E}}\otimes \id_{F}+\id_{E}\otimes R^{\nabla^{F}}$ and conclude that
$$\ch(\nabla^{E\oplus F})=\ch(\nabla^{E})+\ch(\nabla^{F})\ ,\quad \ch(\nabla^{E\otimes F})=\ch(\nabla^{E})\wedge \ch(\nabla^{F})\ .$$
Further show that
$$\ch_{2i}(\nabla^{E^{*}})=(-1)^{i} \ch_{2i}(\nabla^{E})\ , \quad  \ch(\nabla^{\Hom(E,F)})=\ch(\nabla^{E^{*}})\wedge \ch(\nabla^{F})\ .$$ 
 \end{prob}

The identity
$$\det(1+X)=\exp(\Tr(\log(1+X)))=\exp(\sum_{n=1}^{\infty} (-1)^{n-1} \frac{\Tr(X^{n})}{n})$$ implies
$$c(\nabla)=\exp\left(\sum_{n=1}^{\infty} (-1)^{n-1} (n-1)! \ch_{2n}(\nabla)\right)\ .$$
Therefore $c_{i}(\nabla)$ can be expressed as a rational polynomial in the $\ch_{k}$ for $k\le i$.
\begin{prob}
Show that one can reverse this and express
$\ch_{2i}(\nabla)$ as a rational polynomial in the $c_{k}(\nabla)$ for $k\le i$.
For example, in degree $2$ we have 
$$c_{1}(\nabla)=\ch_{2}(\nabla)\ ,$$
Determine  these expressions in degree $4$ and $6$.
Give a formula for
$c_{2}(\nabla^{E\otimes F})$ in terms of the Chern forms of $\nabla^{E}$ and $\nabla^{F}$.
\end{prob}

We now consider the dependence of characteristic forms on the connection. The answer is given in terms of transgression.
Consider two connections $\nabla,\nabla^{\prime}$ on a bundle
$E$. Then we can make a connection
$\tilde \nabla$ on $\pr^{*}_{M}E\to \R\times M$ such that
$\nabla=\tilde \nabla_{|\{0\}\times M}$ and $\nabla^{\prime}=\tilde \nabla_{|\{1\}\times M}$.
For example,
take $$\tilde \nabla=t\pr_{M}^{*}\nabla^{\prime}+(1-t) \pr_{M}^{*}\nabla\ .$$
We say that $\tilde \nabla$ is a path between $\nabla$ and $\nabla^{\prime}$.

Let $\omega$ be a characteristic form of degree $n$.
Then we define the transgression form
$$\tilde \omega(\nabla^{\prime},\nabla):=\int_{[0,1]\times M/M} \omega(\tilde\nabla)\in \Omega^{n-1}(M;\C)\ .$$
\begin{lem}
The class $$[\tilde \omega(\nabla^{\prime},\nabla )]\in \Omega^{n-1}(M;\C)/\im(d)$$
is independent of the choice of the path $\tilde \nabla$ between $\nabla$ and $\nabla^{\prime}$.
We have
\begin{equation}\label{apr1240}d\tilde \omega(\nabla^{\prime},\nabla)=\omega(\nabla^{\prime})-\omega(\nabla)\end{equation} and
\begin{equation}\label{apr1241}[\tilde \omega(\nabla^{\prime\prime},\nabla^{\prime})]+[\tilde \omega(\nabla^{\prime},\nabla)]=[\tilde\omega(\nabla^{\prime\prime},\nabla)]\end{equation}
\end{lem}
\proof
The relation \eqref{apr1240} follows immediately from Stokes' theorem.
We consider the hyperplane $A^{2}=\{x_0+x_1+x_2=1\}\subset \R^{3}$. Its intersection with the positive quadrant is the standard simplex $\Delta^{2}$.  
Assume that we are given connections $\nabla_i$ on $E$ for $i=0,1,2$. Then we consider a connection $\hat \nabla$ on $\pr_M^{*}E\to A^{2}\times M$ such that 
$\nabla_i:=\hat \nabla_{\{x_i=1\}\times M}$. We can further arrange that the restrictions
$\hat \nabla_{\{x_i=0\}\times M}$ coincide with previously given interpolations.
 By Stokes' theorem
\begin{equation}\label{apr2101}\tilde \omega(\nabla_1,\nabla_{0})+\tilde \omega(\nabla_2,\nabla_{1})-\tilde \omega(\nabla_2,\nabla_0)=d \int_{\Delta^{2}} \tilde \omega(\hat \nabla)\ .\end{equation}
This implies \eqref{apr1241} provided we have shown the first assertion. For this we use \eqref{apr2101}
where we take $\nabla_1=\nabla_2$ and the constant interpolation $\pr_M^{*}\nabla$
between them.  Note that in this case $\tilde \omega(\nabla_{2},\nabla_{1})=0$. \hB

\begin{ddd}
The map $$(\nabla^{\prime},\nabla)\mapsto 
[\tilde \omega(\nabla^{\prime},\nabla)]\in \Omega^{n-1}(M;\C)/\im(d)$$
is called the transgression of $\omega$.
\end{ddd}

\begin{kor}
The cohomology class of $\omega(\nabla)$ only depends on the bundle $E$.
\end{kor}
\begin{ddd}
We define the  characteristic class associated to the characteristic form  $\omega$ such that it maps the bundle $E\to M$ to the class
$$\omega(E):=[\omega(\nabla)]\in H^{n}_{dR}(M;\C)\ ,$$
where $\nabla$ is any choice of connection on $E$.
\end{ddd}

\begin{prob}
Show that
$$[\tilde \ch(\nabla_{1}\otimes \nabla_{1}^{\prime}, \nabla_{2}\otimes \nabla_{2}^{\prime})]=[\tilde \ch(\nabla_{1},\nabla_{2})\wedge \ch(\nabla^{\prime}_{1})] +  [\ch(\nabla_{2})\wedge \tilde \ch(\nabla_{1}^{\prime}, \nabla^{\prime}_{2})]\ .$$
\end{prob}

\subsection{Characteristic forms, classes and transgression - examples}

We start with some examples in which Chern classes are calculated explicitly.

\begin{ex}\label{apr2102}{\em
The following problem is important since it shows that Chern classes are non-trivial. In the axiomatic approach (e.g. in \cite{MR0440554}) its solution is used to normalize the Chern classes.
\begin{prob}
Show that the class $c_1(L)$ is non-trivial, where $L\to \C\P^{n}$ is the tautological bundle.
\end{prob}\proof Use Example \ref{apr1260} and calculate $\langle c_{1}(L),[\C\P^{1}]\rangle$. \hB

Consider a complex vector bundle $E\to M$ and
let $\P(E)\to M$ be the associated bundle of projective spaces.
Recall that a  point $p\in \P(E)_{m}$ is a line in $E_{m}$. The tautological bundle $L\to \P(E)$ is the 
subbundle $L\subset p^{*}E$ given by
$$L:= \{(H,e)\in p^{*}(E)\subset \P(E)\times E\: |\:e\in H\}\ .$$
\begin{prob}  \label{apr3002} (Leray-Hirsch theorem)
Show that $(c_{1}(L)_{\P}^{i})_{i=0,\dots,\dim(E)-1}$ forms a basis of
$H^{*}(\P(E);\C)$ as a $H^{*}(M;\C)$-module.
\end{prob}
\proof
Use Example \ref{apr2102}  and the fact that
$H^{*}(\C\P^{n};\C)\cong \C[c]/(c^{n+1})$ for a generator $c\in H^{2}(\C\P^{n};\C)$ to show that
the classes $(c_{1}(L)_{\P}^{i})_{i=0,\dots,\dim(E)-1}$ restrict to a basis of  $H^{*}(\P(E)_{m};\C)$ for all $m\in M$.  Then argue by induction over the cells of the basis $M$. \hB 
 
 It is known that, as a ring,  $H^{*}(\C\P^{n};\C)\cong \C[c]/(c^{n+1})$ with $c:=-c_1(L)$.
\begin{prob}
Calculate the Chern classes of $T\C\P^{n}$.
\end{prob}
\proof
To this end construct a sequence
$$0\to L^{*}\otimes T^{*}\C\P^{n}\to \C\P^{n}\times \C^{n+1}\to L^{*} \to 0$$
and deduce that
$$c(T\C\P^{n})=(1+c)^{n+1}\ .$$ \hB 
}
\end{ex}

\begin{ex}{\em 
We reconsider the Poincar\'e bundle $P\to J(M)\times M$ (see \ref{apr1204}).
We have a canonical isomorphism $a: H_1(M;\C) \stackrel{\sim}{\to} H^{1}_{dR}(J(M)^{0};\C)$.
 Let $x\in H_1(M;\C)$ be represented by a homotopy class $\gamma\in \pi_1(M,m)$. Then $\frac{1}{2\pi i}\rho^{-1}(\gamma)d\rho(\gamma)\in \Omega^{1}_{cl}(J(M))$ represents $a(x)$. We consider the map $a$ as a class
$$\tilde a\in H_{dR}^{1}(M;\C)\otimes_{\C} H_{dR}^{1}(J(M)^{0};\C) \subset H_{dR}^{2}(J(M)^{0}\times M;\C)$$
\begin{lem}
We have
$c_1(P)= \tilde a$.
\end{lem}
\proof
We use the notation introduced in \ref{apr1204}.
The curvature of $\nabla^P$ lifts to the $\pi$-invariant form $-d\alpha$.
By construction the curvature represents the class $-\theta\in H^{1}_{dR}(M;\C)\otimes_{\C} \Omega_{cl}(J(M)^{0},\C)$. Its total cohomology class is thus $- 2\pi i\: a$. \hB
We observe that $c_1(P)$ is non-zero if $H^{1}_{dR}(M;\C)\not=0$.

}
\end{ex}

We now investigate the transgression of a characteristic form between two flat connections and how it gives rise to a  secondary characteristic class.
If $\nabla$ and $\nabla^{\prime}$ are flat, then
$d\tilde \omega(\nabla^{\prime},\nabla)=0$ by Problem \ref{apr1601}.
\begin{ddd}\label{apr1290}
The secondary characteristic class associated to $\omega$ maps 
the pair of flat connections $\nabla^{\prime},\nabla$ to the cohomology class
$$[\tilde \omega(\nabla^{\prime},\nabla])\in H^{n-1}_{dR}(M;\C)\ .$$
\end{ddd}

Note that $\omega(\nabla,\nabla)=0$.

\begin{ex}
{
\em 
We consider the trivial bundle
$S^{1}\times \C$ with trivial connection $d$.
For $\alpha\in \Omega^{1}(S^{1},\C)$ we consider the connection
$d+\alpha$. 
\begin{prob}\label{apr1291}
Calculate
$\tilde c_1(d+\alpha,d)$.
\end{prob}
The result of this calculation is
$$\tilde c_1(d+\alpha,d)=[-\frac{1}{2\pi i} \alpha]\in H^{1}_{dR}(S^{1};\C)\ .$$
This class is non-zero in general. It determines the holonomy of $d+\alpha$ by
\begin{equation}\label{apr1280}\hol_{d+\alpha}(S^{1})=\exp(-\int_{S^{1}}  \alpha)=\exp(-\langle [\alpha],[S^{1}]\rangle )=\exp(2\pi i \langle \tilde c_1(d+\alpha,d),[S^{1}]) .\end{equation} \hB

  By naturality, the calculation determines the class $\tilde c_{1}(\nabla^{\prime},\nabla)$ for any pair of flat connections on a line bundle $L\to M$. We have
\begin{equation}\label{apr12200}\exp(2\pi i \langle \tilde c_{1}(\nabla^{\prime},\nabla),[\gamma]\rangle)=\hol_{\nabla^{\prime}}(\gamma)\hol_{\nabla}(\gamma)^{-1}\ .\end{equation}
Like the holonomy the class $\tilde c_{1}(\nabla^{\prime},\nabla)$ can vary continously with the flat connections in a non-trivial way. On the other hand we have the following.

\begin{lem}\label{apr1294}
If $\nabla$ and $\nabla^{\prime}$ belong to the same path component of the space of flat connections, then $\tilde c_{n}(\nabla^{\prime},\nabla)=0$ for all $n\ge 2$. A similar statement holds true for all characteristic forms of degree $\ge 3$.

\end{lem}
\proof
We only consider the case of $c_n$.
We consider a connection $\tilde \nabla$ on $\pr_{M}^{*}E\to \R\times M$ such that
$\tilde \nabla_{|\{t\}\times M}$ is flat for all $t$. Then
$R^{\tilde \nabla}=dt\wedge \iota_{\partial_{t}}R^{\tilde \nabla}$. This yields the following formula for the total Chern class
$$c(\tilde \nabla)=1-dt\wedge \frac{1}{2\pi i} \Tr R^{\tilde \nabla}\ .$$ The vanishing of the higher degree components implies the assertion. \hB
The classes $\tilde c_{n}(\nabla^{\prime},\nabla)$ can thus be used to show that two flat connections belong to different path components. Examples for this will be given later, see Examples \ref{apr1602}.

}
\end{ex}

\begin{ex}\label{apr12100}{\em We consider a vector bundle $E\to M$.
For an automorphism  $F\in \Aut(E)$ of $E$ we define the form
$$\omega(F):=[\tilde \omega(F^{*}\nabla,\nabla)]\in \Omega^{n-1}(M;\C)/\im(d)\ .$$
\begin{lem}
We have
$d\omega(F)=0$ and therefore
$\omega(F)\in H^{n-1}_{dR}(M;\C)$.
The class $\omega(F)$ is independent of the choice of $\nabla$ and only depends on the homotopy class of $F$.
It satisfies
\begin{equation}\label{apr1250}\omega(F\circ F^{\prime})=\omega(F)+\omega(F^{\prime})\ .\end{equation}
\end{lem}
\proof
The first assertion follows from $\omega(\nabla)=\omega((F,\id)^{*}\nabla)$.
The identity \eqref{apr1250} follows from \eqref{apr1241}. 
A homotopy can be understood as an element $\tilde F\in \Aut(\pr_M^{*}E)$,
where $\pr_M:\R\times M\to M$ denotes the projection. Then $\omega(\tilde F)$ is closed. 
Hence its restriction to $\{t\}\times M$ does not depend on $t\in \R$. \hB

For every characteristic form of degree $n$ we have defined a group homomorphism
$$\omega:\pi_0(\Aut(E))\to H_{dR}^{n-1}(M;\C)\ .$$

Here is an alternative interpretation of the class $\omega(F)$. We consider the action of $\Z$ on
$\pr_{M}^{*}E=\R\times E\to \R\times M$ given by
$n(t,e):=(t+1,F^{-1}(e))$, $n(t,m):=(t+1,m)$, 
and let
$E(F)\to S^{1}\times M$ be the quotient.
\begin{ddd} The complex vector bundle
$E(F)\to S^{1}\times M$ is called the suspension of $E$ with respect to $F$.
\end{ddd}
We then have
$$\omega(F)=\int_{S^{1}\times M/M}\omega(E(F))\ .$$
\begin{prob}
Prove this!
\end{prob}
 
Let $G$ be a Lie group and $\rho:G\to GL(n,\C)$ be a representation. 
Let $G\times \C^{n}\to G$  be the $n$-dimensional trivial bundle.  Then we have a tautological element $F_{\rho}\in \Aut(G\times \C^{n})$
given by $F_\rho(g,v):=(g,\rho(g)v)$.
We therefore get a class
$$\omega(F_{\rho})\in H^{n-1}_{dR}(G;\C)\ .$$
Recall that a cohomology class $x\in H^{*}_{dR}(G;\C)$ is called primitive, if
$$\mu^{*}x=x\cup 1+1\cup x\ ,$$
where $\mu:G\times G\to G$ is the multiplication map.
\begin{lem}
$\omega(F_{\rho})$ is primitive.
\end{lem}
\proof
On $G\times G\times \C^{n}\to G\times G$ we have
the automorphisms $F_i$ given by 
$F_0(h,g,v):=(h,g,\rho(h)v)$\ ,
$F_1(h,g,v) :=(h,g,\rho(hg)v)$, and 
$F_2(h,g,v):=(h,g,\rho(g)v)$. Note that
$F_0\circ F_2=F_1$. We have
$F_0=\pr_0^{*} F_\rho$, $F_1=\mu^{*} F_\rho$ and $F_2=\pr_1
^{*}F_0$, where $\pr_0,\pr_1,\mu:G\times G\to G$ are the projections and multiplication.
We thus have
$\mu^{*}\omega(F_\rho)=\pr_0^{*}\omega(F_\rho)+\pr_1^{*}\omega(F_\rho)$ as required. \hB

If we apply this construction to the Chern classes $c_{k}$ and the standard representation $\id$ of the group $U(n)$ we get primitive classes $c_{2k-1}=c_k(F_\id)\in H_{dR}^{2k-1}(U(n);\C)$. It is known that
\begin{equation}\label{apr2301}H^{*}_{dR}(U(n);\C)=\Lambda_\C(c_{1},\dots,c_{2n-1})\ .\end{equation}
We can calculate the forms $\ch_{2n}(F_\id)= \tilde \ch_{2n}((F_\id,\id)^{*}\nabla,\nabla)$ explicitly. We consider the Mauer-Cartan form $g^{-1}dg\in \Omega^{1}(U(n))\otimes \Mat(n,\C)$.
\begin{lem}\label{apr2601}
We have
$$\ch_{2n}(F_\id)=\frac{(-1)^{n-1}(n-1)!}{(2\pi i)^{n}(2n-1)!} \Tr (g^{-1}dg)^{2n-1}\ .$$
\end{lem}
\proof Explicit calculation following the definitions.
\begin{prob}\label{apr1271}
Do this calculation!
\end{prob}
}
\end{ex}

\begin{prob}
Calculate $H^{*}(U(n);\Z)$ inductively using the Serre spectral sequences
associated to the fibrations
$$ U(n-1)\to U(n)\to S^{2n-1}\ .$$
Deduce \eqref{apr2301}.
\end{prob}

\subsection{Metrics and unitarity}

We consider a complex vector bundle $E\to M$ with a hermitean metric $h$.
Let $\nabla$ be a connection on $E$.
\begin{prob} Show that
there is a unique connection $\nabla^{*}$ (called the adjoint of $\nabla$ w.r.t $h$), which
is characterized by 
$$d h(\phi,\psi)= h(\nabla \phi,\psi)+h(\phi,\nabla^{*}\psi)\ ,\quad  \forall \phi,\psi\in \Omega^{0}(M,E)\ .$$
\end{prob}
For $\alpha\in \Omega^{1}(M,\End(E))$ we have
$$(\nabla+\alpha)^{*}=\nabla^{*}-\alpha^{*}\ , $$
where here $*$ only acts on the endomorphism part.

\begin{ddd}
A connection on a metrized bundle $(E,h)$ is called unitary if $\nabla^{*}=\nabla$.
\end{ddd}

We write $\bu(E)\subset \End(E)$ for the subbundle of antihermitean endomorphisms.
If a connection $\nabla$ is unitary, then so is $\nabla+\alpha$ for all $\alpha\in \Omega^{1}(M,\bu(V))$. For any connection $\nabla$ define
$$\omega:=\nabla^{*}-\nabla\in \Omega^{1}(M,\End(E))\ .$$
Then
$$\nabla^{u}:=\nabla+\frac{1}{2} \omega$$ is unitary. It is called the symmetrization of $\nabla$.
\begin{prob} Show this result. If $A\in \Omega^{0}(M,\End(E))$ is invertible and symmetric, then we define $h_{A}(\phi,\psi):=h(A\phi,\psi)$. Calculate
$\nabla^{*_{h}}-\nabla^{*_{h_{A}}}$. Assume that $A$ is scalar and calculate $R^{\nabla^{*_{h}}}-R^{\nabla^{*_{h_{A}}}}$ in this case.
\end{prob}

\begin{ddd}
A connection is called unitarizable if it is unitary for some choice of hermitean metric.
\end{ddd}

There are various obstructions against unitarizablility as the following exercise shows. 
\begin{prob}\label{apr1270}
We consider the trivial bundle $\R^{2}\times \C^{2}\to \R^{2}$ with connections
$$\nabla:=d+\left(\begin{array}{cc}0&1\\0&0\end{array}
\right) dx\ ,\quad \nabla^{\prime}:=d+\left(\begin{array}{cc}0&1\\0&0\end{array}
\right) ydx$$
Show that $\nabla$ is unitarizable, while $\nabla^{\prime}$ is not.
 The connection $\nabla$ descends to the quotient
$\R^{2}/\Z^{2}$, but this descent is not unitarizable. 
\end{prob}


\begin{prob}
Show that
$$R^{\nabla^{*}}=-(R^{\nabla})^{*}\ .$$
Hence, if $\nabla$ is unitary, then 
$R^{\nabla}\in \Omega^{2}(M,\bu(E))$.
Conclude that for unitary connections $\ch(\nabla)$ and $c(\nabla)$ are real forms.
\end{prob}

Assume now that $\nabla$ is flat. Then $\nabla^{*}$ is flat, too. For a characteristic form
$\omega$ of degree $n$ we therefore can consider the associated secondary class Def. \ref{apr1290}
$$\tilde \omega(\nabla^{*},\nabla)\in H_{dR}^{n-1}(M;\C)\
 .$$
A priori this is a characteristic class for metrized bundles.

\begin{lem}
The class
$\tilde \omega(\nabla^{*},\nabla)\in H_{dR}^{n-1}(M,\C)$ does not depend on the choice of
metric.
\end{lem}
\proof If $\tilde h$ is a metric on the flat bundle
$\pr_M^{*}E\to \R\times M$, then  $\tilde \omega((\pr_M^{*}\nabla)^{*},\pr_M^{*}\nabla)$ is closed. Hence its restriction to $\{t\}\times M$ is independent of $t$. 
This shows the assertion, since every two metrics can be connected. 
\hB 

\begin{ddd}
We write
$\check \omega(\nabla)\in H_{dR}^{n-1}(M,\C)$ for the
characteristic class for flat connections obtained from
$\omega$.
\end{ddd}

If $\nabla$ is unitarizable, then $\check \omega(\nabla)=0$.
If $\deg(\omega)\ge 3$, then it follows from Lemma \ref{apr1294} that 
$\check{\omega}(\nabla)$ only depends on deformation class of the flat connection $\nabla$.
This is the rigidity result \cite{MR1867006}.

\begin{ex}
{\em Let $\lambda\in \C\setminus\{0\}$.
We consider the action of $\Z$ on
$\R\times \C\to \R$ by $n(t,z):=(t+n,\lambda^{n} z)$.
It preserves the trivial connection.
Hence, by taking the quotient, we get a flat bundle
$(E_\lambda,\nabla_\lambda)\to S^{1}$.
We want to calculate
$\check c_1(\nabla_\lambda)\in H_{dR}^{1}(S^{1};\C)$.
The holonomy of $\nabla_\lambda$ along  $S^{1}$ is multiplication $\lambda^{-1}$.
The holonomy of $\nabla_\lambda^{*}$ is thus $\bar \lambda^{-1} $.
Note that $E_\lambda$ is trivial.
Hence
we have
$$\check c_1 (\nabla_\lambda)=\tilde c_1(\nabla_\lambda^{*},\nabla^{triv})+\tilde c_1(\nabla^{triv},\nabla_\lambda)\ .$$
We now use
 \eqref{apr1280} and get
$$ \exp(2\pi i \check  c_1 (\nabla_\lambda),[S^{1}]\rangle)=
 |\lambda|^{-2}\ .$$
 It follows that
 $$ \langle \check c_1 (\nabla_\lambda),[S^{1}]\rangle= \frac{i}{\pi} \log|\lambda| \ .$$
This calculation shows that $\check  c_{1}(\nabla)$ varies continuously with the flat connection $\nabla$.
In contrast, by Lemma \ref{apr1294} the classes $\check c_{n}(\nabla)$ for $n\ge 2$ only depend on the path components of $\nabla$. 
}
\end{ex}

The following is a very non-trivial example.
\begin{ex}\label{apr1602}
{\em We consider a number ring $R$ and an embedding $\sigma:R\to \C$. 
We fix $k,n\in \nat$ and  consider a manifold $M$ together with an $n$-equivalence
$f:M\to BGL(R,k)$. Since $BGL(R,k)$ is a countable $CW$-complex its skeleta can be  realized 
by smooth manifolds which gives the existence of $M$. The map $f$ classifies a
$GL(R,k)$-bundle $\tilde M\to M$, and we let
$E:=\tilde M\times_{GL(R,k),\sigma} \C^{k}\to M$ be the associated bundle with flat connection 
$\nabla$.  The class
$$\check \ch_{2j}(\nabla)\in H^{2n-1}(M;\C)$$
is called the Kamber-Tondeur class (compare \cite{MR1303026}). In \cite{MR0387496}  Borel has shown the following:
If $\sigma$ is  real and $j$ is even, then the class
$\check \ch_{2j}(\nabla)$ is non-zero for sufficiently large $k,n$.
If $\sigma$ is a complex embedding, then given $j$ the class  $\check \ch_{2j}(\nabla)$ is non-zero for sufficiently large $k,n$. A proof is beyond the scope of this course.

The non-triviality of $\check \ch_{6}(\nabla)$ shows that $\nabla$ can not be connected with any adjoint connection $\nabla^{*}$ by  a path of flat connections on $E$.  }
\end{ex}

\subsection{Integrality}

Let \begin{equation}\label{apr3001}\epsilon_{\C}:H^{n}(M;\Z)\to H^{n}_{dR}(M;\C)\end{equation}
be the map induced by the inclusion $\Z\subset \C$ and the de Rham isomorphism.

\begin{ddd}\label{apr2402}
A class $x\in H^{n}_{dR}(M;\C)$ is called integral if it belongs to the image of $\epsilon_{\C}$.
\end{ddd}
Equivalently, a class $x\in  H^{n}_{dR}(M;\C)$ is integral if and only if
$$\langle x,z\rangle \in \Z$$
for all smooth cycles $z\in Z_{n}(S_\infty(M))$.
We shall use the fact that the subset of integral de Rham cohomology classes is
closed under pull-back along smooth maps and products. All this follows from the fact that \eqref{apr3001}
comes from a natural transformation of multiplicative cohomology theories.

\begin{ddd}
A characteristic form for complex vector bundes $\omega$ is called integral if
$\omega(E) $ is integral 
for all complex vector bundles $E\to M$.
\end{ddd}

\begin{prop}\label{apr1295}
The Chern classes $c_{n}$ are integral.
\end{prop}
\proof
We first check integrality of $c_{1}(H)$ for line bundles $H$.
We use Example \ref{apr1260} in order to see that
$c_{1}(L)$ is integral for the tautological bundle $L\to \C\P^{n}$.
Given a line bundle $H$ on a manifold $M$ there exists $n$ and a map $f:M\to \C\P^{n}$ such that
$H\cong f^{*}L$. Hence $c_{1}(H)=f^{*}c_{1}(L)$ is integral.

We now discuss $c_{n}$. If the bundle $E$ has a decomposition $E\cong H_{1}\oplus \dots\oplus H_{k}$ into line bundles, then $$c(E)=\prod_{i=1}^{k} c(H_{i})=\prod_{i=1}^{k} (1+c_{1}(H_{i}))$$
has integral homogeneous components. Let now $E_{0}:=E\to M$ be general of dimension $r$. 
In this case we argue by the splitting principle which goes as follows. We can choose a decomposition $\pr_{1}^{*}E_{0}\cong L_{1}\oplus E_{1}$, where
$\pr_{1}:\P(E_{0})\to M$ is the projective bundle of $E$ and $L_{1}\to \P(E_{0})$ is the canonical bundle. 
We apply the same construction to $E_{1}\to \P(E_{0})$ and then inductively. We obtain a bundle
$q:F(E)\to M$ and a decomposition $q^{*}E\cong H_{1}\oplus\dots \oplus H_{r}$, where 
$q=\pr_{1}\circ \pr_{2}\circ \dots\circ \pr_{r-1}$ and e.g. $H_{1}=(\pr_{2}\circ \dots\circ \pr_{r-1})^{*}
L_{1}$. Hence $q^{*}c(E)$ has integral components. It is known that the integral cohomology of $F(E)$ is a free module over the integral cohomology ring of  $M$. This follows from a Leray-Hisch argument
similar to \eqref{apr3002}. We finally use the following assertion:
\begin{prob}
Show that if $x\in H^{n}_{dR}(M,\C)$ and $q^{*}x$ is integral, then so is $x$. 
\end{prob}
\hB

It follows that integral polynomials in the $c_{i}$ are integral.

\begin{ex}\label{jun1102}{\em
The following exercise gives a direct integral interpretation of the highest non-trivial Chern class.
Let $E\to M$ be a complex vector bundle of dimension $n$   and $S(E)\to M$ be its sphere bundle. We consider the associated Serre  spectral sequence. Its second term has two rows
$$E_2^{*,0}\cong H^{*}(M;\Z) \ ,\quad E_2^{*,2n-1}\cong H^{*}(M;\Z)\ori_{S^{2n-1}}\ .$$ The only non-trivial differential is
$d_{2n-1}$. The corresponding edge sequence is called Gysin sequence.
We define the Euler class  
$$\chi:=d_{2n-1}(\ori_{S^{2n-1}})\in H^{2n}(M;\Z)\ .$$
\begin{prob}
Show that $c_n(E)=\chi$.
\end{prob}
\proof
Show that this formula is compatible with sums and reduce to the one-dimensional case.\hB

}\end{ex}

\begin{lem}\label{may0204}
The Newton classes $s_{n}:=n!\ch_{n}$ are integral.
\end{lem}
\proof
For line bundles  $H$ we have
$\ch(H)=\exp(c_{1}(H))$. Hence
$s_{n}(H)$ is integral by \ref{apr1295}. If $E$ decomposes into line bundles $E\cong H_{1}\oplus \dots\oplus H_{k}$, then $s_{n}(E)=\sum_{i=1}^{k} s_{n}(H_{i})$ is integral. The general case is reduced to the decomposable case by the splitting principle as in the proof of \ref{apr1295}.
\hB

If the characteristic form $\omega$ is integral, then
$\omega(F)\in H^{n-1}_{dR}(M)$ (see Example \ref{apr12100}) is integral for every
$F\in \Aut(E)$. This allows to define an absolute invariant of flat connections $\nabla$ on trivializable bundles $E$. Let $\phi,\phi^{\prime}$ be trivializations. Then we have
$F(\phi)=\phi^{\prime}$ for a suitable automorphism $F$. 
Hence by \eqref{apr1241}
$$[\omega(\nabla,\nabla^{\phi})]=[\omega(\nabla,\nabla^{\phi^{\prime}})]+\omega(F)\ ,$$
where $\nabla^{\phi}$ denotes the connection induced by the trivialization $\phi$ as in Problem \ref{apr1202}.
It follows that the class
$$[\tilde \omega(\nabla,\nabla^{\phi})]\in \frac{H^{n-1}_{dR}(M;\C)}{\im \epsilon_{\C}}$$
is independent of the choice of the trivialization $\phi$.

\begin{ddd}\label{apr2220}
We define the Chern-Simons invariant of a flat connection $\nabla$ on a trivializable bundle $E$ associated to the characteristic form $\omega$ as
$$\tilde \omega(\nabla):=[\tilde \omega(\nabla,\nabla^{\phi})]\in 
 \frac{H^{n-1}_{dR}(M;\C)}{\im \epsilon_{\C}}\ ,$$
 where $\phi$ is some choice of trivialization of $E$.
 If $M$ is compact oriented of dimension $n-1$, then we set
$$\cs_{\omega}(\nabla):=[\langle \tilde \omega(\nabla,\nabla^{\phi}),[M]\rangle]\in \C/\Z\ .$$
\end{ddd}

 \begin{prob}
Calculate $\tilde c_{1}(\nabla)$.
\end{prob}
\proof
We have by \eqref{apr12200}
$$\exp(2\pi i \langle \tilde c_{1}(\nabla),[\gamma]\rangle )=\hol_{\nabla}(\gamma)\ .$$


 \begin{ex}\label{aug1940}{\em

Let $\omega$ be an integral characteristic form of degree $n$.
Let $B$ be a space with a complex vector bundle $V^{\delta}\to B$ with structure group reduced to $GL(k,\C^{\delta})$. Below we refer to $V^{\delta}$ as a flat bundle. By $V\to B$ we denote the associated bundle with structure group $GL(k,\C)$.
We assume that $V$ is trivializable.
 
If $f:M\to B$ is a continuous map from a manifold, then the flat structure on $f^{*}V$ provides a flat connection which we denote, by abuse of notation, by $f^{*}\nabla^{V^{\delta}}$. It is characterized uniquely by
$f^{*}\nabla^{V^{\delta}}(f^{*}\phi)=0$, where $\phi:U\to V^{\delta}$, $U\subseteq B$ is a local section of $V^{\delta}$.

  Let $f:M\to B$ represent an oriented bordism class in $\MSO_{n-1}(B)$. In particular $M$ is closed oriented of dimension $n-1$.  We consider the Chern Simons invariant \ref{apr2220}
$$ \cs^{V^{\delta}}_{\omega}(f):= \cs_{\omega}(f^{*}\nabla^{V^{\delta}})\in \C/\Z\ .$$
\begin{lem}\label{aug1901} The Chern-Simons invariant
$\cs_{\omega}(f)$ only depends on the bordism class of $f$.
It gives a homomorphism $\cs^{V^{\delta}}_{\omega}: \MSO_{n-1}(B)\to \C/\Z$.
\end{lem}
\proof
Additivity under disjoint union is clear. 
If $F:Z\to B$ is a zero bordism, then by Stokes' Theorem
$$\langle \tilde \omega(f^{*}\nabla^{V^{\delta}},\nabla^{triv}),[M]\rangle=\int_{M} \tilde \omega(f^{*}\nabla^{V^{\delta}},\nabla^{triv})
=\int_{Z} d \tilde \omega(F^{*}\nabla^{V^{\delta}},\nabla^{triv})=0\ .$$
This shows bordism invariance.  \hB

The orientation $\kappa:\MSO\to H\Z$ induces a natural transformation of homology groups
$\kappa:\MSO_{n-1}(B)\to H_{n-1}(B)$. Geometrically, if $f:M\to B$ represents $[f]\in \MSO_{n-1}(B)$, then $\kappa([f])=f_{*}[M]$. Later in \ref{jun0101} we will observe that
$\cs_{\omega}$ factorizes over this transformation.

The following is a consequence of a generalization  of \ref{apr1294} from Chern classes to $\omega$. 
\begin{prob}
If $\deg(\omega)\ge 3$, then
the Chern-Simons invariant $\cs^{V^{\delta}}_\omega$ only depends on the deformation class of $V^{\delta}$.
\end{prob}

Let $\lambda\in \C$ determine the character $\Z\to \C$ and therefore a flat bundle $V^{\delta}_{\lambda}\to S^{1}$. 
\begin{prob}Calculate the composition
$$ \Z\cong \pi_{1}(S^{1})\stackrel{can}{\to} \MSO_1(S^{1})\stackrel{\cs_{c_{1}}^{V^{\delta}_{\lambda}}}{\to} \C/\Z\ .$$
Use the result to conclude that  the homomorphism $can$ is an isomorphism.
\end{prob}
}
\end{ex}

\begin{ex}{\em 
Trivializability of the bundle $V\to B$ underlying $V^{\delta}$ is an annoying condition.
One can extend the range of the definition of the Chern-Simons invariant $\cs_{2n-1}^{V}$ as follows. We are going to define an invariant of flat bundles on $M$ which extend as bundles to some zero bordism:
Assume that
$\nabla$ is a flat connection on $E\to M$, and that $Z$ is a zero-bordism with
an extension $F\to Z$ of $E$. We choose some extension $\tilde \nabla$ of $\nabla$ to $Z$, not necessarily flat. Then we define  
\begin{equation}\label{aug1902}\cs_{\omega}(\nabla):=[\int_{Z}\omega(\tilde \nabla)]\in \C/\Z\ .\end{equation}
\begin{prob}
Show that this is a well-defined invariant of the flat connection $\nabla$.
\end{prob}
The common domain of both definitions of $\cs_{\omega}$ are flat connections on bundles $E\to M$ 
where $E$ is trivializable and $M$ is zero-bordant. 
\begin{prob}
Show that both constructions coincide on this common domain.
\end{prob} 
Let $(L,\nabla^{L})$ be the tautological line bundle of $\C\P^{n}$. For $k\in \nat $ we consider the power $H:=L^{\otimes k}$. 
Let $\pi:M\to \C\P^{n}$ be the unit sphere bundle of  $H$.
The bundle
$\pi^{*}H\to M$ is canonically trivialized.
Let $\nabla^{H,triv}$ be the associated trivial connection. The bundle $\pi^{*}L$ acquires a canonical flat connection $\nabla$ characterized as follows. For a local section $\psi$ of $\pi^{*}L$ we have
$\nabla \psi=0$ if and only if 
$\nabla^{H,triv}\psi^{\otimes k}=0$.
\begin{prob}\label{jul3130}
Calculate $\cs_{c_1^{n+1}}(\nabla)\in \C/\Z$.
\end{prob}
\proof
Let $q:Z\to \C\P^{n}$ be the disc bundle of $H$. Then $q^{*}L$ is an extension of $\pi^{*}L$ across $Z$. We let $\tilde \nabla$ be any extension of $\nabla$ over $Z$. Then we have
$$\int_Z c_1^{n+1}(\tilde \nabla)=\int_Z c_1^{n+1}(q^{*}\nabla^{L})+\int_Z d\widetilde{c_1^{n+1}}(\tilde \nabla, q^{*}\nabla^{L})=\int_M \widetilde{c_1^{n+1}}(\nabla, \pi^{*}\nabla^{L})\ .$$
Here we use Stokes theorem and that the integral over $Z/\C\P^{n}$ of a form pulled back from $\C\P^{n}$ vanishes. 
We now must calculate the transgression explicitly.
Let $\alpha\in \Omega^{1}(M;\C)$ be defined by $$\nabla+\alpha=\pi^{*}\nabla^{L}\ .$$
Then we form the connection 
$$\hat \nabla:=\pr_M^{*} \nabla +t\alpha$$
on $\pr_M^{*}\pi^{*}L\to \R\times M$. Its curvature is given by 
$$R^{\hat \nabla}=dt\wedge \alpha+t\pr_{\C\P^{n}}^{*} R^{\nabla^{L}}\ ,$$
where we use that
$\pi^{*}R^{\nabla^{L}}=[\nabla,\alpha]$.
We have
$$\widetilde{c_1^{n+1}}(\nabla, \pi^{*}\nabla^{L})=\frac{(-1)^{n+1}}{(2\pi i)^{n+1}}(n+1)\int_0^{1} t^{n} dt\wedge \alpha\wedge \pi^{*}(R^{\nabla^{L}})^{n}=\frac{-1}{2\pi i} \alpha\wedge \pi^{*}c_1(\nabla^{L})^{n}    \ .$$  
We now calculate $\int_{M/\C\P^{n}}\alpha$. We fix a base point and
identify the fibre of $H$ with $\C$. The fibre of $M$ is then identified with $U(1)$, and a typical parallel local section of $\pi^{*}L$ is given by a branch of $\psi(u):= u^{\frac{1}{k}}$. We get
$$\frac{du}{ku} \psi(u)=(d\psi)(u)=(\pi^{*}\nabla^{L}\psi)(u)=((\nabla+\alpha)\psi)(u)=\alpha \psi(u)$$
and conclude that $\alpha=\frac{du}{ku}$.
It follows that
$$\int_{M/\C\P^{n}}\alpha=\frac{2\pi i}{k}\ .$$
We conclude that
$$\cs_{c_1^{n+1}}(\nabla)=[\int_M \widetilde{c_1^{n+1}}(\nabla, \pi^{*}\nabla^{L})]=[-\frac{1}{k}\int_{\C\P^{n}} c_1(\nabla^{L})^{n}]
=[-\frac{1}{k}] \ .$$ 
\hB

\begin{prob}
Show that $\cs_{c_1^{n+1}}(\nabla^{\otimes r})=[-\frac{r}{k}]$.
\end{prob}

 
We can now modify the construction of the  bordism invariant \ref{aug1901} as follows.
Let $V^{\delta}\to B$ be flat and its underlying continuous bundle $V$  be classified   by $v:B\to BGL(k,\C)$.
Previously $\cs_{\omega}^{V^{\delta}}:\MSO_{n-1}(B)\to \C/\Z$ was defined under the condition that $v$ is homotopic to a constant map.  Using the construction \eqref{aug1902} we can now extend the definition to get a homomorphism
$$\cs_{\omega}^{V^{\delta}}:\ker\left(v_*:MSO_{n-1}(B)\to MSO_{n-1}(BGL(k,\C))\right)$$
by setting
$$\cs_{\omega}^{V^{\delta}}(f):=\cs_{\omega} (f^{*}\nabla)\ .$$
In \ref{apr2510} we will further  extend this homomorphism to all
of $MSO_{n-1}(B)$ using differential cohomology.
 
Consider the map 
 $v:BSL(k,\C^{\delta})\to BGL(k,\C)$ and observe that we then have
$0=v_{*}:MSO_{3}(BSL(k,\C^{\delta}))\to MSO_{3}(BGL(k,\C))$.  
Let $V^{\delta}\to BSL(k,\C^{\delta})$ be associated to the standard representation of $SL(k,\C^{\delta})$. 
\begin{prob}
Show that $\cs^{V^{\delta}}_{c_{1}^{2}}:MSO_{3}(BSL(k,\C^{\delta}))\to \C/\Z$ vanishes.
Furthermore show that
$\cs^{V^{\delta}}_{c_{2}}:MSO_{3}(BSL(k,\C^{\delta}))\to \C/\Z$
is non-trivial by calculating examples.
\end{prob}
\proof Calculate  
$\cs_{c_{2}}(f^{*}\nabla)-\cs_{c_{2}}((f^{*}\nabla)^{*})$ for suitable $f:M\to BSL(k,\C)$ for some metric on $f^{*}V$,
compare this with the Kamber-Tondeur class and use Borel's result that the latter generates $H^{3}(BSL(k,\C^{\delta});\R)$. \hB 

The Chern-Simons invariants for Seifert manifolds with maps to $BSp(1)^{\delta}$ have been determined in D. Auckly \cite{MR1277058}. More calculations can be found in \cite{MR1266279}, \cite{MR1054574}.
Here is a special case. 
Let $M\to S^{1}$ be a $2$-torus bundle. Then we have a sequence
\begin{equation}\label{may0301}0\to \pi_1(T^{2})\to \pi_1(M)\to \pi_1(S^{1})\to 0\ .\end{equation}
We consider generators $A,B\in \pi_1(T^{2})$ and an element $T\in \pi_1(M)$ which maps to a generator of $\pi_1(S^{1})$. We consider a representation
$\rho:\pi_1(M)\to Sp(1)$ with $\rho(T)=J$ (quaternionic notation) and $\rho(A)=\exp(2\pi i \phi)$, $\rho(B)=\exp(2\pi i \psi)$ for some $\phi,\psi\in \R$. Let $(V\to M,\nabla)$ be the associated two-dimensional flat bundle.
\begin{prob} 
Calculate 
$\tilde c_2(\nabla)$ in terms of the data $\psi,\phi$ and the extension \eqref{may0301}.
Discuss, under which conditions $\rho$ exists.
\end{prob}
}\end{ex}

\begin{ex}{\em
 Let $\omega$ be an integral characteristic class of degree $n$. Then
we can define an invariant of connections on $n-1$-dimensional manifolds.
Assume that $\nabla$ is a connection on a bundle $E\to M$ over an oriented closed $n-1$-dimensional manifold $M$ which extends to a bundle $F\to Z$ on an oriented  zero bordism $Z$ of $M$.
Then we can choose an extension  $\hat \nabla$ of   $\nabla$ to $F$ with product structure near to the boundary and define
$$\cs_\omega(\nabla):= [\int_{Z} \omega(\hat \nabla)]\in \C/\Z\ .$$

\begin{prob}\label{aug1930}
Show that 
$\cs_\omega(\nabla)$ is independent of the choice of the zero bordism and extension. 
 \end{prob}
 \proof
 Two choices can be glued along $M$.
 The difference of the corresponding two integrals is integral by the integrality of $\omega$. See also the proof of  \ref{aug1903}.
 \hB 
 The assumption of a product structure simplifies the argument for \ref{aug1930} but is not necessary.
 \begin{prob}
 Verify this assertion.
 \end{prob}
 \proof
 See \ref{may2130} for a more general case. \hB 

Let $\omega$ be an integral characteristic form of degree $n$.
We consider a hermitean bundle $(E,h)$ on a manifold $M$ which extends to a zero bordism of $M$.
Let $\cA(E)^{u,0}$ denote the space of unitary connections on $E\to M$ with trivial determinant (see \ref{may0602}). This is  an affine space over $\Omega^{1}(M,\su(E))$, (compare \ref{may0601}), where
$\su(E)$ denotes the anti-hermitean trace-free bundle endomorphisms.
Observe that if $M$ is closed oriented of dimension $n-1$, then we get a function $\cs_{\omega}:\cA(E)^{u,0}\to \R/\Z$, the Chern-Simons functional. 
\begin{prob}
\begin{itemize}
\item
Calculate $d\cs_{\omega}$ in general.
\item Assume that $n=3$ and
characterize the critical points of the Chern-Simons functional
$\cs_{c_{2}}$.
\end{itemize}\end{prob}
\proof Notice that $d\cs_{\omega}(\nabla)$ is a linear map
$\Omega^{1}(M,\su(E))\to \R$.
\hB 
 Let $P\to B$ be a $U(1)$-principal bundle, $E\to B$ the associated line bundle, and $\nabla$ be a connection on $E$. For  $\kappa\in \Z$ we let $p:M\to B$ be the unit sphere bundle of $$E^{k}:=\underbrace{E\otimes \dots \otimes E}_{k\times}\ .$$ The bundle $p^{*}E^{k}\to M$
is canonically trivialized. Let $q:Z\to B$ be the unit disc bundle and $\tilde \nabla$ be any connection on $q^{*}E^{k}$ which extends the trivial connection on $M$.
\begin{prob}\label{may0201}
Show that for all $n\ge 1$
$$\left[\int_{Z/B} c_1(\tilde \nabla)^{n}\right]= -\left[\kappa^{n} c_1^{n-1}(\nabla)\right] \in \Omega^{2n-2}(B;\C)/\im(d) \ .$$
\end{prob}

 If $V\to M$ is a real vector bundle with connection $\nabla$, then we can form
the complex vector bundle $V\otimes_{\R}\C$ with induced connection $\nabla\otimes \C$.
\begin{ddd}
If $V\to M$ is a real vector bundle, then we define its Pontrjagin classes by
$$p_{k}(V):=(-1)^{k}c_{2k}(V\otimes \C)\ .$$ If $\nabla$ is a connection, then we define the Pontrjagin form
$p_{k}(\nabla):=(-1)^{k}c_{2k}(\nabla\otimes \C)$.
\end{ddd}
The Pontrjagin forms are integral characteristic forms for real vector bundles.
 
 It is known that   $\pi_{3}(\MSO)=0$. Given a compact oriented Riemannian three manifold  $M$
 we thus can find a zero bordism
$Z$.  We choose a Riemannian metric on $Z$ extending the metric on $M$ with a product structure near $\partial Z=M$ and let $\nabla$ be  the Levi-Civita connection on $TZ$.
\begin{ddd}\label{apr2231} We define the Chern-Simons invariant of $(M,g)$ by
$$\CS(M,g):=[\int_{Z} p_{1}(\nabla)]\in \C/\Z .$$
\end{ddd}

\begin{lem}\label{aug1903}
$\CS(M,g)$ is well-defined independently of the choice of the Riemannian metric and the zero bordism.
\end{lem}
\proof
Let $Z^{\prime}$ be another zero bordism. We equip $Z^{\prime}$ with a metric extending $g$. Then we can obtain a compact Riemannian manifold $-Z\cup_{M} Z^{\prime}$. We have
$$\int_{Z^{\prime}} p_{1}(\nabla)-\int_{Z} p_{1}(\nabla)=\int_{-Z\cup_{M} Z^{\prime}}p_{1}(\nabla)\in \Z$$
by the integrality of $c_{2}$. 
Hence $\CS(M,g)$ only depends on the Riemannian manifold $M$. \hB
Note that one can improve the Chern-Simons invariant using the following additional integrality of $p_1$ of tangent bundles given by the signature theorem (see \cite{MR1215720}).
$$\langle p_1(TN),[N]\rangle =3\sign(N)\in 3\Z$$
for every oriented closed  four manifold $N$. Hence we can define, using the notation as above,
\begin{equation}\label{may0603}\CS_{refined}(M,g):=[\frac{1}{3} \int_Z p_1(\nabla)]\ .\end{equation}

\begin{prob}
Show that the assumption of a product structure is important here.
\end{prob}
\proof
The product structures makes sure that the restriction of the Levi-Civita
connection to the boundary induces the Levi-Civita connection of the boundary. Without a product structure
the corresponding difference is measured by the second fundamental form.
\hB 

\begin{prob}
We consider the lens space $L^{3}_k$.
Calculate $\CS(L^{3}_{k},g)$ and $\CS_{refined}(L^{3}_k,g)$, where $g$ is the metric induced from the round metric on $S^{3}$.
\end{prob}
The result is
$$\CS(L^{3}_{k},g)=[\frac{1}{k}]\ .$$
Here is a trick for the calculation.
The $U(1)$-principal bundle structure on
$p:L^{3}_k\to \C\P^{1}$ gives a decomposition $TL^{3}_k\cong p^{*}T\C\P^{1}\oplus H$, where $H
$ is trivialized by the $U(1)$-action. We let
$\nabla^{\prime}:=p^{*}\nabla^{T\C\P^{1}}\oplus \nabla^{H}$  be an adapted connection.
 We extend $\nabla^{\prime}$ to the disc bundle $q:D\to \C\P^{1}$ with boundary $L^{3}_k$ in the form $q^{*}\nabla^{T\C\P^{1}}\oplus \hat \nabla$, where $\hat \nabla$ extends $\nabla^{H}$.
 Then we have
\begin{equation}\label{aug1910}\CS(L^{3}_k,g)=[\int_{D} p_1(\nabla^{\prime})]+[\int_{L^{3}_k}\tilde p_1(\nabla^{g},\nabla^{\prime})]\ .\end{equation}
Note that $\hat \nabla$ is not trivial since the trivialization of $H$ does not extend to the disc bundle. We rather have $\int_{D } c_1^{2}(\hat\nabla)=-k^{2}\in \Z$ independently of the extension, see Problem \ref{may0201}.  Together with the addivity of the Pontrjagin form 
$$p_1(q^{*}\nabla^{T\C\P^{1}}\oplus \hat \nabla)=q^{*}c_1^{2}( \nabla^{T\C\P^{1}})+c_1^{2}(\hat \nabla^{H})$$ we get $[\int_{D} p_1(\nabla^{\prime})]=[-k^{2}]=0$.
In order to
calculate the second integral in \eqref{aug1910}
we go to the $k$-fold covering $S^{3}\to L^{3}_k$. 
The family of connections interpolating between $\nabla^{g}$ and $\nabla^{\prime}$  pulls back to a family interpolating between the Levi-Civita connection on $TS^{3}$ associated to the round metric and the lift of $\nabla^{\prime}$.
We can extend the  round metric from $S^{3}$ (considered as equator of $S^{4}$) as round metric to the upper half of $S^{4}$. Its  Pontrjagin form vanishes. On the other hand we can extend the lift of  $\nabla^{\prime}$ similarly as above preserving the decomposition. 
We get
$$0=\langle p_1(TS^{4}),[S^{4}]\rangle =k \int_{L^{3}_k}\tilde p_1(\nabla^{g},\nabla^{\prime}) - 1\ ,$$
hence
$$\int_{L^{3}_k}\tilde p_{1}(\nabla^{g},\nabla^{\prime})=\frac{1}{k}\ .$$
For the refined version we get
$$\CS_{refined}(L^{3}_k,g)=[\frac{1-k^{3}}{3k}]\ .$$
\hB

}
\end{ex} 
 
 \begin{ex}\label{apr2341}{\em
 For a Riemannian manifold $(M,g)$ with Levi-Civita connection $\nabla$
 the Pontrjagin form $p_{i}(g):=p_{i}(\nabla)$ only depends on the Weyl tensor (see \cite{MR0263098} for an argument). In particular, we have $p_{i}(g)=p_{i}(g^{\prime})$ if $g$ and $g^{\prime}$ are conformally equivalent.
 
 At the moment we can define, generalizing \ref{apr2231}, a Chern-Simon invariant
 $\CS(M,g)$ for $4n-1$-dimensional manifolds for which there exists a zero bordism $(Z,h)$ by
 $$\CS(M,g)=[\int_{Z} p_{n}(h)]\in \C/\Z\ .$$
 $\CS(M,g)$ is a conformal invariant of $M$. It is, for example, an obstruction against finding a locally conformally flat zero bordism. 
 \begin{kor}
$(L^{3}_{k},g)$ does not bound totally umbilically a locally  conformally flat manifold. 
\end{kor}
Below, in \ref{apr2340} we will drop the assumption that $M$ is zero bordant.
}\end{ex}

\subsection{Integral refinement}

Let $\omega$ be an integral characteristic form of degree $n$ and $E\to M$ be a vector bundle.
Then we have $$\omega(E)\in \im(\epsilon_\C:H^{n}(M;\Z)\to H^{n}_{dR}(M;\C))\ .$$
This suggests the following definition.

\begin{ddd}
An integral refinement $\omega^{\Z}$ of $\omega$ associates to each bundle $E\to M$
a class $\omega^{\Z}(E)\in H^{n}(M;\Z)$ with
$\epsilon_\C(\omega^{\Z}E))=\omega(E)$ such that for every pull-back diagram 
\begin{equation}\label{apr1910}\xymatrix{E^{\prime}\ar[d]\ar[r]^{F}&E\ar[d]\\M^{\prime}\ar[r]^{f}&M}\end{equation}
we have
$$f^{*}\omega^{\Z}(E)=\omega^{\Z}(E^{\prime})\ .$$
\end{ddd}

\begin{theorem}\label{apr2001}
An integral characteristic class for complex vector bundles has a unique integral refinement.
\end{theorem}
\proof
The functor which associates to each manifold $M$ the set 
of isomorphism classes of complex vector bundles on $M$ is represented  by the space
$$BU:=\bigsqcup_{n\ge 0} BU(n)\ .$$ The calculation of the integral cohomology of this space is  
a basic result in algebraic topology, see \eqref{apr1911}. It could be used to give a short conceptual proof of the theorem. For pedagogical reasons in order to present a basic technique
in the field of differential cohomolgy we give another argument which only  invests the following consequence:

\begin{kor}
The integral cohomology of $BU$ is concentrated in even degrees and torsion-free.
\end{kor}

We have a Bockstein sequence
$$\dots \to H^{k}(BU;\C)\to H^{k}(BU;\C/\Z)\to H^{k+1}(BU;\Z)\to H^{k+1}(BU;\C)\to \dots\ .$$
We see  that $H^{odd}(BU;\C/\Z)\cong 0$, $H^{odd}(BU;\C)\cong 0$ and  $H^{ev}(BU;\Z)\to H^{ev}(BU;\C)$ is injective.   

Recall that a map of spaces $f:X\to Y$ is called $n+1$-connected if it induces a bijection on $\pi_{0}$ and
the homotopy fibre of $f$ at each point of $Y$ has trivial homotopy groups up to degree $n$. Note that this implies (use Serre's spectral sequence) that the pull back
$f^{*}:H^{k}(Y;A)\stackrel{\sim}{\to} H^{k}(X,A)$ is an isomorphism for all $k\le n$ and every abelian group $A$.

We assume that $\omega$ is of degree $n$ and $M$ is $m$-dimensional. We define $r:=\max(n,m)+1$. If  $E\to M$ is a vector bundle,
then there exists a manifold $N$ with a vector bundle  $F\to N$ classified by  an $r$-connected map $u:N\to BU$
and a map $h:M\to N$ such that $E\cong h^{*}F$. For $N$ we can take, for example, an approximation of an $r$-skeleton of $BU$. 
 The map $u$ then induces an isomorphism in complex and integral cohomology in degree $\le r$.
If $n$ is odd, then $\omega(F)=0$ and hence $\omega(E)=h^{*}\omega(F)=0$.
If $n$ is odd, then we can take $\omega^{\Z}=0$.

From now on assume that $n$ is even.  Then there exists a unique class $w\in H^{n}(N;\Z)$ such that $\epsilon_\C(w)=\omega(F)$. We define $\omega^{\Z}(E):=h^{*}w\in H^{n}(M;\Z)$.
We must show that $\omega^{\Z}(E)$ is well-defined independently of the choices of $h$ and
$F\to N$.  

Assume that $F^{\prime}\to N^{\prime}$, $h^{\prime}:M\to N^{\prime}$ constitutes another choice. Then we produce a diagram
$$\xymatrix{&N\ar[dr]^{g}\ar[drr]^{u}&&\\M\ar[ur]^{h}\ar[dr]^{h^{\prime}}&&N^{\prime\prime}\ar[r]^{u^{\prime\prime}}&BU\\&N^{\prime}\ar[ur]^{g^{\prime}}\ar[urr]^{u^{\prime}}&&}$$
which commutes up to homotopy, where $u^{\prime\prime}$ is an $r$-connected map, as well.
We conclude that $g^{*}w^{\prime\prime}=w$ and $g^{\prime*}w^{\prime\prime}=w^{\prime}$ and therefore $h^{*}w=h^{\prime*}w^{\prime}$. 

This finishes the verification that $\omega^{\Z}$ is well-defined.
Finally we show that $\omega^{\Z}$ is natural.

Given a diagram \eqref{apr1910} we consider the extension
\begin{equation}\label{apr1910nnn}\xymatrix{E^{\prime}\ar[d]\ar[r]^{F}&E\ar[d]\ar[r]&F\ar[d]\\M^{\prime}\ar[r]^{f}&M\ar[r]^{h}&N}\ .\end{equation}
We can assume that $N\to BU$ is a $r^{\prime}$-connected map.
Then by construction
$\omega^{\Z}(E^{\prime})=f^{*}\omega^{\Z}(E)$.
\hB

If $u:N\to BU$ is $n+1$-connected, then we  have an isomorphism  $u^{*}:H^{n}(BU;\Z)\to H^{n}(N;\Z)$. Given an integral characteristic form $\omega$ of degree $n$ we thus get a  class $\omega^{\Z}\in H^{n}(BU;\Z)$ such that $u^{*}\omega^{\Z}=\omega^{\Z}(F)$.
The argument above shows that $\omega^{\Z}$ is well-defined. The class
$\omega^{\Z}$ is called the universal characteristic class associated to the integral characteristic form $\omega$.

We apply this to the Chern forms $c_n$ and get integral classes $c_n^{\Z}\in H^{2n}(BU;\Z)$.
We have a decomposition
$$H^{*}(BU;\Z)=\prod_{m\ge 0} H^{*}(BU(m);\Z)\ .$$
We let
$$c_n=(c_{n,1},c_{n,2},\dots,)$$
be the corresponding decomposition of the Chern classes.
Then it is known that the integral cohomology of the classifying spaces are given as polynomial rings by 
\begin{equation}\label{apr1911}H^{*}(BU(m);\Z)\cong \Z[c_{1,m},\dots,c_{n,m}]\ .\end{equation}

\begin{prob}
Use the Gysin sequence of the sphere bundles $S^{2n-1}\to BU(n-1)\to BU(n)$ in order to show
\eqref{apr1911} by induction w.r.t. $n$.
\end{prob}
In general the integral refinement of an integral characteristic form contains more information than the associated complex valued cohomology class.
\begin{ex}{\em
We consider a $U(1)$ principal bundle
$\pi:M\to \C\P^{n}$ with first Chern class $pc\in H^{2}(\C\P^{n};\Z)$ for a prime $p\in \nat$.
\begin{prob}
Show that  $\pi_1(M)\cong \Z/p\Z$.   
\end{prob} 
We have a non-canonical isomorphism
$J(M)\cong \Z/p\Z$. Consider the Poincar\'e bundle $P\to J(M)\times M$ (see Example \ref{apr1204}). For $\rho\in J(M)$
we consider the restriction $P_\rho:=P_{| \{\rho\}\times M}\to M$.

\begin{prob}
Show that
$\rho\mapsto c_1^{\Z}(P_\rho)$
provides an isomorphism of abelian groups
$$\Z/p\Z\cong J(M)\stackrel{\sim}{\to} H^{2}(M;\Z)\ .$$
\end{prob}
\proof
We identify $M$ with the unit sphere bundle  of $L^{p}$.
Then we observe that $\pi^{*}L^{p}$ is canonicaly trivialized.
Hence we can define a unique flat connection $\nabla^{\pi^{*}L}$
which induces the trivial connection on $\pi^{*}L^{p}$. This connection is characterized by
$\nabla^{\pi^{*}L}\phi=0$ for local sections $\phi$ such that $\nabla^{triv}\phi^{\otimes p}=0$.
We further observe $(\pi^{*}L,\nabla^{\pi^{*}L})=P_{\rho_{0}}$
for some generator $\rho_{0}\in J(M)$. 
We further observe that $c^{\Z}_{1}(\pi^{*}L)=\pi^{*}c^{\Z}_{1}(L)$ generates $H^{2}(M;\Z)$.
We get the isomorphism
$$J(M)\ni x=\rho_{0}^{k}\mapsto c_{1}^{\Z}(\pi^{*}L^{k})=c_{1}^{\Z}(P_{\rho_{0}^{k}}) \in H^{2}(M;\Z)\ .$$ 
\hB

Note that $\epsilon_{\C}(c_1^{\Z}(P_\rho))=0$.

 Let $E\to M$ be a vector bundle with connection.
Then we can form its determinant $\det(E)\to M$ (a line bundle) which has an induced connection $\nabla^{\det(E)}$.

\begin{prob}\label{may0602}
Show that
$c_1(\nabla)=c_1(\nabla^{\det(E)})$ and
$c_1^{\Z}(E)=c_1^{\Z}(\det(E))$.
\end{prob}
The integral total Chern class is defined by 
$$c^{\Z}(E)=1+c^{\Z}_{1}(E) +c_{2}^{\Z}(E)+\dots\ .$$
Recall the Newton classes from  Lemma \ref{may0204}.

\begin{prob}
Show that  
$$s_{4}^{\Z}(E)=c_1^{\Z}(E)^{2}-2c_2^{\Z}(E)\ .$$
Show that for two complex vector bundles $E,F$ we have the relation 
$$c^{\Z}(E\oplus F)=c^{\Z}(E)\cup c^{\Z}(F)$$
of total Chern classes. 
\end{prob} 

}\end{ex}

\begin{ex}{\em
For a real vector bundle $V\to M$ we define the integral Pontrjagin classes by
$$p_{i}^{\Z}(V):=(-1)^{i}c^{\Z}_{2i}(V\otimes \C)\in H^{4i}(M;\Z)\ .$$ The integral total Pontrjagin class is defined by $$p^{\Z}(V):=1+p_{1}^{\Z}(V)+p_{2}^{\Z}(V)+\dots\ .$$
\begin{prob} Show for 
a real vector bundle $V\to M$ we have $2c^{\Z}_{i}(V\otimes \C)=0$ for all odd $i\ge 1$.
Show further that for two real bundles $V,W$ we have
$$p^{\Z}(V\oplus W)=p^{\Z}(V)\cup p^{\Z}(W)+\mbox{$2$-torsion}\ .$$ 
\end{prob}
See \cite[Ch. 15]{MR0440554} for more information.
\begin{prob}
If $E_{\R}$ is the underlying real vector bundle of a complex vector bundle, then
$c^{\Z}_{j}(E_{\R}\otimes \C)=0$ for all odd $j\ge 1$. 
Conclude that
$$p^{\Z}(E_\R\oplus F_\R)=p^{\Z}(E_\R)\cup p^{\Z}(F_\R)$$ and
calculate $p(TP\C^{n}_{\R})$ explicitly.
\end{prob}

\begin{prob}
Calculate
$c^{\Z}_{1}(L\otimes \C)$, where $L\to \R\P^{\infty}$ is the universal real line bundle.
\end{prob}
 
\begin{prob}
We let $M_{r}\to S^{4}$ be the $SU(2)$-principal bundle classified by 
$c^{\Z}_2(M_{r})=r\ori_{S^{4}}\in H^{4}(S^{4};\Z)$, $r\in \nat$.
Calculate $H^{4}(M_{r};\Z)$ and $p_{1}^{\Z}(TM_{r})$.
 \end{prob}
}
\end{ex}

\section{Smooth Deligne cohomology}

\subsection{Recollections on sheaf theory}

We assume basic knowledge of sheaf theory.
We consider the category of smooth manifolds $\Mf$ as a site with the topology given by open coverings. It comes with a structure sheaf $C^{\infty}$ which associates to every manifold $M$ its algebra of complex-valued smooth functions. By $\Ab$   we denote the category of abelian groups. 
It is an abelian category. We let $\Ch(\cA)$ denote  the category of chain complexes of objects of an abelian category $\cA$.

 An abelian Lie group $A$ represents a sheaf of abelian groups
$$\underline{A}\in \Sh_{\Ab}(\Mf)\ .$$ If $A$ is not discrete, then $A^{\delta}$ denotes $A$ with the discrete topology.
We have a natural map of sheaves
$\underline{A^{\delta}}\to \underline{A}$ which is far from being an isomorphism.

\begin{prob}
Show that an exact sequence 
$$0\to A\to B\to C\to 0$$
of abelian Lie groups induces an exact sequence of sheaves
\begin{equation}\label{aug1920}0\to \underline{A}\to \underline{B}\to \underline{C}\to 0\ .\end{equation}
\end{prob} 
\proof
The main point is to show that 
$ \underline{B}\to \underline{C}$ is surjective. Use that
$B\to C$ has local sections.
\hB 
Note that it is important to consider the sequence \eqref{aug1920} in sheaves. 
As a sequence of presheaves it is not exact in general.
\begin{prob}
Discuss this assertion in the example of the exponential sequence 
$$0\to \Z\to \R\to \R/\Z\to 0$$
\end{prob}

The difference of the notions of exactness in the categories of sheaves and presheaves is the source of  sheaf cohomology. In particular we let $H^{*}(M;\underline{A})$
denote the cohomology of $M$ with coefficients in $A$, where $H^{*}(M;\dots)$
  is the higher derived image of the functor $$\Gamma(M,\dots):\Sh_{\Ab}(\Mf)\to \Ab$$
  of evaluation at $M$. These cohomology groups can be calculated using $\Gamma(M,\dots)$-acyclic resolutions.
  
  If $A$ is an abelian group, then we can define the cohomology of $M$ with coefficients in $A$ using simplicial or homotopy theoretic means, e.g. by
  $$H^{n}(M;A):=[M,K(A,n)]$$ or
  $$H^{n}(M;A):=H^{n}(\Hom(C_{*}(M),A))\ ,$$
  where
 $K(A,n)$ is the  $n$'th Eilenberg-MacLane space of $A$ and $C_{*}(M)$ is the simplicial chain complex of $M$. We shall use the fact that there are canonical isomorphisms between these definitions and
 $$H^{n}(M;A):= H^{n}(M;\underline{A}^{\delta})\ .$$
  
 \begin{prob}\label{may0610}
 Show that the three versions of the cohomology of $M$ with coefficients in a discrete abelian group $A$ mentioned above
are canonically isomorphic.
\end{prob}

\begin{prob}
Let $\cF$ be a sheaf on $\Mf$ and consider the presheaf
$\cH^{k}(\cF)$ which associates to $M\in \Mf$ the cohomology group
$H^{k}(M,\cF)$. Let $\widetilde{\cH^{k}(\cF)}$ be its sheafification and show that
$$\widetilde{\cH^{k}(\cF)}\cong \left\{\begin{array}{cc}
\cF&k=0\\
0&k\ge 0
\end{array}\right. \ .
$$
\end{prob}

 For applications we must know examples of $\Gamma(M,\dots)$-acyclic sheaves. 
\begin{lem}\label{apr2104}
A sheaf $\cF$ of   $C^{\infty}$-modules is $\Gamma(M,\dots)$-acyclic.
\end{lem}
\proof
The existence of smooth partitions of unity shows that $\cF$ is fine.
A fine sheaf is $\Gamma(M;\dots)$-acyclic. In the following we give the details of the argument.
 
 To an open covering $\cU$ of $M$ we associate the \v{C}ech complex $\check{C}(\cU;\cF)$.
We let
$\check{C}(M;\cF)$ be the colimit of  \v{C}ech complexes over a cofinal system of open coverings.
For every $n\ge 0$ the functor $\Sh_{\Ab}\to \Ab$, $\cF\to  \check{C}^{n}(M;\cF)$, is exact and
$H^{0}(\check{C}^{n}(M;\cF))=\Gamma(M,\cF)$. We thus get a $\delta$-functor of  \v{C}ech-cohomology
$(H^{*}(\check{C}^{n}(M;\dots)),\delta)$ extending $\Gamma(M,\dots)$. Uniqueness of such $\delta$-functor extensions implies that
$H^{*}(\check{C}^{n}(M;\cF))\cong H^{*}(M;\cF)$. Hence we can calculate
$H^{*}(M;\cF)$ as  \v{C}ech cohomology.

For the Lemma it thus suffices to see that the  \v{C}ech complex of a sheaf $\cF$ of $C^{\infty}$-modules is exact. 
We use the $C^{\infty}$-module structure through the following.
Let $V\subset U$ be an inclusion of open subsets, $x\in \cF(V)$, 
and $\chi\in C_{c}^{\infty}(V)$. The sheaf axioms imply that there exists an extension by zero $(\chi x)_{0}\in \cF(U)$ of $\chi x$.
\begin{prob}Prove this!
\end{prob}
 Let $\cU=(U_{i})_{i\in I}$ be a covering and $(\chi_{i})_{\in I}$ be a partition of unity. 
Then we define a homotopy
$h:\check{C}^{*}(\cU;\cF)\to \check{C}^{*-1}(\cU;\cF)$ by
$$h(x)_{j_{0},\dots,j_{n-1}}=(-1)^{n}\sum_{i\in I} (\chi_{i} x_{j_{0},\dots,j_{n-1},i})_{0}\ .$$
\begin{prob}
Check that $\partial h+h\partial=\id$.
\end{prob}
It follows that
$H^{\ge 1}(\check{C}^{*}(\cU;\cF))=0$ for every covering $\cU$, and thus
$H^{\ge 1}(\check{C}(M;\cF))=0$ by the exactness of the colimit.
\hB

  \begin{ex}{\em
 The complex valued de Rham complex can be interpreted as a sheaf of  complexes $\Omega_{\C}\in \Ch(\Sh_{\Ab}(\Mf))$. The sheaves $\Omega_{\C}^{n}$ are sheaves of $C^{\infty}$-modules an therefore $\Gamma(M;\dots)$-acyclic by Lemma \ref{apr2104}. By the Poincar\'e Lemma
 $\Omega_{\C}$ resolves the sheaf $\underline{\C}^{\delta}$. Hence we have a canonical isomorphism $$H_{dR}(M;\C)\cong H^{*}(M;\underline{\C}^{\delta})\ .$$

\begin{prob}
Show that for all $n>0$ we have
$H^{n}(M; \underline{\R} )=0$.
\end{prob}

\begin{prob}
Calculate and compare
$H^{*}(T^{n};\underline{U(1)})$ and $H^{*}(T^{n};\underline{U(1)}^{\delta})$.
\end{prob}

\begin{prob}
Let
$\Omega^{k}_{cl,\C}$ be the sheaf of closed complex
$k$-forms. Calculate
$H^{*}(M,\Omega^{k}_{cl,\C})$.
\end{prob}
}
 \end{ex}

Up to this point we have discussed the cohomology of sheaves. We now introduce the notion of cohomology for complexes of sheaves called hypercohomology.

\begin{ddd}
If $\cF\in \Ch^+(\Sh_{\Ab}(\Mf))$ is a lower bounded complex of sheaves of abelian groups, then by
$H^{*}(M;\cF)$ we mean the hypercohomology of $\cF$.
Let $\cF\to \cC$ be a quasi-isomorphism of lower bounded  complexes of sheaves such that $\cC$ consists of injective sheaves, then by definition
$$H^{*}(M;\cF):=H^{*}(\cC(M))\ .$$
\end{ddd}

\begin{prob}
Show that the hyper cohomology of $\cF$ is well-defined (independent of the resolution), and that it suffices to resolve by  complexes  of $\Gamma(M;\dots)$-acyclic sheaves.
\end{prob}

\begin{prob}
Show that the notion of hypercohomology extends the definition of cohomology for sheaves
if we consider a sheaf as a complex of sheaves concentrated in degree $0$.
\end{prob}

\begin{prob}
Calculate the hypercohomology groups of the complexes  
$$\underline{\Z}\to \underline{\C^{\delta}}\ , \quad  \underline{\Z}\to \underline{\C}$$
in terms of homotopy theoretic data.
\end{prob}

From now on we will omit the prefix ``hyper''.

\begin{prob}\label{may0615}
For $p\ge 1$ consider the complex of sheaves on $\Mf$
$$\cK(p)\quad  :\quad 0\to \underline{\C^{*}}\stackrel{d\log}{\to} \Omega^{1}_{\C}\to \dots\to \Omega^{p}_{\C}\to 0\ .$$
Calculate $H^{*}(M;\cK(p))$  in terms of homotopy theoretic and differential geometric data of $M$.
\end{prob}

Let $$\cC:C^{0}\to C^{1}\to C^{2}\to \dots $$ be a complex  of sheaves (or more general a complex in some abelian category).
Then we define an increasing  filtration by subcomplexes 
\begin{eqnarray*}F^{0}\cC&:& Z^{0}\\
F^{1}\cC&:&C^{0}\to Z^{1}\\
&\dots&\\
F^{k}\cC&:&C^{0}\to \dots \to C^{k-1}\to Z^{k}\end{eqnarray*}
We have quasi-isomorphisms
$$F^{i}\cC/F^{i-1}\cC\cong H^{i}(\cC)[-i]\ .$$
\begin{ddd}
The  spectral sequence associated  to this filtration
$$(E_r,d_r)\Rightarrow H^{*}(M,\cC)$$
with $$E_2^{p,q}=H^{p}(M,H^{q}(\cC))$$ is called the hypercohomology spectral sequence.
\end{ddd}
\begin{prob}
Discuss the hypercohomology spectral sequence for the complexes $\cK(p)$ introduced in \ref{may0615}.
\end{prob}

\subsection{Deligne cohomology}

A map of cochain complexes  $f:C\to D$ of abelian groups (or more general, of objects in an abelian category like sheaves) can be extended to an exact triangle
$$\dots \Cone(f)[-1]\to C\to D\to \Cone(f)\to\dots\ .$$
Here
$$\Cone(f)^{i}=C^{i+1}\oplus D^{i}\ , \partial(x\oplus y)=(-\partial x\oplus (\partial y-f(x))\ .$$
The triangle  induces a long exact sequence in cohomology
$$\to H^{n-1}(\Cone(f))\to H^{n}(C)\to H^{n}(D)\to H^{n}(\Cone(f))\to \ .$$

We let $\sigma^{\ge n}:\Ch\to \Ch$ be the stupid truncation functor given by 
$$(\sigma^{\ge n}C)^{k}:=\left\{\begin{array}{cc}C^{k}&k\ge n\\
0&k<n
\end{array}
\right\}\ .$$
Note that
$$H^{k}(\sigma^{\ge n} C)=\left\{\begin{array}{cc} 0&k<n\\ Z^{n}(C) &k=n\\
H^{k}(C)&k>n\end{array}\right\} .$$

In a similar manner we define
$$(\sigma^{<n}C)^{k}:=\left\{\begin{array}{cc}C^{k}&k< n\\
0&k\ge n
\end{array}
\right\}\ .$$
Note that
$$H^{k}(\sigma^{< n} C)=\left\{\begin{array}{cc} 0&k\ge n\\ C^{n-1}/\im(d) &k=n-1\\
H^{k}(C)&k<n-1\end{array}\right\} .$$

All this applies  to complexes of sheaves.

We now come to Deligne cohomology.
We consider the natural map
$$\underline{\Z}\oplus \sigma^{\ge n}\Omega_{\C}\to \Omega_{\C} ,\quad (z,x)\mapsto z-x$$
 using the not written inclusions $\underline{\Z}\to \underline{\C}$ and $\sigma^{\ge n}\Omega_{\C}\to \Omega_{\C}$.
\begin{ddd}\label{may0201nnn}
We define the $n$'th Deligne complex $\cD(n)\in \Ch(\Sh_{\Ab}(\Mf))$ by
$$\cD(n):=\Cone\left(\underline{\Z}\oplus \sigma^{\ge n}\Omega_{\C}\to \Omega_{\C} \right)[-1]\ .$$
 We define the $n$'th  Deligne cohomology by 
$$\hat H_{Del}^{n}(M;\Z):=H^n(M;\cD(n))\ .$$
\end{ddd}
Note that we take different complexes $\cD(n)$ for different $n$.
The Deligne complex  fits into
an exact triangle
$$\dots \to \cD(n) \to \underline{\Z}\oplus \sigma^{\ge n}\Omega_{\C}\to \Omega_{\C}\to \cD(n)[1]\to\dots\ .$$
The interesting piece of the associated long exact sequence is
$$\to H^{n-1}(M;\Z)\to H^{n-1}_{dR}(M;\C )\to \hat H^{n}_{Del}(M;\Z)\stackrel{(I,R)}{\longrightarrow} H^{n}(M,\Z)\oplus \Omega^{n}_{cl}(M;\C)\to H_{dR}^{n}(M;\C)\to  $$

\begin{prob}\label{apr2211nnn} Show that $\hat H_{Del}^{n}(M;\Z)$ naturally fits 
into an exact sequence
\begin{equation}\label{apr2211}H^{n-1}(M;\Z)\to H_{dR}^{n-1}(M; \C) \stackrel{a}{\to} \hat H_{Del}^{n}(M;\Z)\stackrel{(I,R)}{\to}  H^{n}(M;\Z)\times_{H^{n}(M;\C)}\Omega^{n}_{cl}(M;\C)\to 0\ .\end{equation} 
Calculate
$H^{i}(M;\cD(n))$ for all $i\not=n$.
The result is
$$\hat H^{i}(M;\cD(n))\cong \left\{\begin{array}{cc}
H^{i-1}(M;\C/\Z )&i<n\\
H^{i}(M;\Z)&i>n\end{array}\right. \ .$$
Show that
$$H^{n-1}(M;\C/\Z )\cong \ker(R:\hat H^{n}_{Del}(M;\Z)\to \Omega_{cl}^{n}(M))\ ,$$
and that the composition
$$  H^{n-1}(M;\C/\Z )\to  \hat H^{n}_{Del}(M;\Z)\stackrel{I}{\to} H^{n}(M,\Z)$$
is the  negative of the Bockstein operator associated to the sequence of coefficients  $$0\to \Z\to \C\to \C/\Z\to 0\ .$$
\end{prob}
\proof
 We have the quasi-isomorphisms $$\Cone(\underline{\Z}\to \Omega_{\C})[-1]\cong \Cone(\underline{\Z}\to \underline{\C^{\delta}})[-1]\cong \underline{\C/\Z}^{\delta}[-1]\ .$$
We consider the short exact sequence of complexes of sheaves
\begin{equation}\label{may1702}0\to \Cone(\underline{\Z}\to \Omega_{\C})[-1]\to \cD(n)\to \Cone(\sigma^{\ge n}\Omega_{\C}\to 0)[-1]\to 0
\end{equation}
which is induced by the natural inclusion and projection.
Since $$ \Cone(\sigma^{\ge n}\Omega_{\C}\to 0)[-1]\cong \sigma^{\ge n}\Omega_{\C}$$ and
$$H^{i}(M;\sigma^{\ge n}\Omega_{\C})=0 \:\:\forall i<n\ , \quad H^{n}(M; \sigma^{\ge n}\Omega_{\C})=\Omega^{n}_{cl}(M,\C)$$ we get from the long exact cohomology sequence associated to \eqref{may1702}
that
$$H^{i-1}(M;\C/\Z)\cong H^{i}(M;\cD(n))$$
for all $i<n$ and 
\begin{equation}\label{may1701}0\to H^{n-1}(M;\C/\Z)\to \hat H^{n}_{Del}(M;\Z)\stackrel{R}{\to}  \Omega^{n}_{cl}(M,\C)
 \ .\end{equation}
For $i>n$ we obtain from the long exact sequence for the cone $\cD(n)$ that
$$I:\hat H^{i}(M;\cD(n))\to H^{i}(M;\Z)$$ is
an isomorphism. Finally, the assertion for $i=n$  follows  from exact sequence of the cone $\cD(n)$, too.

The long exact sequence of the cone $\Cone(\underline{\Z}\to \Omega_{\C})[-1]$
is the Bockstein sequence   shifted down by $-1$. The first map in \eqref{may1702} induces a map from the Bockstein sequence to the long exact sequence of the cone $\cD(n)$.
In particular we get
$$\xymatrix{H^{n-1}(M;\Z)\ar@{=}[d]\ar[r]& H^{n-1}(M;\C )\ar[r]\ar@{=}[d]&H^{n-1}(M;\C/\Z) \ar[r]^{-\beta}\ar[d]& H^{n}(M;\Z)\ar[r]\ar[d]^{\id\oplus 0}& H^{n}(M;\C )\ar@{=}[d]\\ 
H^{n-1}(M;\Z)\ar[r]&H^{n-1}(M;\C )\ar[r]&\hat H^{n}_{Del}(M;\Z)\ar[r]^{I\oplus R}&H^{n}(M;\Z)\oplus \Omega_{cl}^{n}(M)\ar[r]& H^{n}(M;\C )}\ ,
$$
where $\beta$ is the Bockstein operator. The minus sign comes from the shift.
 \hB 

It turns out to be useful to have different representations of Deligne cohomology. 
We define
 \begin{equation}\label{jun1301}\cE(n):=\Cone(\underline{\Z} \to \sigma^{<n} \Omega_{\C})[-1]\ ,\end{equation} i.e.  
$$\cE(n)\quad : \quad 0\to \underline{\Z}\to \Omega_{\C}^{0}\stackrel{-d}{\to} \dots \stackrel{-d}{\to} \Omega^{n-1}_{\C}\to 0$$
where $\underline{\Z}$ sits in degree $0$ and hence $\Omega^{i}_{\C}$ in degree $i+1$. This complex sits
 in  the short exact sequence of complexes of sheaves
$$0\to  \Cone(\sigma^{\ge n}\Omega_{\C}\stackrel{=}{\to} \sigma^{\ge n}\Omega_{\C})[-1]  \to   \cD(n) \to  \cE(n)\to 0\ .$$
The left map is the natural inclusion and the right map the projection.   Since the cone of an isomorphism is acyclic  the projection map is a quasi-isomorphism.  
We conclude that
\begin{equation}\label{may1703}\hat H^{n}_{Del}(M;\Z)\cong H^{n}(M;\cE(n))\ .\end{equation}
We now consider the exact sequence
$$0\to \sigma^{<n}\Omega_{\C}[-1]\to \cE(n)\to \underline{\Z}\to 0\ .$$
It induces a long exact sequence in cohomology. Its most interesting piece is
$$H^{n-1}(M;\Z)\to \Omega^{n-1}(M;\C)/\im(d)\stackrel{a}{\to} H^{n}(M,\cE(n))\stackrel{I}{\to} H^{n}(M;\Z)\to 0\ .$$
The notation $I$ is used since this map corresponds to $I$ under the isomorphism \eqref{may1703}.\begin{ddd}
We define 
$$a: \Omega^{n-1}(M;\C)/\im(d) \to \hat H^{n}_{Del}(M;\Z)$$
to be the map given by
$$\Omega^{n-1}(M;\C)/\im(d)\stackrel{a}{\to}  H^{n}(M,\cE(n))\stackrel{\eqref{may1703}}{\cong}  \hat H^{n}_{Del}(M;\Z)$$
\end{ddd}

\begin{kor}
We have an exact sequence 
\begin{equation}\label{may1710}H^{n-1}(M;\Z)\to \Omega^{n-1}(M;\C)/\im(d)\stackrel{a}{\to} \hat H_{Del}^{n}(M;\Z)\stackrel{I}{\to} H^{n}(M;\Z)\to 0\ .\end{equation}
\end{kor}


\begin{prob}
Show that 
$a: \Omega^{n-1}(M;\C)/\im(d)\to \hat H_{Del}^{n}(M;\Z)$ extends the map
$a:H^{n-1}(M;\C)\to  \hat H_{Del}^{n}(M;\Z)$ in \eqref{apr2211}
and that
$R\circ a=d$.
\end{prob}
\proof
The first assertion follows from the commutativity
$$\xymatrix{\Omega_{\C}[-1]\ar[r]\ar[d]&\cD(n)\ar[d]\\
 \sigma^{<n}\Omega_{\C}[-1]\ar[r]&\cE(n)}\ .$$
The upper horizontal map induces the map $a:H^{n-1}(M;\C)\to \hat H^{n}_{Del}(M;\Z)$, while the 
lower horizontal map induces $a: \Omega^{n-1}(M;\C)/\im(d)\to \hat H_{Del}^{n}(M;\Z)$.

For the second assertion we consider the following web of exact sequences
$$\xymatrix{&0\ar[d]&0\ar[d]&0\ar[d]&\\0\ar[r]&\Cone(0\to \sigma^{\ge n}\Omega_{\C})[-1]\ar[r]\ar[d]&\Cone(\Z\to \Omega_{\C})[-1]\ar[r]\ar[d]&\Cone(\Z\to \sigma^{<n}\Omega_{\C})[-1]\ar@{=}[d]\ar[r]&0\\0\ar[r]&\Cone(\sigma^{\ge n}\Omega_{\C}\to \sigma^{\ge n}\Omega_{\C})[-1]\ar[d]\ar[r]&\cD(n)\ar[r]^{\pi}\ar[d]^{R}&\cE(n)\ar[d]\ar[r]&0\\0\ar[r]&\Cone(  \sigma^{\ge n}\Omega_{\C}\to 0)[-1]\ar@{=}[r]\ar[d]&\Cone(  \sigma^{\ge n}\Omega_{\C}\to 0)[-1]\ar[d]\ar[r]&0\ar[d]\ar[r]&0\\&0&0&0&}\ .$$
Let $\delta_{v}$ and $\delta_{h}$ denote the boundary maps associated to the left vertical and the upper horizontal exact sequence.
Then we know from homological algebra (see Exercise \ref{may2301}) that
$$\delta_{v}\circ R=-\delta_{h}\circ \pi\ .$$
Let now $$\omega\in \Omega^{n-1}(M,\C)/\im(d)\cong  H^{n}(M,\sigma^{<n}\Omega_{\C}[-1])\ .$$
Then in $H^{n+1}(M,\Cone(0\to \sigma^{\ge n}\Omega_{\C})[-1])$ we have
$$[0\oplus R(\omega)]=\delta_{v}(R(\pi^{-1}(a(\omega))))=-\delta_{h}(0\oplus \omega)= [0\oplus d\omega]\ .$$ 
We read off that $R(a(\omega))=d\omega$.
\hB

Let us collect all this information in the following commuting diagram, called the differential cohomology diagram.
\begin{prop}
The Deligne cohomology fits into the differential cohomology diagram
$$\xymatrix{&\Omega^{n-1}(M;\C)/\im(d)\ar[dr]^{a}\ar[rr]^{d}&&\Omega^{n}_{cl}(M;\C)\ar[dr]&\\\textcolor{red}{H^{n-1}_{dR}(M;\C)}\ar[ur]\textcolor{red}{\ar[dr]}&&\hat H_{Del}^{n}(M;\Z)\ar[ur]^{R}\ar@{->>}[dr]^{I}&&\textcolor{red}{H^{n}_{dR}(M;\C)}\\&
\textcolor{red}{
H^{n-1}(M;\underline{\C/\Z}^{\delta})}\ar[rr]^{-Bockstein}\ar@{^{(}->}[ru]&&\textcolor{red}{H^{n}(M;\Z)}\ar[ur]&}
$$
where the diagonal compositions are exact and the part marked in red is a segment of the long exact Bockstein sequence.
\end{prop}

\begin{ex}\label{may2301}
{\em
Here is an exercise in homological algebra.
We consider a web of short exact complexes
$$\xymatrix{&0\ar[d]&0\ar[d]&0\ar[d]&\\
0\ar[r]&A\ar[r]\ar[d]&D\ar[r]\ar[d]&G\ar[r]\ar[d]\ar[r]&0\\
0\ar[r]&B\ar[r]\ar[d]&E\ar[r]\ar[d]&\ar[r]H\ar[d]&0\\
0\ar[r]&C\ar[r]\ar[d]&F\ar[r]\ar[d]&I\ar[r]\ar[d]&0\\
&0&0&0&}$$
of cochain complexes.
Let $e\in H^{n}(E)$ be such that its image in $H^{n}(I)$ vanishes.
Then we can lift the image of $e$ in $H^{n}(H)$ to a class of $H^{n}(G)$ and apply the boundary operator  $\delta_{uh}$ of the upper horizontal sequence to get a well-defined class
$$U(e)\in \frac{H^{n+1}(A)}{\delta_{uh}\delta_{rv}H^{n-1 }(I)}\ .$$ Similarly we can lift the image of $e$
in $H^{n}(F)$ to $H^{n}(C)$ and apply the boundary operator $\delta_{lv}$ of the left vertical sequence
to get a class $$V(e)\in \frac{H^{n+1}(A)}{\delta_{lv}\delta_{lh}H^{n-1}(I)}\ .$$
 \begin{prob}
Show that
$$[U(e)]=[-V(e)]\in \frac{H^{n+1}(A)}{\delta_{uh}\delta_{rv}H^{n-1}(I)+\delta_{lv}\delta_{lh}H^{n-1}(I)}\ .$$
\end{prob}}
\end{ex}

\begin{prob}
Relate the structure of $\hat H^{n}_{Del}(M;\Z)$ with the hypercohomology spectral sequence of the complex of sheaves $\cE(n)$.
\end{prob}

Deligne cohomology $\hat H_{del}^{n}(\dots;\Z)$ is a contravariant functor from $\Mf$ to
$\Ab$. It is not homotopy invariant, but its deviation from homotopy invariance is measured by the homotopy formula.
Let $i_{t}:M\to [0,1]\times M$ be the inclusion determined by $t\in [0,1]$.
\begin{prop}\label{apr2201}
If $\hat x\in \hat H^{n}_{Del}([0,1]\times M;\Z)$, then
$$i_{1}^{*}\hat x-i_{0}^{*}\hat x=a(\int_{[0,1]\times M/M} R(\hat x))\ $$
\end{prop}
\proof
By homotopy invariance of integral cohomology   we know that
there exists $y\in H^{n}_{dR}(M;\Z)$ such that $I(\hat x)=\pr_{M}^{*}y$.
We can choose a lift $\hat y\in \hat H^{n}_{Del}(M;\Z)$ such that $I(\hat y)=y$.
Then $I(\pr_{M}^{*}\hat y-\hat x)=0$ and therefore $\hat x=\pr_{M}^{*}\hat y+a(\omega)$ for some
$\omega\in \Omega^{n-1}([0,1]\times M;\C)$.  Note that $R(\hat x)=\pr_{M}^{*}R(\hat y)+d\omega$.
Using $\pr_{M}\circ i_{0}=\pr_{M}\circ i_{1}$, Stokes' theorem, and that $\int_{[0,1]\times M/M}\circ  \pr_{M}^{*}=0$ we get
$$i_{1}^{*}\hat x-i_{0}^{*}\hat x=a(i_{1}^{*} \omega-i_{0}^{*}\omega)=a(\int_{[0,1]\times M/M} d\omega)=a(\int_{[0,1]\times M/M} R(\hat x))\ .$$
\hB

Let $M=U\cup V$ be a decomposition into open submanifolds.
\begin{prob}
Show that there exists a Mayer-Vietotis sequence of the form
$$\hspace{-1.5cm}\dots\to H^{n-2}(U\cap V;\C/\Z)\to\hat H^{n}_{Del}(M;\Z)\to \hat H_{Del}(U;\Z)\oplus \hat H_{Del}(V;\Z)\to \hat H_{Del}(U\cap V;\Z)\to H^{n+1}(M;\Z)\to$$
which extends to the left and right by the Mayer-Vietoris sequences of
$H^{*}(\dots;\C/\Z)$ and $H^{*}(\dots;\Z)$.
\end{prob}
\proof
It is clear that sheaf cohomology $H^{*}(M;\cD(n))$ has a Mayer-Vietoris sequence.
The rest follows from the calculations in  Proposition \ref{apr2211nnn}. \hB

\begin{ex}{\em
We have canonical isomorphisms
$\hat H_{Del}^{0}(\pt;\Z)\cong \Z$ and $\ev:\hat H_{Del}^{1}(\pt,\Z)\cong \C/\Z$.
Let $M$ be any smooth manifold.
We define a map
$$\phi:\hat H^{1}_{Del}(M;\Z)\to \C/\Z^{M}$$
by
$\phi(x)(m)=\ev(m^{*}x)$, where
we consider $m$ as a map $m:\pt\to M$.
\begin{prob}\label{may0106}
Show that $\phi(x)\in C^{\infty}(M;\C/\Z)$ and
$ d\phi(x)= R(x)$.
Show that
$\phi:\hat H^{1}_{Del}(M;\Z)\to C^{\infty}(M,\C/\Z)$ is an isomorphism of groups.
\end{prob}
\proof
We use the Five Lemma and that
$C^{\infty}(M,\C/\Z)$ fits into an exact sequence
$$H^{0}(M;\Z)\to \Omega^{0}(M;\C) \to C^{\infty}(M,\C/\Z)\to H^{1}(M;\Z) \to 0$$
which is compatible with \eqref{may1710} by $\phi$. \hB 
 
\begin{prob}\label{apr2205}
Calculate $\hat H_{Del}^{*}(S^{1};\Z)$.
Show that there is a unique class $\hat e\in \hat H^{1}_{Del}(S^{1};\Z)$
with $R(\hat e)=dt$ (here $t$ is a coordinate given by $\R\to \R/\Z\cong S^{1}$) which vanishes after restriction to the base point $\{t=0\}\in S^{1}$.Calculate the restriction
$\hat e_{|\{t\}}$ for every  point $\{t\}\in S^{1}$. Show that the class $\hat e$ is primitive with respect to the 
group structure of $S^{1}$.
\end{prob}

\begin{prob}\label{apr2206}
Use \ref{apr2205} in
order to define a map $\ev:\hat H_{Del}^{2}(M;\Z)\to C^{\infty}(M^{S^{1}};\C/\Z)$ for every manifold $M$.
Here we understand smooth functions from the free loop space $M^{S^{1}}$  in the diffeological sense.
\end{prob}
 }
 
 \end{ex}
 
 \begin{ex}
 
{\em 
Let $G$ be a simple simply-connected compact Lie group. }
\begin{prob}
Show that there is a unqiue (up to sign)
biinvariant class $x\in \hat H^{3}_{Del}(G;\Z)$ such that $I(x)\in H^{3}(G;\Z)$
is a  generator. Show that this class is primitive.
Discuss the non-simply connected case.
\end{prob}

{\em 
We  consider a connected semisimple compact Lie group $G$ and a class
$x\in H^{3}(G;\Z)$. There is a unique biinvariant form $\lambda\in \Omega^{3}(G)$ 
which is (automatically) closed and represents the image of $x$ in real cohomology.
Since $H^{2}(G;\R)=0$ there exists a unique differential refinement
$\hat x\in \hat H^{3}_{Del}(M;\Z)$ such that $I(\hat x)=x$ and $R(\hat x)=\lambda$.
Let $T\subseteq G$ be the maximal torus. We identify
$H^{2}(T;\C/\Z)\cong \Hom(\Lambda^{2}\pi_{1}(T,1),\C/\Z)$.
}
\begin{prob}
Show that $\lambda_{|T}=0$ and calculate the class
$$\hat x_{|T}\in H^{2}(T;\C/\Z)\cong \Hom(\Lambda^{2}\pi_{1}(T,1),\C/\Z)\ .$$
\end{prob}
\proof
First observe that there exists a uniquely determined invariant bilinear form $B$ on $\gaaa$ such that
$\lambda_{|T_{1}G}(X,Y,Z)\cong B([X,Y],Z)$ for $X,Y,Z\in \gaaa$.
This implies $\lambda_{|T}=0$ since $\taaa$ is abelian.
Then calculate the class
$\hat x_{|T}\in H^{2}(T;\C/\Z)$ either by describing it as a homomorphism
$\Lambda^{2}\pi_{1}(T,1)\to\C/\Z$ or as an isomorphism class of central extensions \begin{equation}\label{feb2013-1}
0\to \C/\Z\to \widehat{\pi_{1}(T,1)}\to  \pi_{1}(T,1)\to 0\ .\end{equation}
\hB
I do not have an easy solution for the second part of  this exercise. If the class $x$ is determined by a central extension of the loop group of $K$, then the central extension \eqref{feb2013-1} is obtained by restriction if we view $\pi_{1}(T,1)$ as a subgroup of the loop group in the canonical manner.
\end{ex}

\begin{ex}\label{apr2230}{\em 
If $M$ is a $n-1$-dimensional connected closed oriented
 manifold, then we have a canonical identification
 $$\ev:\hat H^{n}_{Del}(M;\Z)\cong H^{n-1}_{dR}(M;\C)/\im(\epsilon_{\C})\cong \C/\Z\ .$$
 This follows from the exact sequence \eqref{apr2211}. Explicitly, the identification is
 given by
 $\hat x\mapsto [\int_{M}\omega]$,
 where $\omega\in \Omega^{n-1}(M;\C)$ is such that
 $a(\omega)=\hat x$.
 Indeed, such form exists  by \eqref{may1710} since $I(\hat x)=0$   for dimensional reasons, and $\omega$ is well-defined up to integral forms.  Therefore the class in $\C/\Z$ of the integral does not depend on the choice.
 }
 \end{ex}

Recall the complexes introduced in \ref{may0615}. 
 We consider the map
\begin{equation}\label{may0611}\cE(p)\to \cK(p-1)[-1]\end{equation}
given by
$$\begin{array}{cc}0&\mbox{degree zero}\\
2\pi i \id_{\Omega_{\C^{q}}}&q\ge 1\\
\Omega^{0}_{\C}(M) \in f\mapsto\exp(2\pi i f)\in \underline{\C^{*}}(M)&\mbox{degree one}
 \end{array}$$
 \begin{prob}\label{may0616}
 Show that \eqref{may0611} is a quasi-isomorphism so that
 $$\hat H_{Del}^{p}(M;\Z)\cong H^{p-1}(M;\cK(p-1))\ .$$
\end{prob}

\subsection{Differential refinements of integral characteristic classes}

Let $\omega$ be an integral characteristic form for complex vector bundles of degree $n$. It has a unique integral refinement $\omega^{\Z}$ by Theorem \ref{apr2001}.

\begin{ddd}
A differential refinement of $\omega$ associates to every vector bundle with connection $(E,\nabla)$ on $M$ a class $\hat \omega(\nabla)\in \hat H_{Del}(M;\Z)$ such that
$$R(\hat \omega(\nabla))=\omega(\nabla)\ ,\quad I(\omega( \nabla))=\omega^{\Z}(E)$$
and for every map $f:M^{\prime}\to M$ we have
$f^{*}\hat \omega( \nabla)=\hat \omega(f^{*}\nabla)$.
\end{ddd}
A longer, but in some cases clearer notation for evaluation of the differential refinement on $(E,\nabla)$ would be $\hat \omega(E,\nabla)$.
Note that $$\ker(R)\cap \ker(I)=H^{n-1}_{dR}(M;\C)/\im(\epsilon_\C)$$ is non-trivial. Therefore the differential refinement $\hat \omega(\nabla)$ can potentially contain finer information than the pair of the form $\omega(\nabla)$ and $\omega^{\Z}(E)$.

The homotopy formula  Proposition \ref{apr2201} determines how $\hat \omega(\nabla)$ depends on the connection $\nabla$.

\begin{lem}\label{apr2202}
If $\nabla$ and $\nabla^{\prime}$ are two connections on the same bundle, then we 
  have
$$\hat \omega(\nabla^{\prime}) -\hat \omega(\nabla)=a(\tilde \omega(\nabla^{\prime},\nabla))\ .$$
\end{lem}
\begin{prob}
Prove this Lemma.
\end{prob}

\begin{theorem}\label{may2007}
An integral characteristic form of degree $n$ admits a unique differential refinement.
\end{theorem}
\proof
We can assume that $n$ is even.
Let $\omega$ be an integral characteristic form and $\omega^{\Z}$ be its unique
integral refinement by Theorem \ref{apr2001}.
Let $E\to M$ be a vector bundle with connection $\nabla$.
Assume that $n=\dim(M)$. Then we choose an $n+1$-connected map 
$u:N\to BU$ classifying $F\to N$ and a connection $\nabla^{F}$, a map $f:M\to N$ and an isomorphism $E\cong f^{*}F$. Since $H^{n-1}_{dR}(N;\C)=0$ we are forced to define
$\hat \omega(\nabla^{F})\in \hat H^{n}_{Del}(M;\Z)$ uniquely such that $R(\hat \omega(\nabla^{F}))=\omega(\nabla^{F})$ and $I(\hat \omega(\nabla^{F}))=\omega^{\Z}(F)$. By naturality we are forced to define
$\hat \omega(f^{*}\nabla^{F})=f^{*}\hat \omega(\nabla^{F})$. By Lemma \ref{apr2202}
we are forced to define
$$\hat \omega(\nabla)=f^{*}\hat \omega(\nabla^{F})+a(\tilde \omega(\nabla,f^{*}\nabla^{F}))\ .$$
This already shows the uniqueness clause.
It remains to show that $\hat\omega$ is well-defined and natural.
We argue as in the proof of Theorem \ref{apr2001} using the same notation.
We must show that
$$f^{*}\hat \omega(\nabla^{F})+a(\tilde \omega(\nabla,f^{*}\nabla^{F}))-
f^{\prime*}\hat \omega(\nabla^{F^{\prime}})-a(\tilde \omega(\nabla,f^{\prime *}\nabla^{F^{\prime}}))=0\ .$$
We have
$$\hat \omega(\nabla^{F})=g^{*}\hat \omega(\nabla^{F^{\prime\prime}})+a(\tilde\omega(\nabla^{F },g^{*}\nabla^{F^{\prime\prime}}))$$
and
$$\hat \omega(\nabla^{F^{\prime}})=g^{\prime *}\hat \omega(\nabla^{F^{\prime\prime}})+a(\tilde\omega(\nabla^{F^{\prime}},g^{\prime *}\nabla^{F^{\prime\prime}}))\ .$$
Hence we must see that
\begin{multline*}f^{*}g^{*}\hat \omega(\nabla^{F^{\prime\prime}})+f^{*}a(\tilde\omega(\nabla^{F },g^{*}\nabla^{F^{\prime\prime}}))-f^{\prime *}g^{\prime *}\hat \omega(\nabla^{F^{\prime\prime}})-f^{\prime*}a(\tilde\omega(\nabla^{F^{\prime}},g^{\prime *}\nabla^{F^{\prime\prime}}))\\
+a(\tilde \omega(\nabla,f^{*}\nabla^{F}))-a(\tilde \omega(\nabla,f^{\prime *}\nabla^{F^{\prime}}))=0 \end{multline*}
This reduces to 
$$f^{*}g^{*}\hat \omega(\nabla^{F^{\prime\prime}})-f^{\prime *}g^{\prime *}\hat \omega(\nabla^{F^{\prime\prime}})-\tilde \omega(f^{*}g^{*}\nabla^{F^{\prime\prime}}
,f^{\prime}g^{\prime}\nabla^{F^{\prime\prime}})=0$$
which holds true by the homotopy formula,  Proposition \ref{apr2201}, and the additivity \eqref{apr1241} of the transgression.

Naturality of $\hat \omega$ is now easy and left to the reader.
\hB

\begin{ex}{\em
In these examples we discuss some general properties of differential refinements of integral characteristic forms.

We consider the differential refinement $\hat \omega$ of an integral characteristic form $\omega$ of degree $n$. Let $Z$ be a compact oriented $n$-manifold with boundary $M$ and $E\to Z$ be a
complex vector bundle with connection $\nabla$.
\begin{prob}\label{may2130}
Show that
$$\ev(\hat \omega(\nabla_{|M}))=[\int_{Z} \omega(\nabla)]\in \C/\Z\ .$$
\end{prob}
\proof 
Since $H^{n}(Z;\Z)=0$ we
 have $\hat \omega(\nabla)=a(\alpha)$ for some form $\alpha\in \Omega^{n-1}(Z,\C)$.
 Note that $d\alpha=\omega(\nabla)$.
 It follows by Stokes theorem that
$$ \ev(\hat \omega(\nabla_{|M}))=[\int_{M}\alpha_{|M}]=[\int_{Z}d\alpha]=[\int_{Z}\omega(\nabla)]\ .$$
\hB

\begin{prob}\label{apr2203}
If $\omega$ has degree $n\ge 1$, $E$ is trivialized, and  $\nabla$ is the trivial connection, then
$\hat \omega(\nabla)=0$. 
\end{prob}
\proof
Indeed, we can assume that $n\ge 2$. A trivial bundle with trivial connection can be obtained as pull-back from a point. We now use that $\hat H_{Del}^{\ge 2}(\pt;\Z)=0$. 
\begin{prob}\label{jun1201}
Let $E\to M$ be a $k$-dimensional vector bundle with connection $\nabla$.
Show $\hat c_n(\nabla)=0$ for all $n>k$.
\end{prob}
\proof
One can show by adapting the proof of Theorem \ref{may2007} that an integral characteristic form for $k$-dimensional vector bundles has a unique differential refinement. Since for $n>k$ the restriction of $c_n$ to $k$-dimensional bundles vanishes, so does its differential refinement. 
\hB 
 Let $E\to M$ be a vector bundle, $\nabla$ a connection on $E$ and $h$ be a hermitean metric. \begin{prob}
 Show that $\overline{\hat  c_{n}(\nabla)}=\hat c_n(\nabla^{*})$.
\end{prob}
\proof
Verify this one the level of curvatures and then use uniqueness of differential extensions. \hB 

Let $E\to M$ be a vector bundle with connection $\nabla^{E}$ and $(\det(E),\nabla^{\det(E)})\in \Line_\nabla(M)$
be its determinant.
\begin{prob}
Show that $\hat c_1(\nabla^{E})=\hat c_1(\nabla^{\det(E)})$.
\end{prob}
\proof
Use the uniqueness of differential refinements. \hB

}
\end{ex}

\begin{ex}{\em

In the following examples 
we consider $\hat c_{1}$ and the classification of line bundles with connection using $\hat H^{2}_{Del}(\dots,\Z)$.

We have (see Exercise \ref{apr2205}) an isomorphism
$$\hat H_{Del}^{2}(S^{1};\Z)\cong H^{1}(S^{1}; \C/\Z )\cong \C/\Z\ .$$ 
We consider the connection $\nabla=d+\alpha$ on the trivial one-dimensional bundle on $S^{1}$
where $\alpha\in \Omega^{1}(S^{1};\C)$. Then we have
$$\hat c_{1}(\nabla)=[-\frac{1}{2\pi i}\int_{S^{1}}\alpha]\in \C/\Z\ .$$ In order to see this note
that have by \ref{apr2203} and \ref{apr2202}  
$$\hat c_{1}(\nabla)=a(\tilde c_{1}(\nabla,d))\ .$$
We now use \ref{apr1291}.
\begin{prob}
Conclude that for a line bundle with connection $(L,\nabla^{L})$ on a manifold $M$ we have
$$\ev(\hat c_{1}(\nabla))(\gamma)=\hol_{\nabla}(\gamma)\ , \quad \gamma\in M^{S^{1}} \ ,$$
where $\ev$ is the evaluation map found in \ref{apr2206} 
\end{prob}

Let $\Line_{\nabla}(M)$ be the group of isomorphism classes of line bundles with connection under the tensor product operation. 
\begin{prob}
Verify the existence of inverses in $\Line_{\nabla}$. Show that
$$\hat c_{1}:\Line_{\nabla}(M)\to \hat H^{2}_{Del}(M;\Z)$$
is a natural isomorphism.
\end{prob}
\proof
Note that $BU(1)$ classifies line bundles. We know that
$BU(1)\cong K(\Z,2)$ and $H^{*}(BU(1);\Z)\cong \Z[c_1^{\Z}]$. Here the universal first Chern  class $c_1^{\Z}$ is such that
if $L\to M$ is classified by $l:M\to BU(1)$, then $c_1^{\Z}(L)=l^{*}c_1^{\Z}$.
Every  class $x\in H^{2}(M;\Z)$ can be written  as $x=l^{*}c_1^{\Z}$ for some map $l:M\to BU(1)$ and therefore is the first Chern class of a line bundle.
Furthermore, the tensor product of line bundles induces the $h$-space structure on $BU(1)$ and $c_1^{\Z}$ is primitive. This has the effect that $c_1^{\Z}(L\otimes L^{\prime})=c_1^{\Z}(L)+c_1^{\Z}(L^{\prime})$.

Given $(L,\nabla^{L})\in \Line_{\nabla}(M)$ we choose a bundle $H$ such that
$L\otimes H^{-1}$ is trivializable, e.g. such that $c^{\Z}_{1}(H)=-c^{\Z}_{1}(L)$. 

 Let $\nabla^{H}$ be any connection on $H$. Then we can find an $\alpha\in \Omega^{1}(M;\C)$ such that
$c_{1}(\nabla^{L})+c_{1}(\nabla^{H})=d\alpha$. Then
$(H,\nabla^{H}+2\pi i \alpha)$ is the inverse of $(L,\nabla^{L})$.

Note that $\hat c_1(L,\nabla^{L})\in \hat H^{2}_{Del}(M;\Z)$ is characterized completely by
its holonomy function $\ev(\hat c_1(L,\nabla^{L}))\in C^{\infty}(M^{S^{1}};\C/\Z)$ (see Exercise \ref{apr2206}). 
Since the holonomy of $(L,\nabla^{L})\otimes (H, \nabla^{H})$ is the product of holonomies of the factors
we see that $\hat c_{1}$ is additive. 

Assume that $\hat c_{1}(L,\nabla^{L})=0$. Then
$(L,\nabla)$ has trivial holonomy along every path and thus can be trivialized (including the connection).

On the other hand, let
$x\in \hat H^{2}_{Del}(M;\Z)$ be given. Then we choose a line bundle
$(L,\nabla^{L})$ such that $c^{\Z}_{1}(L)=I(x)$. We can further adjust $\nabla$ such that
$R^{\nabla}=-2\pi i R(x)$.
Then
$x-\hat c_{1}(L,\nabla)=a(\alpha) $ for some
$\alpha\in \Omega_{cl}^{1}(M;\C)$.
We consider the connection $\nabla^{\prime}:=d-2\pi i\alpha$ on the trivial bundle $L^{triv}\to M$.
Then
$\hat c_{1}((L,\nabla)\otimes (L^{triv},\nabla^{\prime}))=x$.
\hB 
  
Let $\cU$ be a covering of $M$. Recall the definition \ref{may0615} of the complex of sheaves $\cK(1)$.
\begin{prob}
Show that
a \v{C}ech cocycle $c\in C^{1}(\cU,\cK(1))$
can naturally be identified with the glueing data for a 
line bundle with connection. Use this to construct the isomorphism 
$$\hat c_{1}:\Line_{\nabla}(M)\stackrel{\sim}{\to} H^{1}(M;\cK(1))$$
explicitly  on the level of \v{C}ech cohomology. Verify compatibility with the isomorphism \ref{may0616}.
\end{prob}
\proof
Note that
$$\cK(1): \underline{\C^{*}}\stackrel{d\log}{\to} \Omega^{1}\to 0\ .$$
Let $\cU:=(U_\alpha)$ be a covering.
Then a one-cocycle is given by
$x:=(\omega_\alpha,g_{\alpha,\beta})$, where $\omega_\alpha\in \Omega^{1}(U_\alpha,\C)$ and $g_{\alpha,\beta}\in C^{\infty}(U_\alpha\cap U_\beta;\C^{*})$.
The relation $\delta x=0$ is equivalent to
$\omega_\beta-\omega_\alpha=d\log g_{\beta,\alpha}$ on $U_\alpha\cap U_\beta$  and
$g_{\alpha,\beta}g_{\beta,\gamma}=g_{\alpha,\gamma}$ on the triple intersections
$U_\alpha\cap U_\beta\cap U_\gamma$.
This is exactly the cocycle condition for a line bundle locally trivialized by sections $s_\alpha$ on $U_\alpha$ such that $s_\alpha=g_{\alpha,\beta} s_\beta$ on $U_\alpha\cap U_\beta$ and a connection such that
$\nabla \log s_\alpha=\omega_\alpha$. \hB

This has a higher-degree analog. \v{C}ech cocycles for classes in
$H^{2}(M;\cK(2))\cong \hat H^{3}_{Del}(M;\Z)$ correspond to Hitchin's descent data for geometric gerbes with band $\underline{\C^{*}}$.
The group $\hat H_{Del}^{3}(M;\Z)$ classifies isomorphisms classes of geometric gerbes.
See \cite{MR1876068} and \cite{MR2362847} for more information.

 }\end{ex}

\begin{ex}\label{apr2501}{\em
In the following examples we generalize the Chern-Simons invariants from relative to absolute invariants.

Let $\nabla$ be a flat connection on a trivializable bundle $E\to M$.
Then $$\hat \omega(\nabla)\in \frac{ H^{n-1}_{dR}(M;\C)}{\im(\epsilon_{\C})}=\ker(R)\cap \ker(I)\ .$$
With this identification
$\hat \omega(\nabla)$ coincides with the Chern-Simons invariant of the flat connection $\nabla$ introduced in Definition \ref{apr2220}. 
Indeed, if $\nabla^{triv}$ is a trivial connection, then
$\hat \omega(\nabla)=a(\tilde \omega(\nabla,\nabla^{triv}))$
and $\tilde \omega(\nabla)=\tilde \omega(\nabla,\nabla^{triv})$. 
More generally, if $\nabla$ is just flat, then we get
$$\hat \omega(\nabla)\in H^{n-1}(M; \C/\Z )\ .$$
This is the generalization of the Chern-Simons invariant dropping the trivializability
condition.
 
 Let $\nabla$ on $E\to M$ have finite holonomy. 
\begin{prob}
Show that
$\hat \omega(\nabla)\in H^{n-1}(M, \C/\Z )$
is a torsion class (determine and fill in the missing details of the argument
below).
\end{prob}
\proof
There exists a reduction of structure groups of $E$  to a finite group $H$.
Let $f:M\to BH$ a classifying map so that $E\cong f^{*}F$ for some bundle $F\to BH$.
Since higher-degree rational cohomology of $BH$ is trivial we have an isomorphism
$$B:H^{n-1}(BH; \C/\Z )\cong H^{n}(BH;\Z)$$
from the Bockstein sequence.  We now observe that
$$\hat \omega(\nabla)=f^{*}B^{-1}(\omega(F))\ .$$\hB 
}
\end{ex}

\begin{ex}\label{apr2510}{\em
We generalize the theory developed in \ref{aug1940}.
Let $\omega$ be an integral characteristic form of degree $n$.
We consider a space $B$ with a flat bundle $V\to B$ (a complex vector bundle with structure group $GL(n,\C^{\delta}))$. If $f:M\to B$ is a map, then $f^{*}V\to M$ becomes a complex vector bundle with a flat connection $\nabla$. If $M$ is closed, oriented, $n-1$-dimensional, then we define
$$\cs^{V}_\omega(f):=\hat \omega(\nabla)\in \hat H^{n}_{Del}(M;\Z)\cong  \C/\Z\ .$$
\begin{prob}
Show that $\cs^{V}_\omega(f)$ only depends on the oriented bordism class of $f$ and that we get a homomorphism
$$\cs^{V}_\omega:\MSO_{n-1}(B)\to \C/\Z\ .$$
\end{prob}
\proof Use \ref{may2130}. \hB

 \begin{lem}\label{jun0101}

There exists a class
$u\in H^{n-1}(B;\C/\Z)\cong \Hom(H_{n-1}(B;\Z);\C/\Z)$
such that
$$\cs_\omega^{V}([f:M\to B])=u(f_*[M])\ .$$
In particular, $\cs_\omega^{V}$ factorizes over the natural transformation
$\kappa:\MSO_{n-1}(B)\to H_{n-1}(B;\Z)$.
\end{lem}
\proof 
Let $g:\tilde B\to B$ be an $n+1$-connected approximation of $B$ and $\tilde V:=g^{*}V$.
Then we have   isomorphisms $g_*:\MSO_{n-1}(\tilde B)\stackrel{\sim}{\to}  \MSO_{n-1}( B)$
and $g^{*}:H^{n-1}(B;\C/\Z)\stackrel{\sim}{\to}  H^{n-1}(B;\C/\Z)$, and
  $$\cs_{\omega}^{V}\circ g_*=\cs_{\omega}^{\tilde V}\ .$$

From now on we can assume that $B$ is a smooth manifold and $V$ has a flat connection $\nabla^{V}$. We have $\hat \omega(\nabla^{V})\in H^{n-1}(B;\C/\Z)\subseteq \hat H_{Del}(M;\Z)$ and
$$\cs_\omega^{V}(f)=\ev(f^{*}\hat \omega(\nabla^{V}))=\hat \omega(\nabla^{V})(f_*[M])\ .$$
\hB

On $B\Z/k\Z$ we consider the flat line bundle $V$ given by the character $\Z/k\Z\to U(1)$, 
$[1]\mapsto \exp(2\pi i k^{-1})$.
Let $L^{2n-1}:=S^{2n-1}/\mu_{k}$, where the group of $k$th roots of unity $\mu_{k}$ acts by multiplication on
$S^{2n-1}\subset \C^{2n}$.
We consider the canonical class $[f:L_{k}^{2n-1}\to B\Z/k\Z]\in \MSO_{2n-1}(B\Z/k\Z)$ 
\begin{prob}
Calculate
$\cs^{V}_{c^{n}_{1}}(f)$!
\end{prob}
}
\end{ex}

 \begin{ex}{\em
Let $(M,g)$ be a closed oriented connected Riemannian three-manifold.
Then we can form the class $-\hat c_{2}(\nabla)\in \hat H^{4}_{Del}(M;\Z)$, where $\nabla$ is the Levi-Civita connection on $TM\otimes \C$. It corresponds to
the Chern-Simons invariant $\CS(M,g)$ (Definition \ref{apr2231}) under
the identification $\hat H^{4}_{Del}(M;\Z)\cong \C/\Z$ (Example \ref{apr2230}).
\begin{prob}
Show this assertion.
\end{prob}\proof
Use \ref{may2130}. \hB 
 
We consider the Hopf fibration
$S^{7}\to S^{4}$. We realize $S^{7}$ as the unit sphere in $\bH^{2}\cong \R^{8}$. The group of unit quaternions
$Sp(1)$ acts on $S^{7}$ by right multiplication. This turns the Hopf fibration into an $Sp(1)$-principal bundle. We define a connection such that the horizontal distribution is the orthogonal complement
of the vertical bundle in the round geometry of $S^{7}$. Let $E\to S^{4}$ be the complex vector bundle with induced connection $\nabla$ associated to the representation of $Sp(1)$ on $\C^{2}\cong \bH$
by left multiplication. We consider the class
$$\hat c_{2}(\nabla)\in \hat H^{4}_{Del}(S^{4};\Z)\ .$$
For $r>0$ let $S(r)\subset S^{4}$ be the distance sphere centered at the north pole.
\begin{prob}
Calculate
$$\ev(\hat c_{2}(\nabla)_{|S(r)})\in \C/\Z\ .$$
\end{prob}
\proof
We first show that
$c_{2}(\nabla)=\vol_{S^{4}}$ (the normalized volume form).
Then by \ref{may2130}
$$\ev(\hat c_{2}(\nabla)_{|S(r)})=[\vol(B(r))]\in \C/\Z\ .$$ \hB 
}
\end{ex}

\begin{ex}\label{apr2340}{\em 
\begin{ddd}
If $\nabla$ is a connection on a real vector bundle, then we define the differential lift of the Pontrjagin form by 
$$\hat p_{i}(\nabla):=(-1)^{k}\hat c_{2k}(\nabla\otimes \C)\in \hat H^{4n}_{Del}(M;\Z)
\ .$$
\end{ddd}
If $(M,g)$ is a Riemannian manifold and $\nabla$ is the Levi-Civita connection on $TM$, then
we set  $\hat p_{k}(g):=\hat p_{k}(\nabla)$. 
We can now extend the Definition of the Chern-Simons invariant to the non-bounding case.
\begin{ddd} 
If $(M,g)$ is a closed oriented connected $4n-1$-dimensional Riemannian manifold, then we define
$$\CS(M,g):=\hat p_{n}(g)\in \hat H^{4n}_{Del}(M;\Z)\cong \C/\Z\ .$$
\end{ddd}
\begin{prob}
Verify, that this extends \ref{apr2341}.
\end{prob}
\begin{lem}
$\CS(M,g)$ is a conformal invariant.
It vanishes e.g. if $M$ bounds (with product structure) a locally conformally flat manifold.
\end{lem}
\proof
We consider a metric  on $\R\times M$ of the form 
$\tilde g=f(dt^{2}+  g)$, where $f\in C^{\infty}(\R\times M)$. 
Then we have  $p_{n}(\tilde g)=0$.
It follows from the homotopy formula that 
$\hat p_{n}(\tilde g)_{|\{t\}\times M}$ is independent of $t$.
\begin{prob}
Let $\Z/p\Z$ act on $S^{4n-1}\subset \C^{2n}$ diagonally by
$[1]\mapsto (\zeta^{q_{1}},\dots,\zeta^{q_{2n}})$, where $\zeta=\exp(2\pi i p^{-1})$
is a primitive root of unity and the numbers $g.g.T(q_{1},\dots,q_{2n},p)=1$ and set
$L^{4n-1}_{p}(q_{1},\dots,q_{2n}):=S^{4n-1}/(\Z/p\Z)$ with the induced round metric $g$.
Calculate
$$\CS(L^{4n-1}_{p}(q_{1},\dots,q_{2n}),g)\in \C/\Z\ .$$
\end{prob}
}

\end{ex}

\begin{ex}{\em
 If $V\to M$ is a real vector bundle with connection $\nabla$, then we write $\nabla\otimes \C$ for the induced connection on the complexification $V\otimes \C$.
We have 
$c_{2n+1}(\nabla\otimes V)=0$ for all $n\ge 0$.
It follows that
$$\hat c_{2n+1}(\nabla\otimes \C)\in H^{4n+1}(M;\C /\Z)\subseteq \hat H_{Del}^{4n+2}(M;\Z)\ .$$ In particular, this class is independent of the connection and only depends on the real bundle $V$.

\begin{prob}\label{may0211}
Calculate
$\beta(\hat c_1(\nabla \otimes \C))\in H^{2}(\R\P^{n};\Z)\cong \Z/2\Z$ for the canonical bundle $(V,\nabla)$ of $\R\P^{n}$, where $\beta:H^{1}(\R\P^{n};\Z/2\Z)\to H^{2}(\R\P^{n};\Z)$ is the Bockstein.
\end{prob}
\proof
$\beta(\hat c_1(\nabla \otimes \C))$ 
is the generator. Indeed, we have an embedding $\R\P^{n}\hookrightarrow \C\P^{n}$
such that $V\otimes \C$ is the restriction of the canonical bundle $L\to \C\P^{n}$.
We have
$\beta(\hat c_1(\nabla \otimes \C))=c_1^{\Z}(L)_{|\R\P^{n}}$. The right-hand side is known to restrict to a generator.
\hB 

We have a natural map
$i:H^{*}(M;\Z/2\Z)\to H^{*}(M;\C/\Z)\to \hat H^{*}_{Del}(M;\Z)$
induced by $\Z/2\Z\hookrightarrow \C/\Z$, $[1]\to[\frac{1}{2}]$.
With this notation
$\hat c_1(\nabla\otimes \C)=i(w_1(V))$.
}\end{ex}

\begin{ex}
{\em 
We consider the bottom of the Whitehead tower of $BO$:
$$\xymatrix{ 
BFivebrane=BO\langle9\rangle\ar[d]&\\
BString=BO\langle8\rangle\ar[d]^{\pi_{8}}\ar[r]^{\frac{p_{2}}{6}}&K(\Z,8)\\
BSpin=BO\langle4\rangle\ar[d]^{\pi_{4}}\ar[r]^{\frac{p_{1}}{2}}&K(\Z,4)\\
BSO=BO\langle2\rangle\ar[d]\ar[r]^{w_{2}}&K(\Z/2\Z,2)\\
BO\ar[r]^{w_{1}}&K(\Z/2\Z,1)
}$$
In order to understand this discuss the following topological exercise.
\begin{prob}
We have the Pontrjagin classes $p_{i}\in H^{4i}(BO;\Z)$.
Let $\pi_{k}:BO\langle k\rangle\to BO$ be the projection.
Observe that there are unique classes $\frac{p_{1}}{2}\in H^{4}(BSpin;\Z)$ and
$\frac{p_{2}}{6}\in H^{8}(BString;\Z)$ (these symbols are names!) such that
$\pi_{4}^{*}p_{1}=2\frac{p_{1}}{2}$ and $\pi_{8}^{*}p_{2}=6\frac{p_{2}}{6}$. 
\end{prob}
\begin{ddd}
A (stable) $O\langle k-1\rangle$-structure on
a real vector bundle $V\to M$ is a lift of the classifying map $v:M\to BO$
 of the stabilization of $V$ to $\tilde v:M\to BO\langle k\rangle$.
 \end{ddd}
 We use special names (oriented, spin, string and fivebrane) in the first few cases. 
\begin{prob}\label{may2020}
\begin{enumerate}
\item
Let $(V,\nabla)$ be a real vector bundle with spin structure. Show that
$\hat c_{1}(\nabla\otimes \C)=0$.
\item Show that the characteristic form $\frac{1}{2}p_{1}$ for real spin vector bundles with connection
is integral and has unique differential refinement $\widehat{\frac{p_{1}}{2}}$.
 Show that $\widehat{\frac{p_{1}}{2}}$ is additive under direct sum.
\item Show that the characteristic form $\frac{1}{6}p_{2}$ for real string vector bundles with connection
is integral and has unique differential refinement $\widehat{\frac{p_{2}}{6}}$.
 Show that $\widehat{\frac{p_{2}}{6}}$ is additive under direct sum.
\end{enumerate}
\end{prob}
}
\end{ex}

\begin{ex}\label{may2004}{\em
We refer to \cite{MR1393940} for more details on the following material.
Let $G$ be a Lie group with Lie algebra denoted by  $\gaaa$
  and $\pi:P\to M$ be a $G$-principal bundle. By definition this means that
$P$ admits a free right $G$-action such that $M\cong P/G$ and  $\pi$ has local sections. A connection on $P$ is a $G$-invariant decomposition $TP\cong T^{v}\pi\oplus T^{H}\pi$, where $T^{v}\pi:=\ker(d\pi)$ is the vertical bundle generated by the fundamental vector fields $A^{\sharp}\in \cX(P)$, $A\in \gaaa$, of the $G$-action and
$T^{H}\pi$ is a complementary bundle. Equivalently, a connection determines and is determined through
$$T^{h}\pi=\ker(\omega)$$ by a 
form $\omega\in \Omega^{1}(P)\otimes \gaaa$ which satisfies
\begin{enumerate}
\item $R_{g}^{*}\omega=\Ad(g)^{-1}\omega$
\item $\omega(A^{\sharp}(p))=A$ for all $p\in P$ and $A\in \gaaa$.
\end{enumerate}
The curvature of the connection $\omega$ is defined by
$$R^{\omega}:=d\omega+\frac{1}{2}[\omega,\omega]\ .$$
We define the adjoint bundle $\Ad(P):=P\times_{G,\Ad}\gaaa$.
Then we can consider the curvature as a horizontal invariant two-forms
$$R^{\omega}\in (\Omega^{2}(P)\otimes \gaaa)^{hor, G}\cong \Omega^{2}(M,\Ad(P))\ .$$
It satisfies the Bianchi identity
$$\nabla^{\omega}R^{\omega}:=d\Omega+[\omega,\Omega]=0\ .$$

Let us now assume that the  Lie group $G$ has finitely many components. By $I^{*}(G)\subseteq S^{*}(\gaaa_{\C}^{*})$ we denote the algebra  of  $G$-invariant complex-valued symmetric polynomials on $\gaaa$. An element $\phi\in I^{k}(G)$ induces a characteristic form
$$\cw(\phi):(P\to M,\omega)\mapsto \cw(\phi)(\omega)\in \Omega_{cl}^{2k}(M;\C)$$
for $G$-principal bundles with connection
by
$$\cw(\phi)(\omega)=\phi(\frac{R^{\omega}}{2\pi i})\ ,$$
where we apply $\phi$ fibrewise.
 \begin{prob}
Show  the Bianchi-identity $\nabla^{\omega} R^{\omega}=0$ and use it in order to verify that
$\cw(\phi)(\omega)$ is  closed.
\end{prob}
We get a characteristic class for $G$-principal bundles
$$\cw(\phi)(P):=[\cw(\phi)(\omega)]\in H^{2k}(M;\C)$$
(for some choice of connection $\omega$ on $P$) and therefore a universal class
$$\cw(\phi)\in H^{2k}(BG;\C)\ .$$
\begin{prob}
Use transgression in order to show that
$\cw(\phi)$ is well-defined. 
\end{prob}
 
\begin{ddd}
The homomorphism of graded rings  
 $$\cw:I^{*}(G)\to H^{2*}(BG;\C)$$ constructed above is called the Weil homomorphism.
 \end{ddd}
   In fact,  for compact $G$ the map $\cw$ is an isomorphism  by  a theorem of H. Cartan.
\begin{prob}\label{jun1103}
Prove the theorem of Cartan that for compact Lie groups the Weil homomorphism is an isomorphism.
 \end{prob}
\proof
Here are the main steps. Compare e.g. \cite{MR0500997} for a complete argument.

  Let $T\subseteq G$ be a maximal torus and $W:=W(G,T)$ be the Weyl group. 
  
 Show that the restriction induces an isomorphism
$$I^{*}(G)\stackrel{\sim}{\to} S^{*}(\taaa^{*}_{\C})^{W}\ .$$
The main point here is that for every $A\in \gaaa$ there exists $g\in G$
such that $\Ad(g)(A)\in \taaa$.
Further we identify 
$$S^{*}(\taaa^{*}_{\C})^{W}\stackrel{\sim}{\to}H^{2*}(BT;\C)^{W}$$
by an explicit   calculation.

Then we have the diagram
$$\xymatrix{I (G)\ar[d]^{\cong}\ar[r]^{\cw}&H (BG;\Z)\ar[d]^{i^{*}}\\
S^{*}(\taaa^{*}_{\C})^{W}\ar[r]^{ \cong}&H(BT;\Z)}$$
It suffices to show that
$i^{*}:H (BG;\Z)\to H(BT;\Z)$ is injective.  
We have a fibration $$B/T\to BT\to BG\ .$$
It suffices to  find an element  
$\alpha\in H^{\dim(G/T)}(BT;\C)$ such that
$f_!(\alpha)=1$. Indeed, then we have
$$\beta=f_!(\alpha\cup i^{*}\beta)\ ,\quad \forall \beta\in H^{*}(BG;\C)\ .$$

 Let $\Delta^{+}(\gaaa,\taaa)\subset  \Hom(T,\C^{*})$ be the set of positive roots. For every
$\alpha\in \Delta^{+}(\gaaa,\taaa)$ let $L_\alpha\to BT$ be the associated Line bundle. One checks that
$$ \alpha:=c\prod_{\alpha\in \Delta^{+}(\gaaa,\taaa)} c_1(L_\alpha)$$
does the job for an appropriate choice of the normalization $c$.
\hB

As before we say that
$\phi\in I^{k}(G)$ is integral if there exists a class
$z\in H^{2k}(BG;\Z)$ such that $\epsilon_{\C}(z)=\cw(\phi)$.
Unlike the case $G=U(n)$ we can not expect that $z$ is uniquely determined by 
$\phi$. We define the graded ring
$\tilde I^{*}(G)$ as the pull-back
$$\xymatrix{\tilde I^{*}(G)\ar[r] \ar[d]&I^{*}(G)\ar[d]^{\cw}\\
H^{2*}(BG;\Z)\ar[r]^{\epsilon_{\C}}&H^{2*}(BG;\C)}\ .$$
For $\tilde \phi\in \tilde I^{*}(G)$ we let $\phi\in I^{*}(G)$ and
$\phi_{\Z}\in H^{2k}(BG;\Z)$ denote the underlying invariant polynomial and integral
cohomology class. The following theorem is due to 
Cheeger-Simons \cite[Thm. 2.2 ]{MR827262}.

\begin{theorem}\label{may2010} 
Assume that $G$ is a Lie group with finitely many components.
For every 
$\tilde \phi \in \tilde I^{k}(G)$ there exists a unique $\hat H^{2k}_{Del}(\dots;\Z)$-valued characteristic class for $G$-principal bundles with connection
$$\hat \cw(\tilde \phi):(P\to M,\omega) \mapsto \hat \cw(\tilde \phi)(\omega)\in \hat H^{2k}_{Del}(M;\Z)$$
such that
$$R(\hat \cw(\tilde \phi)(\omega))=\cw(\phi)(\omega)\ ,\quad I(\hat \cw(\tilde \phi)(\omega))=\tilde \phi_{\Z}(P)\ .$$
\end{theorem}
\proof
Observe that $H^{odd}(BG;\C)=0$.
We can now proceed as in the proof of  Theorem \ref{may2007}.
Let $\tilde \phi\in \tilde I^{k}(G)$.
If $N\to BG$ is a $2k+1$-connected approximation, then $H^{2k-1}(N;\C)=0$. If 
$(P\to N,\omega)$ is a $G$-principal bundle with connection, then the class
$$\hat \cw(\tilde \phi)(\omega)\in \hat H^{2k}_{Del}(M;\Z)$$
is uniquely determined by
$$R(\hat \cw(\tilde \phi)(\omega))=\cw(\phi)(\omega)\ ,\quad I(\hat \cw(\tilde \phi)(\omega))=\tilde \phi_{\Z}(P)\ .$$
\hB 

\begin{prob}
Find out where we have used the assumption that $G$ has finitely many connected components.
\end{prob}
\proof
We use that $H^{odd}(BG;\C)=0$. This is not true without this assumption. For example, we have
  $$H^{1}(BGL(n;\Q^{\delta});\C)\cong
\Hom(\Q^{*},\C)\cong \prod_{p\in \nat,prime} \C\ .$$  

Assume now that $G$ is a Lie group with finitely many components.
Since we consider the cohomology with coefficients in a rational vector space we can reduce to the case that $G$ is connected.
Now one has to go through the structure theory of connected Lie groups which says that it can be reduced to product of simple  groups by iteratively factoring out   normal abelian subgroups.  A simple group can further be reduced to its maximal compact Lie group. For the compact group we apply the surjectivity part of Cartan's theorem, or rather its substep, the injectivity of $H^{*}(BG;\C)\to H^{*}(BT;\C)$
for the restriction to a maximal torus. Vanishing of the odd cohomology is then preserved if one argues backwards through the
previous reduction steps.
\hB

\begin{prob}
View the assertions in \ref{may2020} as special cases of Theorem \ref{may2010}.
\end{prob}
}\end{ex}

\begin{ex}
{\em 
Let $\tilde \phi=(\phi,\phi_{\Z})\in \tilde I^{*}(G)$ and assume that $\omega$ has finite holonomy $H\subseteq G$. The bundle  $P\to M$ has an $H$-reduction classified by a map $h:M\to BH$. Consider further the  map
$$B:H^{*}(BG;\Z)\to H^{*}(BH;\Z)\cong H^{*-1}(BH;\C/\Z)\ .$$
The following is \cite[Prop. 2.10]{MR827262}.
\begin{prob}
Show that $\hat \cw(\tilde \phi)(\omega)=h^{*}B(\phi_{\Z})$.
\end{prob}

 We consider  a Lie group $G$ with finitely many component and its discrete version $G^{\delta}$.
\begin{prob}
Show that for every $\tilde \phi\in \tilde I^{k}(G)$ there exists an element
$\hat \phi\in H^{2k-1}(BG^{\delta};\C/\Z)$ which represents the 
universal class $\hat \cw(\tilde \phi)$ for $G$-bundles with flat connection.
\end{prob}
}\end{ex}

\begin{ex}{\em
In this example we study the Euler class. 
We identify
$$\so(2n)_{\C}\cong \Lambda^{2}\C\ ,\quad X \cong \frac{1}{2}\sum_{i,j=1}^{2n} \langle e_{i},Xe_{j}\rangle e^{i}\wedge e^{j}$$
and define  the Pfaffian $\Pf\in S^{n}(\so(2n)_{\C}^{*})$
by the identity $$\Pf(X)\vol_{\R^{n}}=\underbrace{X\wedge\dots \wedge X}_{n\times}\quad \forall X\in \so(2n)_{\C}\ .$$
We further introduce the normalized version 

\begin{equation}\label{jun0611}\bar \Pf(A):=\Pf(-i A)\ . \end{equation}
Let $S(\xi_{m})\to BSO(m)$ be the unit sphere bundle of the universal bundle $\xi_{m}\to BSO(m)$.
We consider the associated Leray-Serre spectral sequence. The differential
$d_{2}:E^{0,m-1}_{2}\to E_{2}^{m,0}$ is a map
$$d:H^{m-1}(S^{m-1};\Z)\to H^{m}(BSO(m);\Z)\ ,$$ and we define the universal Euler class 
$\chi\in H^{m}(BSO(m);\Z)$ by  $\chi:=d_{2}(\ori_{S^{m-1}})$.
\begin{prob}
\begin{enumerate}
\item Show that $\bar \Pf\in I^{n}(SO(2n))$.
\item
Show that $2\chi=0$ if $m$ is odd.
\item Show that $\cw(\bar \Pf)=\epsilon_\C(\chi)$ if $m$ is even.
\end{enumerate}
\end{prob}
\proof
For 2. use the fibrewise antipodal map.

For 3. we argue as follows. The map $U(n)\to SO(2n)$ induces an isomorphism of maximal tori and therefore an injection
  $i^{*}:H^{*}(BSO(2n);\C)\to H^{*}(BU(n) ;\C)$. It therefore  suffices to show the identity as an identity of characteristic classes for  $n$-dimensional complex  vector bundles.  
This has been done in Example \ref{jun1102}. \hB

We define $$\tilde \Pf:=\left\{ \begin{array}{cc}(\chi,\bar \Pf) &m \:\:\mbox{even}\\
(\chi,0)&m\:\:\mbox{odd}\end{array}\right\}\in  \tilde I^{n}(SO(m)) $$ and $$\hat \chi:=\hat \cw( \tilde \Pf)\ .$$ 
If $V\to M$ is a real oriented $m$-dimensional vector bundle with connection $\nabla$, then we write
$$\hat \chi(\nabla):=\hat \chi(\omega)\ ,$$
where $\omega$ is the corresponding connection on the oriented frame bundle $SO(V)\to M$.
\begin{ddd} The characteristic class $\hat \chi$  is the  differentially refined version of the Euler class. 
 \end{ddd} 
Let $V\to M$ be a real oriented euclidean vector bundle of dimension $m$ with a flat euclidean connection $\nabla^{V}$.
Let $\beta$ be the fibrewise volume $m-1$-form (induced by the metric) on the sphere bundle $\pi:S(V)\to M$ normalized to have integral one.
Using the connection we extend this to a form $\beta\in \Omega^{m-1}(S(V);\C)$. Flatness of the connection implies $d\beta=0$. The segment  of  the homological Gysin sequence
$$\dots H_{0}(M;\Z)\to H_{m}(S(V);\Z)\stackrel{\pi_{*}}{\to} H_{m}(M;\Z)\to 0$$
shows:
\begin{enumerate}
\item For every $x\in H_{m}(M,\Z)$ we can find a preimage $y\in H_{m}(S(V);\Z)$
such that $\pi_{*}(y)=x$.
\item The class  of the evaluation 
$f(\nabla)(x):=[\langle [\beta],y \rangle]\in \C/\Z$
does not depend on the choice of the preimage.
\end{enumerate}
We therefore get a class
$$f(\nabla)\in \Hom(H_{m-1}(M);\C/\Z)\cong H^{m-1}(M;\C/\Z)\ .$$
\begin{prob}
Show that
$$a(f)=\hat \chi(\nabla)\ .$$
\end{prob}
Assume in addition that $M$ is closed, oriented and $m-1$-dimensional.
Choose a triangulation $(\sigma_i:\Delta^{m-1}\to M)_{i\in I}$ of $M$.
Let $(m_j)_{j\in J}$ be the set of vertices. For every $j\in J$ choose
a unit vector $v_j\in V_{m_j}$. Let $b_i$ denote the barycenter of the simplex $\sigma_i$ .
Parallel transport along 
straight lines in $\sigma_i$
of the vectors  associated to the vertices of $\sigma_i$ gives 
collections of unit vectors $u_i:=(w_i(0),\dots,w_i(m-1))\in V_{b_i}$.
We assume that these collections $u_i$ are bases for all $i\in I$.
The basis $u_i$ spans a geodesic simplex in the unit sphere $S(V_{b_i})$
of oriented volume $\vol(u_i)$. The following is due to Cheeger-Simons \cite[Thm. 8. 14]{MR827262}.
\begin{prob}
Show that
$$\ev(\hat \chi(V))=[\sum_{i\in I} \vol(u_i)]\in \C/\Z\ .$$
\end{prob}

    If $(E,\nabla)$ is a  complex vector bundle  with connection of dimension $n$, then we let $(E_{|\R},\nabla_{|\R})$ denote the underlying oriented real vector bundle with connection.
  \begin{prob} Show that for a complex vector bundle
  $E\to M$  with connection $\nabla$ we have
  $\hat c_{n}(\nabla)=\hat \chi(\nabla_{|\R})$.
 \end{prob}
\proof
Since we consider characteristic classes of complex vector bundles it suffices to show that the underlying curvatures coincide.
We check this first in the case $n=1$.  At a fixed point $x\in M$ we can   identify $E_x\cong \C$  and assume      that $R^{\nabla}=i  \alpha$ for some $\alpha\in \Lambda^{2}T^{*}_xM$. Then we have
$$R^{\nabla_\R}=\left(\begin{array}{cc}0&-\alpha\\
\alpha&0\end{array}\right)$$
with respect to the standard basis $\{e_1,e_2\}$ of $\C_{|\R}\cong \R^{2}$. This implies at this point
$$  \cw(\bar \Pf)(R^{\nabla_{|\R}})=\Pf(-i \frac{R^{\nabla_{|\R}}}{2\pi i})=-\frac{ \alpha}{2\pi }\ .$$
On the other hand,
$$c_1(\nabla)=\frac{-R^{\nabla}}{2\pi i}=-\frac{ \alpha}{2\pi  } \ .
$$
We extend this to the general case using diagonalization. \hB 
 }\end{ex}

\begin{ex}{\em
Let $P\to M$ be a $G$-principal bundle with connection $\omega$.
Let $\tilde \phi=(\phi_\Z,\phi)\in \tilde I^{n} (G)$.
\begin{prob}\label{jun1105}
Show that there exists a  natural (under pull-back of principal bundles  with connection) form $\theta(\omega)\in \Omega^{2n-1}(P;\C)$ such that
$d\theta(\omega)=\phi(R^{\omega})$.
Show that $$\pi^{*}\hat \cw(\tilde \phi)(\omega)=a(\theta(\omega))\ .$$
\end{prob}
\proof
Note that $\pi^{*} P\to P$ is canonically trivialized and has the trivial connection $\omega_0$. We define
$\theta:=\widetilde{\cw(\phi)}(\omega,\omega_0)$
using transgression along the linear path from $\omega_0$ to
$\pi^{*}\omega$. The second assertion is a consequence of the homotopy formula \ref{apr2201}. \hB

We consider $G:=SU(2)$ and $\tilde \phi \in I(SU(2))$ such that
$\phi_\Z=c_2$. 
\begin{prob}
Give an explicit  formula for
$\theta(\omega)$ in \ref{jun1105}.
\end{prob}

}
\end{ex}

\subsection{Multiplicative structure}

Let $R\in \Sh_{\Rings}(\Mf)$ be a sheaf of commutative rings. Then the sheaf cohomology
$H^{*}(M;R)$ becomes a graded commutative ring.  Here are the basic steps to see this.
The notion of sheaf of rings can be formalized
using the symmetric monoidal structure on $\Sh_{\Ab}(\Mf)$. The evaluation
$\Gamma(M;\dots)$ is lax symmetric monoidal. Using the fact that the category of sheaves contains sufficiently many flat sheaves we see that these properties descend to the derived category $D^{\ge 0}(\Sh_{\Ab}(M))$ so that the derived functor of evaluation preserves rings, too. 
Differential graded algebras (abbreviated by dga's) are ring objects in $\Ch(\Ab)$ and in particular  present rings in the derived category $D^{\ge 0}(\Sh_{\Ab}(M))$.

The sheaf of de Rham complexes $\Omega_{\C}$ is a sheaf of dga's which resolves the sheaf of rings
$\underline{\C^{\delta}}$. The wedge product of forms induces a product on $H_{dR}(M;\C)$.
The de Rham isomorphism
$\Rham:H^{*}(M;\underline{\C^{\delta}})\stackrel{\sim}{\to} H_{dR}(M;\C)$
is induced by the quasi isomorphism of sheaves of dga's $\underline{\C}^{\delta}\to \Omega_{\C}$
and is hence multiplicative. Moreover, the map
$\epsilon_{\C}:H^{*}(M;\Z)\to H_{dR}^{*}(M;\C)$ is multiplicative
since it is induced by the composition
$\underline{\Z}\to \underline{\C^{\delta}}\to \Omega_{\C}$
of maps of sheaves of dga's. This leaves open the question why the multiplicative structures on
integral cohomology defined using sheaf theory or simplicial cohomology coincide.  An argument will be given in   Lemma \ref{may2909}.

\begin{ddd}\label{may0210}
A product on Deligne cohomology is the datum of a  graded commutative ring structure  (denoted by $\cup$) on
$\hat H^{*}_{Del}(M;\Z)$ for every manifold $M$ such that
\begin{enumerate}
\item $f^{*}:\hat H^{*}_{Del}(M;\Z)\to \hat H^{*}_{Del}(M^{\prime};\Z)$ is a homomorphism of rings for every smooth map $f:M^{\prime}\to M$,
\item $R:\hat H^{*}_{Del}(M;\Z)\to \Omega^{*}_{cl}(M;\C)$ is multiplicative for all $M$,
\item $I:\hat H^{*}_{Del}(M;\Z)\to H^{*}(M;\Z)$ is multiplicative for all $M$,
\item and $a(\alpha)\cup x=a(\alpha\cup R(x))$ for all $\alpha\in \Omega^{*}(M;\C)/\im(d)$ and $x\in \hat H^{*}_{Del}(M;\Z)$.
\end{enumerate}
\end{ddd}

\begin{prop}\label{aug2001}
There exists a unique product on Deligne cohomology.
\end{prop}
\proof
We first show existence.
We shall use the complexes $\cE(p)$ defined in \eqref{jun1301} in order to represent Deligne cohomology, see \eqref{may1703}. We 
construct products
\begin{equation}\label{may0102}\cup:\cE(p)\otimes \cE(q)\to \cE(p+q)\end{equation}
by
$$
x\cup y:=\left\{\begin{array}{cc} xy& \deg(x)=0\:\mbox{or}\: \deg(y)=0\\
x\wedge (-dy)&\deg(x)>0\: \mbox{and}\: \deg(y)=q>0\\
0&\mbox{otherwise}\end{array}\right.
$$
\begin{prob}
Show that this is a morphism of complexes.
Furthermore verify associativity.
\end{prob}
\proof
We let $\partial$ denote the differential of the complexes, while $d$ is the de Rham differential.
$\deg(x)$ denotes the degree of $x$ as an element in $\cE(p)$, not the degree as a form. 
We show that the $\cup$-product is a morphism of complexes. We consider the case that $p,q>0$.
There are many cases:
\begin{enumerate}
\item
If $\deg(x)=0$ and $\deg(y)=0$, then we have in $\Omega^{0}(M;\C)=\cE(p+q)^{1}$ that
$\partial (x\cup y)=xy$ and
$\partial x\cup y+x\cup \partial y=xy$.
\item If $0<\deg(y)<q$, then in $\Omega^{\deg(y)}(M,\C)$ we have
$\partial (x\cup y)=-xdy$ and $\partial x\cup y+x\cup \partial y=-xdy$.
\item If $\deg(y)=q>0$, then in $\Omega^{q}(M,\C)$ we have
$\partial (x\cup y)=-xdy$ and $\partial x\cup y+x\cup \partial y=-xdy$.
\item If $\deg(x)>0$ and $0<\deg(y)<q$, then we have
$\partial (x\cup y)=0$ and $\partial x\cup y +(-1)^{\deg(x)}x\cup \partial y=0$.
\item If $\deg(x)>0$ and $\deg(y)=q>0$, then
$\partial (x\cup y)=-d(x\wedge -dy))=dx\wedge dy$ and
$\partial x\cup y+(-1)^{\deg(x)} x\wedge dy=dx\wedge dy$.
\end{enumerate}
Let us check associativity in  some cases We consider the product
$\cE(p)\otimes \cE(q)\otimes \cE(r)\to \cE(p+q+r)$.
\begin{enumerate}
\item If $\deg(x)>0$ and $0<\deg(y)<q$, and $\deg(z)=r$.
Then we have
$(x\cup y)\cup z=0$ and $x\cup (y\cup z)=0$.
\item  If $\deg(x)>0$ and $\deg(y)=q$, and $\deg(z)=r$,
then $(x\cup y)\cup z=xdy\wedge dz$ and
$x\cup (y\cup z)=x\wedge dy\wedge dz$.
If $\deg(x)=\deg(y)=0$, then
$(x\cup y)\cup z=xyz$ and $x\cup (y\cup z)=xyz$.
\end{enumerate} 
The remaining cases are similar.
\hB

We define
$H:\cE(p)\otimes \cE(q)\to \cE(p+q) [-1]$ by
$$H(x\otimes y)=\left\{\begin{array}{cc} 0&\mbox{$\deg(x)=0$ or $\deg(y)=0$}\\
(-1)^{\deg(x)}x\wedge y &\mbox{otherwise}
\end{array}\right.
$$
Let $s:\cE(p)\otimes \cE(q)\to \cE(q)\otimes \cE(p)$ be the symmetry in the symmetric tensor category graded abelian groups. 
\begin{prob}
Show that
$H$ is a homotopy between
$\cup$ and $\cup\circ s$.
\end{prob}
\proof
Again we consider several cases.
\begin{enumerate}
\item Assume that
$0<\deg(x)<p$ and $0<\deg(y)<q$.
Then we have
\begin{eqnarray*}\partial H(x\otimes y)+H(\partial (x\otimes y))&=&-d(-1)^{\deg x} (x\wedge y)+H(-dx\otimes y+(-1)^{\deg(x)} x\otimes (-dy))\\
&=&-(-1)^{\deg(x)}dx\wedge y +x\wedge dy \\&&+(-1)^{\deg(x)}dx\wedge y -x\wedge dy  \\&=&0\\
 (-1)^{\deg(x)\deg(y)} y\cup x-x\cup y&=&0\\
\end{eqnarray*}
\item If $0<\deg(x)<p$ and $\deg(y)=q>0$, then 
\begin{eqnarray*}\partial H(x\otimes y)+H(\partial (x\otimes y))&=&-d(-1)^{\deg x} (x\wedge y)+H(-dx\otimes y )\\
&=&-(-1)^{\deg(x)}dx\wedge y +x\wedge dy \\&&+(-1)^{\deg(x)}dx\wedge y   \\&=&x\wedge dy\\
 (-1)^{\deg(x)\deg(y)} y\cup x-x\cup y&=&x\wedge dy
 \end{eqnarray*}
\item For $\deg(x)=0=\deg(y)$ we have
\begin{eqnarray*}
\partial H(x\otimes y)+H(\partial (x\otimes y))&=0\\
  (-1)^{\deg(x)\deg(y)} y\cup x-x\cup y&=&yx-xy\\
&=&0
\end{eqnarray*}
If $\deg(x)=p$ and $0<\deg(y)<q$ then
\begin{eqnarray*}\partial H(x\otimes y)+H(\partial (x\otimes y))&=&-d(-1)^{\deg x} (x\wedge y)+H((-1)^{\deg(x)} x\otimes (-dy))\\
&=&-(-1)^{\deg(x)}dx\wedge y +x\wedge dy \\&&  -x\wedge dy  \\&=&-(-1)^{\deg(x)}dx\wedge y\\
 (-1)^{\deg(x)\deg(y)} y\cup x-x\cup y&=&-(-1)^{\deg(x)\deg(y)} y\wedge dx\\
&=& -(-1)^{\deg(x)} dx\wedge y
 \end{eqnarray*}
\end{enumerate}
\hB 
We conclude that \eqref{may0102} induces a product
$$\cup:\hat H^{p}_{Del}(M;\Z)\otimes \hat H^{q}_{Del}(M;\Z)\to \hat H^{p+q}_{Del}(M;\Z)$$
which is natural, associative and graded commutative.

We now show that this product is compatible with the structure maps.
First of all we observe that $R$ is induced by the map
$$\tilde R :\cE(p)\to \Omega^{p}_{\C}\ ,$$
$$y\mapsto -dy\ , \deg(y)=p\ , \quad y\mapsto 0 \:\mbox{if $\deg(y)<p$}\ .$$
Indeed, if $\deg(y)=p$, then
$((0\oplus -dy)\oplus y)\in \cD(p)$ is a lift, and
$R$ is induced by the projection onto the $-dy$-component.
We see that
$$\tilde R(x\cup y)=R(x)\wedge R(y)\ .$$
Simlarly,
$I$ is induced by the projection
$$\tilde I:\cE(p)\to \underline{\Z}$$
$$y\mapsto y\ , \deg(y)=0\ , \quad y\mapsto 0 \:\mbox{if $\deg(y)>0$}\ .$$
Indeed, if $\deg(y)=0$, then
$((y\oplus 0)\oplus 0)\in \cD(p)$ is a lift and
 $I$ is induced by the projection onto the first component.
\begin{prob}
Show the compatibility with $a$.
\end{prob}
We now show uniqueness. Assume that
$\cup^{\prime}$ is a second product.
Then we consider $$B:=\cup^{\prime}-\cup:\hat H^{p}_{Del}(\dots;\Z)\otimes \hat H^{q}_{Del}(\dots;\Z)=\hat H^{p+q}(\dots;\Z)\ .$$
This is a bilinear natural transformation.
\begin{prob}
Show that $B$ factorizes over a bilinear transformation
$$H^{p}(\dots;\Z)\otimes H^{q}(\dots;\Z)\to H^{n-1}(\dots;\C/\Z)\ .$$
\end{prob}
We now argue that such a transformation is necessarily zero.
A natural transformation of functors
$$H^{p}(\dots;\Z)\times H^{q}(\dots;\Z)\to H^{n-1}(\dots;\C/\Z)$$ is represented by a map of Eilenberg-Mac Lane spaces
$$K(\Z,p)\times K(\Z,q)\to K(\C/\Z,p+q-1)\ .$$
If the transformation is bilinear, then this map factorizes over
$$K(\Z,p)\wedge K(\Z,q)\to K(\C/\Z,p+q-1)\ .$$
Since this smash product  is $p+q-1$-connected this latter  map is homotopic to a constant map.
We conclude that $B=0$. It follows $\cup=\cup^{\prime}$. 
\hB

\begin{ex}\label{jun1420}{\em
We identify $J(S^{1})\cong S^{1}$.
Consider the Poincar\'e bundle $P\to S^{1}\times S^{1}$
with its canonical connection $\nabla^{P}$ described in \ref{apr1206}.
 \begin{prob}
Show that in $\in \hat H_{Del}^{2}(S^{1}\times S^{1};\Z)$
$$ \hat c_{1}(\nabla^{P})=\pr_{1}^{*}\hat e\cup \pr_{2}^{*}\hat e \ .$$
\end{prob}
\proof
Define
$$\delta:=\hat c_{1}(\nabla^{P})-\pr_{1}^{*}\hat e\cup \pr_{2}^{*}\in  \hat H_{Del}^{2}(S^{1}\times S^{1};\Z)\ .$$
We first see  that $R(\delta)=0$. Hence  $$\delta\in H^{1}( S^{1} \times S^{1};\C/\Z)\cong \Hom(H_{1}( S^{1} \times S^{1};\Z),\C/\Z)\ .$$ Therefore we must calculate
the holonomy of $\delta$ along the two basis cycles $ S^{1} \times\{1\}$ and $\{1\}\times S^{1}$.
The holonomies of both terms vanish separately. \hB

}\end{ex}

\begin{ex}\label{may2001}{\em
We consider two classes
$x,y\in \hat H^{1}_{Del}(S^{1};\Z)$ and want to calculate
$x\cup y\in \hat H^{2}_{Del} (S^{1};\Z)\cong \C/\Z$.
Note that we can write
$x=n\hat e+a(f)$ and $y=m\hat e+a(g)$ for some $f,g\in C^{\infty}(S^{1};\C)$.
These functions are determined up to integral constants by the differential equations
$$R(x)-ndt=df\ , [f(0)]=[x_{|\{0\}}]\ ,\quad R(y)-mdt=dg\ , [g(0)]=[y_{|\{0\}}]\ .$$
Then
$$x\cup y=nm\hat e\cup \hat e + a(( mf-ng) dt+f\wedge dg)\ .$$
This gives 
$$\ev(x\cup y)=nm\ev(\hat e\cup \hat e)+ [\int_{S^{1}} (fm-ng) dt+f\wedge dg ]\ .$$
We must calculate $\ev(\hat e\cup \hat e)$. Let $\diag :S^{1}\to S^{1}\times S^{1}$
be the diagonal. By Example \ref{jun1420} we have
$$\ev(\hat e\cup \hat e)=\ev(\diag^{*} \hat c_1(\nabla^{P}))=\hol_\nabla(\diag)\ .$$
Since $\diag$ is bordant to the union of the cycles $S^{1} \times\{1\}$ and
$\{1\}\times S^{1}$ on which $\nabla^{P}$ is trivial, we can calculate the holonomy
along $\diag$ by a curvature integral over the bordism, a triangle covering half of the torus.
We get
$$\ev(\hat e\cup \hat e)=[\frac{1}{2}]\ .$$


The final formula is 
$$\ev(x\cup y)=\left[\frac{nm}{2}+ \int_{S^{1}} (fm-ng) dt+f\wedge dg \right]\ .$$
}
  \end{ex}

\begin{ex}{\em 
Note that $H^{*}(M;\C/\Z)$ is an $H^{*}(M;\Z)$-module in the natural way.
Let $i:H^{*}(M;\C/\Z)\to \hat H^{*+1}_{Del}(M;\Z)$ be the inclusion. The following identity extends \ref{may0210}, 4.
\begin{prob}\label{may0212}
Show that
$i(x)\cup y=i(x\cup I(y))$.
\end{prob}
\proof
First observe that the formula holds true if $x=a(\alpha)$.
Furthermore observe that $i(x)\cup y$ only depends on $I(y)$.
Both assertions  follow from Definition \ref{may0210}, 4. From this conclude that if $\deg(x)=n$ and $H^{n}(M;\Z)$ is torsion free, then the formula holds true.  For the general case observe that
there exists a smooth map  $f:M\to\tilde M$ and classes
$\tilde x\in H^{n}(\tilde M;\C/\Z)$ and $\tilde y\in H^{*}(\tilde M;\Z)$ such that
$H^{n}(\tilde M;\Z)$ is torsion-free and $f^{*}\tilde x=x$ and $f^{*}\tilde y=I(y)$.
Since the identity in question holds true for $\tilde x$ and any differential lift of $\tilde y$, it also holds true
for $x$ by naturality of the product.
\hB

}

\end{ex}

\begin{ex}\label{may3003}{\em
Let $H^{*}(M;\Z)\ni z\mapsto \overline{z}\in H^{*}(M;\Z/2\Z)$
denote  the mod-$2$ reduction. We have a Steenrood square
$$\Sq^{2n}:H^{2n+1}(M;\Z/2\Z)\to  H^{4n+1}(M;\Z/2\Z)\ .$$
Finally we have a natural map
$$i:H^{4n+1}(M;\Z/2\Z)\stackrel{\Z/2\Z\to \C/\Z}{\longrightarrow} H^{4n+1}(M;\C/\Z) \ .$$
The following has first been observed in \cite{MR2509714}. Our argument below differs from Gomi's and seems to be much simpler.

\begin{prob}\label{jun0501} 
Show that for $\hat x\in \hat H_{Del}^{2n+1}(M;\Z)$ we have
$\hat x\cup \hat x=i(\Sq^{2n}(\overline{I(\hat x)}))$.
\end{prob}
\proof 
We first discuss the case that $n\ge 1$.
We can assume that the classifying map $x$ of $I(\hat x)$ has a factorization
$$\xymatrix{&N\ar[d]^{ y }\\M\ar@{.>}[ur]^{f}\ar[r]^{ x }&K(\Z,2n+1)}\ ,$$
where $y$ is at least  $\max(4n+2,\dim(M))$-connected. We can assume that $f$ is the inclusion of a submanifold and that there exists a class
$ \hat y\in \hat H^{2n+1}_{Del}(N;\Z)$ such that
$f^{*} \hat y= \hat x$. To see this proceed as follows:
\begin{enumerate}
\item Replace $N$ by $N\times \R^{k}$, $ y $ by $  y \circ \pr_{N}$, and $f$ by $(f,j)$, where $j:M\to \R^{k}$ is an embedding.
\item Choose some class $ \hat y_{0}\in \hat H^{2n+1}_{Del}(N;\Z)$  with $I(\hat  y_{0})$ classified by $y$ and a form
$\alpha\in \Omega^{2n}(M;\C)$ such that $f^{*} \hat y_{0}=\hat x+a(\alpha)$.
\item Since $f$ is an embedding of a submanifold we can choose a form $\beta\in \Omega^{2n}(N;\C)$ such that $f^{*}\beta=\alpha$.
\item We set $ \hat y:= \hat  y_{0}+a(\beta)$.
\end{enumerate}
If $n\ge 1$, then the Bockstein operator induces an isomorphism
$$B:H^{4n+1}(K(\Z,2n+1);\C/\Z)\stackrel{\sim}{\to} H^{4n+2}(K(\Z,2n+1);\Z)\ , $$
and if $u\in H^{2n+1}(K(\Z,2n+1);\Z)$ is the universal element, then we have
$$i(\Sq^{n}(\bar u))=B^{-1}(u\cup u)\ .$$
The same holds true on $N$. Since $R(\hat  y\cup  \hat y)=R(\hat y)\cup R(\hat y)=0$ we have
$\hat y\cup  \hat y\in H^{4n+1}(N;\C/\Z)$ so that
$$B(\hat y\cup  \hat y)=I(\hat y\cup \hat y)=B(i(\Sq^{n}(\overline{I(\hat y)} )))\ .$$
Since $B$ is injective we get the result in the case $n\ge 1$.

If $n=1$, then the universal calculation is $\hat e\cup \hat e=a(\frac{dt}{2})$.
The assertion now follows from $\Sq^{0}=\id$ and $i(\overline{\ori_{S^{1}}})=a(\frac{dt}{2})$.
\hB
}
\end{ex}

\begin{ex}{\em
The following examples concern the calculation of products of degree-one classes.

We consider the maps
$\pr_{1},\pr_{2},\mu:T^{2}\to S^{1}$.
\begin{prob}
Show that
$$\pr_{1}^{*}\hat e\cup \pr_{2}^{*}\hat e
\cup \mu^{*}\hat e=0\ .$$
\end{prob}
\proof
Use that $\hat e$ is primitive and Gomi's result \ref{jun0501}. \hB
 
Let $f,g\in C^{\infty}(M;\C/\Z)$ and $x,y\in H^{1}_{Del}(M;\Z)$ be the corresponding classes (see  \ref{may0106}). We are interested in the product $x\cup y$.
\begin{prob}
For every map $\gamma:S^{1}\to M$ calculate
$\ev(\gamma^{*}(x\cup y))\in \C/\Z$.
\end{prob}
\proof
If $f=\exp(2\pi i h)$ for some function $h\in C^{\infty}(M;\C)$ (i.e. $I(x)=0$), then
$$\ev(x\cup y)=[\int_{S^{1}} \gamma^{*} h g^{-1}dg]\ .$$
For the general case we can proceed as follows.
We consider $(f,g):M\to T^{2}$ and observe that
$x\cup y=(f,g)^{*}(\pr_{1}^{*}\hat e\cup \pr_{2}^{*}\hat e)$. 
By exercise \ref{jun1420} we must calculate $\hol_{(f,g)^{*}\nabla^{P}}(\gamma)$.
  Or we write $\ev(\gamma^{*}(x\cup y))=\gamma^{*}x\cup \gamma^{*}y$ and use Problem 
\ref{may2001}.
  \hB 
  
 We consider the manifold
$\C^{*}\setminus \{1\}$ and the functions
$$[\frac{1}{2\pi i}\ln z],[\frac{1}{2\pi i}\ln(1-z)]\in C^{\infty}(\C^{*}\setminus \{1\},\C/\Z)$$
as elements
$f,g\in \hat H^{1}_{Del}(\C^{*}\setminus \{1\};\Z)$.
\begin{prob}\label{may2002}
Calculate $f\cup g$.
\end{prob}
\proof
The result is
$$f\cup g=0\ .$$
First observe by calculation that
$R(f\cup g)=0$. Furthermore, we have $I(f\cup g)=0$. It follows that
$f\cup g=a([\alpha])$ for some
$[\alpha]\in H^{1}(\C^{*}\setminus \{1\};\C)/H^{1}(\C^{*}\setminus \{1\};\Z)$.
We choose a basis $u,v\in H_{1}(\C^{*}\setminus \{1\};\Z)$ given by small counterclockwise circles around $0$ and $1$. This basis induces an identification
$H^{1}(\C^{*}\setminus \{1\};\C)/H^{1}(\C^{*}\setminus \{1\};\Z)\cong (\C/\Z)^{2}$.
We thus have to calculate
$\int_{u}\alpha$ and $\int_{v} \alpha$. 
We now use the result of Example \ref{may2001}. \hB
}\end{ex}




%

%

%
%
 
\begin{ex}{\em
In the next examples we calculate products involving higher-dimensional classes.
The calculation is usually easy if one of the classes is topologically trivial. The general case is
usually difficult and requires some tricks.

Let $M$ be a closed oriented surface,
$(L,\nabla)\in \Line_{\nabla}(M)$ be a line bundle with connection on $M$ and
$u\in C^{\infty}(M;\C/\Z)$.
Then we have
classes
$x:=\hat c_{1}(L,\nabla)\in \hat H^{2}_{Del}(M;\Z)$ and $y:=[u]\in \hat H^{1}_{Del}(M;\Z)$ (see  \ref{may0106}).
\begin{prob}
Calculate
$\ev(x\cup y)\in \C/\Z$.
\end{prob}
\proof
Assume that $L$ is trivial $\nabla^{L}=d+2\pi i\alpha$.
Then
$\hat c_{1}(L)=a(-\alpha)$ and
$$\ev(x\cup y)=[-\int_{M} \alpha\wedge  du]\ .$$
The general case is complicated. \hB

%
%
%


\begin{prob}\label{may0111}
On the two-torus $T^{2}$ with affine coordinates $s,t$ calculate the product
$\hat c_{1}(L,\nabla)\cup [u]$ for
$(L,\nabla)\in \Line_{\nabla}(T^{2})$ and
$u\in C^{\infty}(T^{2},\C/\Z)$.
\end{prob}
\proof
We first let $(L,\nabla)$ be such that
$\hat c_{1}(L,\nabla)=k\pr_{1}^{*}\hat e\cup \pr_{2}^{*}\hat e$.
We further assume that
$u:T^{2}\to \C/\Z$ is a homomorphism so that
$[u]=m\pr_{1}^{*} \hat e+n\pr_{2}^{*}\hat e$.
Then
$$\hat c_{1}(L,\nabla)\cup [u]= -mk\pr_{1}^{*}(\hat e\cup \hat e)\cup \pr_{2}^{*}\hat e+n k\pr_{1}^{*}\hat e\cup \pr_{2}^{*}(\hat e\cup \hat e)=a([\frac{k(n+m)}{2}ds\wedge dt])\ .$$
The general case now follows by pertubation with forms. \hB

Let $M$ be an oriented connected closed surface and $y\in H^{1}(M;\Z)$. Then
there exists a map
$f:M\to T^{2}$ of degree one and a class
$u\in H^{1}(T^{2};\Z)$ such that $f^{*}u=y$.
Without loss of generality we can assume that $y$ is primitive. Then there exists
a dual element $x\in H^{1}(M;\Z)$ such that
$\langle x\cup y,[M]\rangle=1$ and $(x,y)$ span a hyperbolic summand of the first
cohomology. The classes $x,y$ can be represented by closed integral forms
$\alpha,\beta$. The map
$f:M\to T^{2}$ can then be obtained as the period map associated to $x,y$.
A line bundle bundle $L\to M$ can be written as $L\cong f^{*}H$ for some line bundle
$H\to T^{2}$. 

We can now first calculate
$\ev(\hat c_1(H,\nabla^{H})\cup \hat u)$ for some function $u:T^{2}\to \C/\Z$ representing $u$.
The pairing
$\hat c_1(L,\nabla^{L})\cup \hat y$ is then obtained from $\ev(\hat c_1(H,\nabla^{H})\cup \hat u)$ and a perturbation by forms.

}
\end{ex}

\begin{ex}{\em 
In the following examples we consider the duality pairing in Deligne cohomology induced by the product and evaluation.
We consider a closed oriented $n-1$-dimensional manifold.
We define a pairing
$$\langle\dots,\dots\rangle\colon \hat H^{p}_{Del}(M;\Z)\otimes \hat H^{n-p}_{Del}(M;\Z)\to \C/\Z$$
by
$$\langle x,y\rangle:=\ev(x\cup y)\ .$$
We have the following differential version of Poincar\'e duality. 
\begin{prob}
Show that this pairing is non-degenerated in the sense that
$\langle x,y\rangle=0$ for all $y\in \hat H^{n-p}_{Del}(M;\Z)$ implies $x=0$.
\end{prob}
\proof  See e.g.  \cite{MR2283960} for details. \hB

\begin{prob}
Let $M$ be a two-dimensional and orientable. Show that
$$\cup:\hat H^{1}_{Del}(M;\Z)\otimes \hat H^{1}_{Del}(M;\Z)\to \hat H^{2}_{Del}(M;\Z)$$
is surjective.
\end{prob}
}
\end{ex}





\begin{ex}{\em 

In the following we discuss the compatibility of differential refinements of integral characteristic forms with products. 
Let $\omega$ and $\omega^{\prime}$ be integral characteristic forms for complex vector bundles of degree $n$ and $m$.
Then $\omega\wedge \omega^{\prime}$ is a characteristic form of degree $n+m$.

\begin{prob}
Show that
$$\widehat{\omega\wedge \omega^{\prime}}=\hat \omega\cup \hat \omega^{\prime}\ .$$
Conclude that
$\hat s_4=\hat c_1-2\hat c_2$.
\end{prob}
\proof 
We first consider
the difference
$(\omega\wedge \omega^{\prime})^{\Z}-\omega^{\Z}\cup  \omega^{\prime,\Z}$.
This would be an integral refinement of zero and hence vanishes by the uniqueness of integral refinements.
We now consider the difference
$\widehat{\omega\wedge \omega^{\prime}}-\hat \omega\cup \hat \omega^{\prime}$. This would be a 
differential refinement of zero and again vanishes, this time by the uniqueness of differential refinements.
\hB

\begin{prob}
Show that $\hat \cw$  (see Theorem \ref{may2010}) is additive and multiplicative. More precisely, for a $G$-principal bundle with connection $(P\to M,\omega)$ and
$\tilde\phi,\tilde \psi\in \tilde I^{*}(G)$ we have e.g.
$$\hat \cw(\tilde \phi)(\omega)\cup \hat \cw(\tilde \psi)(\omega)=\hat \cw(\tilde\phi\tilde \psi)(\omega) \ .$$
\end{prob}

Let $(E,\nabla)\cong  (E_1,\nabla_1)\oplus (E_2,\nabla_2)$ be a decomposition of a complex vector bundle with connection.
\begin{prob}
Show that
$$\hat c(\nabla)=\hat c(\nabla_1)\cup \hat c(\nabla_2)\ .$$
\end{prob} 
 
\begin{prob}
Conclude that for a real vector bundle $V\to M$ with connection $\nabla$ we have for all $n\ge 0$ that 
$$2\hat c_{2n+1}(\nabla\otimes \C)=0\ .$$
Show by example that for every $n$ there exists a real bundle $(V,\nabla)$ such that
 $\hat c_{2n+1}(\nabla\otimes \C)\not=0$.
\end{prob}
\proof
Use the splitting principle, \ref{may0211} and \ref{may0212}. \hB

We conclude that for a decomposition of real vector bundles
$$(V,\nabla)\cong (V_1,\nabla_1)\oplus (V_2,\nabla_2)$$ we have
$$\hat p(\nabla)-\hat p(\nabla_1)\cup \hat p(\nabla_2)=\mbox{$2$-torsion}$$
and that the r.h.s is non-trivial in general.

\begin{prob}\label{jun1202}
Let $n_{i}=\dim(V_{i})$ be even for $i=1,2$. Show that
$$\hat p_{\frac{n_{1}}{2}}(\nabla_{1})\cup \hat p_{\frac{n_{2}}{2}}(\nabla_{2})=\hat p_{\frac{n_{1}+n_{2}}{2}}(\nabla)\ .$$
\end{prob}
\proof Use \ref{jun1201} in order to exclude disturbing terms. \hB

\begin{prob}
If $(V,\nabla^{V})$ and $(W,\nabla^{W})$ are real oriented vector bundles with connection, then
we have
$$\hat \chi(\nabla^{V})\cup \hat \chi(\nabla^{W})=\hat \chi(\nabla^{V}\oplus \nabla^{W})\ .$$
\end{prob}
\proof
Show that $\chi(V)\cup \chi(W)=\chi(V\oplus W)$.
Then extend this to the differential refinements using unicity considerations as in similar cases. \hB

\begin{prob}
Show the identity $\hat \chi^{2}=\hat p_{n}$
of differential characteristic classes for real oriented $2n$-dimensional vector bundles with connection.
\end{prob}
\proof
We must show that the corresponding elements in
$\tilde I(SO(2n))$ coincide.  The main point is the equality
$\chi_{\Z}=p_{n}^{\Z}$.
One can use the splitting principle and \ref{jun1202} in order to reduce to the case $n=1$. 
\hB

}
\end{ex}

\subsection{Cheeger-Simons differential characters}

The goal of this section is to relate the Deligne cohomology $\hat H_{Del}^{*}(M;\Z)$ with the predating classical definition of the group of Cheeger-Simons differential characters \cite{MR827262}. 
We let $\sing_{\infty}(M)$ be the simplicial complex of smooth simplices of $M$ and
$C_{*}(\sing_{\infty}(M))$ be the associated chain complex. Integration over simplices induces the de Rham map
$$\Rham:\Omega^{*}_{\C}(M)\to \Hom(C_{*}(\sing_{\infty}(M));\C)\ , \quad\omega\mapsto\{z\mapsto \int_{z}\omega\} \ .$$
A version of the de Rham Lemma says that this map induces an isomorphism
$$\Rham:H^{*}_{dR}(M;\C)\stackrel{\sim}{\to} H^{*}(M;\C)$$
between de Rham cohomology and the smooth singular cohomology.

If $c=\sum_{\sigma\in C_{n}(\sing_{\infty}(M))} n_{\sigma}\sigma$ is a smooth chain, then we define its support by 
$$|c|:=\bigcup_{\{\sigma\in C_{n}(\sing_{\infty}(M))| n_{\sigma}\not=0 \}} \sigma(\Delta^{n})\subseteq M\ .$$ 
It is a compact subset of $M$ which lookes like something at most $n$-dimensional. More precisely, we have:
 \begin{lem}\label{aug0770}
There exists an open neighbourhood $U\subseteq M$ of $|c|$ such that
$H^{k}(U;\Z)=0$ for all $k>n$.
\end{lem}
\proof The argument can be found in \cite[Fact 2.1]{MR2365651}. \hB

We let $Z_{n-1}^{\infty}(M):=Z_{n-1}(C_{*}(\sing_{\infty}(M)))$ be the group of smooth $n-1$-cycles on $M$.  For $z\in Z_{n-1}^{\infty}(M) $ we choose an open neighbourhood $U$
of $|z| $ such that $H^{n}_{dR}(M;\C)=0$. For $\hat x\in \hat H^{n}_{Del}(M;\Z)$ we can find a form $\omega\in \Omega^{n-1}(U;\C)$ such that $a(\omega)=\hat x_{|U}$.
  
\begin{ddd}
We define
$\ev_{\hat x}\in \Hom(Z^{\infty}_{n-1}(M),\C/\Z)$ by 
$$\ev_{\hat x}(z):=[\int_{z}\omega]\ .$$
 \end{ddd}

\begin{lem}\label{apr2210}
This evaluation is well-defined. It induces a homomorphism
 $$\ev:\hat H^{n}_{Del}(M;\Z)\to \Hom(Z^{\infty}_{n-1}(M),\C/\Z)\ .$$
 If $z=\delta c$ for some $c\in C_{n}(\sing_{\infty}(M))$, then
 $\ev_{\hat x}(z)=[\int_{c} R(\hat x)]$.
\end{lem}
\proof
In order to prove  well-definedness the main observation is that if
 $\omega^{\prime}$ is another choice, then the difference $\omega^{\prime}-\omega$ is integral. This implies $[\int_{z}(\omega-\omega^{\prime})]=0$. It is easy to see that we get a homomorphism.
For the last assertion note that $R(\hat x)_{|U}=d\omega$. By Stokes theorem
$\int_{z}\omega=\int_{c} R(\hat x)$. \hB

 \begin{ddd}
A Cheeger-Simons differential character of degree $n>0$ is a  
homomorphism
$$\phi:Z^{\infty}_{n-1}(M)\to \C/\Z$$
such that there exists a form $R(\phi)\in \Omega^{n}_{cl}(M,\C)$ with
$$\phi(\delta c)=[\int_{c} R(\phi)]\ .$$ 
We let
$\hat H^{n}_{CS}(M;\Z)$ denote the group of Cheeger-Simons differential characters of degree $n>0$.
We define $\hat H^{0}_{CS}(M;\Z):=\underline{\Z}(M)$.
\end{ddd}

In Lemma \ref{apr2210}  we have constructed a natural  map 
$$\ev:\hat H_{Del}^{n}(M;\Z)\to \hat H_{CS}^{n}(M;\Z)\ .$$
\begin{lem}
The group of Cheeger-Simons characters fits into the differential cohomology diagram
$$\xymatrix{&\Omega^{n-1}(M;\C)/\im(d)\ar[dr]^{a}\ar[rr]^{d}&&\Omega^{n}_{cl}(M;\C)\ar[dr]&\\\textcolor{red}{H^{n-1}_{dR}(M;\C)}\ar[ur]\textcolor{red}{\ar[dr]}&&\hat H_{CS}^{n}(M;\Z)\ar[ur]^{R}\ar@{->>}[dr]^{I}&&\textcolor{red}{H^{n}_{dR}(M;\C)}\\&
\textcolor{red}{
H^{n-1}(M; \C/\Z )}\ar[rr]^{-Bockstein}\ar@{^{(}->}[ru]&&\textcolor{red}{H^{n}(M;\Z)}\ar[ur]&}
$$
\end{lem}
\proof
We first construct the structure maps.
If $\phi\in \hat H^{n}_{CS}(M;\Z)$, then we first observe that the form $R(\phi)$ is uniquely determined. To this end evaluate $\phi$ on the cycles given by the boundaries of small $n$-simplices.
This gives the map $R$. Using that $Z^{\infty}_{k-1}(M)$ is a free $\Z$-module we choose a lift
$\hat \phi:Z^{\infty}_{k-1}(M)\to \C$.
Then  the cocycle $d \hat \phi -\Rham(R(\phi)) \in Z^{n}(C^{*}(M;\Z))$ represents
$I(\phi)\in H^{n}(M;\Z)$. Finally, for $\alpha\in \Omega^{n-1}(M;\C)$ we define
$a(\alpha)(z):=\Rham(\alpha)(z)$. 
\begin{prob}
Verify the uniqueness of $R(\phi)$ in detail (see \cite{MR827262} for an argument).
Show that $I$ is well-defined.
Further verify the  properties of a differential cohomology diagram.
Show that $\ev$ is compatible with the maps $R,a,I$.
\end{prob}

We have an exact sequence 
$$H^{n-1}(M;\Z)\to \Omega^{n-1}(M, \C)/\im(d)\stackrel{a}{\to} \hat H_{CS}^{n}(M;\Z)\stackrel{(I,R)}{\to}  H^{n}(M;\Z)\times_{H_{dR}^{n}(M;\C)}\Omega^{n}_{cl}(M;\C)\to 0$$ 
which is compatible via $\ev$ with \eqref{apr2211}.
We can thus apply the Five Lemma  and see
\begin{kor}
$$\ev:\hat H^{n}_{Del}(M;\Z)\to \hat H^{n}_{CS}(M;\Z)$$
is an isomorphism.
\end{kor}

\begin{ex}{\em 
Let $M$ be closed, connected, oriented and $n-1$-dimensional. The identification
$\hat H^{n}_{Del}(M;\Z)\cong \C/\Z$ (\ref{apr2230}) is given by 
$\hat x\mapsto \ev_{\hat x}(\{M\})$, where $\{M\}\in Z_{n-1}^{\infty}(M)$ is some representative of the fundamental cycle $[M]$ of $M$.
}\end{ex}
 
Note that in \cite{MR827262} also a product for differential characters was defined. By the uniqueness of products \ref{aug2001} the evaluation  $\ev:\hat H^{n}_{Del}(M;\Z)\to \hat H^{n}_{CS}(M;\Z)$ is automatically multiplicative.

\subsection{Integration}

We now consider a proper submersion  $\pi:W\to M$ of dimension $n$. 

\begin{prob}
Observe that a proper submersion is the same as a locally trivial fibre bundle with closed fibres.
\end{prob}
An orientation of $\pi$ is an orientation of the vertical bundle $T^{v}\pi:=\ker(d\pi)$.  If $\pi$ is oriented, 
then we have an integration of forms $$\int_{W/M}:\Omega^{*}(W)\to \Omega^{*-n}(M)$$
such that $d\int_{W/M}\omega=\int_{W/M}d\omega$. It induces an integration
$$\pi^{Rham}_{!}:H_{dR}^{*}(W;\C)\to H_{dR}^{*-n}(M;\C)\ ,\quad \pi^{Rham}_{!}[\omega]:=[\int_{W/M}\omega]$$ in de Rham cohomology.
Moreover, we have an integration
$$\pi_{!}:H^{*}(W;\Z)\to H^{*-n}(M;\Z)$$
in integral cohomology which is compatible with $\pi^{Rham}_{!}$ via the map $\epsilon_{\C}$.
All these integrations are natural for pull-back diagrams
\begin{equation}\label{may0701}\xymatrix{W^{\prime}\ar[d]^{\pi^{\prime}}\ar[r]^{F}&W\ar[d]^{\pi}\\
M^{\prime}\ar[r]^{f}&M}\end{equation}
and functorial for compositions $$W\stackrel{\pi}{\to} M\stackrel{\kappa}{\to} N$$
of oriented proper submersions.

\begin{ex}{\em 
We assume that the fibres $F$ of $\pi$ are connected and
 consider the Serre spectral sequence associated to the bundle $\pi$ for integral or complex cohomology. Every cohomology class $x\in H^{*}(W;\Z)$ has a symbol  $\sigma(x)  \in E_\infty^{*-n,n}\subseteq E_2^{*-n,n}\cong H^{*-n}(M;\Z) $, where the last identification uses the orientation of the fibres. The integration is given by $\pi_!(x):=\sigma(x)$.
\begin{prob}
 Use this description of the integration in order to prove that $\pi_!$ is compatible with $\Rham$ and $\epsilon_\C$, and verify the the remaining assertions about pull-back and composition made above.
\end{prob}}
\end{ex}

We now define the notion of an integration structure   for the differential cohomology $\hat H^{*}_{Del}(\dots,\Z)$.
\begin{ddd}
An integration  for $\hat H^{*}_{Del}(\dots;\Z)$ is given by maps
$$\hat \pi_{!}:\hat H^{*}_{Del}(W;\Z)\to H_{Del}^{*-n}(M;\Z)$$ for all proper oriented submersions $\pi:W\to M$
such that
$$\xymatrix{\Omega^{*-1}(W,\C)/\im(d)\ar[d]^{\int_{W/M}}\ar[r]^{a}&\hat H^{*}_{Del}(W;\Z)\ar[d]^{\hat \pi_{!}}\ar[r]^{I}\ar@/^1cm/[rr]^{R}&H^{*}(W;\Z)\ar[d]^{\pi_{!}}&\Omega_{cl}^{*}(W;\C)\ar[d]^{\int_{W/M}}\\\Omega^{*-n-1}(M;\C)/\im(d)\ar[r]^{a}&\hat H^{*-n}_{Del}(M;\Z)\ar[r]^{I}\ar@/_1cm/[rr]_{R}&H^{*-n}(M;\Z)&\Omega^{*-n}_{cl}(M;\C) }$$ commutes, where $n=\dim(W)-\dim(M)$. Furthermore, for every diagram
\eqref{may0701} we have
$$\hat \pi_{!} \circ F^{*}=f^{*}\circ \hat \pi_{!}\ .$$
\end{ddd}

\begin{theorem}
There exists a unique integration structure for $\hat H^{*}_{Del}(\dots,\Z)$. It is functorial for compositions.
\end{theorem}
\proof
We construct the integration structure in the picture of Cheeger-Simons characters.
Let $\phi\in \hat H^{k}_{CS}(W;\Z)$. Then we define $\hat \pi_{!}(\phi)\in \hat H^{k-n}_{CS}(W;\Z)$ as follows. Let $z\in Z^{\infty}_{k-n-1}(M)$. We choose an open neighbourhood of $U$ of $|z|$  which is  homotopy equivalent to a $k-n-1$-dimensional complex. Then $\pi^{-1}(U)$ is equivalent to a 
$k-1$-dimensional complex. In particular there exists a form $\alpha\in \Omega^{k-1}(\pi^{-1}(U);\C)$ such that $a(\alpha)=\phi_{|\pi^{-1}(U)}$. By naturality of the integration and its compatibility with $a$ 
we are forced to define
$$\hat \pi_{!}(\phi)(z):=[\int_{z}
 \int_{\pi^{-1}(U)/U} \alpha]\ .$$
\begin{prob}
Show that this construction is compatible with $a$, $R$, pull-backs, and that it is functorial.
Further show that the integration is compatible with Mayer-Vietoris sequences induced by decompsitions of the base.
\end{prob}
In order to show compatibility with $I$ we argue as follows. The integration $\hat \pi_{!}$
induces an integration $\pi_{!}^{\prime}$ in integral cohomology which is compatible with pull-back along diagrams \eqref{may0701} and complexification $\epsilon_{\C}$. It coincides with $\hat \pi_{!}$ if $H^{*}(M;\Z)$ is torsion-free, in particular if $M$ is a point.
\begin{prob}
Show that there is at most one functorial integration for integral cohomology which is 
compatible with pull-back along diagrams \eqref{may0701}, Mayer-Vietoris sequences associated to decompositions of the base, and which coincides with  
$\pi_{!}$ for $M=\pt$.
\end{prob}
By this exercise we have  $\pi_{!}^{\prime}=\pi_{!}$.
\hB 

Alternatively, the existence of an integration for $\hat H^{*}_{Del}(\dots,\Z)$ follows from \ref{jun0740}.

\begin{prob}
Show the projection formula
$$\hat \pi_{!}(\pi^{*}x\cup y)=x\cup \hat \pi_{!}(y)\ , \quad x\in \hat H_{Del}^{*}(M;\Z)\ ,\quad y\in \hat H_{Del}^{*}(W;\Z)\ .$$
\end{prob}
\proof
This is not easy. See \cite{MR2179587} or \cite{MR2674652} for arguments based on different models. Alternatively, this follows from \ref{jun1401}, 2.\hB 

\begin{ex}{\em 
Let $(P,\nabla^{P})$ be the Poincar\'e bundle over $J(S^{1})\times S^{1}$.
\begin{prob}
Calculate the integrals of $\hat c_{1}(\nabla^{P})$ over the two projections.
\end{prob}
\proof
Note that these integrals are classes in $\hat H^{1}_{Del}(\dots;\Z)$ and thus $\C/\Z$-valued functions.
\hB 
Consider now a surface $M$ of genus $g$ and again the Poincar\'e bundle $(P,\nabla^{P})$ over $J(M)\times M$.
\begin{prob}
Calculate the integral of $\hat c_{1}(\nabla^{P})^{g+1}$ along the projection $J(M)\times M\to M$.
\end{prob}
\proof
Note that the result can be interpreted as an element of $\Line_{\nabla}(M)$.
The question is to characterize this isomorphism class of line bundles. \hB 
}
\end{ex}

\begin{ex}{\em
Let $p:W\to B$ be an oriented bundle of compact manifolds with boundary $\partial W\to B$. Let $q:\partial W\to B$ denote the restriction of $p$ to the boundary.
We consider   $x\in \hat H_{Del}^{n}(W;\Z)$.
\begin{prob}
Show that $$q_{!}(x)=a(\int_{W/B} R(x))\ .$$
\end{prob}
}\end{ex}

\section{Differential extensions of generalized cohomology theories}\label{aug2040}

\subsection{The $\infty$-categorical black-box}

In the following we will use some $\infty$-categorial language as a black box. General references are \cite{MR2522659} and \cite{highalg}. 
Our philosophy is to work, as much as possible, in a model independent way. We try to hightlight the places where models are required for calculations. 
The $\infty$-categorical black box provides the set-up and   the space where our objects live. It ensures correctness of the constructions. 

The input and output of calculations
belong to the classical world, and the dependence of the output from the input does not depend on the details of interior machine. In particular, as  we try to emphasize,  this relation can be revealed in praxis just using the 
formal properties of the $\infty$-categorical machine.

First of all, we use the word $\infty$-category as an abbreviation for $(\infty,1)$-category.

If $C$ is a category, then its nerve $N(C)$ is an $\infty$-category.

If $\bC$ is an $\infty$-category and $W$ is a collection of morphisms, then we can form
its localization $\iota:\bC\to \bC[W^{-1}]$. It is characterized by an obvious universal property that
$$\Fun_{W^{-1}}(\bC,\bD)\cong \Fun(\bC[W^{-1}],\bD)\ ,$$
where $\Fun_{W^{-1}}(\bC,\bD)$ denotes the full subcategory of functors which map morphisms from $W$ to equivalences.
If $X$ is an object of $\bC$, then usually we use the same symbol in order to denote  the corresponding object in $\bC[W^{-1}]$. If it is confusing not to distinguish, then we use the notation  ${}_\infty X$ or $\iota(X)$ for the latter.

Let $\sSet$, $\Top$,  $\Sp$ and $\Ch$ denote the categories of simplicial sets, topological spaces, spectra or chain complexes.
Then $\Nerve(\sSet)[W^{-1}]$, $\Nerve(\Top)[W^{-1}]$, $\Nerve(\Sp)[W^{-1}]$ or $\Nerve(\Ch)[W^{-1}]$ denote the localizations of the associated $\infty$-categories at weak equivalences, stable equivalences or quasi-isomorphisms, respectively. All these $\infty$-categories have symmetric monoidal versions, the first with respect to the cartesian product, the second with respect to $\wedge$, and the last with respect to $\otimes$.

For two objects $X,Y\in \bC$ we let $\map(X,Y)\in \Nerve(\sSet)[W^{-1}]$ be the mapping space.
The categories $\Nerve(\Sp)[W^{-1}]$ or $\Nerve(\Ch)[W^{-1}]$ are stable
and therefore have mapping spectra which will be denoted by 
$\Map(X,Y)$. Note that $\map(X,Y)\cong \Omega^{\infty} \Map(X,Y)$. 

These $\infty$-categories a complete and cocomplete. They are hence tensored and cotensored over $\Nerve(\sSet)[W^{-1}]$, or   the equivalent $\infty$-category $\Nerve(\Top)[W^{-1}]$.
We will write the tensor of an object $C$ with a space $X$ by $C\otimes X$, and the cotensor will be denoted by $C^{X}$.  The stable $\infty$-categories $\Nerve(\Ch)[W^{-1}]$ and $\Nerve(\Sp)[W^{-1}]$ are tensored and cotensored over spectra  $\Nerve(\Sp)[W^{-1}]$, and we use a similar  notation for these structures with a subscript $s$ added. Note that for a space $X$ we have $C\otimes X\cong C\otimes_{s}\Sigma^{\infty}_{+}X$ and $C^{X}\cong C^{{}_{s}\Sigma^{\infty}_{+}X}$.

We will often use the following argument.
Let
$\Phi$ and $\Psi$ be  colimit preserving functors $ \Fun(\Nerve(\Top)[W^{-1}],\bC)$.
\begin{lem}\label{may1510}
An transformation (equivalence)
$\Phi(\pt)\to \Psi(\pt)$ extends naturally to a transformation (equivalence)
$\Phi\to \Psi$. 
There are corresponding versions for contravariant functors which map colimits to limits.
\end{lem}
\proof
Indeed, for every space $X\in \Nerve(\Top)[W^{-1}]$ we can form its category of simplices
$\Simp(X)$ and get a natural equivalence
$$\colim_{\Nerve(\Simp(X))} \pt\stackrel{\sim}{\to} X\ .$$  
We obtain
$$\Phi(X)\cong \colim_{\Nerve(\Simp(X))} \Phi(\pt)\to  \colim_{\Nerve(\Simp(X))} \Psi(\pt)\cong \Psi(X)\ .$$
\hB

\begin{ex}{\em 
 \begin{prob}\label{jun2502} Show that the tensor and cotensor for spectra can be written in terms of the smash product and internal mapping object $\Map$  as follows:
$$E\otimes X\cong E\wedge \Sigma^{\infty}_+X\ , \quad E^{X}\cong \Map( \Sigma_+^{\infty}X,E)\ , \quad X\in \Nerve(\sSet)[W^{-1}]\ , E\in \Nerve(\Sp)[W^{-1}]\ .$$
\end{prob}
\proof
We discuss the cotensor.
Observe that both sides  map colimits in the $X$-variable to limits and are canonically equivalent for $X=\pt$.
 Lemma \ref{may1510} gives the assertion. \hB

\begin{prob}\label{jun2504}
Show that the tensor and cotensor structure with $\Nerve(\sSet)[W^{-1}]$ on $\Nerve(\Ch)[W^{-1}]$ is given in terms of the chain complex $C_*(X)$ of $X\in \Nerve(\sSet)[W^{-1}]$ by
$$A\otimes X=A\otimes C_*(X)\ , \quad A^{X}\cong \underline{\Map}(C_*(X),A)\ , \quad A\in 
\Nerve(\Ch)[W^{-1}]\ ,$$ where $\underline{\Map}$ denotes the mapping chain complex.
\end{prob}
\proof
We discuss the tensor. Both sides coincide for $X=*$ and preserve colimits. \hB

For a chain complex $A\in \Nerve(\Ch)[W^{-1}]$ let $H_k(A)$ denote the homology group in degree $n$.

\begin{prob}\label{jun2503}
Show that for a space $X\in \Nerve(\sSet)[W^{-1}]$, and integer  $k\ge 0$ and a chain complex $A\in \Nerve(\Ch)[W^{-1}]$ we have
\begin{equation}\label{jun2501}H_{0}(A)\cong \pi_{0}(\map(\Z,A))\end{equation} and  more generally
$$H_{k}(A^{X})\cong \pi_{k}(\map(C_{*}(X),A))\ .$$
\end{prob} 
\proof
We assume the first equivalence \eqref{jun2501} and deduce the second.  We get
\begin{eqnarray*}
H_{k}(A^{X})&\cong& \pi_{0}(\map(\Z[k],A^{X}))\\
&\cong&\pi_{0}(\map(\Z[k]\otimes X,A))\\
&\stackrel{\ref{jun2504}}{\cong}&\pi_{0}(\map(C_{*}(X)[k],A))\\
&\cong&\pi_{k} (\map(C_{*}(X),A))\ .
\end{eqnarray*}
In order to show \eqref{jun2501} it seems to be necessary to use a model for the $\infty$-category
$\Nerve(\Ch)[W^{-1}]$. \hB }
\end{ex}

Let  $\bC$ be an $\infty$-category. 
\begin{ddd}
We call  $\Sm(\bC):=\Fun(\Nerve(\Mf^{op}),\bC)$ the $\infty$-category of smooth objects in $\bC$.
\end{ddd}

Recall that $\Mf$ has a topology given by open coverings.
\begin{ddd} We say that a smooth object $X\in \Sm(\bC)$ satisfies descent, if for every open covering $\cU$ of $M$ the augmentation map
$$X(M)\to \lim_{\Nerve(\Delta)} X(U^{\bullet}) $$
is an equivalence in $\bC$, where $U^{\bullet}\in \Mf^{\Delta^{op}}$ is the nerve of $\cU$.
\end{ddd}
For the following we make the technical assumption that $\bC$ is  presentable. This holds true for all examples considered here.
We let $\Sm^{desc}(\bC)\subseteq \Sm(\bC)$ denote the full subcategory of objects satisfying descent. Then by \cite[6.2.2.7]{MR2522659} we have an adjunction
$$L:\Sm(\bC)\leftrightarrows \Sm^{desc}(\bC)\ .$$
The functor $L$ is called sheafification.
In general it seems very difficult to evaluate sheafifications explicitly. See Example \ref{jul1702} for   a method  which can be applied in some interesting cases.

   \begin{ex}\label{jul1702}{\em
Let $f:X\to Y$ be a morphism in $\Sm(\bC)$. Sometimes we need a criterion ensuring  that
$L(f):L(X)\to L(Y)$ in $\Sm^{desc}(\bC)$  is an equivalence (see Lemma \ref{jun3005}).   We define the functor
$\cL:\Sm(\bC)\to \Sm(\bC)$ which on objects acts as 
$$\cL(X)(M):=\colim_{\cU}\lim_{\Nerve(\Delta)} X(\cU^{\bullet})\ ,$$
where the first colimit is over the filtered system of open coverings of $M$. 
\begin{prob}
Describe this  precisely as a functor between $\infty$-categories.
\end{prob}
 \proof
We consider the functor $\Mf^{op}\to\Filt$
which associates to every manifold $M$ the filtered partially ordered 
set of open coverings of $M$ indexed by the points of $M$. Then we let $\tilde \Mf$ be the transport category of pairs
$(M,\cU)$ of manifolds and coverings  whose morphisms $(M,\cU)\to (M^{\prime},\cU^{\prime})$ are smooth maps $f:M\to M^{\prime}$ such that $ U_{m}\subseteq f^{-1}(U^{\prime}_{f(m)}) $ for all $m$.  We have the \v{C}ech nerve functor
$\tilde \Mf\to \Mf^{\Delta^{op}}$, $(M,\cU)\mapsto \cU^{\bullet}$. Precomposition with this functor gives the functor $\tilde \cL:\Sm(\bC)\to \Fun(\Nerve(\tilde\Mf)^{op},\Fun(\Nerve(\Delta),\bC))$.
Finally we compose with $\lim_{\Delta}$ and get
$$\Sm(\bC)\to  \Fun(\Nerve(\tilde\Mf)^{op},\bC)\ .$$
We have an adjunction
$$\Phi_{*}:\Fun(\Nerve(\tilde \Mf)^{op},\bC)\leftrightarrows \Fun(\Nerve(\Mf)^{op},\bC):\Phi^{*}\ ,$$
where $\Phi:\tilde \Mf\to \Mf$ forgets the coverings. Here $\Phi_{*}$ is the Kan-extension functor.
We define $\cL:=\Phi_{*} \circ \tilde \cL$.
By the pointwise formula for the Kan extension functor we have
$$\cL(X)(M)\cong \colim_{(( N,\cU),M\to N)\in \tilde \Mf^{op}_{/M}}\lim_{\Nerve(\Delta)} X(\cU^{\bullet})\ .$$
Now
$$\tilde \Mf^{op}_{/M}=(\tilde \Mf_{M/})^{op}$$
contains the category of coverings of $M$ as a cofinal subcategory. 
Hence it suffices to take the colimit over the coverings of $M$ and we get 
$$\cL(X)(M)=\colim_{\cU} \lim_{\Nerve(\Delta )} X(\cU^{\bullet})\ .$$
 \hB

There is a natural morphism $\id\to \cL$.
We define
$$\cL^{\infty}:=\colim (\id\to \cL\to \cL^{2}\to  \cL^{3}\to \dots):\Sm(\bC)\to \Sm(\bC)\ .$$
By construction the natural morphism  $\cL^{\infty}\to (\cL^{\infty})^{2}$ is an equivalence.
We let $$\Sm^{\cL^{\infty}}(\bC)\subseteq \Sm(\bC)$$  be the essential image of $\cL^{\infty}$. By the recognition principle
\cite[Prop. 5.2.7.4]{MR2522659}, since condition 3. is satisfied, this is a localization. If $F$ satisfies descent, then we observe that $F\to \cL^{\infty} F$ is an equivalence.
Hence we have a sequence of localizations
$$\Sm^{desc}(\bC)\subseteq \Sm^{\cL^{\infty}}(\bC)\subseteq \Sm(\bC)\ .$$ In particular, the natural transformation
\begin{equation}\label{jun3001}L\to L\circ \cL^{\infty}\end{equation} is an equivalence.
\textcolor{red}{\begin{prob}\label{aug0760}
Do we have
$\Sm^{desc}(\bC)\cong \Sm^{\cL^{\infty}}(\bC)¬¨‚Ä†?$
\end{prob}
}

\begin{lem} \label{jun3005}
If $f:X\to Y$ is a morphism in $\Sm(\bC)$  such that $\cL(f):\cL(X)\to \cL(Y)$ is an equivalence,
then $L(f):L(X)\to L(X)$ is an equivalence.
\end{lem}
\proof
The condition implies that
$\cL^{\infty}(f)$ is an equivalence. Now by \eqref{jun3001}
we have an equivalence $L(\cL^{\infty}(f))\cong L(f)$ so that $L(f)$ is an equivalence, too. 
\hB

Here is an application of Lemma \ref{jun3005}. 
Let $C(\Omega_{cl}^{n})\in \Sm(\Nerve(\Ch)[W^{-1}])$ be the functor which associates to every manifold $M$ the complex which is concentrated in degree zero and given there by $\Omega^{n}_{cl}(M)$.
Since $(\sigma^{\ge n} \Omega)[n]$ satisfies descent the natural inclusion map $C(\Omega^{n}_{cl})\to (\sigma^{\ge n} \Omega)[n]$ extends to $i:L(C(\Omega^{n}_{cl}))\to (\sigma^{\ge n} \Omega)[n]$.
\begin{prob}\label{jun3006}
Show that $i$ is an equivalence.
\end{prob}
\proof
We check that
$\cL( C(\Omega^{n}_{cl}))\to \cL(  (\sigma^{\ge n} \Omega)[n])\cong (\sigma^{\ge n} \Omega)[n] $ is an equivalence and apply Lemma \ref{jun3005}.
Indeed, if $\cU$ is a good covering of $M$, then
$$\check{C}(\cU ,C(\Omega^{n}_{cl}))\to \check{C}(\cU ,(\sigma^{\ge n} \Omega)[n] )$$
is an equivalence. Since good coverings are cofinal in all coverings we get the desired result. 
\hB

}
\end{ex}



We have a functor
$$t:\Nerve(\Mf)\to \Nerve(\Top)[W^{-1}]$$ 
which associates to smooth manifold its underlying topological space.
It  satisfies codescent
in the sense that for every open covering $\cU$ of a manifold $M$ we have  an equivalence
$$\colim_{\Nerve(\Delta)^{op}} t(U^{\bullet})\to t(M)\ .$$

We now assume that the $\infty$-category $\bC$ is cotensored  over $\Nerve(\Top)[W^{-1}]$.   

\begin{ddd}
We define the smooth function object functor 
$$\Funk:\bC\to \Sm^{desc}(\bC)$$ by
$$\Funk(X)(M):= X^{t(M)}\ .$$
\end{ddd}
More formally we should write the functor as the adjoint of the composition
$$
\bC \times \Nerve(\Mf)^{op}\stackrel{\id\times  t}{\to} \bC\times  \Nerve(\Top)[W^{-1}]^{op}\stackrel{\dots^{\dots}}{\to}\bC\ .$$
Since $t$ satisfies codescent we see that $\Funk(X)$ indeed satisfies descent.

\begin{ex}{\em 
We consider a spectrum $E\in \Nerve(\Sp)[W^{-1}]$.
\begin{lem}
We have
$$\pi_{k}(\Funk(E)(M))=E^{-k}(M)\ .$$
\end{lem}
\proof

We have
\begin{eqnarray*}
\pi_{k}(\Funk(E)(M))&\stackrel{def}{=}&\pi_k (E^{t(M)})\\
&\stackrel{\ref{jun2502}\cong }{\to}&\pi_{k}(\Map( \Sigma_{+}^{\infty}t(M),E)\\
&\cong& \pi_{0}(\Map( \Sigma_{+}^{\infty+k}t(M),E)\\
&\stackrel{def}{=}&E^{-k}(M)\ .
\end{eqnarray*}
\hB
}\end{ex}

\begin{ex}{\em
We consider a complex of sheaves $A\in \Ch(\Sh_\Ab(\Mf))$ as an object of
$\Sm(\Nerve(\Ch))$.

\begin{prob}\label{jun2801}
Show that $A\in \Sm(\Nerve(\Ch))$   satisfies descent. Furthermore, show by example that its image  ${}_\infty A\in   \Sm(\Nerve(\Ch)[W^{-1}])$ may not satisfy descent.
Finally show that if $A$ is a complex of fine sheaves, then
${}_\infty A$ satisfies descent.
\end{prob}
\proof
Descent in $\Sm(\Nerve(\Ch))$ is equivalent to the degree-wise sheaf condition.
 
In the case of $\Sm(\Nerve(\Ch)[W^{-1}])$,
 for an open covering $\cU$ of $M$,   the desired limit is represented by the \v{C}ech complex (see Problem \ref{jun2401} for an argument)
$$\check{C}(\cU,  A)\cong  \lim_{\Nerve(\Delta)} {}_\infty A(U^{\bullet})\ .$$  
Let $A=\underline{\Z}$. Then ${}_\infty A$   does not satisfy descent.
To this end consider a good (i.e. all multiple intersections are empty or contractible) covering of a manifold with higher-degree integral cohomology which is calculated by   $\check{C}(\cU,  \underline{\Z} )$.
Of course, ${}_\infty \underline{\Z} (M)$ lives in degree $0$.

Since the \v{C}ech complex of a fine sheaf is acyclic (Lemma \ref{apr2104}) we get descent in this case.
\hB
}
\end{ex}

 \begin{ex}\label{may2905} 
{\em 
 Let $*$  denote the final category. It has one object and one morphism. There is a natural
functor $p:\Mf\to *$. We have an equivalence $\Fun(N(*),\bC)\cong \bC$.
\begin{prob}\label{jul1701}
Let $X\in \bC$ and $p^{*}X\in \Sm(\bC)$ be the constant smooth object.
Show that $L(p^{*}X)(*)\cong X$.
\end{prob}
\proof
We have a natural morphism
$$c:p^{*}X\to \Funk(X)\ .$$
We use homotopy invariance of
$p^{*}X$ and $\Funk(X)$ to show that the induced map
$\cL(p^{*}X)\to \cL(\Funk(X))$
is an equivalence.
This implies by \ref{jun3005} that
$L(p^{*}X)\to L(\Funk(X))$
is an equivalence.
Since $\Funk(X)$ satisfies descent 
we have $L(\Funk(X))(*)\cong \Funk(X)(*)\cong X$.
\hB

\begin{ddd}
A smooth object $X\in \Sm^{desc}(\bC)$ is called constant if there exists an object $X_{*}\in \bC$ such that
$L(p^{*}X_{*})\cong X$.
\end{ddd}
It is clear by \ref{jul1701} that 
$X_*$ must be the evaluation of $X$ on the manifold $*$.
The following reproduces a result of Dugger, \cite{MR1870515}.
\begin{prob}\label{may2908}
Show that a smooth space $X\in \Sm^{desc}(\Nerve(\Top)[W^{-1}])$ is constant if and only if $
X=\Funk(X_{*})$ for some space $X_{*}\in \Nerve(\Top)[W^{-1}]$.
\end{prob}
\proof
 There is a canonical map $L(p^{*}X(*))\to \Funk(X(*))$. Show that it is an equivalence.
 Show that both $L(p^{*}X(*))$ and $\Funk(X(*))$, are homotopy invariant (see Problem  \ref{jun2511} for the first).
 Then  use coverings and descent in order to reduce to the problem to show that
 $L(p^{*}X(*))(U)\to \Funk(X(*))(U)$ is an equivalence for all contractible $U$.
Finally, this holds true since by homotopy invariance both  evaluations are equivalent to $X(*)$.
 \hB
}
\end{ex}

\begin{ex}{\em
The following exercise provides a tool for the explicit computation of some limits in $\Nerve(\Ch)[W^{-1}]$.
We consider a cosimplicial chain complex $A\in \Fun(\Delta ,\Ch)$ and its version
${}_{\infty}A\in \Fun(\Nerve(\Delta),\Nerve(\Ch)[W^{-1}])$.  
We let $\tot(A)\in \Ch[W^{-1}]$ be the chain complex obtained as total complex of the double complex
associated to $A$ by normalizing the simplicial direction.
\begin{prob}\label{jun2401}
Show that there is an equivalence
$$\lim_{\Nerve(\Delta)} {}_{\infty}A\cong  \tot(A)\ .$$
\end{prob}
\proof 
We first show that
$\tot(A)\cong R\lim_\Delta(A)$.
We use the projective model category on $\Ch$ in order to 
present $\Nerve(\Ch)[W^{-1}]$. In order to present 
$\Fun(\Nerve(\Delta),\Nerve(\Ch))[W^{-1}]$ we then use the injective
model category structure on $\Ch^{\Delta}$. We consider the object
$C_*(\Delta^{\bullet})\in \Ch^{\Delta}$ and form for $A\in \Ch^{\Delta}$ the internal 
$\underline{\Hom}(C_*(\Delta^{\bullet}),A)\in \Ch^{\Delta} $.
We write $\underline{\hom}(A,B):=\lim_\Delta \underline{\Hom}(A,B)$.
By an explicit calculation we observe that
$$\tot(A)= \underline{\hom}(C_*(\Delta^{\bullet}),A)  \ .$$
With the chosen model category  structures the functor $$\underline{\hom}:(\Ch^{\Delta})^{op}\times \Ch^{\Delta}\to \Ch$$ is bi-Quillen.
The projection $\Delta^{\bullet}\to *$ induces the map
$C_*(\Delta^{\bullet})\to \underline{\Z}$, and dually 
 \begin{equation}\label{jul0301}\lim_\Delta A\cong \underline{\hom}( \underline{\Z},A)\to \underline{\hom}(C_*(\Delta^{\bullet}),A)\cong \tot(A)\ . \end{equation}

A level-free lower bounded chain complex is projectively cofibrant in $\Ch$.
Hence $\underline{\Z}$ and $C_*(\Delta)$ are injectively cofibrant in $\Ch^{\Delta}$. 
It follows that \eqref{jul0301} is an equivalence if $A$ is fibrant. Moreover, $\tot$ preserves equivalences.
 Let $A\to R(A)$ be a fibrant replacement. 
We get a chain of equivalences
$$R\lim_\Delta(A)\stackrel{def}{=}  \lim_\Delta R(A)\stackrel{\eqref{jul0301}}{ \cong }
\tot(R(A))\cong \tot(A)\ .$$

The final step is to observe that the natural map
 $$ \Fun(\Nerve(\Delta),\Nerve(\Ch))[W^{-1}]\to \Fun(\Nerve(\Delta),\Nerve(\Ch)[W^{-1}])$$
 induces an equivalence which maps $ A $ to ${}_{\infty}A$ such that
 $$\lim_{\Nerve(\Delta)} {}_{\infty}A\cong R\lim_{\Delta} (A)\ .$$
 This rigidification result uses the fact that the class of  quasi isomorphisms  $W$ is the class of weak equivalences of a combinatorial model category structure on $\Ch$ (see \cite[Prop. 1.3.3.12]{highalg}).
  \hB

The corresponding problem for the colimit of a simplicial chain complex $A\in \Fun(\Nerve(\Delta^{op}),\Nerve(\Ch))$ is much simpler. We define
$$\tot(A)^{n}:=\bigoplus_{p-q=n}A^{p}([q])$$ and the differential
$d=(-1)^{q}d_{A}+\sum_{i=0}^{q}\partial_{i}^{*}$ on the summand $A^{p}([q])$.
\begin{prob}
Show that there is an equivalence
$$\colim_{N(\Delta)}{}_{\infty} A\cong \tot A\ .$$
\end{prob}
\proof
We again use
the equivalence
$$\Fun(\Nerve(\Delta^{op}),\Nerve(\Ch))[W^{-1}]\cong \Fun(\Nerve(\Delta^{op}),\Nerve(\Ch)[W^{-1}])\ .$$
The left-hand side can be obtained from a model category structure on $\Ch^{\Delta^{op}}$. We let $Q$ denote the cofibrant replacement functor. 
Note that $\tot$ preserves quasi-isomorphisms.
Then we have a sequence of maps
$$\tot\stackrel{\sim}{\leftarrow} \tot\circ Q\to\colim_{\Delta^{op}} \circ Q\stackrel{def}{=}L\colim_{\Delta^{op}}\to \colim_{\Delta^{op}} $$
of colimit preserving functors. Now the constant diagram $c(\Z)$ is cofibrant, and 
$\tot(c(\Z))\cong \colim_{\Delta^{op}}(c(\Z))$ by inspection.
Since every object of $\Ch^{\Delta^{op}}$ can be written as a homotopy colimit 
of a diagram of $c(\Z)$ we see that
$\tot\cong L\colim_{\Delta^{op}}$. Finally we identify $L\colim_{\Delta^{op}}$ 
with the colimit in $\Fun(\Nerve(\Delta^{op}),\Nerve(\Ch))[W^{-1}]$ in the $(\infty,1)$-categorial sense by showing that  the latter has required universal property of a left-derived functor.
\hB

}\end{ex}

 \begin{ex}\label{jun1901}{\em
The following is similar to a theory developed in a different context  in \cite{MR1813224}.
Let $\bC$ be a presentable $\infty$-category and $I:=[0,1]\in \Mf$ be the unit interval.
\begin{ddd}
We say that $X\in \Sm(\bC)$ is homotopy invariant, if
$\pr^{*}:X(M)\to X(I\times  M)$ is an equivalence for all manifolds $M$. 
\end{ddd}We let
$\Sm^{I}(\bC)\subseteq \Sm(\bC)$ denote the full subcategory of homotopy invariant objects.
We have an adjunction
\begin{equation}\label{jun2510}\cI:\Sm(\bC)\leftrightarrows \Sm^{I}(\bC)\ .\end{equation}

Let $i_{0},i_{1}:*\to I$ be the inclusions of the endpoints of the interval. 
\begin{ddd}
 We say that $X\in \Sm(\bC)$ is elementary homotopy invariant, if for every manifold  $M$ the morphisms
 $i_{0}^{*},i_{1}^{*}:X(I\times M)\to X(M)$ are equivalent.
 \end{ddd}
 
 \begin{prob}
 Show that elementary homotopy invariance is equivalent to homotopy invariance.
 \end{prob}
 \proof
 Assume that $X$ is homotopy invariant.
 Then the compositions
 $$X(M)\stackrel{\pr^{*}}{\cong} X(I\times M)\stackrel{i_{j}^{*}}{\to} X(M)\ ,\quad  j=0,1$$
 are both equivalent to $\id_{X(M)}$. Since the first map is an equivalence we conclude
 that $i_{0}^{*}$ and $i_{1}^{*}$ are equivalent.
 
Now assume that
$X$ is elementary homotopy invariant.
The composition
$$X(M)\stackrel{\pr^{*}}{\to} X(I\times M)\stackrel{i_{0}^{*}}{\to} X(M)$$
is equivalent to the identity. We must show that
$$X(I\times M)\stackrel{i_{0}^{*}}{\to} X(M) \stackrel{\pr^{*}}{\to} X(I\times M)$$
is also equivalent to the identity. To this end we consider  
the map $\mu:I\times I\to I$, $\mu(s,t):=st$, and the diagram
$$\xymatrix{I\ar[dr]^{i_{0}\times \id_{I}}\ar[r]&{}*\ar[rd]^{i_{0}}&\\
&I\times I\ar[r]^{\mu}&I\\
I\ar[ur]^{i_{1}\times \id_{I}}\ar[rru]_{\id_{I}}&&}\ .$$
We insert this diagram into $X(\dots\times M)$ in order to obtain 
the desired result. 
\hB

 We define a functor
$$\bar \bs:\Sm(\bC)\to \Sm(\bC)$$ as the composition
of
$$\bs :\Sm(\bC)\to \Sm(\Fun(\Nerve(\Delta^{op}),\bC))$$
and
$$\colim_{\Nerve(\Delta^{op})}: \Sm(\Fun(\Nerve(\Delta^{op}),\bC))\to \Sm(\bC)\ ,$$
where $\bs$ is  precomposition by
$$\Mf\to \Fun(\Delta,\Mf), \, \quad M\mapsto \Delta^{\bullet}\times M \ .$$ 
We have a natural map of cosimplicial manifolds
$\pr:\Delta^{\bullet}\times M\to M$ (the target is the constant cosimplicial manifold)  which induces
a morphism $\id\to \bar \bs$.

\begin{prob}
If $X$ is homotopy invariant, then $X\to \bar \bs(X)$ is an equivalence.
\end{prob}
\proof
We show that the map of simplicial objects
$X(M)\to X(\Delta^{\bullet}\times M)$ induced by the projection is an equivalence.
Hence we must show that $\pr^{*}:X(M)\to X(\Delta^{n}\times M)$ is an equivalence for every $M$ and $n\ge 0$. 

Let $i_{(0,\dots,0)}:*\to \Delta^{n}$ be the inclusion of the zero corner.
It suffices to show that $(i_{(0,\dots,0)}\times \id_{M})^{*}$ is an inverse of
$\pr^{*}$. The non-trivial part is the verification that  
$\pr^{*} \circ i_{(0,\dots,0)}^{*}\cong \id_{X(\Delta^{n}\times M)}$.

We let $\mu:I\times \Delta^{n}\to \Delta^{n}$ be the map
$$(s,(0\le t_{0}\le  \dots \le t_{n}\le 1))\mapsto (0\le st_{0}\le  \dots\le st_{n}\le 1)\ .$$
Then
we have a commutative diagram
$$\xymatrix{\Delta^{n}\ar[dr]^{i_{0}\times \id_{\Delta^{n}}}\ar[r]&{}*\ar[rd]^{i_{(0,\dots,0)}}&\\
&I\times \Delta^{n}\ar[r]^{\mu}&\Delta^{n}\\
\Delta^{n}\ar[ur]^{i_{1}\times \id_{\Delta^{n}}}\ar[rru]_{\id_{\Delta^{n}}}&&}\ .$$
If we insert this into $X(\dots\times M)$ and use 
that $X$ is elementary homotopy invariant, we get the desired equivalence. \hB

 \begin{prob} 
 If a morphism $f:X\to Y$ in $\Sm(\bC)$ induces an equivalence $\bar \bs (f):\bar \bs(X)\to \bar \bs(Y)$, then
 $\cI(f):\cI(X)\to \cI
(Y)$ is an equivalence.
\end{prob}
\proof
This is similar to \ref{jun3005}.
We define $$\bar \bs^{\infty}:=\colim (\id\to \bar \bs\to \bar\bs^{2}\to \dots)\ .$$ 
Then $\bar \bs^{\infty }\to (\bar \bs^{\infty})^{2}$ is an equivalence. Use  \cite[Prop. 5.2.7.4]{MR2522659} (Condition 3.) in order to   see that
the essential image $\Sm^{\bar \bs^{\infty}}(\bC)$ of $\bar \bs^{\infty}$
is a localization 
$$\bar \bs^{\infty}:\Sm^{\bar \bs^{\infty}}(\bC)\subseteq \Sm(\bC)\ .$$ 
Since $\bar \bs^{\infty}$ preserves homotopy invariant objects we have
a chain of localizations
$$\Sm^{\cI}(\bC)\subseteq\Sm^{\bar \bs^{\infty}}(\bC)\subseteq \Sm(\bC)$$
and therefore $\cI\circ \bar \bs^{\infty}\cong \cI$.
If $\bar \bs(f)$ is an equivalence, then so is
$\bar \bs^{\infty}(f)$ and hence $\cI\circ \bar \bs^{\infty}(f)\cong \cI(f)$. 
\hB

\textcolor{red}{  \begin{prob}\label{aug0761}Is
 $\Sm^{\cI}(\bC)\cong\Sm^{\bar \bs^{\infty}}(\bC)?$\end{prob}}
 
There seems to be a great simplification in the case where $\bC=\Ch[W^{-1}]$ or $\bC=\sSet[W^{-1}]$.
\begin{prob}\label{jun3001n1}
Assume that  $X\in \Sm(\Ch[W^{-1}])$ or $\bC=\sSet[W^{-1}]$.  Is  $\bar \bs (X)$  homotopy invariant?
\end{prob}
\proof
For chain complexes one can probably use the argument given in \cite[Lemma 2.18]{MR2242284}. For simplicial sets one can construct a map
$$X(I\times \Delta^{\bullet}\times M)\to X(\Delta^{\bullet}\times M)^{\Delta^{1}}$$ 
which gives the desired homotopy of the restrictions to the points $0,1\in I$ after realization.
\hB 
 
If  \ref{jun3001n1} has a positive answer it implies similar statements for $\bC=\Mod(H\Z)$ or $\bC=\Sp[W^{-1}]$.


%
 


\begin{prob}\label{jun2511}
Show that for $X\in \bC$  the object $L(p^{*}\bC)$ is homotopy invariant.
\end{prob}
\proof
 We have a diagram
$$\xymatrix{\Sm^{I,desc}(\bC)\ar[r]\ar[d]&\Sm^{desc}(\bC)\ar[d]\\
\Sm^{I}(\bC)\ar[r]&\Sm(\bC)}$$
of inclusions of full subcategories.
Taking left adjoints we get
$$\xymatrix{\Sm^{I,desc}(\bC)&\Sm^{desc}(\bC)\ar[l]^{\cI_{|\Sm^{desc}(\bC)}}\\
\Sm^{I}(\bC)\ar[u]^{L_{|\Sm^{I}(\bC)}}&\Sm(\bC)\ar[u]^{L}\ar[l]^{\cI}}\ .$$

Note that for $X\in \bC$ we have obviously $p^{*}X\in \Sm^{I}(\bC)$ so that
$$L(p^{*}X)=L(I(p^{*}X))=I(L(p^{*}X))\ ,$$
in other words, $L(p^{*}X)$ is homotopy invariant.
  \hB 
}\end{ex}

\subsection{Eilenberg-MacLane and de Rham}

Let $A\in \Ch$ be a chain complex of complex vector spaces.
It gives rise to a constant sheaf of chain complexes $\underline{A}\in \Ch(\Sh_{\Ab}(\Mf))$.
Recall that $\Omega_{\C}\in \Ch(\Sh_{\Ab}(\Mf))$ denotes the sheaf of complex de Rham complexes.
\begin{ddd} \label{aug1001} The sheaf of differential forms with coeffcients in $A$ is defined by  
$$\Omega A:= \Omega_{\C}\otimes_{\C} \underline{A} \in \Ch(\Sh_{\Ab}(\Mf))\ .$$
\end{ddd}
\begin{prob}
If $M$ is compact, then we have
$$\Omega A(M)=\Omega_{\C}(M)\otimes_{\C} A\ .$$
Show by example that the inclusion
$$\Omega_{\C}(M)\otimes_{\C} A\hookrightarrow \Omega A(M)$$ may be   strict
without the compactness assumption.
\end{prob}

Let $H^{*}(A)$ be the cohomology of the chain complex $A$. We consider $H^{*}(A)$  as a chain complex with trivial differential. 
\begin{prob}
Show that for every manifold $M$ there exists an isomorphism
$$H^{*}(\Omega H^{*}(A)(M))\cong H^{*}(\Omega A(M))\ .$$
\end{prob}
\proof
We can choose a quasi-isomorphism
$H^{*}(A)\to A$.
 Hence we get   a quasi-isomorphism of complexes of sheaves
$\Omega H^{*}(A)\to \Omega A$ which induces the asserted isomorphism on the level of cohomology.
\hB

 Recall that $\Nerve(\Ch)[W^{-1}]$ is a stable $\infty$-category and hence has mapping spectra. We consider the group of integers $\Z$ as an object of $\Nerve(\Ch)[W^{-1}]$ by viewing it as a chain complex concentrated in degree zero.
The  Eilenberg-MacLane spectrum is defined by $$H\Z:=\Map(\Z,\Z)\ .$$ It can be considered as a commutative algebra in $\Nerve(\Sp)[W^{-1}]^{\wedge}$ so that we can form its module category $\Mod(H\Z)$. The homotopy groups of $H\Z$ are given by 
$$\pi_*(H\Z)\cong \left\{\begin{array}{cc}
\Z&*=0\\0&*\not=0\end{array}\right. \ .$$ 
This looks   innocent at a first glance but seems difficult to show without going to a concrete realization of the $\infty$-category $\Nerve(\Ch)[W^{-1}]$.

\begin{fact}
There exists a unique symmetric monoidal equivalence
$$H:\Nerve(\Ch)[W^{-1}]\stackrel{\sim}{\to} \Mod(H\Z)$$
such that $H(\Z)=H\Z$. We call $H$ the Eilenberg-MacLane spectrum functor.
\end{fact}
\proof
The domain and target are stable symmetric monoidal $\infty$-categories. The domain is generated 
as a stable $\infty$-category under colimits by the chain complex $\Z$ in degree zero, and $\Mod(H\Z)$ is generated in a similar way by $H\Z$.
We obtain $H$ as the unique colimit preserving functor which maps
$\Z$ to $H\Z$. See  \cite[Prop. 7.1.2.7]{highalg} for a general argument of this type.
\hB

\begin{lem}
We have a canonical isomorphism
$$H_{*}(A)\cong \pi_{*}(H(A))$$
 \end{lem}
\proof
We can write
$$H(A)\cong \Map(\Z,A)\ .$$
Indeed, the functor $\Map(\Z,\dots)$   preserves colimits, commutes with shifts, and produces $H\Z$ for $A=\Z$ by the definition of $H\Z$.
Using this presentation of $H(A)$ we can calculate the homotopy groups as follows:
\begin{eqnarray*}
\pi_k(H(A))&\cong&\pi_k(\Map(\Z,A))\\
&\cong &\pi_0(\Map(\Z[k],A))\\
&\cong &\pi_0(\map(\Z[k],A)) \\
&\stackrel{\ref{jun2503}}{\cong} &H_k(A) \ .
\end{eqnarray*}
\hB

 For a $H\Z$-module spectrum
 $E\in \Mod(H\Z)$ we form the smooth object   $$\Funk(E)\in \Sm(\Mod(H\Z))$$ 
in $\Mod(H\Z)$.
We now consider a chain complex $A\in \Ch$ of complex vector spaces and the object
$\Omega A\in \Sm(\Nerve(\Ch))$ defined in \ref{aug1001}.
It induces an object
$ \Omega A\in \Sm(\Nerve(\Ch)[W^{-1}])$ and therefore a smooth spectrum
$$H(  \Omega A)\in \Sm(\Mod(H\Z))\ .$$
On the other hand we can form the smooth spectrum
$$\Funk(H(A))\in \Sm(\Mod(H\Z ))\ .$$
Assume for example that $A$ is concentrated in degree $0$.
Then we have
$$\pi_{k}(\Funk(H(A))(M))\cong H(A)^{-k}(M)\cong H^{-k}(M;A)\cong H^{-k}(\Omega A(M))\ .$$
This isomorphism extends to general complexes.
More precisely, we have the following  version of  the de Rham isomorphism.

\begin{prop}\label{may2901}
There exists an equivalence of smooth spectra
$$\Rham: H(  \Omega A)\to \Funk(H(A))\ .$$
\end{prop}
\proof
In the following we give an argument which is modeled on the classical proof using integration over simplices.
An independent, more functorial proof will be given in \ref{may2909}.
Let $C_{*}(M)$ be the smooth singular complex of $M$.
Integration over simplices gives a quasi-isomorphism
$$\int:\Omega A(M)\stackrel{\sim}{\to} \underline{\Map}(C_{*}(M),A)\ ,$$
where  $\underline{\Map}$ denotes the internal mapping object of $\Ch$.
We therefore get an equivalence of smooth spectra
$$H(\int):H( \Omega A(\dots))\stackrel{\sim}{\to} H(  \underline{\Map}(C_{*}(\dots),A))\ .
$$
We now use
\begin{eqnarray*}
H(  \Map(C_{*}(\dots),A))&\cong&
   \Map_{\Mod(H\Z)}(H( C_{*}(\dots)),H(A))\\
 &\stackrel{!}{\cong}&  \Map_{\Mod(H\Z)}( \Sigma^{\infty}_{+}(\dots)\wedge H\Z,H(A))\\
   &\cong&\Map_{\Nerve(\Sp)[W^{-1}]}( \Sigma^{\infty}_{+}(\dots),H(A))\\
   &\cong&\Funk(H(A))(\dots)\ , 
\end{eqnarray*}
where the marked equivalence is given in the following Lemma.

\begin{lem}
There exists an equivalence in $\Fun(\Nerve(\Mf),\Mod(H\Z))$ 
$$H( C_{*}(\dots))\stackrel{\sim}{\to}  H\Z\wedge \Sigma^{\infty}_{+}(\dots) \ .$$
\end{lem}
\proof
We use that $H$ preserves the tensor structure and get
$$H(C_{*}(\dots))\stackrel{\ref{jun2504}}{\cong} H(\Z\otimes \dots)\cong H(\Z)\otimes \dots\stackrel{\ref{jun2502}}{\cong} H\Z\wedge \Sigma^{\infty}_{+}(\dots)\ .$$  \hB 
 
As mentioned above an alternative proof of the de Rham equivalence will be given in \ref{may2909}.

Note that $\Omega A$ is a complex of fine sheaves, i.e. it admits an action of partitions of unity.
By Problem \ref{jun2801} we gave
$$\Omega A\in \Sm^{desc}(\Nerve(\Ch)[W^{-1}])\ .$$
Since $H:\Nerve(\Ch)[W^{-1}]\to \Nerve(\Sp)[W^{-1}]$ is the composition of an equivalence and a right-adjoint, it preserves the descent property.
In particular we have  \begin{equation}\label{jun2802}H(\Omega A)\in \Sm^{desc}(\Nerve(\Ch)[W^{-1}])\ . \end{equation}

\begin{ex}{\em 
We consider the sheaf $\Omega_{cl}^{n}\in \Sh_{\Ab}(\Mf)$ of closed $n$-forms. 
\begin{prob}
Calculate the homotopy and cohomology groups  groups of the evaluations on $M$ of the following presheaves and sheaves derived from $\Omega_{cl}^{n}$ in terms of the de Rham cohomology on $M$.
\begin{enumerate}
\item We let $C(\Omega_{cl}^{n})\in \Sm(\Nerve(\Ch)[W^{-1}])$ be the  presheaf of chain complexes given by  $\Omega_{cl}^{n}$ located in degree $0$.
Calculate $H^{*}(C(\Omega_{cl}^{n}))(M))$.
\item Calculate $H^{*}(L(C(\Omega_{cl}^{n}))(M))$, where $L:\Sm(\Nerve(\Ch)[W^{-1}])\to \Sm^{desc}(\Nerve(\Ch)[W^{-1}])$
  is the sheafification.
\item  
Similarly let $C^{\le 0}(\Omega_{cl}^{n})\in \Sm(\Ch^{\le 0}[W^{-1}])$ be the negatively graded chain complex represented by $\Omega^{n}_{cl}$ and $L:\Sm(\Nerve(\Ch^{\le 0})[W^{-1}])\to \Sm^{desc}(\Nerve(\Ch^{\le 0})[W^{-1}])$. 
Calculate $H^{*}(L(C^{\le 0}(\Omega_{cl}^{n}))(M))$ 
\item Calculate the groups
$H^{*}(\bar \bs C(\Omega^{n}_{cl})(M))$, $H^{*}(\bar \bs C^{\le 0}(\Omega^{n}_{cl})(M))$, 
$H^{*}( L(\bar \bs C(\Omega^{n}_{cl}))(M))$ and  $H^{*}(L (\bar \bs C^{\le 0}(\Omega^{n}_{cl}))(M))$, 
 where $\bs$ is as in \ref{jun1901}.
 \item We consider $\Omega_{cl}^{n}$ as a smoth constant simplicial abelian group
$S(\Omega_{cl}^{n})\in \Sm(\Nerve(\sAb)[W^{-1}])$.
Calculate $\pi_{*}(S(\Omega_{cl}^{n})(M))$.
\item  
Calculate $\pi_{*}(\bar \bs S(\Omega_{cl}^{n})(M))$.
 \end{enumerate}
\end{prob}
\proof
These calculations can be found, with different notation, in \cite{MR2192936}.
In order to understand sheafifications use \ref{jun3006}.
\begin{enumerate}
\item $$H^{*}(C(\Omega_{cl}^{n}))(M))=\left\{\begin{array}{cc}
\Omega^{n}_{cl}(M)&*=0\\
0&*\not=0\end{array}\right.$$
\item $$H^{*}(L(C(\Omega_{cl}^{n}))(M))=\left\{\begin{array}{cc}
\Omega^{n}_{cl}(M)&*=0\\
H_{dR}^{n+*}(M)&*> 0\\
0&*<0\end{array}\right.$$
\item $$H^{*}(L(C^{\le 0}(\Omega_{cl}^{n}))(M))=\left\{\begin{array}{cc}0&*\ge 0\\
 H^{n+*}_{dR}(M)&*<0\end{array}\right.$$
\item $$H^{*}(\bar \bs C(\Omega^{n}_{cl})(M))=H^{*}(\bar \bs C^{\le 0}(\Omega^{n}_{cl})(M))= \left\{\begin{array}{cc}
H^{n+*}_{dR}(M)&*\le 0\\
0&*>0\end{array}\right.$$
  $$H^{*}( L(\bar \bs C(\Omega^{n}_{cl}))(M))=H^{*}_{dR}(M)$$
 $$ H^{*}(L (\bar \bs C^{\le 0}(\Omega^{n}_{cl}))(M))=\left\{\begin{array}{cc}
H_{dR}^{n+*}(M)&*<0\\
0&*\ge 0\end{array}\right.$$
\item $$\pi_{*}(S(\Omega_{cl}^{n})(M))=\left\{\begin{array}{cc}
\Omega^{n}_{cl}(M)&*=0\\
0&*\not=0\end{array}\right.$$
\item $$\pi_{*}(\bar \bs S(\Omega_{cl}^{n})(M))=  \left\{\begin{array}{cc}
H^{n-*}_{dR}(M)&*\ge 0\\
0&*<0\end{array}\right.$$

\end{enumerate}
}\end{ex}

\subsection{Differential function spectra and differential cohomology}

A homotopy theoretic construction of a differential extension of a generalized cohomology has first been given in the foundational paper by Hopkins and Singer \cite{MR2192936}.  In that paper the differential cohomology was defined in terms of the homotopy groups of differential function spaces. The approach in the present paper uses the differential function spectrum. Such objects have already been   built in \cite{MR2192936} from  differential function spaces.   In our present course we give a direct construction of differential function spectra. It emphasizes the stable aspect of differential cohomology. This will be important for the construction of operations like integration or transfer.

Before we start with the construction we define 
the notion of a differential extension of a generalized cohomology theory.
This is important since there are many examples where a differential extension has been
constructed by other methods, e.g. of geometric or analytic nature. The axioms list the properties
which one expects to hold in every model of a differential refinement of a generalized cohomology theory.

We start with a description of the data which goes into the definition.
For an abelian group $A$ we can form the Moore spectrum $MA$.
It is a connective spectrum characterized by  
$$ H\Z\wedge MA \cong  H(A) \ .$$
For a spectrum $E$ we write  $EA:=E\wedge MA$. Note that
$M\Z\cong S$ and that we have a canonical map
$M\Z\to MA$ inducing $E\to EA$. 
\begin{prob}
Show that for every $i\in \Z$ there is an exact sequence
$$0\to  \pi_{i}(E)\otimes A\to \pi_{i}(EA)\to \Tor(\pi_{i-1}(E)\otimes A)  \to 0\ .$$
\end{prob}
\proof
Start with a presentation $$0\to \bigoplus_{\alpha} \Z\to \bigoplus_{\beta} \Z\to A\to 0$$
of the group $A$. Model this by a map between wedges of sphere spectra and obtain $MA$ as cofibre
$$  \bigvee_{\alpha} S\to \bigvee_{\beta} S\to  MA \ .$$
Smash with $E$ and discuss the associated long exact sequence in homotopy.
 \hB

\begin{prob}
Show that $H\C\cong M\C$. Furthermore, if $F$ is a $H\C$-module spectrum, then every map
$c:E\to F$   of spectra  uniquely extends to a map 
$c_{\C}:E\C\to F$ of $H\C$-module spectra.
\end{prob}

The construction of the differential extension of $E^{*}$ will depend on the choice of differential data.
\begin{ddd}
A differential data is a triple
$(E,A,c)$ consisting of
\begin{enumerate}
\item a spectrum $E\in \Nerve(\Sp)[W^{-1}]$,
\item 
a chain complex $A\in \Ch$ of complex vector spaces, and \item a map
$c:E\to H( A)$ in $ \Nerve(\Sp)[W^{-1}]$.\end{enumerate}
We say that
$(E,A,c)$ is strict if $c$ induces an equivalence
$c_{\C}:E\C\to H(  A)$. A morphism
of differential data $(E,A,c)\to (E^{\prime},A^{\prime},c^{\prime})$ is a commutative diagram in  $\Nerve(\Sp)[W^{-1}]$
$$\xymatrix{E\ar[r]^{f}\ar[d]^{c}&E^{\prime}\ar[d]^{c^{\prime}}\\
H(  A)\ar[r]^{H(\phi)}&H(  A^{\prime})}\ ,$$
for some morphisms $f:E\to E^{\prime}$ and $\phi:A\to A^{\prime}$.
\end{ddd}

Given the spectrum $E$ we can take $A:=\pi_{*}E\otimes \C$  and let  $c:E\to H(A)$ be a map uniquely determined up to homotopy such that that it induces $\pi_{*}(E)\to \pi_{*}(A)\cong \pi_{*}(E)\otimes \C$,
$x\mapsto x\otimes 1$. This datum is strict.
\begin{ddd}\label{jul1901}
We call a differential data of this form the canonical differential data.
\end{ddd}

\begin{prob}
Analyse the problem whether one can define  a functor from spectra to differential data
which maps a spectrum to its canonical differential data.
\end{prob}
\proof The answer is ``no``. 
The problem is that the map $c$ is unique up to homotopy, but the homotopy is not unique.
\hB

\begin{ddd}
A differential extension of the cohomology theory $E^{*}$ associated to a differential data $(E,A,c)$ is a tuple $(\hat E^{*},R,I,a)$ consisting of
\begin{enumerate}
\item 
 a functor
$\hat E^{*}:\Mf^{op}\to \{\mbox{$\Z$-graded abelian groups}\}$,
\item a transformation $R:\hat E^{*}\to \Omega A^{*}_{cl}$,
\item a transformation $I:\hat E^{*}\to E^{*}$, and 
\item a transformation $a:\Omega A^{*-1}/\im(d)\to \hat E^{*}$,
\end{enumerate}
such that
\begin{enumerate}
\item $\Rham\circ R=c\circ I$,
\item $R\circ a=d$,
\item  \begin{equation}\label{may2910}E^{*-1}\stackrel{c}{\to}\Omega A^{*-1}/\im(d)\stackrel{a}{\to} \hat E^{*}\stackrel{I}{\to} E^{*}\to 0\end{equation}
is exact. 
\end{enumerate}
 \end{ddd}
 
 There is an obvious notion of a morphism between differential extensions associated to the same differential data. In view of \eqref{may2910}   a morphism is automatically an isomorphism by the Five Lemma.

\begin{ex}\label{aug1002}{\em
Observe that $(\hat H^{*}_{Del},R,I,a)$ is a differential extension of $H\Z^{*}$ associated to the canonical data.
}\end{ex}
Let $(\hat E^{*},R,I,a)$ be a differential extension. The homotopy formula measures the deviation of the functor $\hat E$ from homotopy invariance.
\begin{lem}[Homotopy formula]\label{aug1003}
For $x\in \hat E^{*}([0,1]\times M)$ we have the equality  
$$i_1^{*}x-i_0^{*}x=a(\int_{[0,1]\times M/M} R(x))$$
in  $\hat E^{*}(M)$.
\end{lem}
\proof
The same argument as in Proposition \ref{apr2201}. \hB 

In the following we give the  general construction of differential extensions using differential function spectra.
Recall that the smooth objects $\Funk(E)$, $H( \sigma^{\ge n}\Omega A)$ and $H(\Omega A)$ satisfy descent.
Modelled on the construction of Deligne cohomology we make the following definition.

 \begin{ddd}\label{may2902}
 We define the $n$'th differential function spectrum $$\Diff^{n}(E,A,c)\in \Sm^{desc}(\Nerve(\Sp)[W^{-1}])$$ by
 $$\Diff^{n}(E,A,c):=\Cone\left(\Funk(E)\vee H( \sigma^{\ge n}\Omega A)\stackrel{c-\Rham}{\to} \Funk(H( A))\right)[-1]\ .$$
 \end{ddd}
 
 Recall the definition of the Deligne complex \ref{may0201nnn}. The following problem continues  Example \ref{aug1002}.
 \begin{prob}
 Let $E=H\Z$ and $A=\C$ and $c:H\Z\to H( \C)$ be the canonical map.
 Show that there is a natural equivalence
 $$H(\cD(n))\cong \Diff^{n}(H\Z,\C,c)\ .$$
\end{prob}

We now come back to the general case and define the differential cohomology groups in terms of the homotopy groups of the differential function spectrum.
\begin{ddd}
The differential cohomology functor  associated to
$(E,A,c)$ is defined by 
$$\hat E^{n}:=\pi_{-n}(\Diff^{n}(E,A,c))\ .$$
\end{ddd}
Note that  in degree $n$ we use the   differential function spectrum indexed by $n$. This situation is similar as in the construction of Deligne cohomology, Definition \ref{may0201nnn}.

We now discuss the dependence of the differential function spectrum on the data.
 One can describe the $\infty$-category of strict differential data as a pull-back in $\infty$-categories $\infty\Cat$
$$\xymatrix{\Data^{str} \ar[r]\ar[d]&\Nerve(\Ch)\ar[d]^{H}\\
\Nerve(\Sp)[W^{-1}]\ar[r]^{\dots\wedge H\C}&\Mod(H\C)}\ .$$

\begin{prob}\label{jun2601}
Show that there exists a natural map of monoids in $\Nerve(\sSet)[W^{-1}]$
$$ \Omega \map(E,H(A)) \to \ed(E,A,c) \ .$$
\end{prob}
\proof
Let $x$ be an object of an $\infty$-category $\bC$. Then we can identify the group of automorphisms of $x$ as a pull-back 
$$\xymatrix{\aut(x)\ar[r]\ar[d]&\Fun(B\Z,\C)\ar[d]\\
\bullet\ar[r]^{x}&\bC}$$
where the right vertical map is induced by the base point $\bullet\to B\Z$, and the group structure comes from the cogroup structure of $B\Z$. Therefore
 $$\xymatrix{ \aut(E,A,c)\ar[r]\ar[d]&\Fun(B\Z,\Data^{str})\ar[d]\\
\bullet\ar[r]^{(E,A,c)}& \Data^{str}}\ .$$
We now insert the definition of $\Data^{str}$ as a pull-back, move the pull-back outside the functor
$\Fun(B\Z,\dots)$, and restrict to   $\id_{E}$ and $\id_{A}$ in the lower left and upper right corners. This gives a diagram with
$\Fun(B\Z,\map(E\C,H(A)))$ in the upper left corner and therefore a map 
$$\Omega \map(E\C,H(A))\cong \Fun(B\Z,\map(E\C,H(A)))\to  \aut(E,A,c)\ .$$
Finally we precompose use the canonical map
$\map(E,H(A))\to \map(E\C,H(A))$. \hB 
 
We now show how one can interpret  the construction of the differential function spectrum as a funtor from the $\infty$-category of data to spectra.
Let $\Lambda^{2}_2$ be the category of the shape
$$\bullet \to \bullet\leftarrow \bullet\ .$$
Then we get a natural functor
$$\tilde P^{n}:\Data^{str}\to \Fun(\Nerve(\Lambda^{2}_2),\Sm^{desc}(\Nerve(\Sp)[W^{-1}]))\ , (E,A,c)\mapsto  (\Funk(E)\to \Funk(E\C)\leftarrow H(\sigma^{\ge n}(\Omega A))\ .$$

We have 
$$\Diff^{n}(E,A,c) \cong \lim_{\Nerve(\Lambda^{2}_2)} \tilde  P^{n}(E,A,c)\ .$$ 
 \begin{prob}
Give a precise $\infty$-categorical description of the functor
$P^{n}:\Data^{str}\to\Sm^{desc}(\Nerve(\Sp)[W^{-1}])$ explained above. 
\end{prob}
\proof
 We have a natural transformation in $\Fun(\Lambda^{2}_2,\infty\Cat)$
 $$\xymatrix{\Nerve(\Sp)[W^{-1}]\ar@{:>}[dr]^{\id}\ar[d]^{\Funk}\ar[r]^{\dots\wedge H\C}&\Nerve(\Sp)[W^{-1}]\ar[d]^{\Funk}&\Ch\ar[l]_{H}\ar[d]^{H(\sigma^{\ge n}\Omega_\C\otimes\underline{\dots})}\ar@{:>}[dl]^{\Rham}\\
\Sm^{desc}(\Nerve(\Sp)[W^{-1}])\ar[r]^{\id}&\Sm^{desc}(\Nerve(\Sp)[W^{-1}])&\Sm^{desc}(\Nerve(\Sp)[W^{-1}])\ar[l]_\id}\ .$$
Taking the limit over $\Nerve(\Lambda^{2}_2)$ we get the functor 
$P^{n}:\Data^{str}\to \Sm(\Nerve(\Sp)[W^{-1}])$. \hB

For completeness let us mention that in order to capture non-strict data we can  form the limit of
$$\xymatrix{\Data  \ar[rr]\ar[d]&&\Nerve(\Ch)\ar[d]^{H}\\\Nerve(\Sp)[W^{-1}] &\Nerve(\Sp)[W^{-1}]^{\Nerve([1])}\ar[r]\ar[l]&\Nerve(\Sp)[W^{-1}] } \ .$$

We now calculate the homotopy groups of the differential function spectrum.

\begin{theorem}
For  $i>n$ we have $$\pi_{-i}(\Diff^{n}(E,A,c))\cong  E^{i}(M) .$$
For $i<n$ we have $$\pi_{-i}(\Diff^{n}(E,A,c))\cong  \Cone(c)^{i-1}(M)\ ,$$
and if $(E,A,c)$ is strict, then 
$$\pi_{-i}(\Diff^{n}(E,A,c))\cong E\C/\Z^{i-1}(M)\ .$$
Finally, for $i=n$  we have an exact sequence
\begin{equation}\label{may2810}E^{n-1}(M)\to H(A)^{n-1}(M)\stackrel{a}{\to}\hat E^{n}(M)\stackrel{I\oplus R}{\to}  E^{n}(M)\oplus \Omega A_{cl}^{n}(M)\to H(A)^{n}(M)\to 0\ .\end{equation}
\end{theorem}
\proof
We have a long exact sequence of homotopy groups associated to the cone which for
$i>n$ reads
as
\begin{eqnarray*}E^{i-1}(M)\oplus H(A)^{i-1}(M)\twoheadrightarrow H(A)^{i-1}(M)\hookrightarrow \pi_{-i}(\Diff^{n}(E,A,c))&&\\ \to  E^{i}(M)\oplus H(A)^{i}(M)\twoheadrightarrow H(A)^{i}(M)&&\ .\end{eqnarray*}
This gives the assertions for $i>n$.
A similar reasoning using $\pi_{-n}(H(\sigma^{\ge n}\Omega A)=\Omega A^{n}_{cl}$ gives the case $i=n$.
 
 We have a triangle
$$\to \Cone\left(\Funk(E)\stackrel{c}{\to} \Funk(H(A))\right)[-1]\to \Diff(E,A,c)\to H(\sigma^{\ge n} \Omega A)\to $$
If $(E,A,c)$ is strict, then we have an equivalence $$\Cone\left(\Funk(E)\stackrel{c}{\to} \Funk(H(A))\right)[-1]\cong \Funk(E\C/\Z)[-1]\ .$$ Since for $i<n$ we have
$\pi_{-i}H(\sigma^{\ge n} \Omega A)=0$ we conclude the assertion for $i<n$.
 \hB
  
We now consider the fibre sequence defining the spectrum $\cE(n)$
$$\Cone(H(\sigma^{\ge n}\Omega A)\stackrel{=}{\to} H(\sigma^{\ge n}\Omega A))[-1] \to \Diff^{n}(E)\stackrel{\pi}{\to} \cE(n)\ ,$$
where the first map is the composition of the de Rham equivalence
$$  \Cone(H(\sigma^{\ge n}\Omega A)\to H(\Omega A))[-1]
\stackrel{\sim}{\to}
  \Cone(H(\sigma^{\ge n}\Omega A)\to \Funk(H(A))) $$
  with the two obvious inclusions.
  The map $\pi$ induces an equivalence.
  We can identify
  $$\cE(n)\cong \Cone(\Funk(E)\to H(\sigma^{<n}\Omega A))[-1]\ .$$
  The map $H(\sigma^{<n}\Omega A))\to \cE(n)$ induces a map
  $$a:\Omega A^{n-1}(M)/\im(d) \to \hat E^{n}(M)\ .$$
  The same arguments as for Deligne cohomology show:
  \begin{prop}
  The tuple $(\hat E^{*},R,I,a)$ is a differential extension of $E^{*}$. In detail, 
  the map $a:\Omega A^{n-1}(M)/\im(d) \to \hat E^{n}(M)$ extends the map
  $H(A)^{n-1}(M)\to \hat E^{n}(M)$.  We have
  $$R\circ a=d\ .$$
  The sequence
  \begin{equation}\label{jul1002}E^{n-1}(M)\to \Omega A^{n-1}(M)/\im(d)\to \hat E^{n}(M)\to E^{n}(M)\to 0\end{equation}
  is exact. 
  If $(E,A,c)$ is strict, then 
  we have the differential cohomology diagram
  $$\xymatrix{&\Omega A^{n-1}(M)/\im(d)\ar[dr]^{a}\ar[rr]^{d}&&\Omega A^{n}_{cl}(M)\ar[dr]&\\\textcolor{red}{E\C^{n-1}(M)}\ar[ur]\textcolor{red}{\ar[dr]}&&\hat E^{n}(M)\ar[ur]^{R}\ar@{->>}[dr]^{I}&&\textcolor{red}{E\C^{n}(M)}\\&
\textcolor{red}{
E\C/\Z^{n-1}(M)}\ar[rr]^{-Bockstein}\ar@{^{(}->}[ru]&&\textcolor{red}{E^{n}(M)}\ar[ur]&}
\ .$$
\end{prop}

\begin{ex}{\em
We consider the differential data $(H\Z,\C,H\Z\to H\C)$ in order to define the differential extension $(\widehat{H\Z},R,I,a)$.
\begin{prob}\label{may2911}
Show that there is a unique isomorphism of differential extensions between 
$(\widehat{H\Z},R,I,a)$ and $( \hat H^{*}_{Del},R,I,a)$.
\end{prob}
\proof See \cite{MR2608479} or \cite{MR2365651}.
}
\end{ex}
\begin{ex}{\em Recall the action \ref{jun2601} of $\Omega \map(E,H(A))$ on the datum
$(E,A,c)$.
\begin{prob}
Describe the action  of $\pi_{0}(\Omega \map(E,H(A)))\cong H(A)^{-1}(E)$ on $\hat E^{n}$.
\end{prob}
\proof Let $\phi\in  \Omega \map(E,H(A))  $ and $\hat x\in \hat E^{n}(M)$. 
 We write $\nu_\phi(\hat x)$ for the action of $\phi$ on $\hat x$. Note that
 $\phi(I(\hat x))\in H(A)^{n-1}(M)$. We get
$$\nu_\phi(\hat x)=a(\phi(I(\hat x)))\ .$$
\hB 
}
\end{ex}

Unlike   generalized   cohomology theories, differential cohomology is not homotopy invariant. The deviation from homotopy invariance is measured by the homotopy formula \ref{aug1003}. But it still has an interesting Mayer-Vietoris sequence as a consequence of the descent property of the differential function spectrum.

\begin{prob}Construct the Mayer-Vietoris sequence in differential cohomology.
Assume that the data is strict.
\end{prob}
\proof If $M=U\cup V$ is a decomposition into two open submanifolds, then by the descent property of the differential function spectrum 
we get a pull-back
$$\xymatrix{\Diff^{n}(E,A,c)(M)\ar[r]\ar[d]&\Diff^{n}(E,A,c)(U)\ar[d]\\
\Diff^{n}(E,A,c)(V)\ar[r]&\Diff^{n}(E,A,c)(U\cap V)}$$
in $\Nerve(\Sp)[W^{-1}]$. The interesting segment of the long exact sequence is
$$\dots \to E\C/\Z^{n-2}(U\cap V)\to \hat E^{n}(M)\to \hat E^{n}(U)\oplus \hat E^{n}(V)\to \hat E^{n}(U\cap V)\to E^{n+1}(M)\to \dots\ .$$ 
It extends by the Mayer-Vietoris sequence of $E\C/\Z^{*}$ and $E^{*}$ to the left- and right-hand sides. \hB

\begin{ex}
Let $V$ be some abelian group and $G$ be a monoid acting on $V$. 
\begin{prob} Show that there is a unique extension differential extension of the cohomology theory $H^{*}(\dots;V)$, and that the action of $G$ extends to an action on this differential extension.
\end{prob}
 \end{ex}

 \begin{ex}
 {\em 
 A differential cohomology class $x\in \hat E^{n}(M)$ is essentially the datum of an underlying class
 $I(x)\in E^{n}(M)$ together with form data $R(x)\in \Omega A^{n}_{cl}(M)$ combined in an appropriate homotopical manner. This leads to the following two very suggestive alternative descriptions of differential cohomology suggested to me by Bruce Williams. These models make the nature of the structure maps $R,I,a$ completely clear.
 \begin{prob}\label{jul2101}
 Show that there are natural equivalences between $\Diff^{n}(E,A,c)$ and the pull-backs in $\Sm(\Nerve(\Sp)[W^{-1}])$
$$ \xymatrix{\Diff^{n}(E,A,c)^{\prime}\ar[r]\ar[d]&H(L(\Omega A_{cl}^{n}))\ar[d]\\
 \Funk(E)\ar[r]&\Funk(H(A))}\ ,\quad \xymatrix{\Diff^{n}(E,A,c)^{\prime\prime}\ar[r]\ar[d]&H(L(\Omega A_{cl}^{n-1}/\im(d)))\ar[d]\\
 \Funk(E)\ar[r]&\Funk(H(A))}\ .$$
 \end{prob}
 \proof
 Use \ref{jun3006}. \hB 
 
 Let $\iota:\Mod(H\Z)\to \Nerve(\Sp)[W^{-1}]$ be the forgetful functor.
\textcolor{red}{
\begin{prob}\label{aug0754}
Do we have an equivalence $L\circ \iota\cong \iota\circ L$
of transformations $\Sm(\Mod(H\Z))\to \Sm^{desc}(\Nerve(\Sp)[W^{-1}])$.
\end{prob}}
A positive answer would allow to interchange $H$ and $L$ in right-upper corners of the pul-back diagrams in \ref{jul2101}.

 }
 \end{ex}

\subsection{Differential $K$-theory}
  
One of the most important examples of differential cohomology theories is the differential extension of complex $K$-theory $\bKU^{*}$. It was first  used in string theory and $M$-theory in order to model
the global topology of fields with differential form field strength, see e.g.  \cite{MR1919425}, \cite{MR1769477}. It is also closely related to the local index theory of Dirac operators and can be used to formulate refined index theorems capturing secondary information, see \ref{aug0750}, \cite{MR2602854}. To some extend, the  consideration of differential complex $K$-theory motivated the development of differential cohomology theory as a field of mathematical  research.
Predating the constructions of full differential extension of complex $K$-theory there were geometric constructions of the functor $\bKU\R/\Z$, starting with  Karoubi \cite{MR1076525}, \cite{MR1273841}, \cite{MR913964} and Lott \cite{MR1312690}.
There are various topological, geometric and  analytic models of differential complex $K$-theory, see e.g.   \cite{MR2192936}, \cite{MR2664467}, \cite{MR2550094},  \cite{MR2521641}, \cite{2009PhDT}, \cite{MR2231056}. This variety of models  called for the investigation of uniqueness in \cite{MR2608479}.

Here we present differential complex $K$-theory  as an example of our general construction. 
The main tool to relate it with geometric problems is the cycle map defined in \ref{aug1005}.

We consider a ring spectrum $\bKU$  representing the generalized cohomology theory $\bKU^{*}$, called complex $K$ theory. In order to be able to  connect with geometry we fix an identification
\begin{equation}\label{may2801}\Omega^{\infty}\bKU\cong \Z\times BU\end{equation} as $h$-spaces. 
If $\bKU$ is defined using the classical approach via vector bundles, then this identification is built in into the construction. On the other hand there are more abstract constructions of $\bKU$ e.g. from the  multiplicative formal group law using the Landweber exact functor theorem. In this case the relation with vector bundles is not obvious.

We fix the Bott element $b\in \pi_{2}(\bKU)$ represented by the map
$$S^{2}\cong \C\P^{1}\stackrel{L^{*}}{\to} BU(1) \hookrightarrow BU\to \{0\}\times BU\hookrightarrow \Z\times BU\cong \Omega^{\infty}\bKU\ ,$$
where $L^{*}\to \C\P^{1}$ is the dual of the tautological bundle (see \ref{apr1260}).
Note that the cohomological degree of the Bott element $b$ is $-2$.
The choice of $b\in \pi_{2}(\bKU)$ determines
 an isomorphism of rings
$\pi_{*}(\bKU)\cong \Z[b,b^{-1}]$. We define the chain complex $A:=\C[b,b^{-1}]\in \Ch$ with trivial differential. There is a   map
$$c:\bKU\to H(A)$$ uniquely determined up to a unique homotopy class of homotopies which induces the inclusion
$$\pi_{*}(\bKU)\cong \Z[b,b^{-1}]\hookrightarrow \C[b,b^{-1}]\cong \pi_{*}(H(A))$$
on the level of homotopy groups. 
\begin{prob}
Verify this assertion.
\end{prob}

We let
$(\bKU,A,c)$ be our strict differential data for the differential extension
$\hbKU^{*}$.

\begin{prob}
Calculate as abelian groups $\hbKU^{*}(*)$ and $\hbKU^{*}(S^{1})$. 
\end{prob}

A  complex vector bundle  $V\to M$ represents a class  $\bKU^{0}(M)$. In order to make this precise we have to use the identification
 \eqref{may2801}.
 Note that $\bigsqcup_{n\ge 0} BU(n)$ is the classifying space for complex vector bundles.
Let $$v:M\to \bigsqcup_{n\ge 0} BU(n)$$ be the classifying map of $V$.
\begin{ddd} We define the topological cycle 
 $\cycl(V)\in \bKU^{0}(M)$ of the complex vector bundle $V\to M$ as the homotopy class of the composition
$$M\stackrel{v}{\to}   \bigsqcup_{n\ge 0} BU(n)\to \Z\times BU\cong \Omega^{\infty}\bKU\ .$$
\end{ddd}
Let $\Vect:\Mf^{op}\to \CommMon(\Set)$ be functor which associates to a manifold $M$ the commutative monoid of isomorphism classes of vector bundles with respect to the direct sum. Then we can interpret the topological cycle map
$$\cycl:\Vect\to \bKU^{0}$$
as a map of  smooth monoids $\Sm(\CommMon(\Set))$. This uses the additional precison that \eqref{may2801} is a map of $h$-spaces.

Let $\nabla$ be a connection on the complex vector bundle $V\to M$. Then we can define a Chern character form
$$\ch(\nabla):=\Tr\exp(-\frac{bR^{\nabla}}{2\pi i}) \in \Omega A^{0}_{cl}(M)\ .$$
This is essentially  \eqref{may2803}. We add the factor $b$ of cohomological degree $-2$ in order to get an element of degree zero.
\begin{prob}\label{may2820}
Show  that
$\Rham([\ch(\nabla)])=c(\cycl(V))$
in $H(A)^{0}(M)$.
\end{prob}
\proof
The natural transformation
$\Vect(M)\ni V\mapsto \Rham([\ch(\nabla)])\in H(A)^{0}(M)$
is additive and hence factorizes over a transformation
$\bKU^{0}(M)\to H(A)^{0}(M)$. The restriction of this transformation to compact manifolds
 is thus represented by a map of $h$-spaces $\Omega^{\infty}\bKU\to \Omega^{\infty} H(A)$, i.e. a
primitive cohomology  class in $ H(A)^{0}(\Omega^{\infty}\bKU)$. 
We know that $H^{ev}(\Omega^{\infty} \bKU;\C)^{prim}$ is spanned by the classes $\ch_{2n}$ for $n\ge 0$.
It follows that
$H(A)^{0}(\Omega^{\infty}\bKU)^{prim}$ is spanned by the classes
$b^{n}\ch_{2n}$ for $n\ge 0$. 
We must only check that
$\Rham(\ch(\nabla))=b^{n}\ori_{S^{2n}}\in H(A)^{0}(S^{2n})$
when the (virtual) bundle $V\to S^{2n}$ represents the element in $\tilde \bKU^{0}(S^{2n})$ which corresponds to $b^{n}\in \bKU^{2n}(*)$.
Using multiplicativity of the Chern character it suffices to check this for $n=1$.
The Bott generator is $\cycl(L^{*})-\cycl(1)$.
In this case, by \ref{apr1260}, 
$$\Rham(\ch(\cycl(L^{*}))) -\Rham(\ch(\cycl(1)))=b\ori_{S^{2}}$$
as required.
\hB 

We let $$\Vect^{geom}:\Mf^{op}\to \CommMon$$ be the functor which associates to a manifold $M$ the commutative monoid (under direct sum) of isomorphism classes of vector bundles with connection.
\begin{ddd}\label{aug1005}
A differential refinement of the topological cycle map is a transformation
$$\hcycl:\Vect^{geom}\to \hbKU^{0}$$ of semigroup valued functors 
satisfying
\begin{enumerate}
\item
$R(\hcyc(V,\nabla))=\ch(\nabla)\in \Omega A_{cl}^{0}$
\item $I(\hcyc(V,\nabla))=\cycl(V))\in \bKU^{0}$.
\end{enumerate}
\end{ddd}

\begin{theorem}\label{may2930}
There exists a unique differential refinement of the topological cycle map.
\end{theorem}
\proof
Let $(V,\nabla^{V})$ be a $k$-dimensional vector bundle with connection on a manifold $M$ of dimension $n$.
Then we can find a $n+1$-connected approximation $N\to BU(k)$
with a bundle with connection $(W,\nabla^{W})$ and a smooth map
$h:M\to N$ such that $h^{*}(W,\nabla^{W})\cong (V,\nabla^{V})$.
Observe that $0=h^{*}:H(A)^{-1}(N)\to H(A)^{-1}(M)$. We choose
some element $w\in \hbKU^{0}(N)$ with $R(w)=\ch(\nabla^{W})$ and $
I(w)=\cycl(W)$. Such an element exists by \eqref{may2810} and is uniquely determined up to 
elements in $a(H(A)^{-1}(N))$. By naturality we are forced to define
$$\hcycl(V,\nabla^{V}):=h^{*}(w)\ .$$
If we can show that this class is well-defined then the remaining properties
(naturality, aditivity) follow easily. 

Given a second choice we argue as in the proof of Theorem \ref{apr2001}.
Using the notation introduced there we can assume that
$(W,\nabla^{W})\cong g^{*}(W^{\prime\prime},\nabla^{W^{\prime\prime}})$ and
$(W^{\prime},\nabla^{W^{\prime}})\cong g^{\prime*}(W^{\prime\prime},\nabla^{W^{\prime\prime}})$.

\begin{prob}
Give more details.
\end{prob}
We must show that $h^{*}g^{*} w^{\prime\prime}=h^{\prime*}g^{\prime*} w^{\prime\prime}$.
Let $H:I\times M\to N^{\prime\prime}$ be the homotopy from
$g\circ h$ to $g^{\prime}\circ h^{\prime}$. The bundle
$H^{*}(W^{\prime\prime},\nabla^{W^{\prime\prime}})$ can be glued to a bundle
$(\tilde W,\nabla^{\tilde W})$ on $S^{1}\times M$.
By the homotopy formula
$$h^{*}g^{*} w^{\prime\prime}-h^{\prime*}g^{\prime*} w^{\prime\prime}=a(R(H^{*}(w^{\prime\prime})))=a(\int_{S^{1}\times M/M} \ch(\nabla^{\tilde W}))=0$$
since $$\int_{S^{1}\times M/M} \ch(\nabla^{\tilde W})\in \im(\ch:\bKU^{-1}(M)\to H(A)^{-1}(M))\ .$$
\hB 

Let us add a remark on the role of cyle maps. First of all it is obvious that the cycle map provides a means to write down elements in differential complex $K$-theory. Reflecting on what it means to calculate an element in differential $K$-theory one can get several possibilities. 
One can map the element to some known group. For example, in the odd case, one can restrict to a point and study the resulting element in $\C/\Z$. But often a very satisfactory way to calculate an element in differential $K$-theory  is to exhibit it as a cycle of an explicit vector bundle with  connection. This makes even more sense in view of the following result.

We let $\Vect^{geom}_{virtual}$ denote the smooth set of virtual (i.e. formal differences of) vector bundles.

\begin{ex}{\em 
Let $M$ be a compact manifold. 
\begin{prob} Show that the cycle map $\hcycl: \Vect^{geom}_{virtual}(M)\to \hbKU^{0}(M)$
is surjective.
\end{prob}
\proof
See \cite{MR2521641}. \hB }
\end{ex}

\begin{ex}{\em
For $n\ge 1$ the group $Spin(n+1)$ acts on $S^{n}$. 
Note that by suspension $\tilde \bKU^{n}(S^{n})\cong \Z g$ for some generator $g$.
\begin{prob}Show that for $n\ge 2$ there exists a unique
$x\in \hbKU^{n}(S^{n})$ with $I(x)=g$ which is $Spin(n+1)$-invariant and evaluates trivially on points. What goes wrong in the case $n=1$?\end{prob}
See \cite[Ch. 5.7]{MR2664467} for an argument.

\begin{prob}
Note that $\tilde \bKU^{0}(S^{2n})\cong b^{n} \tilde \bKU^{2n}(S^{2n})\cong b^{n} \Z$.
Calculate the class $b^{n}$.  More precisely find a virtual bundle with connection $(V,\nabla)$ on $S^{2n}$ such that
$\hcycl(V,\nabla^{V})\in \hbKU^{0}(S^{2n})$ is the unique $Spin(2n+1)$-invariant class 
with $I( \hcycl(V,\nabla^{V}))=b^{n}$.
\end{prob}
}
\end{ex}

\begin{ex}{\em
To a map $f:M\to U(k)$ we can associate the suspension bundle $T(f)\to S^{1}\times M$.
It is a $k$-dimensional complex vector bundle which represents a class
$\cycl(T(f))\in \bKU^{0}(S^{1}\times M)$. This class  only depends on the homotopy class $[f]\in [M,U(k)]$ of $f$. The set of homotopy classes
$[M,U(k)]$ has a natural group structure induced by the group structure on $U(k)$ and one can check that the map $[f]\mapsto \cycl(T(f))$ is a homomorphism.

The projection $\pr:S^{1}\times M\to M$ is $\bKU$-oriented and we have an Umkehr map
$\pr_!:\bKU^{0}(S^{1}\times M)\to \bKU^{-1}(M)$. We define  an odd version of a topological cycle map 
$$\cycl^{-1}:[M,U(k)]\to \bKU^{-1}(M)$$
which maps the homotopy class $[f]$ of the  map $f:M\to U(k)$ to the class
$$\cycl^{-1}([f]):=\pr_! (\cycl(T(f)) \in \bKU^{-1}(M)\ .$$
This map preserves the group structures. We now ask whether one can refine this construction to differential $K$-theory.
\begin{prob}\label{aug0753}
Is there a cycle map 
$$\hcycl^{-1}:C^{\infty}(\dots,U(k))\to \hbKU^{-1}(\dots)$$ 
which preserves the group structure. 
$$I(\hcycl^{-1}(f))=\cycl^{-1}(f)\ .$$  
\end{prob}
\proof
The answer is yes for $k=1$ and no for $k\ge 2$.

Every cycle map $\hcycl^{-1}$ is induced by a universal class $\hat u_k\in \hbKU^{-1}(U(k))$.
It preserves the group structure if and only if $\hat u_k$ is primitive, i.e. satisfies  $\mu^{*} \hat u_k=\pr_1^{*}\hat u_k+\pr_2^{*}\hat u_k
$ on $U(k)\times U(k)$.

We first consider the case $k=1$. One can show that there  exists a unique primitive class
$\hat u_1\in \hbKU^{-1}(U(1))$ whose curvature $R(\hat u_1)$ is the normalized invariant volume
form. 

For $k\ge 2$ the universal class can not be primitive. Otherwise,  the higher-degree components of its curvature would be a primitive forms of degree $\ge 2$, but such forms do not  exist.

\hB}
\end{ex}

\begin{ex}{\em
The $n$th Chern class can be considered as a transformation of smooth sets
$$c_{n}:\Vect\to H\Z^{2n}\ .$$
\begin{prob}
Observe that $c_{n}$ extends to a transformation
$$c_{n}:\bKU^{0}\to H\Z^{2n}\ .$$
Furthermore show, that there exists a unique extension
$$\hat c_{n}:\hbKU^{0}\to \hat H^{n}_{Del}(M;\Z)$$
such that
 $$\hat c_n(\hcycl(V,\nabla))=\hat c_{n}(\nabla)$$
 (where $\hat c_{n}$ on the r.h.s. is the differential refinement of $c_{n}$ according to Theorem \ref{may2007})
\end{prob}
\proof
If $M$ is compact, then
we can identify the set $\bKU^{0}(M)$ with the set of stable equivalence classes of vector bundles on $M$. Since the Chern class $c_{n}$ is well-defined on stable equivalence classes we get the factorization.

In order to construct the differential refinement we again use manifold approximations of $BU$.
See \cite{MR2734150} for details. \hB
}\end{ex}

\begin{ex}{\em \label{jun0210}
Consider a multiplicatively closed subset $S\subset \nat$ and the spectrum $\bKU[S^{-1}]$. We have a natural extension 
$$\xymatrix{\bKU\ar[r]^{c}\ar[d]&H(\C[b,b^{-1}])\\\bKU[S^{-1}]\ar@{..>}^{c}[ur]&}
$$
We define
$\hbKU[S^{-1}]$ using the differential data
$(\bKU[S^{-1}],A,c)$.
\begin{prob}
Show that there exists a natural transformation
$\hbKU\to \hbKU[S^{-1}]$
which induces an isomorphism
$\hbKU^{*}(M)\otimes \Z[S^{-1}]\to \hbKU[S^{-1}]^{*}(M)$
for compact manifolds $M$. 
Why do we assume that $M$ is compact?
\end{prob}

For every $k\in S$ we have the Adams operations $\Psi^{k}$ on $\bKU[S^{-1}]$ and
$\Psi^{k}_{H}:\C[b,b^{-1}]\to \C[b,b^{-1}]$, $b\mapsto k^{-1} b$
such that the following diagram commutes canonically.
$$\xymatrix{\bKU[S^{-1}]\ar[r]^{c}\ar[d]^{\Psi^{k}}&H(\C[b,b^{-1}])\ar[d]^{\Psi_{H}^{k}}\\\bKU[S^{-1}]\ar[r]^{c}&H(\C[b,b^{-1}])}\ .$$
Hence for every $k\in S$ we have a map of data
$$\Psi^{k}:(\bKU[S^{-1}],A,c)\to (\bKU[S^{-1}],A,c)$$ which induces 
   differential Adams operation $\hat \Psi^{k}$ on $\hbKU[S^{-1}]$.
\begin{prob}
Verify the relation $\hat \Psi^{k}\hat \Psi^{l}=\hat \Psi^{kl}$.
\end{prob}

More one Adams operations in differential $K$-theory can be found in \cite{MR2740650}.

We consider the torus $T^{3}$ with its action on itself.
\begin{prob}
Calculate
$ \hbKU^{*}[\{5\}^{-1}](T^{3})$
and determine the action of $\hat \Psi^{5}$.
\end{prob}

Let $(L,\nabla)$ be a complex line bundle over a manifold $M$.
\begin{prob}
Show that 
$$\hcycl((L,\nabla)^{\otimes k})=\hat\Psi^{k}(\hcycl(L,\nabla))$$
in $\hbKU[\{k\}^{-1}](M)$.
\end{prob}
}
\end{ex}

\begin{ex}\label{jun0301}
{\em
Recall the notion of a geometric family $\cE$ over $M$ \cite[Def. 2.2]{MR2664467}. It consists of
\begin{enumerate}
\item a proper submersion $f:W\to M$,
\item a Riemannian structure $(T^{h}f,g^{T^{v}f})$ consisting of a horizontal subbundle and a vertical Riemannian metric, 
\item a family of Dirac bundles $V$,
\item an orientation of $T^{v}\pi$.
\end{enumerate}
We say that a geometric family is of degree $0$ or $-1$ if the fibre dimension is even or odd.
Isomorphism classes of geometric families on $M$ form a graded semi-group $\GeomFam^{*}(M)$ under fibrewise disjoint sum.
In  \cite[Def. 2.2]{MR2664467} we have constructed a model  $\hbKU_{\GeomFam}^{*}$ of a differential extension $\hbKU^{*}$ as a group of equivalence classes of geometric families. In particular we have a tautological cycle map
$$\hcycl_{\GeomFam}:\GeomFam^{*}(M)\to \hbKU_{\GeomFam}^{*}(M)\ ,\quad *\in \{0,-1\}$$
which maps a geometric family to its equivalence class. A vector with $V\to M$ with hermitean metric and metric connection
$\nabla$ can naturally be considered as geometric family $\cE$ over $M$, where $f=\id:M\to M$.

The following exercise shows that the differential  cycle map can be extended to geometric families.
\begin{prob}
Show that there is a unique isomorphism of group-valued functors
$$\iota:\hbKU_{\GeomFam}^{0}\cong \hbKU^{0}$$
which is compatible with the structures $a,R,I$.
Show further, that for a hermitean vector bundle $V$ with connection $\nabla$ we have
$$\iota(\hcycl_{\GeomFam}(V,\nabla))=\hcycl(V,\nabla)\ .$$
\end{prob}
\proof
For the first statement we refer to \cite{MR2608479}.
The second follows from the unicity of the cycle map.
\hB

}

\end{ex}

\subsection{Differential Bordism theory}

The complex bordism theory $\bMU^{*}$ plays a fundamental role in stable homotopy theory.
Its coefficient ring carries the universal formal group law. The relation between formal group laws
and complex oriented cohomology theories governs the structure of the stable homotopy category. More specifically, via the Landweber exact functor theorem, complex bordism gives rise to a variety of
complex oriented generalized cohomology theories with Landweber exact formal group laws.  
This correspondence generalizes to the differential case. In \cite{MR2550094} (cf. \ref{aug1021}) we observed that a differential extension of
$\hbMU$ gives rise to differential extensions of complex oriented generalized cohomology theories with Landweber exact formal group laws. In this sense the differential extension of $\bMU$ plays a fundamental role in the theory. 

It turned out that for bordism theories in general there are geometric constructions of 
differential extensions in which additional structures like multiplications or integration are easy  to built in.
This has been used e.g. in \cite{MR2674652} to deliver a bordism model for the differential extension of ordinary integral cohomology in which one has integration and products and a simple verification the projection formula.

The discussion of the cycle map for complex bordism is in may aspects parallel to the cycle map in the holomorphic and algebraic cases where one  associates to Arakelov cycles classes in absolute Hodge cohomology or Arakelov Chow groups. We will not try to review the vast literature in this direction.

In the present course we discuss the differential extension of $\bMU$  because it is fundamental and simple. The constructions for other bordism theories are simple modifications of the case of $\bMU$.

Let $\bMU$ be a spectrum representing the complex bordism theory. 
We let $A_{*}:=\bMU_*\otimes \C$ and fix an equivalence $c:\bMU\C\to H(A)$ so that
 \begin{equation}\label{may2931}\bMU \to \bMU\C\stackrel{c}{\to} H(A) \end{equation}
 induces the complexification map in homotopy.
 We obtain a strict differential data $(\bMU,A,c)$ and therefore a differential extension $(\hbMU,R,I,a)$.
 
In the following we want to  connect $\bMU$ with the geometric picture. To this end we realize
$\bMU$ as the Thom spectrum associated to the map $BU\to BO$.

Cycles for complex cobordism classes are maps between manifolds with complex stable normal bundle. In the following we give some details which are needed to generalize to the geometric situation.  Let $f:W\to M$ be a smooth map of manifolds.
\begin{ddd}\label{may2941}
A representative $\cN$ of the stable normal bundle of $f$ is an exact sequence
$$0\to TW\stackrel{df\oplus \alpha}{\to} f^{*}TM\oplus  \R^{k}\to N\to 0\ .$$
We call $N\to W$ the underlying bundle of $\cN$.
\end{ddd}

The sequence $\cN$ is just the infinitesimal model of an embedding $\iota$ over $M$
$$\xymatrix{W\ar[rr]^{\iota}\ar[dr]^{f}&&M\times \R^{k}\ar[dl]\\
&M&}\ .$$
Given $\iota$ we let $\alpha$ be the second component of $d\iota$.

\begin{ddd}\label{jul2601}
A   cycle for a class in $\bMU(M)$ of degree $n$ is a triple $(f,\cN,J)$ of a proper smooth map
$f:W\to M$ such that $\dim(M)-\dim(W)=n$, a representative of a stable normal bundle $\cN$, and a complex structure $J$ on its the underlying bundle.  
\end{ddd}
There is a natural notion of an isomorphism between cycles. It is very usefulful to
add homotopies of representatives of the stable normal bundle with fixed underlying bundle to the isomorphism relation. This point of view is also adopted in   \cite{MR2550094}.
\begin{ddd}\label{jul3001}We let $\Cycle^{*}_{\bMU}(M)$ denote the graded semigroup of isomorphism classes of cycles with respect to disjoint sum. 
\end{ddd}

If we realize $\bMU$ as the Thom spectrum of the map $BU\to BO$, then
the Thom-Pontrjagin construction gives a map of semigroups
$$\cycl:\Cycle^{*}_\bMU(M)\to \bMU^{*}(M)\ .$$
\begin{prob}
Understand the details.
One can actually put an equivalence relation on $\Cycle^{*}_{\bMU}$
involving stabilization of the representative of the stable normal bundle and
bordism such that $\cycl$ induces an isomorphism on quotients.
\end{prob}
\proof
See  \cite{MR2550094} for details. \hB 

A complex vector bundle has an oriented underlying real vector bundle and therefore a Thom class  for $H\Z$.
Putting the $H\Z$-Thom isomorphisms of the universal bundles $\xi_{n}\to BU(n)$ for all $n\ge 0$ together we get
a Thom isomorphism
$$\Phi:M^{*}(BU)\stackrel{\cong}{\to} M^{*}(\bMU)$$ for every $H\Z$-module spectrum $M$.
Note that on the right-hand side we consider the spectrum cohomology of $\bMU$ with coefficients in $M$. The shifts involved in the construction of $\bMU$ ensure that there is no degree-shift in the Thom isomorphism.
In particular, we have 
the Thom isomorphism $$\Phi:H(A)^{0}(BU)\stackrel{\sim}{\to} H(A)^{0}(\bMU)\ .$$ Under this isomorphism
the spectrum cohomology class
$c:\bMU\stackrel{\eqref{may2931}}{\to} H(A)$ corresponds to a cohomology class
$$u:=\Phi^{-1}(
c) \in H(A)^{0}(BU)\ .$$ 

The space $BU$ classifies stable complex vector bundles.
For a stable complex vector bundle $N\to W$ classified by a map $\rho:W\to BU$ we write
\begin{equation}\label{aug1050}u(N):=\rho^{*} u \in H(A)^{0}(W)\ .\end{equation}

If $(f,\cN,J)$ is a cycle, then $f$ is proper and oriented for ordinary cohomology. Hence  we have an Umkehr  map
$$f_!:H(A)^{*}(W)\to H(A)^{*+n}(M)\ ,$$
(see \ref{aug0501} for more details on Umkehr maps which are denotes there by $I(\iota,\nu)_{!}$).

\begin{prob}
Show that
$c(\cycl(f,M,J))=f_!(u(N))\in H(A)^{n}(M)$.
\end{prob}

We know that $$H\Z^{*}(BU)\cong \Z[c_{1},c_{2},\dots]$$
is a polynomial ring in the universal Chern classes. Hence
we can interpret
 $$u\in H(A)^{0}(BU)\cong A[[c_1,c_2,\dots]]^{0}$$
 as a formal power series in the universal Chern classes with coefficients in $A$. 
 Given a connection $\nabla$ on the complex vector bundle
 $N$ we can define the form
 \begin{equation}\label{aug0310}u(\nabla)=u(c_1(\nabla),\dots )\in \Omega A_{cl}^{0}(W)\end{equation}
by replacing the universal Chern classes by the corresponding Chern forms.
Observe that
$$\Rham(u(\nabla))=u(N)\ .$$

\begin{ex}\label{aug1020}{\em 
The following exercise illustrates the nature of $u$. Note that 
$$A=\pi_*(\bMU\C)\cong \C[[\C\P^{1}],[\C\P^{2}],\dots]$$
is a polynomial ring in  the complex bordism classes $[\C\P^{n}]\in \bMU_{2n}$ of the complex projective spaces.  \begin{prob}
Find an explicit formula for $u$ in terms of the generators $[\C\P^{n}]$ and $c_{k}$.
 \end{prob}
\proof Suitable tools for this calculation can be found in \cite{MR1189136}. \hB 
}\end{ex}

Let $W$ be a  $k$-dimensional manifold. Then we define the complex of distributional forms on $W$ as the topological dual $$\Omega^{*}_{-\infty}(W)=\Omega^{k-*}_c(W,\ori_W)^{\prime}\ ,$$ where $\Omega^{*}_c(W,\ori_W)$ denotes the complex  of Fr\'echet spaces of compactly supported smooth forms with coefficients in the orientation bundle of $W$

We fix a manifold $W$ and consider the sheaf $\Omega_{-\infty}$ of distributional forms on $W$.   Then we define the sheaf $$\Omega A_{-\infty}:=\Omega_{-\infty}\otimes \underline{A}$$
and let $\Omega A_{-\infty}(W)$ denote its complex of global sections.
\begin{prob}
Understand the difference between
$\Omega_{-\infty}(W)\otimes_{\C}A$ and $\Omega A_{-\infty}(W)$.
\end{prob}

 If $f:W\to M$ is an oriented proper smooth map and $n:=\dim(M)-\dim(W)$, then we can define a push-forward
$$f_!:\Omega^{*}_{-\infty}(W)\to  \Omega^{*+n}_{-\infty}(M)$$
as the adjoint of the pull-back
$$f^{*}:\Omega_c(M,\ori_M)\to \Omega_c(W,f^{*}\ori_M )\stackrel{\sim}{\to} \Omega_c(W,\ori_W)\ ,$$ where  the identification $f^{*}\ori_M\cong \ori_W$ is fixed by the orientation of $f$. We  extend the push-forward to the tensor product with $A$ $$f_!:\Omega A^{*}_{-\infty}(W)\to  \Omega A^{*+n}_{-\infty}(M)\ .$$

\begin{prob}
Fill in some details.
 \end{prob}

Let $(f,\cN,J)$ be a $\bMU$-cycle and
$\nabla$ be a  connection
on the underlying complex bundle  of $\cN$.
Then we get a form
$u(\nabla)\in \Omega A^{0}_{cl}(W)\subseteq \Omega A_{-\infty,cl}^{0}(W)$.

 We know that the inclusion
$$\Omega A(M)\to \Omega A_{-\infty}(M)$$ is a quasi-isomorphism.
We use this to extend the de Rham isomorphism to distributional forms
$$\Rham:H^{*}(\Omega A_{-\infty}(M))\stackrel{\sim}{\to} H(A)^{*}(M)$$
(see \cite{MR760450} for details). The following proposition is the bordism analog of  \ref{may2820}.
\begin{prob}
Show that 
$$\Rham(f_! u(\nabla))=c(\cycl(f,\cN,J))\ .$$
\end{prob}

\begin{ddd} \label{aug0510}
A geometric $\bMU$-cycle of degree $n$ is a tuple
$(f,\cN,J,\nabla,\eta)$, where $(f,\cN,J)$ is a $\bMU$-cycle of degree $n$,  $\nabla$ is a connection
on the underlying complex bundle  of $\cN$, and $\eta\in \Omega A^{n-1}_{-\infty}(M)/\im(d)$
is such that
$f_!(u(\nabla))-d\eta\in \Omega A_{cl}^{n}(M)$.
We let
$\Cycle_\bMU^{geom,*}(W)$ denote the graded semigroup of isomorphism classes of geometric $\bMU$-cycles.
\end{ddd}
The sum operation is given by disjoint union of the geometric pieces and the sum of the forms.

If $h:W^{\prime}\to W$, then we can define
$h^{*}(f,\cN,J,\nabla,\eta)$ if $f$ is transverse to $h$.
\begin{prob}
Fill in the details. 
Note that the pull-back of $\eta$ needs an argument.
\end{prob}
\proof
See  \cite{MR2550094} for details.

\begin{prob}\label{may2940}
Show that there exists unique additive differential refinement 
$$\hcycl:\Cycle_{\bMU}^{geom,*}\to \hbMU^{*}(M)$$
such that
$$I\circ \hcycl(f,\cN,J,\nabla,\eta)=\cycl(f,\cN,J)\ , \quad 
R(\hcycl(f,\cN,J,\nabla,\eta))=f_!(u(\nabla))-d\eta$$
which is compatible with the partially defined pull-back.
\end{prob}

\begin{ex}{\em 
Let $(M,g)$ be a closed compact Riemannian manifold. Let $(f,\cN,J)$ be an $\bMU$-cycle of degree $n$. Let us choose a connection $\nabla$ on $N$. Then there exists a unique
$\eta\in \Omega A^{n-1}_{-\infty}(M)$ which is orthogonal to $\ker(d)$ and such that
$f_{!}(u(\nabla))-d\eta$  is harmonic.  
We form
$$\hcycl(f,\cN,J,\nabla,\eta)\in \hbMU^{n}(M)\ .$$
\begin{prob}
Show that
$\hcycl(f,\cN,J,\nabla,\eta)$ does not depend on the choice of the connection $\nabla$.
\end{prob}
\proof 
Use a homotopy argument. \hB

\begin{ddd}
We call  this class 
$$\harm(f,\cN,J):=\hcycl(f,\cN,J,\nabla,\eta)\in \hbMU^{n}(M)$$
the harmonic class associated to the topological cycle 
$(f,\cN,J)$.
\end{ddd}

In the case of $\widehat{H\Z}$ harmonic differential characters have been considered in 
\cite{MR2518994}.

\begin{prob}
Show that
$$\harm:\cycl^{*}_{\bMU}(M)\to  \hbMU^{*}(M)$$
is additive.
\end{prob}

We consider the manifold $(S^{2},g)$ with the standard metric.
Let $\gamma:S^{1}\to S^{2}$ be a smooth embedding.
The normal bundle of $\gamma$ is trivialized. 
We get a representative
$$\cN:0\to TS^{1}\to f^{*}TS^{2}\oplus \R\to S^{1}\to S^{1}\times  \C\to 0$$
of the stable normal bundle with complex structure $J:=i$ on the quotient.
Then
$$\harm(\gamma,\cN,J)\in \hbMU^{1}(S^{2})\ .$$
For $x\in S^{2}$ we get
$$\harm(\gamma,\cN,J)_{|x}\in \hbMU^{1}(*)\cong \C/\Z \ .$$
\begin{prob}
Calculate the element
$\harm(\gamma,\cN,J)_{|x}\in \C/\Z$
for every $x\in S^{2}$.
\end{prob}

We consider the circle $(S^{1},g)$. A point
$x\in S^{1}$ gives naturally rise to  topological
cycle $c(x)\in \Cycle_{\bMU}^{1}(S^{1})$.
\begin{prob}
Describe the dependence of $\hcycl(c(x))$ on $x$.
\end{prob} 
}

\end{ex}

\subsection{Multiplicative structures}\label{jul3101}

In this section we give a general construction of a multiplicative differential extension of a multiplicative generalized cohomology theory.  It was a common belief that this is possible, but the details were an open problem for a while. Multiplicative extensions have been known in special cases
like ordinary cohomology, $\bKU$, $\bMU$ and other bordism theories, or Landweber exact theories.
See \cite{MR2550094} and \cite{2011arXiv1112.4173U}. 
In the present section we restrict attention to the commutative case which captures most examples. There should be an analogous theory for the associative case. 

 The input for the construction is a multiplicative refinement of  differential data. 
\begin{ddd}
A multiplicative differential data is a differential data
$(E,A,c)$ where
\begin{enumerate}
\item  $E\in \CommMon(\Nerve(\Sp)[W^{-1}]^{\wedge})$ is a commutative ring spectrum, \item $A\in \CommMon(\Ch^{\otimes})$ is a commutative DGA over $\C$, and \item $c:E\C\to H(A)$ is a morphism of commutative ring spectra.
\end{enumerate}
A morphism
of differential data $(E,A,c)\to (E^{\prime},A^{\prime},c^{\prime})$ is a commutative diagram in  $\CommMon(\Nerve(\Sp)[W^{-1}]^{\wedge})$
$$\xymatrix{E\ar[r]^{f}\ar[d]^{c}&E^{\prime}\ar[d]^{c^{\prime}}\\
H(  A)\ar[r]^{H(\phi)}&H(  A^{\prime})}\ ,$$
for some morphisms $\phi:E\to E^{\prime}$ in  $\CommMon(\Nerve(\Sp)[W^{-1}]^{\wedge})$ and $f:A\to A^{\prime}$ of DGA's.
\end{ddd}

As before on can
 describe the $\infty$-category of strict multiplicative differential data as a pull-back in $\infty$-categories $\infty\Cat$
$$\xymatrix{\Data^{str,mult} \ar[r]\ar[d]&\Nerve(\CommMon(\Ch)^{\otimes})\ar[d]^{H}\\
\CommMon(\Nerve(\Sp)[W^{-1}]^{\wedge})\ar[r]^{\dots\wedge M\C}&\Mod(H\C)}  \ .$$

Before we turn to the construction of a multiplicative differential cohomology theory in terms of a multiplicative differential function spectrum we first define what we mean by a multiplicative differential extension of a multiplicative cohomology theory. In this way we capture examples constructed by other means.

\begin{ddd}
A multiplicative differential extension of the multiplicative cohomology theory
$E^{*}$ associated to a multiplicative data $(E,A,c)$ is a tuple 
$(\hat E^{*},R,I,a)$ which is a differential extension of $E$ such that $\hat E^{*}$ has values in graded commutative rings, $R$ and $I$ are transformations of ring-valued functors, and for  $x\in \Omega A (M)/\im(d)$ and $y\in \hat E (M)$ we have  $$a(x)\cup y=a(x\wedge R(y))\ .$$
\end{ddd}

Note that the notion of a multiplicative differential extension only depends on a homotopy theoretic part of the data $(E,A,c)$, namely the multiplicative cohomology theory $E^{*}$,  the differential graded algebra $A$, and the transformation of cohomology theories $E^{*}\to H(A)^{*}$ induced by $c$. Therefore we can talk about multiplicative extensions   having only fixed this coarser part datum. This is the way how previous definitions are related to the one given here.

Let $E\in \CommMon(\Nerve(\Sp[W^{-1}])$ be a commutative ring spectrum.
It is a natural question when the canonical strict data (Definition \ref{jul1901})
$(E,A,c)$ with $A:=\pi_*(E)\otimes \C$ can be refined to a multiplictive data.
\begin{ddd}
We say that  $E$ is formal (over $\C$), if there exists an morphism of ring spectra
$c:E\to H(A)$ refining the canonical morphism.
\end{ddd}

Most of our examples are formal:
$H\Z$, $\bMU$ (see \ref{jul1601}), $\bKU$ (see  \ref{jul1902}).
One can also show that the algebraic $K$-theory spectrum of a number ring is formal, see \cite{buta2}.

Observe that $(H\Z,\C,H\Z\to H\C)$ is a multiplicative data and $(\hat H^{*}_{Del},R,I,a)$
is an associated  multiplicative extension.

In the following we give a general construction of multiplicative extensions.
If $A$ is a commutative DGA, then $\Omega A$ becomes a sheaf of commutative DGA's.  \begin{lem}\label{may2909}
The de Rham equivalence refines to an equivalence
$$\Rham:H(\Omega A)\to \Funk(H(A))$$
in $\Sm(\CommMon(\Nerve(\Sp)[W^{-1}]^{\wedge}))$.
\end{lem}
\proof
In principle we can repeat the proof of Proposition \ref{may2901}. The main
point is to observe that the integration map
$$\int:\Omega A(M)\to \underline{\Map}(C_*(M),A)$$
refines to an equivalence in $\CommMon(\Ch[W^{-1}]^{\otimes})$.
This seems to be true (compare \cite{2010arXiv1011.4693A}) but I do not know a good reference
for this precise statement.

Therefore we give an alternative argument which does not involve integration of forms.
We know that $$H(\Omega A)\in \Sm^{desc}(\CommMon(\Nerve(\Sp)[W^{-1}]^{\wedge}))$$ is constant.
We have a natural morphism
$$p^{*}H(\Omega A(*))\to H(\Omega A)$$ in $ \Sm(\CommMon(\Nerve(\Sp)[W^{-1}]^{\wedge}))$,
where $p:\Mf\to *$ is as in \ref{may2905}.
Hence the natural morphism
 \begin{equation}\label{may2906}L(p^{*}(H(A)))=L(p^{*}(H(\Omega A(*))))\to H(\Omega A)\end{equation}
is an equivalence in $\Sm^{desc}(\CommMon(\Nerve(\Sp)[W^{-1}]^{\wedge}))$.
On the other hand, in this category we have the equivalence (see Problem \ref{may2908})
\begin{equation}\label{may2907}L(p^{*}(H(A)))\stackrel{\sim}{\to} \Funk(H(A))\ .\end{equation} We therefore get an equivalence
$$\Rham: H(\Omega A)\xrightarrow{\eqref{may2907}\circ \eqref{may2906}^{-1},\sim}
\Funk(H(A))\ .$$
\hB

\begin{prob}
Verify that the equivalences obtained in \ref{may2901} and \ref{may2909} conincide.
\end{prob}

Note that the cone in Definition \ref{may2902} can be written as a pull-back. 
The presentation of the differential function spectrum can be refined to the multiplicative case by interpreting that pull-back in smooth commutative ring spectra. We will actually refine $$\Diff^{\bullet}(E,A,c):=\bigvee_{n\in \Z} \Diff^{n}(E,A,c)$$ to a commutative ring spectrum.

 Via the sequence of canonical  symmetric monoidal functors
\begin{equation}\label{sep2760}\Nerve(\Set)\stackrel{\iota}{\to}     \Nerve(\sSet)[W^{-1}]\stackrel{\Sigma^{\infty}_+}{\to}\Nerve(\Sp)[W^{-1}] \end{equation} the abelian group $\Z\in \CommMon( \Nerve(\Set )^{\times})$
gives rise to a commutative ring spectrum 
$$\Sigma^{\infty}_+\iota(\Z)\in \CommMon(\Nerve(\Sp)[W^{-1}]^{\wedge})\ .$$

For a commutative ring spectrum  $R$ we define the commutative ring spectrum $$R[z,z^{-1}]:=R\wedge \Sigma^{\infty}_+\iota(\Z)\ .$$ Note that this wedge product is the coproduct in the category of commutative ring spectra.
If one forgets the commutative ring structure, then 
$R[z,z^{-1}]\cong \bigvee_{n\in \Z} R$.

\begin{prob}\label{sep2761}
Show that there is a canonical equivalence  $$\Sigma^{\infty}_+\iota(\Z)\wedge H\Z\to H(\Z[z,z^{-1}])$$
in $\CommMon(\Mod(H\Z)^{\wedge_{H\Z}})$.
\end{prob}
\proof
We first construct a map
$\Sigma^{\infty}_+\iota(\Z)\to  H(\Z[z,z^{-1}])$.
The functor $\Sigma^{\infty}_+$ fits into an adjunction
$$\Sigma^{\infty}_+:\CommMon( \Nerve(\sSet)[W^{-1}]^{\times} )\leftrightarrows \CommMon(\Nerve(\Sp)[W^{-1}]^{\wedge}):\Omega^{\infty}\ . $$
Therefore it suffices to construct a map
$$\iota(\Z)\to \Omega^{\infty}H(\Z[z,z^{-1}])$$
in $\CommMon(\Nerve(\sSet)[W^{-1}]^{\times})$.
We further have an adjunction
$$\pi_0: \CommMon(\Nerve(\sSet)[W^{-1}]^{\times})\leftrightarrows \CommMon(\Nerve(\Set)^{\times}):\iota\ .$$
 The unit of  this adjunction provides the vertical map in the diagram
$$\xymatrix{&\Omega^{\infty }H(\Z[z,z^{-1}])\ar[d]\\\iota(\Z)\ar@{.>}[ur]\ar[r]^{n\mapsto z^{n}}&\iota(\Z[z,z^{-1}])}\ .$$
This vertical map is an equivalence so that we can choose the lift
$\iota(\Z)\to \Omega^{\infty }H(\Z[z,z^{-1}])$ in 
 $\CommMon(\Nerve(\sSet)[W^{-1}]^{\times})$.
 
By extension of coefficients we get a map
$$ \Sigma^{\infty}_+\iota(\Z)\wedge H\Z\to H(\Z[z,z^{-1}])$$
in $\CommMon(\Mod(H\Z)^{\wedge_{H\Z}})$. It induces an isomorphism in homotopy and therefore is an equivalence. 
\hB

Assume that $A$ is a commutative differential graded algebra. Then we can form the tensor product $A[z,z^{-1}]:=A\otimes \Z[z,z^{-1}]$. It is   a coproduct in commutative graded
algebras. Since the equivalence $H$ preserves coproducts we get an equivalence
\begin{equation}\label{sep2701}H(A[z,z^{-1}])\simeq H(A)\wedge_{H\Z} H(\Z[z,z^{-1}]) \stackrel{\ref{sep2761}}{\simeq} H(A)[z,z^{-1}] .\end{equation}

We consider the presheaf of complexes
  $$\sigma^{\ge \bullet} \Omega A:=\bigoplus_{n\in \Z} \sigma^{\ge n} \Omega A\ .$$
If $A$ is a commutative DGA, then $\sigma^{\ge \bullet} \Omega A$ is a commutative DGA as well such that the natural embedding
\begin{equation}\label{sep27022}\bigoplus_{n\in \Z} \sigma^{\ge n} \Omega A\subseteq  \bigoplus_{n\in \Z} \Omega A\cong \Omega A[z,z^{-1}]\end{equation}

preserves the differential graded algebra structures.
In particular, we have a map of presehaves commutative ring spectra
\begin{equation}\label{sep2703}H(\sigma^{\ge \bullet} \Omega A)\to H(\Omega A[z,z^{-1}])\ .\end{equation}

\begin{ddd}\label{aug0602}
We define $\Diff^{\bullet}(E,A,c)\in \Sm(\CommMon(\Nerve(\Sp)[W^{-1}]^{\wedge}))$
as the pull-back
$$\xymatrix{\Diff^{\bullet}(E,A,c)\ar[r]\ar[d]&H(\sigma^{\ge \bullet} \Omega A)\ar[d]^{\Rham}\\
\Funk(E[z,z^{-1}])\ar[r]^{c}&\Funk(H(A)[z,z^{-1}])}\ .$$ 
 \end{ddd}
 The left vertical map involves   \eqref{sep2703}, de Rham equivalence, an the identification \eqref{sep2701}.
Note that
$$\pi_n (\Diff^{\bullet}(E,A,c))=\bigoplus_{k\in \Z} \pi_n(\Diff^{k}(E,A,c))$$
forms a bigraded ring.
In particular, by taking the diagonal, we get the graded ring-valued functor
$$\hat E^{*}:=\bigoplus_{k\ge \Z}   \pi_k(\Diff^{k}(E,A,c))  \ .$$
By construction, the transformations
$$R:\hat E^{*}\to \Omega A_{cl}^{*}\ , \quad I:\hat E^{*}\to E^{*}$$
become transformations between graded commutative ring valued functors.

\begin{prob}
Show that for $x\in \Omega A^{n-1}(M)/\im(d)$ and $y\in \hat E^{m}(M)$ 
$$a(x)\cup y=a(x\cup R(y))\ .$$
\end{prob}

\begin{ex}{\em
If we apply the above construction to the data $(H\Z,\R,c)$, where $c:H\Z\to H\C$ is the canonical map, then we get
the multiplicative structure on $\widehat{H\Z}^{*}$.
\begin{prob}
Show that the canonical isomorphism $\hat H_{Del}^{*}
\cong \widehat{H\Z}$ from Problem \ref{may2911} is 
multiplicative.\end{prob}
}\end{ex}

\begin{ex}{\em
The sphere spectrum $\bS$ is a commutative ring spectrum which comes with a canonical map of ring spectra $c:\bS\to H\C$. We get the differential function ring spectrum $\Diff^{\bullet}(\bS,\C,c)$.
Note that every spectrum is an $\bS$-module.
\begin{prob}\label{jun06101}
Let $(E,A,c)$ be a multiplicative data.
Work out the definition of a module data $(F,B,d)$.
Construct
$\Diff(F,B,d)$ as a $\Diff(E,A,c)$-module.
\end{prob}

\begin{prob}\label{jun0610}
 Show that for every differential data $(E,A,c)$
we get a  $\Diff^{\bullet}(\bS,\C,c)$-module structure on $\Diff^{\bullet}(E,A,c)$.
\end{prob}

}\end{ex}

\begin{ex}{\em
We can fix an isomorphism of commutative ring spectra $$c:\bKU\C\cong H(\C[b,b^{-1}])\ .$$
\begin{prob}\label{jul1902}
Show this!
\end{prob}
\proof
We first consider connective $K$-theory $\bku$. Similarly as in Problem \ref{jul1601}
we produce an equivalence of ring spectra
$$\bku\wedge H(\C)\cong  H(\C[b])\ .$$
It is useful to write the localization $\iota:\Nerve(\Ch)\to \Nerve(\Ch)[W^{-1}]$ explicitly.
So more precisley we should write
$\bku\wedge H(\C)\cong  H(\iota(\C[b]))$.
Then we invert the multiplication by $b$. We get
$$\bKU\cong \bku[b^{-1}]\cong H(\iota(\C[b]))[b^{-1}]\cong H(\iota(\C[b])[b^{-1}])\ .$$ 
Here $\iota(\C[b])[b^{-1}]$ is the inversion of $b$ in the $\infty$-category algebras in over $\iota(\C[b])$. As a final step one has to check that $\iota(\C[b])[b^{-1}]\cong \iota(\C[b,b^{-1}])$. 
\hB

Alternatively one can use the existence of a strictly multiplicative Chern form
for vector bundles with connection. Details can be found in \cite[Sec. 3.6]{buta}.

The datum
$(\bKU,\C[b,b^{-1}],c)$ is a strict multiplicative differential datum.
We therefore get a multiplicative differential extension
$(\hbKU,R,I,a)$.
\begin{prob}
Show that there exists a unique multiplicative differential extension
associated to the datum described above. Observe that the extension is not unique if one
does not require multiplicativity.
\end{prob}
\proof
Compare with \cite{MR2608479}. \hB 

Recall that $\Vect^{geom}$ is a semiring valued functor with product induced by the tensor product.
\begin{prob}
Show that the geometric cycle map constructed in Theorem \ref{may2930}
$$\hcycl:\Vect^{geom}\to \hbKU^{0}$$
is a natural transformation of semiring valued functors.
\end{prob}
\proof
Observe that the deviation from multiplicativity
gives rise to a natural transformation
$$\Vect\times \Vect\to \bKU\C/\Z^{-1}\ .$$
Now use that
$\bKU\C/\Z^{-1}(BU(n)\times BU(m))=0$ for all $n,m\ge 0$. \hB 
}
\end{ex}

\begin{ex}\label{aug2050}{\em
The following exercise concerns aspects of the $\hat \bS^{*}$-module structure found in \ref{jun0610}.
Let $\hat \epsilon:\hat \bS^{*}\to \hbKU^{*}$ be the unit.
We consider the tautological bundle $L\to \C\P^{1}\cong S^{2}$ with its invariant connection $\nabla$.
Then we get a class $\hcycl(L,\nabla)\in \hbKU^{0}(S^{2})$. Note that
$\pi_{1}(\bS)\cong \Z/2\Z\cong \bS^{-1}(*)$.
The generator of this group is the class $\eta$ of the framed manifold $S^{1}$ with the non-bounding  framing coming from the group structure.
 It has a canonical lift
$\hat \eta\in \hat \bS^{-1}(*)$. Let $p_{S^{k}}:S^{k}\to *$ be the projection.
\begin{prob}
Calculate  $\hat \epsilon(\hat\eta)\in  \widehat{\bKU}^{-1}\cong \C/\Z$  and $p_{S^{2}}^{*}\hat \eta\cup  \hcycl(L,\nabla^{L})\in \hbKU^{-1}(S^{2})$.
\end{prob}
\proof
 Similarly, we consider the tautological $\bH$-bundle $E\to \bH\P^{1}\cong S^{4}$ with its invariant connection $\nabla^{E}$.
We get the element
$$\hcycl(L,\nabla)\in \hbKU^{0}(S^{4})\ .$$ The group $SU(2)$ with its right-invariant framing
represents a generator $$\sigma\in \pi_{3}^{s}=\bS^{-3}(*)\cong \Z/24\Z\ .$$
It again lifts canonically to $\hat \sigma\in \hat \bS^{-3}(*)$.
\begin{prob}\label{jul2111}
Calculate $\hat \epsilon(\hat\sigma)\in  \widehat{\bKU}^{-3}\cong \C/\Z$ and $p_{S^{4}}^{*}\hat \sigma\cup  \hcycl(E,\nabla^{E})\in \hbKU^{-3}(S^{4})$.
\end{prob}
\proof
In the following we show how to reduce the calculation of $\hat \epsilon:\hat \bS^{*}(*)\to \hbKU^{*}(*)$ to a topological problem. For $k\ge 1$ we have an isomorphism
$$\bS \C/\Z^{-k-1}(*)\stackrel{\sim}{\to}\hat \bS^{-k}_{flat}(*)\stackrel{I,\sim}{\to}\bS^{-k}(*)\ .$$
Let $u\in \bS^{-k}(*)$, $\hat u\in \hat \bS^{-k}(*)$,  and $u_{\C/\Z}\in   \bS \C/\Z^{-k-1}(*)$ be its preimages under these isomorphisms.
Then $\hat \epsilon(\hat u)=\epsilon_{\C/\Z}(u_{\C/\Z})\in \bKU\C/\Z^{-k-1}(*)\subseteq \hbKU^{k}(*)$.
The map
$$e:\pi_{k}(\bS)\cong \bS^{-k}(*)\to \bKU\C/\Z^{-k-1}(*)\ , \quad  u\mapsto \epsilon_{\C/\Z}(u_{\C/\Z})$$ is the complex variant of Adam's $e$-invariant (cf \cite{2011arXiv1103.4217B} and \ref{jul2110} for more details).  It is known that 
$e(\eta)$ has order $2$ and $e(\sigma)$ has order $12$.
We conclude that
$\hat \epsilon(\hat \eta)$ has order $2$ and $\hat \epsilon(\hat \sigma)$ has order $12$. 
We refer to \ref{jul3140} for the calculations.

 }
\end{ex}
 
  \begin{ex}{\em
Let $S\subset \nat$ be multiplicatively closed. Recall the construction of the Adams operations from Example \ref{jun0210}.
\begin{prob}
Show that there exists a unique multiplicative differential extension of $\hbKU[S^{-1}]$ and that the Adams operations $\Psi^{k}$, $k\in S$, extend to multiplicative operations $\hat\Psi^{k}$ on $\hbKU[S^{-1}]$.
\end{prob}
}
\end{ex}

\begin{ex}{\em
The data $(\bMU,A,\bMU\C\cong H(A))$ is a multiplicative. 
\begin{prob}\label{jul1601}
Show this!
\end{prob}
\proof
The point is to show that the map $c:\bMU\to H(A)$ refines to a morphism of commutative ring spectra. Note that $$A\cong \pi_*(\bMU\wedge H(\C))\cong \pi_0(\map_{\Mod(H(\C))}(H(\C),\bMU\wedge H(\C)))$$ is a free commutative $\C$-algebra in generators $(x_i)_{i\in \nat}$, cf.  \ref{aug1020}. We use the last incarnation of the generators to construct a map of $H\C$-modules
$$M:=H(\C)\langle x_1,\dots\rangle\cong \bigsqcup_{i\in \nat} H(\C)[-\deg(x_i)]\to \bMU\wedge H\C\ .$$
It induces an equivalence of commutative algebras in $\Mod(H(\C))$
$$\Free_{\Mod(H(\C))}(M)\stackrel{\sim}{\to} \bMU\wedge H\C\ .$$
Since $\C$ is a $\Q$-algebra the classical free commutative $\C$-algebra over the graded vector space 
 $\C\langle x_1,\dots\rangle$ (which is $A$) coincides with the  free commutative $\C$-algebra taken in the $\infty$-category $\Mod(\C)$. 
Hence
$$\Free_{\Mod(H(\C))}(M)\cong H(\Free_{\Mod(\C)} \C\langle x_1,\dots\rangle)\cong H(A)\ .$$ \hB

\begin{prob}
Use the arguments  given in the proofs of  \ref{jul1601}   and \ref{jul1902}
to show that a commutative ring spectrum $E$ is formal over $\C$ provided one of the following conditions holds true:
\begin{enumerate}
\item $\pi_{*}(E)\otimes \C$ is a free commutative algebra.
\item $E$ is periodic and $\pi_{*}(E[0,\dots,\infty])\otimes \C$  is a free commutative algebra,
where $E[0,\dots,\infty]$ denotes the connective cover of $E$.
\end{enumerate}
\end{prob}

The multiplicative data $(\bMU,A,c)$ gives rise to a multiplicative extension
$(\hbMU^{*},R,I,c)$.
The graded semigroup
$\Cycle_{\bMU}^{geom}(W)$ has a partially defined product.
More precisely,  $(f,\cN,J,\nabla,\eta)\cup (f^{\prime},\cN^{\prime},J^{\prime},\nabla^{\prime},\eta^{\prime})$
can be defined if $f$ and $f^{\prime}$ are transverse. Details can be found in \cite{MR2550094}.
\begin{prob}
Fill in the details. Show that the differential cycle map
$\hcycl$ obtained in \ref{may2940} is multiplicative.
\end{prob}
\proof
Consider the deviation from multiplicativity.
Use the compatibility of the cycle map with the curvature to show that it induces a transformation
$\bMU^{l}\otimes \bMU^{k}\to \bMU \C/\Z^{k+l-1}$. Then argue that such a transformation is zero
since $\bMU$ is an even spectrum.

\begin{prob}
Assume that $(M,g)$ is a compact K\"ahler manifold or a compact symmetric space.
Show that the harmonic cycle map
$$\harm:\Cycle_{\bMU}^{*}(M)\to \hbMU^{*}(M)$$
is multiplicative.
\end{prob}
\proof
Use formality, or more precisely, that the product of harmonic forms is harmonic.
\hB

 Let $\gamma_{0},\gamma_{1}:S^{1}\to L^{3}_{\Z/p\Z}$ define a link.
Observe that they have canonical refinements to topological cycles
$c(\gamma_{i})\in \Cycle_{\bMU}^{2}(S^{3})$.
Let $\hat c(\gamma_{i})\in \Cycle_{\bMU}^{geom,2}$ be geometric refinements.
\begin{prob}
Calculate
$$\harm(\hat c(\gamma_{0}))\cup \harm(\hat c(\gamma_{1}))\in \bMU^{4}(L^{3}_{\Z/p\Z})\cong \C/\Z\ .$$
  \end{prob}
 
}\end{ex}

\begin{ex}\label{aug1021}{\em
In this example we use the notation $A_{E}:=\pi_{*}(E)\otimes \C$.
We have a morphism of commutative ring spectra $\alpha:\bMU\to \bKU$.
It induces a morphism of rings $A_{\bMU}\to A_{\bKU}$ and a morphism of multiplicative
data
$$\hat \alpha:(\bMU,A_{\bMU},c_{\bMU})\to (\bKU,A_{\bKU},c_{\bKU})\ .$$
 Note that the formal group law over $\pi_{*}(\bKU)$ associated to $\alpha$ is Landweber exact \cite{MR0423332}.
 As a consequence we have the Conner-Floyd theorem \cite{MR0216511},\cite{MR1166518} stating that for a compact manifold $M$ the map $\alpha$ induces an isomorphism
 $$(\bMU(M)^{\cdot}\otimes_{\bMU^{\cdot}} \bKU^{\cdot})^{*}\stackrel{\sim}{\to}\bKU^{*}(M)\ .$$
This extends to differential cohomology. Note that $\hbMU(M)^{\cdot}$ is a $\hbMU^{\cdot}(0)\cong \bMU^{\cdot}$-module.
\begin{prob} Show that if $M$ is compact, then $\hat \alpha$ induces an isomorphism
$$(\hbMU(M)^{\cdot}\otimes_{\bMU^{\cdot}} \bKU^{\cdot})^{*}\stackrel{\sim}{\to}\hbKU^{*}(M)\ .$$
\end{prob}
\proof See \cite{MR2550094} for an argument. \hB 

Let $*\to S^{2}$ be the inclusion of a point. It extends canonically to a cycle $*\in \Cycle_{\bMU}^{2}(S^{2})$. 
\begin{prob}
Calculate $\hat \alpha(\harm(*))\in \hbKU^{2}(S^{2})$. More precisely, find a vector bundle $(V,\nabla)$ with connection such that $b\ \cycl(V,\nabla)=\hat \alpha(\harm(*))$.
\end{prob}
}\end{ex}

\subsection{Relative sites}

In this subsection we introduce some language which is mainly used to capture, in the $\infty$-categorical context, the functoriality of  construction like evaluations with proper support and
integration maps.  If one is solely interested in the latter constructions for a specific map of manifolds, then one should skip most of this subsection and just look how these constructs look like when evaluated on a specific map.

The main application of this theory is the definition of the notion of differential cohomology $\hat E^{*}_{K}(W)$ with support in a closed subset $K\subset W$, or with proper support $\hat E^{*}_{prop/M}(W)$ over a map $f:W\to M$.
Relative versions of $\widehat{H\Z}$ in the Cheeger-Simons and Hopkins-Singer pictures have been considered  and compared in \cite{MR2255012}. See also \cite{MR1976955}.

For a manifold $B$ we let $\Mf/B$ be the site of manifolds over $B$ equipped with the topology induced from $\Mf$.
We write objects in the form $(M\to B)$. A morphism $(M\to B)\to (M^{\prime}\to B)$ is a smooth map $f:M\to M^{\prime}$ which  preserves the structure maps to $B$.
We let $$\Sm_{B}(\bC):=\Fun(\Nerve(\Mf/B)^{op},\bC)$$
denote the $\infty$-catgeory of smooth objects in $\bC$ over $B$.
In the present subsection we generally assume that $\bC$ is a presentable $\infty$-category.

\begin{ex}\label{aug1030}{\em 
Here  is the  typical examples which one should have in mind. Let $V\to B$ be a complex vector bundle. Then we  obtain a smooth abelian group
$\Gamma(V)\in \Sm^{desc}_{B}(\Nerve(\Ab))$ which associates to $f:M\to B$ the space
of sections of the pull-back $f^{*}V\to M$.
Assume further that $V$ has a flat connection $\nabla$. Then we could form the
sheaf of chain complexes
$$\Omega(V,\nabla)\in \Sm^{desc}_{B}(\Nerve(\Ch)[W^{-1}])$$
which associates to $f:M\to B$ the twisted de Rham complex
$\Omega(M,f^{*}V)$ with differential $d\otimes 1+1\otimes f^{*}\nabla$.
}\end{ex}

For a smooth map $f:W\to B$ we have a functor $f_{\sharp}:\Mf/W\to \Mf/B$
given by $$f_{\sharp}(M\to W):=(M\to W\stackrel{f}{\to} B)\ .$$ It induces an adjunction
$$f^{*}:=(f_{\sharp})^{*}:\Sm_{B}(\bC)\leftrightarrows  \Sm_{W}(\bC):f_{*}\ .$$
Note that $f^{*}X(M\to W)\cong X(M\to W\to B)$. 
In order to describe the right-adjoint $f_{*}$ we consider the category $\Mf/(W\to B)$ of diagrams of the form
\begin{equation}\label{jul0110}D:=\xymatrix{V\ar[r]\ar[d]&W\ar[d]\\M\ar[r]&B}\ .\end{equation} 
A morphism $D\to D^{\prime}$ is a pair of smooth maps $V\to V^\prime$ and $M\to M^{\prime}$
which preserve the structure maps.
 We have forgetful functors
$$p:\Mf/(W\to B)\to \Mf/W\ , \quad p(D):=(V\to W)$$
and
$$q:\Mf/(W\to B)\to \Mf/B\ , \quad q(D):=(M\to B)\ .$$
Both, $p$ and $q$ induce adjunctions 
$$p^{*}: \Fun(\Nerve(\Mf/W)^{op} ,\bC)\leftrightarrows \Fun(\Nerve(\Mf/(W\to B))^{op} ,\bC) :p_{*}$$ and
$$q^{*}: \Fun(\Nerve(\Mf/B)^{op} ,\bC)\leftrightarrows \Fun(\Nerve(\Mf/(W\to B))^{op} ,\bC) :q_{*}\ ,$$ 
and we have 
$$f_{*}=q_{*}p^{*}\ .$$
It is now easy to see that
\begin{equation}\label{jul2201}(f_{*}X)(M\to B)=\lim_{(\Mf/(W\to B))/(M\to B)} X(V\to W)\ .\end{equation}

\begin{prob}
Show that $f^{*}$ and $f_{*}$ preserve the subcategories of sheaves.
\end{prob}
\proof
The assertion is clear for $f^{*}$. For $f_{*}$ we can use formula  \eqref{jul2201}. \hB

\begin{ex}{\em
Let $p:B\to *$ be the canonical map. It gives rise to a functor
$p^{*}:\Sm(\bC)\to \Sm_{B}(\bC)$. For example, for $\Omega_{\C}\in \Sm(\Nerve(\Ch)[W^{-1}])$
we get $p^{*}\Omega(M\to B)=\Omega(B)$.  
Continuing example \ref{aug1030}, we have  natural isomorphisms
$$f^{*}\Gamma(V)\cong \Gamma(f^{*}V)\ , \quad f^{*}\Omega(V,\nabla)\cong \Omega(f^{*}V,f^{*}\nabla)\ .$$
If $X\in \Sm(\bC)$, the we often use the notation $X_{B}:=p^{*}X\in \Sm_{B}(\bC)$.
}\end{ex}

 If $f:W\to  B$ is a submersion, then we have a functor
 $f^{\sharp}:\Mf/B\to  \Mf/W$ given by 
 $f^{\sharp}(M\to B):=M\times_{B}W\to W$.
\begin{prob}
Show that for a submersion $f$ we have
$f_{*}(X):=(f^{\sharp})^{*}$. In particular, $f_{*}$ is exact.
\end{prob}
\proof  In this case
$(\Mf/(W\to B))/(M\to B)$
has a final object 
$$\xymatrix{M\times_{B}W\ar[r]\ar[d]&W\ar[d]\\M\ar[r]&B}\ .$$ \hB

We now assume that $\bC$ is pointed by the zero object $*$, i.e. an object which is both, final and initial. Our examples are stable $\infty$-categories, but also $\Nerve(\Ch)$.
We have a functor
$$\fibre:\Fun(\Nerve([1]),\bC)\to \bC$$ which takes
 the fibre $\fibre(f)$ of a map $f:C\to D$ at $*\to D$.
 
 \begin{ex}{\em
 Let $\iota:\Nerve(\Ch)\to \Nerve(\Ch)[W^{-1}]$ be the localization and
 $A\to B$ an object of $\Fun(\Nerve([1]),\Nerve(\Ch))$.  
 \begin{prob}\label{jul2501}
 Show that there is a natural morphism
 $$\phi:\iota(\fibre(A\to B))\to\fibre(\iota(A)\to \iota(B))\ .$$
 Analyse when it is an equivalence.
 \end{prob}
 \proof
 Since $\iota$ preserves the zero object we get
 the natural map $\phi$ by the universal property of the pull-back defining the kernel.
 It is an equivalence for example, if $A\to B$ is level wise surjective.
 
 Note 
$\fibre(A\to B)$ is given by the complex of level wise kernels of the map $A\to B$.
On the other hand,
$\fibre(\iota(A)\to \iota(B))$ is represented by 
the cone $\Cone(A\to B)[-1]$.
The natural map $\phi$  is now represented by the canonical inclusion
$\fibre(A\to B)\to \Cone(A\to B)[-1]$.
If $A\to B$ is level wise surjective, then we get a map of long exact sequences which shows that the map  $\phi$ is a quasi-isomorphism.
\hB 
 }
 \end{ex}
Using these constructs
 we can define  relative versions of smooth  objects.
Let $K\subset M$ be closed and $X\in \Sm(\bC)$.
\begin{ddd}\label{jul0203}
We define the evaluation of $X$ with support in $K$ by
$$X_K(M)=X(M,M\setminus K):=\fibre(X(M)\to X(M\setminus K))\ .$$
\end{ddd}

\begin{ex}{\em 
We get 
$$(\Omega^{n}_{\C})_{K}(M)\cong \{\omega\in \Omega^{n}(M,\C)\:|\:\omega_{|M\setminus K}=0\}=\{\omega\in \Omega^{n}(M,\C)\:|\:\supp(\omega)\subseteq K\}\ .$$
If we consider $\Omega_{\C}\in \Sm(\Nerve(\Ch)[W^{-1}])$, then
$(\Omega_{\C})_{K}$ is not so simple since we must take the kernel in the $\infty$-categorical sense.
We can not simply take all forms supported in $K$.
See \ref{jul2401} for similar discussion.
}\end{ex}
Definition \ref{jul0203} gives the evaluation of a functor $X_{K}$ on objects $X_{K}(M)$. In order to fully construct this functor we proceed as follows.
Consider the category $\Mf_{rel}$ of pairs $(M,U)$ of manifolds and open subsets $U\subseteq M$. A morphism $(M,U)\to (M^{\prime},U^{\prime})$ is a smooth map $f:M\to M^{\prime}$ such that
$ f(U) \subseteq  U^{\prime} $.
\begin{prob}\label{jul0201}
Define  the extension
$$\cE:\Sm(\bC)\to \Fun(\Nerve(\Mf_{rel})^{op},\bC)$$
which on objects acts as $\cE(X)(M,U)=X(M,U)$.
\end{prob}
\proof
We have a forgetful functor $r:\Mf_{rel}\to \Mf$, $(M,U)\mapsto M$.
Furthermore, we have a functor
$$d:\Fun(\Nerve(\Mf_{rel})^{op},\bC)\to \Fun(\Nerve(\Mf_{rel})^{op}, \Fun(\Nerve([1]),\bC))$$
 which maps $X$ to the diagram
 $X(M,\emptyset)\to X(U,\emptyset)$. Then we can write
 $$\cE:=\fibre\circ d\circ r^{*}:\Sm(\bC)\to \Fun(\Nerve(\Mf_{rel})^{op},\bC)\ .$$
 \hB

\begin{prob}
Assume that $X$ satisfies descent, $K\subseteq M$ is closed, $V\subseteq M$ is open and contains $K$. Then we have excision:
The natural map
$X_K(M)\to X_K(V)$ is an equivalence.
\end{prob}
 
\begin{ddd}
We define the compactly supported evaluation
of $X$ by
$$X_c(M):=\colim_K X_K(M)\ ,$$
where the colimit is taken over all compact subsets $K\subseteq M$.
\end{ddd}

\begin{prob}
Define a compactly supported extension functor
 $$\cE_c:\Sm(\bC)\to  \Sm(\bC)$$
 such that
 $$\cE_c(X)(M)\cong X_c(M)\ .$$
 \end{prob}
 
 \begin{ex}{\em
We get 
$$\cE_{c}(\Omega^{n}_{\C})(M)=\{\omega\in \Omega^{n}(M,\C)\:|\:\mbox{$\supp(\omega)$ is compact}\}\ .$$}
 \end{ex}

 \begin{ex}\label{aug1045}{\em 
 We write $\pi_{*}(\Funk(E)_{c})=:E^{-*}_{c}$ for the compactly supported $E$-cohomology functor. This notation coincides with the common usage, but conflicts the use of the symbol $c$ otherwise in the present paper.  The group $E^{-*}_{c}(M)$ is   in general not the group of cohomology classes in $E^{-*}(M)$ which vanish outside of some compact subset.
\begin{prob} Show that
for every $k\in \Z$ there exists a natural isomorphism
\begin{equation}\label{jun0401}E^{k}_{c}(\R^{n})\cong E^{k-n}(*)\ .\end{equation}
\end{prob}
}
\end{ex}

%

We now define the push-forward with proper support along a smooth map $W\to B$.
We let $\tilde \Mf/(W\to B)$ be the  category of objects
$(D,K)$ with   $D$ a diagram as in \eqref{jul0110} and  $K\subset V$ closed such that the induced map $K\to M$ is proper.
 Morphisms   $f:D\to D^{\prime}$ are structure preserving maps such that
$ f(M\setminus K)\subseteq M^{\prime}\setminus K^{\prime}$.  
We define
$$\tilde p^{*}: \Fun(\Nerve(\Mf/W)^{op},\bC)\to \Fun(\Nerve(\tilde \Mf/(W\to B))^{op},\bC)$$ in a natural way so that on objects
$$\tilde p^{*}(X)(D,K):=X_{K}(V\to W)\ .$$
\begin{prob}
Make this precise.
\end{prob}
\proof Similar to \ref{jul0201}. \hB

Let $r:\tilde \Mf/(W\to B)\to  \Mf/(W\to B)$ be the functor which maps $(D,K)$ to $D$.
Then we have an adjunction
$$r_\flat:\Fun(\Nerve(\tilde \Mf/(W\to B))^{op},\bC)\leftrightarrows \Fun(\Nerve( \Mf/(W\to B))^{op},\bC):r^{*}\ .$$
Explicitly, 
$$(r_\flat X)(D)=\colim_K X(D,K)\ .$$

 \begin{ddd}
We define the push-forward with proper support by 
$$f_{!}:=  q_{*}\circ r_\flat\circ \tilde p^{*}\ .$$
 For $X\in \Sm_B(\bC)$ and $(p:M\to B)\in \Mf/B$
 we write $$X_{prop/p}(M\to B)=X_{prop/B}(M\to B):=((\id_{B})_! \id_{B}^{*}X)(M\to B)\ .$$
  \end{ddd}
  \begin{ex} \label{jul0204}{\em
Assume that $f:W\to B$ is an oriented submersion and $n=\dim(B)-\dim(W)$.
Consider $\Omega_\C\in \Sm_B(\Nerve(\Ch))$.
Then $$(f_!f^{*}(\Omega_\C)_B)(M\to B)=\Omega_{\C,prop/M}(M\times_BW)$$
is the space of forms on $M\times_BW$ with proper support over $M$.
Note that $M\times_BW\to M$ is naturally an oriented submersion.
\begin{prob}
Show that integration over the fibre
defines a transformation
$$\int_{W/B}: f_!f^{*} (\Omega_\C)_{B}\to(\Omega_\C[n])_{B}\ .$$
\end{prob}

}
\end{ex}

\begin{ex}{\em

Let $\iota:\Nerve(\Ch)\to \Nerve(\Ch)[W^{-1}]$ be the localization.
Let $\Omega\in \Sm^{desc}_B(\Nerve(\Ch))$ be a complex of soft sheaves.
We have a natural map (see \ref{jul2501})
$$\iota(f_!f^{*}  \Omega)\to f_!f^{*} \iota(\Omega)\ .$$
\begin{prob}\label{jul2401}
Show that
the natural map
$$\iota (f_!f^{*}  \Omega)\to f_!f^{*} \iota(\Omega)$$
is an equivalence.
\end{prob}
\proof
We fix $M\to B$. Then we must show that
$$\iota(\Omega_{prop/M}(M\times_BW\to M))\to \iota(\Omega)_{prop/M}(M\times_BW\to M)$$
is an equivalence. 
We have 
$$\Omega_{prop/M}(M\times_BW\to M)=\colim_{K} \Omega_{K}(M\times_BW\to M)$$
where $K\subseteq M\times_B W$ is proper over $M$. 
On the other hand
$$ \iota(\Omega)_{prop/M}(M\times_BW)=\colim_{K} \iota(\Omega)_K(M\times_BW\to M)\ .$$
The filtered colimit commutes with the limit defining the support.
We therefore get
and
$$\Omega_{prop/M}(M\times_BW\to M)=\fibre(\Omega(M\times_BW\to M)\to \colim_K  \Omega((M\times_BW\setminus K)\to M))$$ and
$$ \iota(\Omega)_{prop/M}(M\times_BW\to M)=\fibre(\iota(\Omega(M\times_BW \to M))\to \colim_K \iota(\Omega(M\times_BW\setminus K)\to M))\ .$$
Since a filtered colimit commutes with cohomology we see that
$$\colim_K \iota(\Omega((M\times_BW\setminus K)\to M))\cong \iota (\colim_K \Omega((M\times_BW\setminus K)\to M))\ .$$
Since $\Omega$ is soft the map
$\Omega(M\times_BW)\to \colim_K  \Omega(M\times_BW\setminus K)$ is surjective.
This implies by \ref{jul2501} that its kernel also represents the kernel taken after localization.
\hB 
 }
 \end{ex}

\begin{ex}{\em 
We consider   a differential data $(E,A,c)$.  Define
$\Diff(E,A,c)^{n}_{K}(W)$ as in \ref{jul0203}.

\begin{prob} \label{jul2202}Assume that the data $(E,A,c)$ is strict.
Make the long exact sequence of the pair $(W,W\setminus K)$ in differential $E$-cohomology explicit.
\end{prob}
\proof For every $n$ we have sequence
$$\dots \to E\C/\Z^{n-2}(W)\to E\C/\Z^{n-2}(W\setminus K)\to \hat E^{n}_{K}(W)\to \hat E^{n}(W)\to \hat E^{n}(W\setminus K)\to
E^{n+1}(K)\to \dots$$
which continues with the long exact pair sequences of $E\C/\Z^{*}$ and $E^{*}$ to the left and the right.
\hB 

\begin{prob}\label{jul1003}
Derive the compactly supported version of the exact sequence \eqref{jul1002}.
\end{prob}
\proof
We have
$$\to E_{c}^{k-1}(M)\to \Omega A_{c}^{k-1}/\im(d_{|\Omega A_{c}^{k-2}})\to\hat E^{k}_{c}(M)\to E_{c}^{k}
\to 0\ .$$
To see this consider the filtered colimit  over the compact subsets of $M$ of
the pull-back diagrams defining $\Diff^{k}(E,W,c)_{K}$.
Commute the colimit inside and use
\ref{jul2501} to commute the colimit with the localization on the de Rham side. \hB

\begin{prob}\label{aug1302}
Assume that the data $(E,A,c)$ is multiplicative. Let $K,K^{\prime}\subseteq W$ be closed.
Show that $\hat E_{K}^{*}(M)$ is an $\hat E^{*}(W)$-module.
Refine this structure to a product
$$\hat E_{K}^{*}(W)\otimes \hat E_{K^{\prime}}^{*}(W)\to \hat E_{K\cap K^{\prime}}^{*}(W)\ .$$
\end{prob}
}\end{ex}
\begin{ex}{\em
We consider the compact subset $\{0\}\subset \R$ and calculate
$\widehat{H\Z}^{*}_{\{0\}}(\R)$ using the exact sequences. We write out the interesting pieces:
$$\C/\Z\to \C/\Z\oplus \C/\Z\to \widehat{H\Z}^{2}_{\{0\}}(\R)\to 0 $$
$$0\to \widehat{H\Z}^{1}_{\{0\}}(\R)\to \Omega^{0}(\R,\C)/\Z \to \Omega^{0}(\R\setminus\{0\},\C)/\Z\oplus \Z$$
$$ 0\to \widehat{H\Z}^{0}_{\{0\}}(\R)\to \Z\to \Z\oplus \Z$$
We conclude that
$$\widehat{H\Z}^{k}_{\{0\}}(\R)\cong \left\{\begin{array}{cc}
0&k\not=2\\
\C/\Z&k=2\end{array} \right.
\ .$$  

\begin{prob} \label{jul1004}
Calculate $\widehat{H\Z}^{k}_{[0,1]}(\R)$.
\end{prob}
\proof
Use again the exact sequences. The result is
$$\widehat{H\Z}^{k}_{[0,1]}(\R)\cong \left\{\begin{array}{cc}
 \C/\Z&k=2\\
 0\to \Omega^{0}_\flat((0,1),\C)\to\widehat{H\Z}^{1}_{[0,1]}(\R) \to \Z\to 0&k=1\\
 0&k\not\in\{1,2\} 
 \end{array} \right.
 \ .$$
 Here $\Omega^{0}_\flat((0,1),\C)\subset \Omega^{0}((0,1),\C) $ are those functions on the interval which extend smoothly by zero to $\R$. 
\hB 

\begin{prob}
Calculate
$\widehat{H\Z}^{*}_{\{1\}}(S^{1})$.
\end{prob}
 \proof
 We write out the relevant sequences.
 
 $$\C/\Z\stackrel{=}{\to} \C/\Z\to \widehat{H\Z}^{2}_{\{1\}}(S^{1})\to \C/\Z\to 0 $$
$$0\to \widehat{H\Z}^{1}_{\{1\}}(S^{1})\to  \widehat{H\Z}^{1}(S^{1}) \to \Omega^{0}(S^{1}\setminus\{1\},\C)/\Z $$
$$ 0\to \widehat{H\Z}^{0}_{\{0\}}(S^{1})\to \Z\stackrel{=}{\to} \Z$$
We conclude that
$$\widehat{H\Z}^{k}_{\{1\}}(S^{1})\cong \left\{\begin{array}{cc}
0&k\not=2\\
\C/\Z&k=2\end{array} \right.
\ .$$ \hB
The coincidence $\widehat{H\Z}^{*}_{\{0\}}(\R)\cong \widehat{H\Z}^{*}_{\{1\}}(S^{1})$ is of course expected by excision.
}\end{ex}

{\begin{ex}
We calculate $\widehat{H\Z}_c^{*}(\R)$ using \ref{jul1003}.
We again write out the relevant sequences
$$0\to \C/\Z\to  \widehat{H\Z}_c^{2}(\R)\to 0$$
\begin{equation}\label{jul1001}0\to  \Omega^{0}_c(\R)\to \widehat{H\Z}_c^{1}(\R)\to \Z\to 0\end{equation}
$$ 0\to \widehat{H\Z}_c^{0}(\R)\to 0$$
We conclude that 
$$\widehat{H\Z}^{k}_{c}(\R^{1})\cong \left\{\begin{array}{cc}
\C/\Z&k=2\\
\eqref{jul1001}&k=1\\
0&k\not\in{1,2}\end{array} \right. .$$
Alternatively we could use \ref{jul1004}
\end{ex}}

\begin{ex}{\em 
Let $W$ be a smooth oriented $n$-manifold with boundary $\partial W$. 
\begin{prob}
Calculate $\widehat{H\Z}^{n}_{\partial W}(W)$, $\widehat{H\Z}_{\partial W}^{n+1}(W)$ and $\widehat{H\Z}_{c}^{n+1}(W\setminus \partial W)$.
\end{prob}
\proof
For the calculation we use appropriate exact sequences, e.g. the pair sequence obtained in  \ref{jul2202}. 
We get $\widehat{H\Z}_{\partial W}^{n}(W)\cong \C/\Z$ and $\widehat{H\Z}_{\partial W}^{n+1}(W)\cong 0$ and
$\widehat{H\Z}_{c}^{n+1}(W\setminus \partial W)\cong \C/\Z$. \hB

\begin{prob}
If $N\subseteq M$ is a closed embedded submanifold of non-zero codimension, then we have for every $k\in \Z$ a natural isomorphism
$\hat E_{N}^{k}(M)\cong E\C/\Z^{k-1}(M,M\setminus N)$.
\end{prob}
\proof
Use  \ref{jul2202}. \hB 

\begin{prob}
If $N\subseteq M$ is an embedded submanifold of  codimension zero, then we have for every $k\in \Z$ a natural  injection
$ \hat E_{c}^{k}(N\setminus \partial N)\subseteq \hat E_{N}^{k}(M)$
\end{prob}
\proof
Note that we can not expect an isomorphism since
the elements of  $\hat E_{c}^{k}(N\setminus \partial N)$ have curvature compactly supported in $N\setminus \partial N$ while the curvatures of the elements in $\hat E_{N}^{k}(M)$ can be smoothly extended by zero to $M$, but need not be supported properly in the interior of $N$. We have seen this effect already in the calculation \ref{jul1004}.\hB 
}

\end{ex}




\begin{ex}{\em
We have the de Rham equivalence in $\Sm^{desc}(\Nerve(\Sp)[W^{-1}])$
$$\Rham:H(\Omega A)\to \Funk(H(A))$$
which as an eqivalence of sheaves of ring spectra if $A$ is a commutative dga.
\begin{prob}
Show that it naturally induces
an equivalence
$$\Rham:H(\Omega A_{prop/M})\to \Funk (H(A))_{prop/M}$$
in $\Sm(\Mod(H\Z))$ (or $\Sm(\Mod(H(\Omega A)))$ if $A$ is a commutative dga).
\end{prob}
\proof Use the naturality of $\Rham$. \hB 
}\end{ex}


%

If the differential data $(E,A,c)$ is multiplicative we can  interpret the construction
of $\Diff^{\bullet}(E,A,c)_{prop/M}$ in $p^{*}\Diff^{\bullet}(E,A,c)$-module spectra. In particular we get
a $p^{*}\hat E^{*}$-module structure on $\hat E^{*}_{prop/M}$, where we define
$$\hat E^{*}_{prop/M}:=\pi_*(\Diff^{\bullet}(E,A,c)_{prop/M})\ .$$
Again this notation is sloppy in the sense note earlier in \ref{aug1045}.

\begin{prob}
Analyse the difference between  
$\hat E^{*}_{prop/M}=\pi_*(\Diff^{\bullet}(E,A,c)_{prop/M})$ and
the group  $\pi_*(\Diff^{\bullet}(E,A,c))_{prop/M}$ of differential cohomology classes having proper support over $M$. 
\end{prob}
\proof
$\hat E^{*}_{prop/M}$
is a subgroup of $\hat E^{*}_{M}$ consisting of those elements which vanish after restriction to the complement of a sufficiently large subset which is proper over $M$.
$\hat E^{*}_{prop/M}$ maps to $\pi_*(\Diff^{\bullet}(E,A,c))_{prop/M}$,  but not injectively in general. \hB 

\begin{ex}
{\em  Let $f:W\to M$ be a smooth map. 
 In order to explicitly  construct   differential cohomology classes which are properly supported over $M$ 
  we need versions of the cycle maps which respect supports. We first consider the case of bordism.
   \begin{ddd} 
A geometric $\bMU$-cycle  (cf. Def. \ref{aug0510}) $(g,\cN,J,\nabla,\eta)$ on $W$
is  properly supported over $M$ if the composition $f\circ g$ is proper and
$\eta$ has proper support over $M$. 

We let
$\Cycle_{\bMU,prop/M}^{geom,*}(W)$ denote the graded semigroup of isomorphism classes of geometric $\bMU$-cycles.
\end{ddd}

\begin{prob}\label{aug0330}
Show that there exists natural refinement of the cycle map  $\hcycl$  introduced in \ref{may2940} to a cycle map  
$$\hcycl:\Cycle_{\bMU,prop/M}^{geom,*}(W)\to \hbMU^{*}_{prop/M}(W)\ .$$
\end{prob}
\proof
 
Let $c:=(g:E\to W,\cN,J,\nabla,\eta)$ be a geometric $\bMU$-cycle which is properly supported over $M$. Then
there exists a closed subset $K\subseteq W$ which is proper over $M$ such that
$g(E)$  and $\supp(\eta)$ are contained in the interior of $K$. 
The restriction of the cycle $c$ to $W\setminus K$ is the zero cycle.
We can consider
$(g,\cN,J)$ as a relative cycle for the pair $(W,W\setminus K)$ and therefore get a class
$$y:=\cycl(g,\cN,J)\in \bMU^{*}(W,W\setminus K)\cong \bMU_{K}^{*}(W)\ .$$

We have an exact sequence (cf. \ref{jul2202})
$$\dots \bMU\C/\Z^{*-2}(W\setminus K)\to  \hbMU^{*}_{K}(W)\to \hbMU^{*}(W)\to \hbMU^{*}(W\setminus K)\ .$$
Since $\hcycl(x)_{|W\setminus K}=0$ we conclude that there exists a
class $\hat y\in \hbMU_{K}(W)$ which maps to $\hcycl(c)\in \hbMU^{*}(W)$. The class
$\hat y\in  \hbMU_{K}(W)$ is unique up to  elements which come from
$ \bMU\C/\Z^{*-2}(W\setminus K)$. The image of 
$\bMU\C/\Z^{*-2}(W\setminus K)\to  \hbMU^{*}_{K}(W)$ is detected by the composition with $I$.
Hence we can choose $\hat y$ uniquely such that $I(\hat y)=y$.
We define
$$\hcycl(c) :=\hat y\ .$$
\hB

A cycle for a differential $K$-theory class on $W$ which is properly supported over $M$ is given by 
the following data $(V_{0},V_{1},\nabla_{0},\nabla_{1}, \phi,K)$, where
$(V_{i},\nabla_{i})$, $i=0,1$ are complex vector bundles  with connection on $W$,
$K\subseteq W$ is closed and proper over $M$, and $\phi:(V_{0})_{|W\setminus K}\to (V_{1})_{W\setminus K}$ is a connection-preserving vector bundle isomorphism.

\begin{prob}
Show that there exists a unique natural and additive way to associate to the data 
$(V_{0},V_{1},\nabla_{0},\nabla_{1}, \phi,K)$ a class
$$\hcycl(V_{0},V_{1},\nabla_{0},\nabla_{1}, \phi,K)\in \hbKU_{K}(W)\ .$$
\end{prob}

}
\end{ex}

\subsection{Thom classes}

Let $E\in \CommMon(\Sp[W^{-1}])$ be a commutative ring spectrum which represents the multiplicative cohomology theory $E^{*}$.
Assume   that $f:W\to M$ is a real vector bundle of dimension $n$.
A Thom class of $W$ is a class
$\nu\in E^{n}_{prop/M}(W)$ whose restriction to each fibre $W_{m}$, $m\in M$ maps to $\pm 1$ under the isomorphism
\begin{equation}\label{jul2520}E_{prop/m}^{n}(m)\cong   E_{c}^{n}(\R^{n})\stackrel{\eqref{jun0401}\sim }{\to} E^{0}(*)\ .\end{equation} 
The first isomorphism depends on an identification of the fibre $W_{m}$ with $\R^{n}$.
If the characteristic of the ring $\pi_*(E)$ is not $2$, 
and if we fix a Thom class on $W\to M$, then the bundle $W$ acquires an induced ordinary orientation. An oriented frame of the fibre $W_m$, i.e  an orientation-preserving isomorphism  $W_{m}\cong \R^{n}$, is characterized by 
the requirement that
$\nu$ goes to $1$ under \eqref{jul2520}.

\begin{ddd}\label{jun0620}
A differential Thom class   on a real $n$-dimensional vector bundle $W\to M$ is a class
$\hat \nu\in \hat E^{n}_{prop/M}(W)$ such that $I(\hat \nu)=:\nu\in  E^{n}_{prop/M}(W)$
is a Thom class. We say that $\hat \nu$ refines $\nu$. We define
$$\Td(\hat \nu):=\int_{W/M} R(\hat \nu)\in \Omega A_{cl}^{0}(W) ,  \quad  \Td(\nu):=\Rham(\Td(\hat \nu))\in H(A)^{0}(W)^{\times}$$ 
\end{ddd}

\begin{prob}
Show that every Thom class   has a differential refinement.
\end{prob}

\begin{prob}\label{jul2510}
 If $\hat \nu$ is a differential Thom class on $W\to M$, then show that 
$$ \Td(\hat \nu)=1+\sum_{j}  e^{j}\omega_j$$ with homogeneous elements 
$e^{j}\in E^{<0}$ and $\omega_j\in \Omega_{cl}^{-\deg(e_j)}(M,\C)$.
Conclude that $\Td(\hat \nu)$ is a unit.
\end{prob} 
\proof 
Note that $\sum_{j}  e^{j}\omega_j$ is nilpotent.
\hB

The topological Thom class gives a Thom isomorphism
$$\Phi_\nu:=\nu\cup f^{*}(\dots):E^{*}(M)\to E_{prop/M}^{*+n}(W)$$
of $E^{*}(M)$-modules.
\begin{prob}
Assume that $(E,A,c)$ is strict.
Show that the differential  version
$$\hat \nu\cup f^{*}(\dots):\hat E^{*}(M)\to \hat E_{prop/M}^{*+n}(W)$$
 is an injective morphism of $\hat E^{*}(M)$ modules.
\end{prob}
\proof  
Let $x\in \hat E^{k}(M)$ and 
assume that
$\hat \nu\cup f^{*}(x)=0$. Then
$R(\hat \nu\cup  f^{*}(x))=R(\hat \nu)\wedge f^{*}R(x)=0$.
We conclude that
$\int_{W/M} R(\hat \nu)\wedge R(x)=0$.
Since $\int_{W/M} R(\hat \nu)$ is a unit by \ref{jul2510} we see that
$R(x)=0$. Hence $x$ is flat. For a flat class $x$ the class  $\hat \nu\cup f^{*}(x)$ is the image
of $x$ under the Thom isomorphism   $E^{*}\C/\Z^{k-1}(M)\to E^{*}\C/\Z_{prop/M}^{k-1+n}(W)$.
Therefore $x=0$. 
\hB

\begin{ddd}\label{aug0850}
We say that two differential Thom classes $\hat \nu^{\prime},\hat \nu$ on a bundle $N\to W$ are homotopic
if there exists a differential Thom class $\tilde{\hat{\nu}}$ on the pull-back $\tilde N:=\pr^{*}N\to I\times W$
which connects  $\hat \nu^{\prime}$ and $\hat \nu$ and satisfies $\Td(\tilde {\hat {\nu}})=\pr^{*}\Td(\hat \nu).$
\end{ddd}
Note that because of the condition on the $\Td$-forms homotopy in the sense of differential Thom classes is a stronger condition than just homotopy as differential cohomology classes. In particular,  if $\hat \nu$ and $\hat \nu^{\prime}$ are homotopic, then we have the equality
$\Td(\hat \nu)=\Td(\hat \nu^{\prime})$.  
 \begin{prob} \label{aug0610}
 We assume that $(E,A,c)$ is strict. 
 Show that homotopy classes of differential  Thom classes refining a given underlying topological Thom  class $\nu$ and with fixed form  $\Td(\hat \nu)$
  are classified by  the group
 $$ \frac{H(A)^{-1}(W)}{   \Td(\nu)\cup c(E^{-1}(W))}\ .$$
 \end{prob}
 \proof
 Let $\hat \nu$ and $\hat \nu^{\prime}$ be two differential Thom classes on a bundle $\pi:N\to W$
refining the same underlying topological Thom class $\nu$  and such that we have an equality of  integrals
$$\Td(\hat \nu):=\int_{N/W}R(\hat \nu^{\prime})=\int_{N/W}R(\hat \nu)\ .$$
Then there exists a form $\eta\in \Omega A^{n-1}_{prop/W}(N)/\im(d)$ uniquely determined up to
$c(E^{n-1}_{prop/W}(N))$ such that
$\hat \nu^{\prime}-\hat \nu=a(\eta)$. 
We calculate that $d\int_{N/W} \eta=0$. 
By the Thom isomorphism in cohomology we have $$\int_{N/W}c(E^{n-1}_{prop/W}(N))=  \Td(\nu)\cup c(E^{-1}(W))\ .$$
We therefore define the difference class
\begin{equation}\label{aug0710}\delta(\hat \nu^{\prime},\hat \nu):=\Rham(\int_{N/W} \eta)\in \frac{H(A)^{-1}(W)}{   \Td(\nu)\cup c(E^{-1}(W))}\ .\end{equation}
We claim that $\delta(\hat \nu^{\prime},\hat \nu)$ exactly detects whether the two Thom classes $\hat \nu$ and $\hat \nu^{\prime}$  are homotopic.

First we show that if $\delta(\hat \nu^{\prime},\hat \nu)=0$, then $\hat \nu^{\prime}$ and $\hat \nu$ are homotopic.
Assume that $\delta(\hat \nu^{\prime},\hat \nu)=0$.
Then we can find an element $y\in E^{-1}(W)$ such that
$$\Rham(\int_{N/W} \eta)=\Td(\nu)\cup c(y)\ .$$ We can find a closed form $\beta\in \Omega A^{n-1}_{prop/W,cl}(N)/\im(d)$ such that $\Rham(\beta)=c(z)$ for $z:=\nu\cup \pi^{*}y$ 
and $\int \beta=\int \eta$. Indeed, if $\beta_{0}$ is some representative of the class $c(z)$, then
we have
$\int_{N/W} \eta-\int_{N/W}\beta_{0}=d\alpha$. We then set
$$\beta:=\beta_{0}+ d(R(\hat \nu)\wedge \pi^{*}(\alpha\cup \Td(\hat \nu)^{-1}))\ .$$
We now have $\hat \nu^{\prime}-\hat \nu=a(\eta)=a(\eta-\beta)$.
We consider the homotopy
$$\tilde{\hat {\nu}} :=\hat \nu+t(\eta-\beta)\ .$$
We must check the curvature condition. 
We have
$$R(\tilde{\hat {\nu}})=\pr^{*}R(\hat \nu)+dt\wedge (\eta-\beta) + t\wedge d \eta\ .$$
Since $\int_{N/W} d\eta=0$ and $\int_{N/W} (\eta-\beta)=0$ we conclude that
$$\Td(\tilde{\hat {\nu}})=\pr^{*} \Td(\hat \nu)\ .$$

We assume now that $\hat \nu^{\prime}$ and $\hat \nu$ are connected by a homotopy
$\tilde {\hat {\nu}}$.  By the homotopy formula
$$\hat \nu^{\prime}-\hat \nu=a(\int_{I\times N/N} R(\tilde {\hat {\nu}}))\ .$$
Hence we get
$$\delta(\hat \nu^{\prime},\hat \nu)=\int_{I\times N/W} R(\tilde {\hat {\nu}})=\int_{I\times W/W} \pr^{*}R(\hat {\nu})
=0\ .$$
 
Finally, let $y\in H(A)^{-1}(W)$ and $\hat \nu$ be given. Then we define
$$\hat \nu^{\prime}:=\hat \nu+a(R(\hat \nu)\wedge \pi^{*}(\Td(\hat \nu)^{-1} \cup y))\ .$$
With this choice
$$\delta(\hat \nu^{\prime},\hat \nu)=[y]\in \frac{H(A)^{-1}(W)}{   \Td(\nu)\cup c(E^{-1}(W))}\ .$$
\hB

The solution of Problem \ref{aug0610} allows to state the following  {\em two-out-of-three}
principle. Assume that $(E,A,c)$ is strict. Let $N=N_{0}\oplus N_{1}$ be a decomposition of the bundle $N\to W$ and $p_{i}:N\to N_{i}$ be the projections.   Then we consider the relation
\begin{equation}\label{aug0650}
\hat \nu=p_{0}^{*}\hat \nu_{0}\cup p_{1}^{*}\hat \nu_{1}
\end{equation}
between differential Thom classes on $N$ and the summands $N_{i}$. 
\begin{prob}\label{aug0670}

Show that two of these classes determine the third uniquely up to homotopy such that
\eqref{aug0650}
holds true.
\end{prob}
\proof
Assume that $\hat \nu$ and $\hat \nu_{1}$ are given.
The equality
$$\Td(\hat \nu_{1})=\frac{\Td(\hat \nu)}{\Td(\hat \nu_{0})}$$  fixes the $\Td$-form of $\hat \nu_{1}$.
If we choose any differential Thom class $\hat \nu_{1}^{\prime}$ with this $\Td$-form, then we have a class
$$\delta:=\delta(\hat \nu,p_{0}^{*}\hat \nu_{0}\cup p_{1}^{*}\hat \nu^{\prime}_{1})\in  \frac{H(A)^{-1}(W)}{   \Td(\nu)\cup c(E^{-1}(W))}\ .$$
We are forced to define
$$\hat \nu_{1}:=\hat \nu_{1}^{\prime}+a(R(\hat \nu_{1}^{\prime})\wedge \pr^{*}(\Td(\hat \nu_{1})^{-1}\cup \delta ))\ .$$

\hB

We now discuss the special case of the Thom isomorphism for trivial bundles.
In particular we relate it with the suspension isomorphism on the one hand, and with the integration of differential forms on the other. Let $\pr: \R^{n}\to *$ be the projection. 
Then we can construct  a commutative diagram in
$\Nerve(\Sp)[W^{-1}]$.  
 It turns out to be useful not to hide the localization $\iota:\Nerve(\Ch)\to \Nerve(\Ch)[W^{-1}]$.
\begin{equation}\label{jun0601}\hspace{-1cm}\xymatrix{\pr_{!}\pr^{*}\Funk(E)\ar[dd]\ar[r]^{c}&\pr_{!}\pr^{*}\Funk(H(A))\ar[dd]&\pr_{!}\pr^{*}H(\Omega A)\ar[l]^{\Rham}&\pr_{!}\pr^{*}H(\iota(\sigma^{\ge k}\Omega A))\ar[l]\\&&H(\iota(\pr_{!}\pr^{*}\Omega A ))\ar[u]_{(ii)}\ar[d]^{(i)}&H(\iota(\pr_{!}\pr^{*}\sigma^{\ge k} \Omega A ))\ar[l]\ar[u]_{(ii)}\ar[d]^{(i)}\\\Funk(E[-n])\ar[r]^{c}&\Funk(H(A)[-n])&H(\iota(\Omega A[-n]))\ar[l]^{\Rham}&H(\iota(\sigma^{\ge k-n} \Omega A[-n]))\ar[l]}
\end{equation}

The two  vertical maps in \eqref{jun0601} marked by $(i)$  are induced by the integration of forms
$\int_{\R^{n}}$ (see \ref{jul0204}). Note that we can not define the integration on
$\pr_{!}\pr^{*} \iota(\Omega A)$ directly since the elements of its evaluation  are not really properly supported forms.

   The arrows marked by $(ii)$ are equivalences, similar to \ref{jul2401}.

The lower right square and the upper right square commute.
 We explain the construction of the left vertical maps in the   evaluation at $M$. For $r>0 $ let $B^{c}(r)\subset \R^{n}$ be the complement of the closed $r$-ball centered at zero. Then we have  a canonical identification
$$\Funk_{prop/M}(E)(M\times \R^{n})\cong \colim_{r} \Map(\Sigma^{\infty} (M_{+}\wedge (\R^{n}/B^{c}(r))),E)$$
If we identify $S^{n}$ with the one-point compactification of $\R^{n}$, then we get a compatible system of natural maps $S^{n}\to \R^{n}/B^{c}(r)$ for all $r>0$ and therefore a map
$$\colim_{r} \Map(\Sigma^{\infty} (M_{+}\wedge (\R^{n}/B^{c}(r))),E)\to \Map(\Sigma^{\infty}M_{+}\wedge \Sigma^{\infty}S^{n},E)\cong \Map(\Sigma^{\infty} M_{+},E[-n]) \ ,$$
were $$E[-n]:=\Map(\Sigma^{\infty} S^{n},E)\ .$$ The composition of these maps
gives the arrow $$\Funk_{prop/M}(E)(M\times \R^{n})\to \Funk(E[-n])(M)$$ and, similarly, 
$$\Funk_{prop/M}(H(A))(M\times \R^{n})\to \Funk(H(A)[-n])(M)\ .$$ The left square in \eqref{jun0601}  commutes.

\begin{prob}\label{jun0611n}
The middle square in \eqref{jun0601} commutes.
If $A$ is a commutative dga, then we can interpret this square in $H(\Omega A(M))$-modules.
\end{prob}
\proof
Note that $A$ is a $\C$-module. Moreover,
the vertical arrow and the composition $(i)\circ (ii)^{-1}$ are isomorphisms whose inverse is given by the multiplication with a representative of the    Thom class $H \Omega_{c,cl}^{n}(\R^{n},\C)$ and its image under $\Rham$. 
The square with the vertical arrows inverted thus commutes since $\Rham$ is multiplicative and preserves the Thom classes.
\hB 

If one takes the product of these diagrams \eqref{jun0601} over all $k$ and uses the fact that the 
diagram of ring spectra $\Diff^{\bullet}(E,A,c)$ acts naturally on all entries one can refine this product to a diagram in $\Mod(\Diff^{\bullet}(E,A,c))$.

\begin{prob}\label{aug0630} Show that
the diagram  \eqref{jun0601} induces a desuspension map
$$\desusp:\pr_{!}\pr^{*}\Diff^{\bullet}(E,A,c)  \to \Diff^{\bullet-n}(E,A,c) $$
(of $\Diff^{\bullet}(E,A,c) $-modules). \end{prob}
\proof
We must interchange the order of $\pr_{!}\pr^{*}$ and the finite limit defining $\Diff^{\bullet}(E,A,c)$. This works since $\pr_{!}$ involves filtered colimits, only. 
\hB 

\begin{kor}\label{aug0901}
We have a desuspension map
$$\desusp:E^{*}_{prop/M}(M\times \R^{n})\to E^{*-n}(M)$$
of $E^{*}(M)$-modules.
\end{kor}

In order to verify functorial properties of differential integration maps we need the following statement
about desuspension in stages.
\begin{prob}\label{aug0902}
We have a refinement of the desuspension map
$$\desusp:\hat E^{*}_{prop/M}(M\times \R^{l}\times \R^{n})\to \hat E^{*-n}_{prop/M}(M\times \R^{l})\ ,$$
and the diagram
$$\xymatrix{\hat E_{prop/M}(M\times \R^{l}\times \R^{n})\ar[rr]\ar[dr]^{\desusp}&&\hat E^{*-n-l} (M)\\&\hat E_{prop/M}^{*-n}(M\times \R^{l})\ar[ur]^{\desusp}&}$$
commutes.
\end{prob}
\proof
Write out the diagram \eqref{jun0601} as a composition of two  similar, appropriately adapted diagrams. 
\hB


\begin{ex}{\em 
Let $\hat \nu\in \widehat{\bS}^{n}_{c }(\R^{n})$ be a differential Thom class in differential stable cohomotopy refining the canonical class in $\bS_{c }^{n}(\R^{n})$. Exterior multiplication by $\hat \nu$ (see \ref{jun0610}) induces a map
$$\hat \nu\times {\dots}:\hat E^{*}(M)\to \hat E_{prop/M}( M\times \R^{n})\ .$$
\begin{prob}
Show that on $\hat E^{*}(M)$ we have the identity
$$\desusp\circ (\hat \nu\times \dots)=\id\ .$$
\end{prob}

}\end{ex}
\begin{ex}\label{jul2610}{\em
In this example we discuss the geometric construction of differential $\bMU$-Thom classes for complex vector bundles with connection.
Let $p:V\to M$ be a  real vector bundle of dimension $2n$ with a complex structure $I$. We consider the zero section $0_{V}:M\to V$.   We have a canonical inclusion $\iota:p^{*}V\to TV$ by the linear structure and an isomorphism 
$TM\oplus V\stackrel{d0_{V}+ 0_{V}^{*}\iota}{\to} 0_{V}^{*}TV$. In particular we get an exact sequence
$$\cN_{0_{V}} : 0\to TM\stackrel{d0_{V}}{\to} 0_{V}^{*} TV\to V\to 0$$
which turns $V$ into the normal bundle of $0_V$. 
We therefore get  an $\bMU$-cycle $(0_{V},\cN_{0_{V}},I)$ of degree $2n$ (cf. Definition \ref{jul2601}) and a class $\nu_{\bMU}:=\cycl(0_{V},\cN_{0_{V}},I)\in \bMU^{2n}_{prop/M}(V)$
(see the proof of \ref{aug0330} for the support condition).
\begin{prob}
\begin{enumerate}
\item Show that $\nu_{\bMU}$ is a $\bMU$-Thom class
of $V$.
\item
Show that $\nu_{\bMU}$ can be refined to a class in $\bMU^{2n}_{0_{V}}(V)$.  \end{enumerate}\end{prob}

Let now $\nabla$ be a connection on $V$ preserving the complex structure $I$. 
Then we can define the form $u(\nabla)\in \Omega A_{cl}^{0}(M)$ as in \eqref{aug0310}.

\begin{prob}\label{aug0335}
Show that there exists a differential Thom class $\hat \nu\in \hbMU^{2n}_{prop/M}(V)$ such that
$$\Td(\hat \nu)=u(\nabla)\ .$$
\end{prob}
\proof
Let $\eta\in  \Omega A_{-\infty}^{2n-1}(V)_{prop/M}$ be some distributional form such that
$0_{V,!}(u(\nabla))-d\eta=0$.
Then we have by \ref{aug0330} a differential Thom class
$$\hat \nu_{0}:=\cycl(0_{V},\cN_{0_{V}},I,\nabla,\eta)\in \hbMU^{2n}_{prop/M}(V)\ .$$
Even without having a solution of Problem \ref{aug0330} it is clear
that there exists a differential Thom class  $\hat \nu_{0}\in  \hbMU^{2n}_{prop/M}(V)$ with
$\Td(\hat \nu_{0})=u(\nabla)-d\int_{V/M} \eta$ since this form represents the correct cohomology class.
Now $d\int_{V/M} \eta$ is smooth and exact. Hence there exists a smooth form
$\mu\in \Omega A^{-1}(M)$ such that $d\mu=d\int_{V/M}\eta$. We can further choose a form
$\tilde \mu\in \Omega A^{2n-1}_{prop/M}(V)$ such that $\int_{V/M}\tilde \mu=\mu$.
We define a corrected Thom form
$\hat \nu:=\hat \nu_{0}+a(\tilde \mu)$ which satisfies
$$\Td(\hat \nu)=u(\nabla)\ .$$ \hB

Note that the homotopy class of the differential Thom class $\hat \nu\in  \hbMU^{2n}_{prop/M}(V)$ is not uniquely determined by the condition 
$\Td(\hat \nu)=u(\nabla)$. In fact,  by the solution of Problem \ref{aug0610} the set of such Thom classes  forms a torsor over $$\frac{H(A)^{-1}(M)}{  u(V)\cup c(\bMU^{-1}(M))}\ ,$$
see \eqref{aug1050} for the characteristic class $u(V)$.  One can fix this ambiguity by requiring naturality.

\begin{prob}\label{aug0720}
Show that there is a unique way to associate a homotopy class of differential Thom classes
$\hat \nu(\nabla^{V})\in \hbMU^{2n}_{prop/M}(V)$ to a $n$-dimensional complex vector bundle with connection $(V,\nabla^{V})$ on a manifold $M$ which is natural under pull-back and such that $\Td(\hat \nu(\nabla^{V}))=u(\nabla^{V})$.
\end{prob}
\proof
We use the same technique which was already successful in the proofs of \ref{may2007}, \ref{may2010} and \ref{may2930}. For the notation see also \ref{apr2001}.  
The classifying space for  $n$-dimensional complex vector bundles is $BU(n)$.
Let $m:=\dim(M)$. We can find a factorization
$$\xymatrix{V\ar[d]\ar[r]^{H}&W\ar[d]\ar[r]&\xi_{n} \ar[d]\\M\ar[r]^{h}&N\ar[r]^{u}&BU(n)}$$
of 
the classifying map of the bundle $V\to M$,  where
$u:N\to BU(n)$ is $m+1$-connected  and $\xi_{n}\to BU(n)$ is the universal $n$-dimensional complex vector bundle. We can assume that $W$ has a connection $\nabla^{W}$ which pulls back to $\nabla^{V}$. 

The odd degree cohomology of $N$ is concentrated in degrees which exceed the dimension of $M$ so that the pull-back
$$h^{*}:\frac{H(A)^{-1}(N)}{   u(W)\cup c(\bMU^{-1}(N))}\to \frac{H(A)^{-1}(M)}{   u(V)\cup c(\bMU^{-1}(M))}$$
vanishes.
If we choose $\hat \nu(W)$ such that $\Td(\hat \nu(W))=u(\nabla^{W})$, then the pull-back
$ H^{*}\hat \nu(W)$
satisfies $\Td(\hat \nu(\nabla^{V}))=u(\nabla^{V})$ as required and is independent of the choice of
$\hat \nu(W)$. We are forced to define
$$\hat \nu(\nabla^{V}):=H^{*}\hat \nu(W)\ .$$
It remains to show that $\hat \nu(\nabla^{V})$
is well-defined independently of the choices, and to verify that our construction of the differential Thom class is natural.
For well-definedness we argue as before. Two choices can be  related with a third by a diagram
 $$\xymatrix{&N\ar[dr]^{g}\ar[drr]^{u}&&\\M\ar[ur]^{h}\ar[dr]^{h^{\prime}}&&N^{\prime\prime}\ar[r]^{u^{\prime\prime}}&BU(n)\\&N^{\prime}\ar[ur]^{g^{\prime}}\ar[urr]^{u^{\prime}}&&}$$
where $g\circ h$ and $g^{\prime}\circ h^{\prime}$ are homotopic. Let $H:I\times M\to N^{\prime\prime}$ be such a homotopy. We can lift this homotopy to an identification of vector bundles $\tilde H:\pr_{M}^{*}V\to H^{*}W^{\prime\prime}$. We can further assume that under this identification
$\tilde H^{*}\nabla^{W^{\prime\prime}}=\pr_{M}^{*}\nabla^{V}$. 
To this end we just have to make sure that $H$ is an embedding. We can reach this situation by modifying $N^{\prime\prime}$. The connection  $\nabla^{W^{\prime\prime}}$ is the  defined by an extension of $\pr_{M}^{*}\nabla^{V}$ on the image of $H$. With these choices
 $$\int_{I\times M/M} u(H^{*}\nabla^{W^{\prime\prime}})=0\ .$$
 The  homotopy formula gives 
$$H^{\prime\prime *}\hat \nu(W^{\prime\prime})-H^{*}\hat \nu(W^{\prime\prime})=a(\int_{I\times V/V} R(\tilde H^{*}\hat \nu(W^{\prime\prime})))\ .$$
The difference class \eqref{aug0710} between the two Thom classes is now given by
$$\delta(H^{\prime\prime *}\hat \nu(W^{\prime\prime}),H^{*}\hat \nu(W^{\prime\prime}))=\left[\int_{V/M} 
\int_{I\times V/V} R(\tilde H^{*}\hat \nu(W^{\prime\prime}))\right]
\in \frac{H(A)^{-1}(M)}{   u(V)\cup c(\bMU^{-1}(M))}\ .
 $$
 But by Fubini
$$ \left[\int_{V/M} 
\int_{I\times V/V} R(\tilde H^{*}\hat \nu(W^{\prime\prime}))\right]=\left[\int_{I\times M/M} u(H^{*}\nabla^{W^{\prime\prime}})
\right]=0
$$ 
 

Naturality is now easy to check. \hB

Note that the argument for Problem \ref{aug0720} only depends on the fact that
the rational cohomology of the classifying space $BU(n)$ is concentrated in even degrees.
So a similar argument shows e.g. :

\begin{kor}\label{aug0721}
There is a unique way to associate a homotopy class of differential Thom classes
$\hat \nu(\nabla^{V})\in \hbKU^{n}_{prop/M}(V)$ to a $n$-dimensional  $Spin^{c}$ vector bundle with $Spin^{c}$-connection $(V,\tilde \nabla^{V})$ on a manifold $M$ which is natural under pull-back and such that $\Td(\hat \nu(\tilde \nabla^{V}))=\Td(\tilde \nabla^{V})^{-1}$ (see \ref{jul3110} for notation).
\end{kor}

}
\end{ex}

\begin{ex}{\em 
In the following example we  consider cannibalistic classes. They arise from the non-compatibility of Thom isomorphisms with cohomology operations. The exercises should clarify how the differential
 refinement of the theory of cannibalistic classes works.  We consider the examples of the Adams operations on complex $K$-theory and the Chern character between complex $K$-theory and ordinary cohomology. 

We consider a real vector bundle $V\to M$ of dimension $n$.
We assume that $V$ has a $Spin^{c}$-structure, i.e.  we have a $Spin^{c}(n)$-principal bundle $Q\to M$ and an identification $Q\times_{Spin^{c}(n)}\R^{n}\cong V$. 
By definition, a $Spin^{c}$-connection on $V$ is a connection $\tilde \nabla$ on $Q$. 
The $Spin^{c}$-structure provides a Thom class $\nu\in \bKU^{n}(M^{V})$ and a Thom isomorphism
$$\Phi_\nu:\bKU^{*}(M)\to \bKU^{*+n}(M^{V})\ .$$ 
This follows from the existence of the Atiyah-Bott-Shapiro orientation \cite{MR0167985}
$$ABS:\mathbf{MSpin^{c}}\to \bKU\ .$$ This orientation is in fact multiplicative (see \cite{MR2122155} for an $E_{\infty}$-version) and
compatible with the complex orientation of $c:\bMU \to \bKU$ given by the multiplicative formal group law $x+_{\bKU}y=x+y+bxy$, i.e. the following diagram commutes
$$\xymatrix{\bMU\ar[rr]\ar[dr]^{c}&&\mathbf{MSpin^{c}}\ar[dl]^{ABS}\\ &\bKU&}\ ,$$
where the upper map is induced by the maps $\beta:U(n)\to Spin^{c}(2n)$, see e.g. \eqref{aug0201}. 
This observation will help in calculations.

We fix $k\in \nat$. We refer to Example \ref{jun0210} for the notation related to Adams operations.
\begin{ddd}
We define the cannibalistic class
$$\rho^{k}(V):=\Phi^{-1}_\nu(\Psi^{k}(\Phi_\nu(1)))\in \bKU[\{k\}^{-1}]^{0}(M)\ .$$
\end{ddd}
Let $\Vect^{Spin^{c}}_\R\in \Sm(\CommMon(\Set))$ be the smooth monoid of real vector bundles with $Spin^{c}$-structure with respect to the sum and $\bKU[\{k\}^{-1}]^{0}(M)^{\times}$ denote the units in the ring $\bKU[\{k\}^{-1}]^{0}(M)$.
\begin{prob}
\begin{enumerate}
\item Show that
$V\mapsto \rho^{k}(V)$ is an exponential (sums go to products) natural transformation
$$\Vect^{Spin^{c}}_\R\to (\bKU[\{k\}^{-1}]^{0})^{\times}$$
in $\Sm(\CommMon(\Set))$.
\item One can twist a $Spin^{c}$-structure by a complex line bundle $L$\footnote{Let $P\to M$ be the $U(1)$-principal bundle associated to $L$. Then the twisted $Spin^{c}$-structure is given by the $Spin^{c}$-principal bundle $Q^{L}:=Q\times_{U(1)} P$. The group $Spin^{c}(n)$ acts on the first factor (note that $U(1)$ is central), and the isomorphism $ Q^{L}\times_{Spin^{c}(n)}\R^{n}\cong V$ is induced from $ Q \times_{Spin^{c}(n)}\R^{n}\cong V$ in the obvious way.}. Let $V^{L}\in \Vect_\R^{Spin^{c}}(M)$ denote the twist of $V$ by $L$.
Calculate $\rho^{k}(V^{L})\cup \rho^{k}(V)^{-1}$.
\item A complex line bundle $E\to M$ has a canonical $Spin^{c}$-structure. Calculate
$\rho^{k}(E)$.
\item Let $V\to B$ be a two-dimensional $Spin^{c}$-bundle. Since the underlying real vector bundle is associated to $SO(2)\cong U(1)$ it has a complex structure. Calculate
$\rho^{k}(V)$.
\end{enumerate}
\end{prob}
\proof
For $1.$ one uses the multiplicativity of the ABS-orientation which relates the Thom classes for
the sum of $Spin^{c}$-bundles with the Thom classes of the product.

For $2.$ one uses that $\nu_{V^{L}}:=\nu_{V}\cup \cycl(L)$. 
We get
$$\rho^{k}(V^{L})=\Phi^{-1}_{\nu_{V^{L}}}(\Psi^{k}(\nu_{V^{L}}))=
\Phi^{-1}_{\nu_{V^{L}}}(\Psi^{k}(\nu_{V})\cycl(L^{k}))=\rho^{k}(V)\cycl(L^{k-1})$$
It follows that
$\rho^{k}(V^{L})\cup \rho^{k}(V)^{-1}=\cycl(L^{k-1})$.

For 3. we first calculate in the universal example  $E\to BU(1)$. Since $\bKU^{*}(BU(1))$ is torsion-free the Chern character is injective. We therefore first calculate 
$\ch(\rho^{k}(E))$. We know that $bz:=0_{E}^{*}\nu_{E}\in \bKU^{2}(BU(1))$ is a coordinate for the multiplicative formal group law. The first Chern class $x\in H^{2}(BU(1);\C)$ of $E$ is the coordinate of the additive formal group law. It follows that
$\exp(bx)-1=\ch(bz)$. 
We apply $\ch\circ 0_{E}^{*}$ to
$$\Psi^{k}(\nu_{E})=\rho^{k}(E)\cup \nu_{E}$$ and get
$$ \Psi^{k}_{H}(\ch(z))=\ch(\rho^{k}(E))\cup \ch(z)\ .$$
 We multiply with $b$ and rearrange the terms:
 $$k^{-1}\Psi_{H}^{k}(\ch(bz))=\ch(\rho^{k}(E))\cup \ch(bz)\ .$$
 We now substitute $z$ by the $x$-variable and carry out the cohomological Adams operation which multiplies $b$ by $k$ in order to obtain
 $$\ch(\rho^{k}(E))=\frac{1}{k}\frac{\exp(kbx)-1}{\exp(bx)-1}=\frac{1+e^{bx}+\dots+e^{(k-1)bx}}{k}\ .$$
 It follows that
 $$\rho^{k}(E)=\frac{\cycl(E^{0})+\cycl(E^{1})+\dots+\cycl(E^{k-1})}{k}\ .$$
 
 We now discuss 4. 
 Note that the $Spin^{c}$-structure  associated to the complex bundle $V$  differs from
 the original one. It can be written in the form $V^{L}$ for a certain line bundle $L\to M$.
 Using the formulas above we have
 $$\rho^{k}(V)=\frac{1+\cycl(V)+\dots+\cycl(V^{k-1})}{k} \cycl(L^{k-1})\ .$$
 In order to understand $L$ we discuss the structure of $Spin^{c}(2)$ in detail.
 
 The identification $SO(2)\cong U(1)$ together with the map $\beta$ in \eqref{aug0201} define a split 
  $s:SO(2)\to Spin^{c}(2)$ of the exact sequence
 $$0\to U(1)\to Spin^{c}(2)\stackrel{\pi}{\to} SO(2)\to 0\ .$$  Let us make this explicit.   We identify $Spin(2)\cong \tilde U(1)$. We can write $Spin^{c}(2)=Spin(2)\times_{\Z/2\Z} \tilde U(1)$.

 Let $t\in SO(2)$. Then take a preimage $\tilde t\in Spin(2)\cong \tilde U(1)$ and 
 define the element $s(t):=[\tilde  t,\tilde t]\in Spin^{c}(2)$.  The split induces a character
 $\chi:Spin^{c}(2)\to \tilde U(1)$,
 $[\tilde a,\tilde b]\mapsto  \tilde b\tilde a^{-1}$.  If $V\to M$ is a $2$-dimensional $Spin^{c}$-bundle, then using $\chi$ we can associate  the  complex bundle $V(\chi)\to M$.
 
 The homomorphism $c:Spin^{c}(2)\to U(1)$, $[\tilde a,\tilde b]\mapsto \tilde b^{2}=\tilde a^{2}(\tilde b\tilde a^{-1})^{2}$ can be written as
 $c=\pi \chi^{2}$. We have $c\circ \beta=\id$. Hence we must take
 $L=V(\chi)^{-1}$.

 \hB 

The map $MSpin^{c}\to MSO\to H\Z$ provides an ordinary orientation
for every $Spin^{c}$-bundle $V\to M$ and hence a Thom class $\nu_\C\in H\C^{n}(M^{V})$.
We get a characteristic class
$$\Td(V):=\left[\Phi^{-1}_{\nu_\C}(\ch(\nu))\right]^{-1}=\left[\Phi^{-1}_{\nu_\C}(\ch(\Phi_\nu(1)))\right]^{-1}\in HA^{0}(M)\ .$$ 
Do not confuse  this characteristic class   of $Spin^{c}$-vector bundles with the class $\Td(\nu)$  (cf. \ref{jun0620}) associated to the  $K$-theory Thom class $\nu$.
They are inverse to each other.

\begin{prob}\label{jul3110}
\begin{enumerate}
\item
Show that $V\to \Td (V)$ is an exponential characteristic class
$$\Td:\Vect^{Spin^{c}}_\R\to HA^{0}\ .$$ 
\item 
By Chern-Weyl theory there exists a natural characteristic form
$$\Td(\tilde \nabla)\in \Omega A^{0}_{cl}(M)\ .$$
Calculate the corresponding invariant polynomial in $I(Spin^{c}(n))$.
\item 
Give a formula for
$\Td(\tilde \nabla)$ in terms of the curvature $R^{\tilde \nabla}$.\end{enumerate}
\end{prob}
\proof
We discuss 2. briefly. We consider the diagram of Lie groups connecting $Spin^{c}(2n)$ with other classical groups.
\begin{equation}\label{aug0201}\xymatrix{&&U(1)\\U(n)\ar[urr]^{\det}\ar[drr]^{\alpha}\ar[r]^{\beta}&Spin^{c}(2n)\ar[dr]^{\pi}\ar[ur]^{c}&\\&&SO(2n)}\ .\end{equation}
Let us explain the map $\beta$. We start with the presentation 
of $Spin(2n)\to SO(2n)$ as  the two-fold connected covering (universal if $n\ge 2$).
Then we have an isomorphism $$Spin^{c}(2n)=Spin(2n)\times_{\Z/2\Z} U(1)$$ (can be taken as a definition here). We have a natural map $U(n)\to SO(2n)$ which induces a surjective map
$\pi_{1}(U(n))\to \pi_{1}(SO(2n))$.  We form the right pull-back-square in the diagram
$$\xymatrix{U(1)\ar[d]^{(\dots)^{2}}&\tilde U(n)\ar[l]_{\sqrt{\det}}\ar[d]^{p}\ar[r]^{a}\ar[r]&Spin(2n)\ar[d]\\U(1)&
U(n)\ar[l]_{\det}\ar[r]&SO(2n)}\ ,$$ where the vertical maps are two-fold coverings.
As indicated, on $\tilde U(n)$ we can form a square-root of the determinant.

We now define the map
$\beta$ such that it maps
$u\in U(n)$ to the class 
$[a(\tilde u),\sqrt{\det}(\tilde u)]\in Spin^{c}(2n)$, where $\tilde u\in \tilde U(n)$ is a lift of
$u$ under $p$. Note that this class is independent of the choice of the lift $\tilde u$.

If we identify
$$I(U(n))\cong \C[[bx_{1},\dots,bx_{n}]]^{\Sigma_{n}}$$
and $I(U(1))\cong \C[[bt]]$, then the image under $\alpha^{*}$ of
$I(SO(2n))$ in $I(U(n))$ is given by $\C[[(bx_{1})^{2},\dots,(bx_{n})^{2}]]^{\Sigma_{n}}$.
Moreover $\det^{*}(bt)=\sum_{i=1}^{n} bx_{i}$. In particular, the restriction  $$\beta^{*}:I(Spin^{c}(2n))\to I(U(n))$$ is injective. 

We consider the universal line bundle $L\to \C\P^{\infty}$ with its Thom classes
$\nu_{\C}$ and $\nu$ for $H\C$ and $\bKU$. Let $0_{L}:\C\P^{\infty}\to L$ be the zero section.
Then $x_{\C}:=   0_{L}^{*}\nu_{\C}$ generates the additive formal group law, while
$x:=  0_{L}^{*}\nu$ generates the multiplicative formal group law.
Hence we have the relation $\ch(b x)=\exp(bx_{\C})-1$.
Let $ Q\in HA^{0}(\C\P^{\infty})$ be the class such that $\ch(b  \nu)=Q b\nu_{\C}$, i.e. $Q=\Td(L)^{-1}$. Then
$\exp(bx_{\C})-1=\ch(b x)=Q bx_{\C}$.
We see that
$$Q=\frac{\exp(bx_{\C})-1}{bx_{\C}}\ .$$
Using 
$$ \frac{e^{x}-1}{x}= \frac{\sinh(x/2)}{x/2}e^{x/2}$$
we conclude that
$$\beta^{*}\Td=e^{-\frac{1}{2}\sum_{i=1}^{n}bx_{i}}\prod_{i=1}^{n} \frac{bx_{i}/2}{\sinh(bx_{i}/2)} \ .$$
If we define $\hA\in I(BSO(2n))$ such that $$\alpha^{*}\hA=\prod_{i=1}^{n} \frac{bx_{i}/2}{\sinh(bx_{i}/2)}\ ,$$
then we see that
$$\Td= \pi^{*}\hA \exp^{-\frac{c^{*}bt}{2}} .$$

\hB

We let $\Vect^{Spin^{c},geom}_\R$ be the smooth monoid of $Spin^{c}$-bundles with $Spin^{c}$-connection with respect to the direct sum. 
 \begin{prob}
\begin{enumerate}
\item Show that there exists a unique  differential refinement
$$\hat \rho^{k}:\Vect^{Spin^{c},geom}_\R \to (\hbKU[\{k\}^{-1}]^{0})^{\times}$$
such that
$$I(\hat \rho^{k}(V,\tilde \nabla))=\rho^{k}(V)\ , \quad R(\hat \rho^{k}(V,\tilde \nabla))=\Psi^{k}_H(\Td(\tilde \nabla))^{-1}\wedge \Td(\tilde \nabla)\ .$$
\item Show that $\hat \rho^{k}$ is exponential.
\item Calculate
$\hat \rho^{k}((V,\tilde \nabla)^{(L,\nabla^{L})})\cup \hat \rho^{k}(V,\tilde \nabla)^{-1}$.
\end{enumerate}
\end{prob}
For some of the solutions see \cite{MR2740650}. \hB

}\end{ex}

\subsection{Orientation and integration - the topological case}\label{aug0501}
 
We first recall the homotopy theoretic construction of the Umkehr map.
For a space $X\in \Top$ and $k\in \Z$ we abbreviate by  $X_+^{k}:=\Sigma^{\infty+k}_+X \in \Nerve(\Sp)[W^{-1}]$  the suspension spectrum of $X$ shifted by $k$. For a real vector bundle $N\to X$ we let $X^{N}\in \Nerve(\Sp)[W^{-1}]$ denote the Thom spectrum of $N$. Note that $X^{X\times \R^{k}}\cong  X_+^{k}$.
  
Let $f:W\to M$ be a proper smooth map. 
For simplicity we assume that   $M$ is connected and all components of  $W$ have the same dimension. We set
 $n:=\dim(M)-\dim(W)$.
For sufficiently large $k$  we can choose a fibrewise embedding
$\iota:W\hookrightarrow M\times \R^{l}$. Its differential provides a representative of the stable normal bundle (see Definition \ref{may2941}) 
$$\cN: 0\to TW \stackrel{d\iota}{\to} f^{*}TM\oplus \R^{l}\to N\to 0\ .$$
We can extend $\iota$ to an open embedding also denoted by 
$\iota:N\to M\times \R^{l}$. If $f:W\to M$ is a submersion, then we can require that this extension
is a map over $M$. This condition will play an important role when we consider the pull-back of representatives of  orientations.

In any case,
it gives rise to a map of Thom spectra
 \begin{equation}\label{may2950}\cl(\iota):M_+^{l }\to W^{N} \end{equation} 
called the clutching map whose description we recall in the following.
We let $\bar N=N\cup\partial N$ be the fibrewise one-point compactification of $N$.
Then we define a map $(M\times \R^{l})_+\to \bar N/\partial N$ which is the inverse $\iota^{-1}$ on $\iota(N)$
 and sends the complement of $\iota(N)$ and the basepoint to the basepoint represented by the contracted boundary $\partial N$.
Stabilization gives the map of spectra $\cl(\iota):M_+^{l}\to W^{N}$.


Let $E$ be a commutative ring spectrum.
For any real vector bundle $N\to W$ over a compact base
 we have an identification
$$E^{*}_{prop/W}(N)\cong E^{*}(W^{N})\ .$$
\begin{prob}
Give an argument.
\end{prob}

Hence we can talk about a  representative $\nu:W^{N}\to E[l+n]$ of a Thom class. 


\begin{ddd}
A representative $(\iota,\nu)$ of an $E$-orientation on the map $W\to M$
is given by the data of an embedding $\iota:N\to M\times \R^{l}$ and a representative of a Thom class
$W^{N}\to E[l+n]$.
\end{ddd}
If $N\to W$ is a real vector bundle, then we have the Thom diagonal
$\diag:W^{N}\to W^{N}\wedge W_{+}$.  
The representative of an $E$-orientation gives rise to a map in $\Mod(E)$  
$$I(\iota,\nu): E \wedge M^{l}_+ \stackrel{\cl(\iota)}{\to}   E \wedge  W^{N} \stackrel{\diag}{\to}E \wedge   W^{N}\wedge W_+ \stackrel{\nu}{\to} E\wedge E[l+n] \wedge  W_+  \stackrel{\mu}{\to}  E[l+n] \wedge W_+\ .$$
\begin{ddd}\label{jul2910}
For every $E$-module spectrum $F\in \Mod(E)$ the map 
$$I(\iota,\nu)_!:=\Map_{\Mod(E)}(I(\iota,\nu),F):\Funk_M(F)(W)\to \Funk_M(F)(M)[n]$$
and the induced map on the level of  homotopy groups
$$\bar I(\iota,\nu)_!:F^{*}(W)\to F^{*+n}(M)$$
are called the Umkehr or integration maps. 
\end{ddd}
The usual notation for the Umkehr map is  $f_!:=\bar I(\iota,\nu)_{!}$. But  in the present section this conflicts with the notation for the proper push-forward of a smooth object so it will be avoided.



 \begin{prob}\label{jul1201}
 In the case that $f:W\to M$ is a proper submersion
give an interpretation of the Umkehr map as a transformation
$$I(\iota,\nu)_!: f_!f^{*} \Funk_M(F)\to \Funk_M(F)[n]$$
\end{prob}
\proof
Since $f$ is a submersion we can and will choose the embedding $\iota:N\to M\times \R^{l}$ as a map over $M$.
If $h:M^{\prime}\to M$ is a manifold over $M$, then there is a functorial way to construct
a pull-back $h^{*}(\iota,\nu)$ of representatives of $E$-orientations. 
We have a pull-back diagram
$$\xymatrix{N^{\prime}\ar[r]^{\lambda}\ar[d]^{\pi^{\prime}}&N\ar[d]\\ W^{\prime}\ar[r]^{H}\ar[d]^{f^{\prime}}&W\ar[d]^{f}\\
M^{\prime}\ar[r]^{h}&M}$$
with $W^{\prime}:=M^{\prime}\times_{M}W$ and $N^{\prime}:=H^{*}N$. We can identify the vector bundle $N^{\prime}\to W^{\prime}$ with the normal bundle of an induced embedding $\iota^{\prime}$ as follows. Let $\iota_{2}:N\to \R^{l}$ be the second component of $\iota$.
We consider the embedding
$$\iota^{\prime}:N^{\prime}\stackrel{(f^{\prime}\circ \pi^{\prime})\times\lambda}{\to} M^{\prime}\times N\stackrel{\id_{M^{\prime}}\times \iota_{2}}{\to} M^{\prime}\times \R^{l}\ .$$  
The corresponding representative of the stable normal bundle is
$$\cN^{\prime}:0\to TW^{\prime}\stackrel{df^{\prime}\oplus H^{*}\alpha}{\to} f^{\prime*}TM^{\prime}\oplus \R^{l}\to N^{\prime}\to 0\ .$$
We define the representative of the Thom class
$$\nu^{\prime}:W^{\prime N^{\prime}}\stackrel{\lambda}{\to} W^{N}\stackrel{\nu}{\to} E[l+n]$$
and let $$h^{*}(\iota,\nu):= (\iota^{\prime},\nu^{\prime})\ .$$

If 
$$\xymatrix{A\ar[rr]^{g}\ar[dr]^{p}&&B\ar[dl]_{q}\\&M&}$$
 is a smooth map over $M$, then there   is  a functorially  induced diagram in $\Mod(E)$ 
$$\xymatrix{ E\wedge A_+^{l}\ar[d]^{I(p^{*}(\iota ,\nu ))}\ar[r]&E\wedge B_+^{l}\ar[d]^{I(q^{*}(\iota,\nu))}\\
 E[l+n]\wedge (A\times_MW)_+\ar[r]&E[l+n]\wedge (B\times_{M}W)_+}\ .$$
 We consider this construction as a functor
 $$\Nerve(\Mf/M)\to \Fun(\Nerve([1])^{op},\Mod(E))\ .$$
 We take the adjoint and compose with
 $\Map_{\Mod(E)}(\dots,F)$. This gives a functor
 $\Nerve(\Mf/M)\times \Nerve([1])\to \Mod(E)$.
 Its evaluation at $0\in [1]$ can be identified with $f_{*}f^{*} \Funk_{M}(F)$, while
 the evaluation at $1\in [1]$ is $\Funk_{M}(F)[n]$.
 Finally, since $f$ is proper we have $f_{!}=f_{*}$. 
  \hB  

\begin{ex}\label{aug1301}{\em 
In the definition of a homotopy between two representatives of topological $E$-orientations we  need the following extension of the construction of the pull-back.
Assume that $M$ decomposes as a product $M=M_{1}\times M_{2}$  and that $h_{1}:M_{1}^{\prime}\to M_{1}$ is a smooth map. For a homotopy, $M_{1}=[0,1]$ and $M_{1}^{\prime}=\{0,1\}$.  We further assume that there exists an open subset
$U\subseteq M_{1}$ containing the image of $h_{1}$ such that the projection
$\pr_{1}\circ f:f^{-1}(U\times M_{2})\to U$
is a submersion. In this case we can choose $\iota$ such that its restriction
 $\iota_{|\pi^{-1}(f^{-1}(U\times M_{2}))}$ is a map over $U$.
Under these assumptions
 we still can define the pull-back
$h^{*}(\iota,\nu)$ by the construction above, where $h=h_{1}\times \id_{M_{2}}$.
}\end{ex}



 
 \begin{ddd}
 We say that two representatives $(\iota_{i},\nu_{i})$, $i=0,1$, of $E$-orientations are homotopic if there exists a representative of an $E$-orientation $ (\tilde \iota,\tilde \nu)$ on $\id\times f:I\times W\to I\times M$ with the property that $\iota:\tilde N\to I\times M\times \R^{l}$ preserves fibres over $I$ in neighbourhoods of the endpoints of the interval, and     
  which pulls back to the representatives of $E$-orientations  $(\iota_{i},\nu_{i})$ at the end-points of the interval $I=[0,1]$. 
  \end{ddd}
  There is a natural operation of stabilization of  the embedding
$\iota:N\to M\times \R^{l}$ to an embedding $\iota^{r}:N\oplus \R^{r}\to M\times \R^{l+r}$.
Note that $W^{N\oplus \R^{r}}\cong \Sigma^{r}W^{N}$.
Let $\nu^{r}:=\Sigma^{r}\nu$ and set $\Sigma^{r}I(\iota,\nu):= (\iota^{r},\Sigma^{r}\nu)$. 
Then we have 
$$I(\iota,\nu)_{!} =I(\Sigma^{r}(\iota,\nu))_{!} : \Funk(F)(W)\to \Funk(F)(M)[n]\ .$$

\begin{ddd}
An $E$-orientation of $f$ is an equivalence class $[\iota,\nu]$ of representatives of $E$-orientations
under the equivalence relation generated by homotopy and stabilization.
\end{ddd}
It is clear that the Umkehr map  in $F$-cohomology
$$\bar I[\iota,\nu]_! :F^{*}(W)\to F^{*+n}(M)$$
as a function of the $E$-orientation
is well-defined

\begin{ex}{\em 
In this example we discuss a prototypical case of an  index theorem. An index theorem is a statement which relates, via cycle maps, a geometric integration with the topological integration map. The prototypical example, which gave the name, is the Atiyah-Singer index theorem which relates the analytic  push-forward of vector bundles in terms of twisted Dirac operators with the topological integration in $\bKU$-theory. 
Let $f:W\to M$ be a proper submersion with representative 
$$\cN_{f}:0\to TW\stackrel{df\oplus \beta}{\to} f^{*}TM\oplus \R^{l}\to N_{f}\to 0$$ 
of the stable normal bundle coming from an embedding $\iota$, and  a complex structure $I_{f}$ on $N_{f}$. The complex structure gives rise to a $\bMU$-Thom class $\nu_{\bMU}\in \bMU^{\dim(N_{f})}(N_{f})_{prop/W}$ as in \ref{jul2610}.
We   let
$\bar I(\iota,\nu)_{!}:\bMU^{*}(W)\to \bMU^{*+n}(M)$ be the integration, where $n:=\dim(M)-\dim(W)$.
 We now define a geometric integration
 $$I:\Cycle_{\bMU}^{*}(W)\to \Cycle_{\bMU}^{*+n}(M)$$ as follows. Let $(g:A\to W,\cN_{g},I_{g})$ be a $\bMU$-cycle on $W$. Then we define
 $$\cN_{f\circ g}:0\to TA\stackrel{dg\oplus \alpha}{\hookrightarrow} g^{*}TW\oplus \R^{k}\stackrel{g^{*}df\oplus g^{*}\beta\oplus \id_{\R^{k}}}{\hookrightarrow} g^{*}f^{*}TM\oplus \R^{l}\oplus \R^{k}\to N_{f\circ g}\to 0$$
 (contract the first to injective maps to one). 
 The bundle $N_{f\circ g}$ has a natural filtration $g^{*}N_{f}\subseteq N_{f\circ g}$ with quotient $N_{g}$.
 We can therefore obtain a well-defined (up to homotopy) identification $N_{f\circ g}\cong g^{*} N_{f}\oplus N_{g}$ by chosing a split.  Note that in view of the definition \ref{jul3001} of $\Cycle_{\bMU}$ it is only the homotopy class of the representative of the stable normal which matters.
 We define the geometric integral of the cycle  $(g:A\to W,\cN_{g},I_{g})$ by $$I(g:A\to W,\cN_{g},I_{g}):=(f\circ g:A\to M,\cN_{f\circ g}, g^{*}I_{f}\oplus I_{g})\ .$$
 
\begin{prob}
Show the index theorem asserting that the following diagram commutes:
$$\xymatrix{\Cycle_{\bMU}^{*}(W)\ar[r]^{cycl}\ar[d]^{I}&\bMU^{*}(W)\ar[d]^{\bar I(\iota,\nu)_{!}}\\
\Cycle_{\bMU}^{*+n}(M)\ar[r]^{\cycl}&\bMU^{*+n}(M)}\ .$$
\end{prob}
\proof
This should follow by a careful analysis of the Thom-Pontrjagin construction. \hB 
}
\end{ex}

\begin{ex}{\em
In the following we contrast index theorems with Riemann-Roch theorems.
A Riemann-Roch theorem is a statement about the compatibility of a natural transformation between cohomology theories and integration maps. In the example below we compare integration in $E$-theory with ordinary cohomology. We assume that the characteristic of the ring $\pi_{*}(E)$ is not $2$. Then 
an $E$-Thom class on a vector $N\to W$ determines an ordinary orientation, hence a Thom class $\nu_{\C}\in H\C^{\dim(N)}(W^{N})$.
\begin{ddd}\label{jul2701}
We define the class
$\Td(\nu)\in H(A)^{0}(W)^{\times}$ uniquely such that 
$$\Td(\nu)\cup \nu_{\C}\cong c(\nu)\ .$$  
 \end{ddd}
The following proposition is an immediate consequence of the defnitions.
\begin{prop}\label{aug0601}
We have the Riemann-Roch theorem 
$$
\xymatrix{E^{k}(W)\ar[d]^{ \bar I(\iota,[\nu])_! }\ar[r]^{c}&H(A)^{k}(W)\ar[d]^{\bar  I(\iota,[c(\nu)])_!}\ar@{=}[r]&H(A)^{k}(W)\ar[d]^{ \bar I(\iota,  [\nu_\C])_! (\Td(\nu)\cup\dots)}\\E^{k+n}(M)\ar[r]^{c}
&H(A)^{k+n}(M)\ar@{=}[r]&H(A)^{k+n}(M)}\ .$$
\end{prop}

 \begin{prob}
Assume that $W\to M$ be a proper holomorphic map between complex manifolds.
Make the Riemann-Roch theorem for complex bordism explicit. In particular consider the cases bundles of curves and surfaces.
 \end{prob}
 \proof
 The exercise consists in the calculation of $\Td(\nu)$.
 }
 \end{ex}
 \subsection{Orientation and integration - the differential  case}

Let $(E,A,c)$ be a multiplicative datum. 
We  assume that $f:W\to M$ is a proper smooth map and set $n:=\dim(M)-\dim(W)$. As in the topological case   we choose an embedding $\iota:W\to M\times \R^{l}$ over $M$ and an extension to an open embedding, also denoted by $\iota:N\to M\times \R^{l}$, of the normal bundle.  We further choose a Thom class $\nu \in E^{k}_{prop/W}(N)$, where $k=l+n$ is the dimension of the normal bundle $\pi:N\to W$.
These structures determine the topological integration map (Definition \ref{jul2910}) 
$$\bar I(\iota,\nu)_{!}: F^{*}(W)\to F^{*+n}(M)$$ for every $E$-module spectrum  $F$. It is given by the composition 
$$\bar I(\iota,\nu)_{!}:F^{*}(W)\stackrel{\pi^{*}}{\to} F^{*}(N)\stackrel{\nu\cup\dots}{\to} F^{*+k}_{prop/W}(W)\stackrel{excision}{\to} F^{*+k}_{prop/M}(M\times \R^{l})\stackrel{\desusp}{\to} F^{*+n}(M)\ .$$

In order to generalize this to the differential case we need the notion of a differential $E$-orientation.
\begin{ddd}
A representative of a differential $E$-orientation refining $(\iota,\nu)$ is a pair $(\iota,\hat \nu)$
where $\hat \nu$ is  a differential Thom class refining $\nu$ (see \ref{jun0620}).
\end{ddd}

For simplicity we only consider the integration for  the differential extension of $E$ itself. The generalization to modules is straight forward, see \ref{jun0701}.

\begin{ddd}
We  define the integration in differential $E$-theory
$$\hat I(\iota,\hat \nu)_{!}:\hat E^{*}(W)\to \hat E^{*+n}(M)$$
associated to the representative of the differential $E$-orientation $(\iota,\hat \nu)$
by the following composition
\begin{equation}\label{aug0903}\hat I(\iota,\hat \nu)_{!}:\hat E^{*}(W)\stackrel{\pi^{*}}{\to} \hat E^{*}(N)\stackrel{\hat \nu\cup\dots}{\to} \hat E^{*+k}_{prop/W}(N)\stackrel{excision}{\to} \hat E^{*+k}_{prop/M}(M\times \R^{l})\stackrel{\desusp}{\to} \hat E^{*+n}(M)\ . \end{equation}
\end{ddd}
Note that we use \ref{aug1302} in the second map.

If $f$ is a proper submersion, then we can choose the embedding $\iota:N\to M\times \R^{l}$ such that it is a map over $M$, i.e. that the diagram
$$\xymatrix{N\ar[d]^{\pi}\ar[rd]^{\iota}&\\W\ar[d]^{f}\ar[r]^{\iota}&M\times \R^{l}\ar[d]^{\pr_{M}}\\
M\ar@{=}[r]&M}$$ 
commutes.
\begin{prob}
Check that
if $\iota:N\to M\times \R^{l}$ is a map over $M$, then the differential integration $\hat I(\iota,\nu)_{!}$ is a morphism of $\hat E^{*}(M)$-modules.
\end{prob}

We define the map of complexes 
$$R(\iota,\hat \nu):\Omega A^{*}(W)\to \Omega A^{*}(M)[n]$$
by
$$R(\iota,\hat \nu)(\alpha):=\int_{M\times \R^{l}/M} (\iota^{-1})^{*}\left[R(\hat \nu) \wedge\pi^{*} \alpha\right]\ , \quad \alpha\in \Omega A(W)\ .$$
The  form $(\iota^{-1})^{*}\left[R(\hat \nu) \wedge\pi^{*} \alpha\right]$ is  first defined on $\iota(N)$. After extension by zero it can be considered as a smooth form on $M\times \R^{l}$ which is properly supported over $M$. So the integral is well-defined. It follows from the closedness of $R(\hat \nu)$ and Stokes' theorem that  $R(\iota,\hat \nu)$ commutes with the differential.
We call $R(\iota,\hat \nu)$ the curvature map of the integration associated to $(\iota,\hat \nu)$.
This is justified by the following exercise.
 \begin{prob}
Check that for $x\in \hat E^{*}(W)$ and $\alpha\in \Omega A(W)$ we have
$$R(\hat I(\iota,\nu)_! (x))=R(\iota,\hat \nu)(R(x))\ , \quad \hat I(\iota,\hat \nu)_!(a(\alpha))=a(R(\iota,\hat \nu)(\alpha))\ .$$
\end{prob}
Sometimes we write the curvature map $R(\hat I_{!})$ as a function of the integration $\hat I_{!}:=\hat I(\iota,\hat \nu)$.
Summing up, the differential integration fits into the   following commuting diagram: $$\xymatrix{\Omega A^{*-1}(W)/\im(d)\ar[d]^{ R(\iota,\hat \nu)}\ar[r]^{a}&\hat E^{*}(W)\ar[d]^{\hat I(\iota,\hat \nu)_{!}}\ar[r]^{I}\ar@/^1cm/[rr]^{R}&E^{*}(W)\ar[d]^{I(\iota,\nu)_{!}}&\Omega A_{cl}^{*}(W)\ar[d]^{ R(\iota,\hat \nu)}\\\Omega A^{*+n-1}(M)/\im(d)\ar[r]^{a}&\hat E^{*-n}(M)\ar[r]^{I}\ar@/_1cm/[rr]_{R}&E^{*+n}(M)&\Omega A^{*+n}_{cl}(M) }\ .$$

Similar constructions of the integration in differential cohomology have been considered in \cite{MR2192936} and \cite{2012arXiv1208.1288F}.

\begin{prob}\label{aug0811}
If the embedding $\iota:N\to M\times \R^{l}$ is a map over $M$, then the curvature map of the integration associated to $(\iota,\hat \nu)$ is given by
$$R(\iota,\hat \nu)(\alpha)=\int_{W/M}\Td(\hat \nu)\wedge \alpha\ ,$$
where $\Td(\hat \nu)$ is defined in \ref{jun0620}.
\end{prob}

We now consider the compatibility of the push-forward with pull-back. We  assume that $f:W\to M$ is a proper submersion and that $\iota:N\to M\times \R^{l}$ is a map over $M$ in order to ensure existence of all pull-backs which are necessary in the construction.
 Let $h:M^{\prime}\to M$ be  a smooth map and consider    
 the bundle $ f^{\prime}:M^{\prime}\times_{M}W\to M^{\prime}$.  We use the notation introduced in the proof of \ref{jul1201}
and define $\hat \nu^{\prime}:=\lambda^{*}\hat \nu$.
Then $h^{*}(\iota,\hat \nu):=(\iota^{\prime},\hat \nu^{\prime})$ is the pulled-back representative of a differential $E$-orientation of $f^{\prime}$.

\begin{prob}
Show that a representative of a differential $E$-orientation $(\iota,\hat\nu)$ on $W\to M$ gives rise to
a transformation
$$\hat I_{!}(\iota,\hat \nu):f_{!}f^{*}  \hat E^{*}_{M}\to \hat E^{*}_{M}\ .$$  
 \end{prob}
 \proof The arguments are similar to \ref{jul1201}. \hB
 This assertion encodes the compatibility of the integration with pull-back diagrams.
 Again we can extend the definition of the pull-back as in \ref{aug1301}.

\begin{ddd}\label{aug1501}
 We say that two representatives $(\iota_{i},\hat \nu_{i})$, $i=0,1$, of differential $E$-orientations  are homotopic if there exists a representative of a differential  $E$-orientation $ (\tilde \iota,\tilde {\hat \nu})$ on $\id\times f:I\times W\to I\times M$  with the property that $\iota:\tilde N\to I\times M\times \R^{l}$ preserves fibres over $I$ in neighbourhoods of the endpoints of the interval,   
 which restricts to $(\iota_{i},\hat \nu_{i})$, $i=0,1$, and whose 
 curvature map satisfies   
 $$R(\tilde \iota,\tilde{\hat \nu})\circ \pr_{W}^{*} = \pr^{*}_M \circ R(\iota_0,\hat \nu_0) \ ,$$
  where   $\pr_M:I\times M\to M$ and $\pr_{W}:I\times W\to W$ are  the projections.
 \end{ddd}
 
 \begin{prob}
 Assume that $\iota:N\to M\times \R^{l}$ is a map over $M$ and that 
 two differential $E$-Thom classes  $\hat \nu_i$, $i=0,1$, are homotopic in the sense of 
 \ref{aug0850}. 
 Check that   the representatives of differential $E$-orientations  $(\iota,\hat \nu_{i})$ are homotopic.
 \end{prob}
 \proof
 Use \ref{aug0811}. \hB 
 

\begin{prob}\label{aug0632} Show that two  homotopic representatives of a differential $E$-orientations induce the same integration map 
$\hat E^{*}(W)\to \hat E^{*+n}(M)$.
\end{prob}
\proof
Use the homotopy formula \ref{aug1003}. \hB 

Next we discuss the stabilization of differential $E$-orientations.  
Recall, that  we can define a stabilization  $(\iota^{r},\Sigma^{r}\nu)$ of a representative of a topological $E$-orientation such that  the topological integration map does not change  on the spectrum level, i.e.  we have 
$I(\iota,\nu)=I(\iota^{r},\Sigma^{r}\nu)$.  
In particular, stabilization does not effect the integration map in cohomology groups.

In the differential case we again stabilize the embedding as before.  In addition
we must define a stabilization of differential $E$-Thom classes. This definition involves the choice of differential Thom classes for trivial bundles. Unlike the topological case there is no natural choice, but if we assume that our data $(E,A,c)$ is strict, then there is a natural choice up to homotopy of Thom classes.  

It suffices to choose a differential Thom class $\hat \nu_{1}\in \hat E^{1}_{c}(\R)$.
We fix this choice   by the condition that
\begin{equation}\label{aug0815}\desusp(\hat \nu_{1})=1\ , \end{equation}
where $\desusp$ is the desuspension map defined in \ref{aug0630}.
Indeed, if $y\in H(A)^{-1}(*)$, then $$\desusp(\hat \nu^{\prime}_{1}+a(R(\hat \nu^{\prime}_{1})\wedge y))=\desusp(\hat \nu^{\prime}_{1})+y\ .$$
We  start with any differential lift   $\hat \nu^{\prime}_{1}$ of the canonical topological Thom class and set $z:=\desusp(\hat \nu^{\prime}_{1})$.
Then $$\hat \nu_{1}:=\hat \nu_{1}^{\prime}- a(R(\hat \nu^{\prime}_{1})\wedge z)$$
satisfies \eqref{aug0815}. Furthermore,  
by \ref{aug0850} it is uniquely determined by this condition up to homotopy of differential Thom classes.
We now define $$\hat \nu_{l}:=\pr_{1}^{*}\hat \nu_{1}\cup \dots\cup \pr_{l}^{*}\hat \nu_{1}\in E^{l}_{c}(\R^{n})\ .$$ Here we use \ref{aug1302} to get compact support of the product.
 By pull-back we then get Thom classes, also denoted by $\hat \nu_{l}$, on trivialized  bundles $M\times \R^{l}$.

We now define the stabilization of the differential Thom class $\hat \nu$ by
$$\Sigma^{l}\hat \nu :=p^{*}\hat \nu\cup q^{*}\hat \nu_{l}\ ,$$
where $p:N\oplus (W\times \R^{l})\to N$ and $q:N\oplus (W\times \R^{l})\to (W\times \R^{l})$ are
the projections. This determines a stabilized representative $(\iota,\Sigma^{l}\hat \nu)$ of a differential $E$-orientation which is natural up to homotopy of Thom classes.
\begin{prob}\label{aug0633}
Show the equality of differential integration maps
$\hat I(\iota^{l},\Sigma^{l}\hat \nu)_{!}=\hat I(\iota,\hat \nu)_{!}$.
\end{prob}
 \proof
This uses \eqref{aug0815},  \ref{aug0901},  and desuspension in stages \ref{aug0902}.  \hB


We assume that the datum $(E,A,c)$ is strict in order to have a well-defined notion of stabilization on the level of homotopy classes of  representatives of differential $E$-orientations.

\begin{ddd}
 On the set of representatives of differential $E$-orientations we consider the equivalence relation generated by stabilization and homotopy. 
A differential $E$-orientation is an equivalence class of representatives of differential $E$-orientations.
\end{ddd}

\begin{kor}
It follows from \ref{aug0632} and \ref{aug0633} that  the integration
$$\hat I(\iota,\hat \nu)_{!}:\hat E^{*}(W)\to \hat E^{*+n}(M)$$
only depends on the equivalence class $[\iota,\hat \nu]$ of differential $E$-orientations.
\end{kor}

\begin{prob}\label{aug0640}
Assume that $f$ is a proper submersion, $(\iota_{0},\hat \nu_{0})$ is a differential $E$-orientation of $f$ where $\iota_{0}:N\to M\times \R^{l}$ is a map over $M$.  Classify differential $E$-orientations refining the  topological $E$-orientation $(\iota_{0},\nu_{0})$  with curvature map $R(\iota_{0},\hat \nu_{0})$.
\end{prob}
\proof We can work with embeddings $\iota:N\to M\times \R^{l}$ which are maps over $M$. Any two become homotopic in this class after suitable stabilization.
We now observe using \ref{aug0811} that a homotopy of differential orientations with fixed
$\iota$  exactly corresponds to a homotopy of differential Thom classes.
It suffices therefore to classify stable homotopy classes of differential Thom classes. Here we use \ref{aug0610}. We finally conclude that
the set of differential $E$-orientations $[\iota,\hat \nu]$  is a torsor over the group $$\frac{H(A)^{-1}(W)}{c(E^{-1}(W))}\ .$$
  \hB 
  
  \begin{prob} \label{aug0870} 
  Drop the assumption that $f$ is a submersion. Fix a curvature map $R(\iota_{0},\hat \nu_{0})$ and
  classify  differential $E$-orientations refining the  topological $E$-orientation $(\iota,\nu)$  with curvature map $R(\iota,\hat \nu)$.
\end{prob}
\proof
\textcolor{red}{I do not know the answer. }
\hB

In \cite[Def. 3.5, Cor. 3.6]{MR2664467}  we introduced a notion of a differential $\bKU$-orientation and discussed the classification of differential $\bKU$-orientations refining a given topological one. The definitions and results  of \cite{MR2664467} differ from the theory presented here. See \ref{aug0701} for more details. A similar remark applies to the comparison with the definitions and the results of \cite{MR2550094} for the complex bordism $\bMU$.



\begin{ex}\label{jun0740}{\em

Consider the differential cohomology theories $\widehat{H\Z}^{*}$ or $\hat \bS^{*}$. 
In this case $\Td(\hat \nu)=1$ for every differential Thom class. Since also $H(A)^{-1}(W)=0$ it follows  that there is a unique up to homotopy  differential Thom class $\hat \nu$ which refines $\nu$. 

  If $f:W\to M$ is a proper submersion with a topological orientation $(\iota,\nu)$, then there exists a preferred differential orientation $(\iota,\hat \nu)$ refining this topological orientation such that $\iota:N\to M\times \R^{l}$ is a map over $M$.  
In particular there exist  preferred differential integration  maps
$\hat I^{H\Z}$ and $\hat I^{\bS}$. \hB

Note that $S^{1}$ has a trivialized normal bundle and is therefore $\bS$-oriented. Hence the projection
$\pr_{M}:S^{1}\times M\to M$ is canonically $\bS$-oriented. We therefore have an integration
$$\hat I_{!}^{F}:\hat F^{*}(S^{1}\times M)\to F^{*-1}(M)\ .$$  The integration over the $S^{1}$-factor played an important role in the uniqueness considerations in  \cite{MR2608479}. In this reference we required
the following additional properties:

\begin{prob}
Show that
\begin{enumerate}
\item If $t:S^{1}\to S^{1}$ is inversion, then $\hat I^{F}_{!}(t^{*}x)=-\hat I_{!}^{F}(x)$.
\item $\hat I^{F}_{!}(\pr_{M}^{*}x)=0$.
\end{enumerate}
\end{prob}

}\end{ex}

\begin{prob}\label{jun1401}
Verify the following assertions:
\begin{enumerate}
 \item If $g:M\to U$ is a second proper submersion, $m=\dim(U)-\dim(M)$ and $(\iota_{f},\hat \nu_{f})$, $(\iota_{g},\hat \nu_{g})$   are differential $E$-orientation    for $f$ and $g$, then there is a natural   construction of a differential $E$-orientation $(\iota_{g\circ f},\hat \nu_{g\circ f})$ 
   such that
$$\hat I(\iota_{g\circ f},\hat \nu_{g\circ f})_{!}=\hat I(\iota_{g},\hat \nu_{g})_{!}\circ \hat I(\iota_{f},\hat \nu_{f})_!:\hat E^{k}(W)\to \hat E^{k+n+m}(U)\ .$$
\item Composition of differential orientations commutes with pull-back. 
\item
We consider a morphism of multiplicative differential data
$\phi:(E,A,c)\to (E^{\prime},A^{\prime},c^{\prime})$. A differential $E$-orientation $(\iota,\hat \nu)$ naturally induces a differential $E^{\prime}$-orientation
$(\iota,\phi_{*}(\hat \nu))$, and we have the equality 
$$\phi_*\circ \hat I(\iota,\hat \nu)_{!}=\hat I(\iota,\phi_{*}(\hat \nu))_{!}\circ   \phi_{*}:\hat E^{k}(W)\to \hat E^{\prime,k+n}(M)\ .$$

\end{enumerate}
\end{prob}
\proof
For 1. one just writes out the composition of the sequences \eqref{aug0903} for $f$ and $g$.
Then one adds commuting cells which connect this composition with corresponding sequence
\eqref{aug0903} for the composition $g\circ f$. Note that the embedding
$\iota_{g\circ f}:N_{g\circ f}:=N_{g}\times_{M}N_{f}\to U\times \R^{l+k}$ is canonically induced by the embeddings
$\iota_{g}$ and $\iota_{f}$, and
$\hat \nu_{g\circ f}=\pr_{N_{f}}^{*}\hat \nu_{f}\cup \pr_{N_{g}}^{*}\hat \nu_{g}$.
 In the argument one should further use \ref{aug0901} and \ref{aug0902}. 

For assertion 2. one should write out an even bigger diagram.

The argument for assertion 3. is clear. \hB

\begin{ex}{\em 

We now introduce the notion of a differentially $E$-oriented zero bordism of a  smooth proper map $f:W\to M$ with a differential $E$-orientation $(\iota,\hat \nu)$. It is given by a smooth proper map
$F:Z\to I\times M$ with a differential $E$-orientation $(\tilde \iota,\tilde {\hat \nu})$
such that
\begin{enumerate}
\item near the  end-points of the interval the composition $\pr_{I}\circ F:Z\to M\to I$ is a submersion and the embedding 
$\tilde \iota:\tilde N\to I\times M\times \R^{l}$ is a map over $I$,
\item the restriction of the differentially $E$-oriented map $F$ to $\{1\}\times M$ is 
 isomorphic to the differentially $E$-oriented map $f$, and 
\item $\{0\}\times M$ does not intersect  the image of $F$.
\end{enumerate}
Observe the similarities with the notion of a homotopy \ref{aug1501}.
Note that there are no restrictions on $F$ over points of $I$ which outside of some neighboorhod of the endpoints of the interval.


\begin{prob}[Bordism formula]\label{jun0839}
  Show  that for $x\in \hat E^{*}(Z)$ we have the bordism formula
$$\hat I(\iota,\hat\nu)_{!} (x_{|W})=a\left(\int_{I\times M/M}R(\tilde \iota,\tilde{\hat \nu}) (R(x))\right)\ .$$
\end{prob}
\proof
This is a consequence of the homotopy formula Prop. \ref{apr2201} applied to
$\hat I(\tilde \iota,\tilde{\hat \nu})_{!}(x)$. Note that
$\hat I(\tilde \iota,\tilde{\hat \nu})_{!}(x)_{|\{0\}\times M}=0$ and
$\hat I(\tilde \iota,\tilde{\hat \nu})_{!}(x)_{|\{1\}\times M}=\hat I(\iota,\hat\nu)_{!} (x_{|W})$.

\hB

}
\end{ex}

\begin{ex} \label{jun0701}
{\em
Let $E\in \CommMon(\Nerve(\Sp)[W^{-1}])$ be a commutative ring spectrum. Then we can form the $\infty$-category $\Mod(E)$ of $E$-module spectra.
We consider a multiplicative differential data $(E,A,c)$. 
If $F$ is an $E$-module, then we can define a notion of a differential $(E,A,c)$-module data $(F,B,d)$ in a natural way so that $\Diff^{\bullet}(F,B,d)$ becomes a $\Diff^{\bullet}(E,A,c)$-module (cf. \ref{jun06101}).

 Let now $f:W\to M$ be a proper smooth map. 
If  $(\iota,\hat \nu)$ is a differential $E$-orientation of $f$, then we obtain a differential $F$-integration $\hat I^{F}(\iota,\hat \nu)_{!}$   in a natural way. 
\begin{prob}
Work out the details.
\end{prob}
If the differential $E$-orientation is clear from the context we often use the shorter notation
$\hat I^{F}_{!}:=\hat I^{F}(\iota,\hat \nu)_{!}$.

 Every spectrum $F$ is an $\bS$ -module.
Assume that $f:W\to M$ is a proper submersion which is  topologically oriented for $\bS$.  Then by \ref{jun0740} there exists a preferred differential orientation $(\iota,\hat \nu)$ where $\iota:N\to M\times  \R^{l}$ is a map over $M$. Hence we get an a preferred differential $F$-integration
$\hat I^{F}_{!}$. The associated curvature map is
$$R(\hat I^{F}_{!})(\alpha)=\int_{W/M}\alpha\ .$$
}\end{ex}

\subsection{Higher $e$-invariants and index theorems}

\begin{ex}[The higher complex $e$-invariant]\label{jul2110}{\em 
In this very long example we provide an explicit calculation of integrals  in differential complex $K$-theory.
From the point of view of topology, what we are going to calculate, is  a higher version of Adam's
$e$-invariant.  But already the restriction of this higher $e$-invariant to a point is interesting. In this case, the calculation made here provides the solution to previous exercises \ref{aug2050}.

A closed $n$-dimensional manifold $M$ with a framing of a representative of the stable normal bundle is oriented for $\bS$. We use the procedure explained in \ref{jun0701} 
 in order to construct an induced differential  integration 
$ \hat I^{  \bKU}_{!}:\hat \bKU^{0}(M)\to \hat \bKU^{-n}(*)$. 
The closed framed $n$-manifold $M$  represents (via Thom-Pontrjagin) a stable homotopy class    $[M]\in \bS_{n}(*)\cong \pi_{n}^{s}$. 
Assume that $n$ is odd. 
Let $$e_{\C}:\pi_{n}^{s} \to \C/\Z$$
be the complex version of the $e$-invariant of Adams defined using the unit $\epsilon$  by 
$$\epsilon_{\C}:\pi_{n}^{s}=\pi_{n}(\bS)\stackrel{\sim}{\to} \pi_{n+1}(\bS\C/\Z)\stackrel{\epsilon}{\to} \pi_{n+1}(\bKU \C/\Z)\cong \C/\Z$$  
(cf. \ref{jul2111}).

\begin{prob}\label{jun0703}
Show that in $\C/\Z$ we have $e_{\C}([M])= \hat I^{  \bKU}_{!}(1)$.
\end{prob}
\proof
See \cite{MR2664467}. \hB

We consider the groups $SO(3)$ and $SU(2)$
as framed manifolds with respect to their left and right-invariant framings and write
$[SO(3)_{l}]$, $[SO(3)_{r}]$, $[SU(2)_{l}]$ and $[SU(2)_{r}]$ for the corresponding 
stable homotopy classes.

\begin{prob}
Calculate the $e_{\C}$-invariant of these stable homotopy  classes using Exercise \ref{jun0703}.
\end{prob}

More generally, let $f:W\to M$ be a proper submersion with a framing $fr$ of the vertical bundle $T^{v}f$.  
This framing determines an $\bS$-orientation of $f$. Note that $T^{v}f$ is a stable complement of the stable normal bundle, and   stable homotopy classes of framings of 
$T^{v}f$ correspond bijectively to stable homotopy classes of   framings of representatives of the stable normal bundle.
 We consider the class $$q(f,fr):=I^{\bS}_{!}(1)\in \bS^{-n}(M)\ ,$$ where $n:=\dim(M)-\dim(W)$. This stable cohomotopy class
is of course difficult to understand, in general.

The $\bS$-orientation  again can be lifted to a differential $\bS$-orientation
and the associated integration $\hat I^{\bS}_!$ does not depend on the choices.
It again induces a differential $\bKU$-integration and therefore a map
$\hat I^{\bKU}_{!}:\hbKU^{*}(W)\to \hbKU^{*+n}(M)$. 

\begin{prob}
Show that $\hat I^{\bKU}_{!}(1)$
 is flat. \end{prob}

\begin{ddd}
We define the higher $e_{\C}$-invariant of $(f,fr)$ by
$$e_{\C}(f,fr):=\hat I^{\bKU}_{!}(1)\in \bKU\C/\Z^{-n-1}(M)\ .$$
\end{ddd}
Below we will often simply write $e_{\C}(f)$ if the framing is clear from the context.

\begin{prob}
Show that $e_{\C}(f)$ only depends on the class $q(f,fr)$.
\end{prob}
As we shall see the higher $e_{\C}$-invariant of $f$ is much easier to calculate than $q(f,fr)$ and can be used to 
give some information on the latter.

We consider the special case where $f:W\to M$ is a $G$-principal bundle. The vertical bundle $T^{v}f$ has a canonical framing given by the fundamental vector fields of the $G$-action.  
\begin{prob}\label{jul3140}
Calculate the class
$$\hat I^{\bKU}_{!}(1)\in \bKU\C/\Z^{-n-1}(M)$$
 in terms of characteristic classes of the principal bundle.
 \end{prob}
 \proof

 Let $T$ be the maximal torus of $G$. We fix an embedding
$\kappa:S^1\to T$. Let $D^2\subset \C$ be the unit disc with the standard action of $S^1$.
Then $q:Z:=W\times_{S^1}D^2\to M$ is a fibre-wise zero-bordism. It fits into a diagram
$$\xymatrix{W\ar[ddr]_{f}\ar[dr]\ar[rr]&&Z\ar[ddl]^q\ar[dl]^s\\&W/S^1\ar[d]^r&\\&M&}\ .$$

The vertical tangent bundle of $q$ can be written in the form
$T^{v}q=s^*T^{v}r\oplus  s^* V_1$, where $V_1:=W\times_{S^1}\C$.
Let $\taaa$ denote the Lie algebra of $T$ and  $\Delta^+\subset \taaa^*_\C$ be  a system of positive roots of the pair $(\gaaa,\taaa)$. Then we have a $T$-equivariant decomposition
$$\gaaa\cong \taaa\oplus\bigoplus_{\alpha\in \Delta^+} \gaaa_\alpha\ ,$$
where $\gaaa_\alpha$ is a complex one-dimensional  representation of $T$ with weight $\alpha$.
Let $\bar \taaa:=\taaa/\kappa_*(\R)$. For a weight $\beta$ of $S^1$  we consider the complex line bundle $V_\beta:=W\times_{S^1,\beta}\C\to W/S^1$.  For a weight $\alpha\in \taaa^*_\C$ of $T$ we get a weight $\kappa^*\alpha$ of $S^1$. 
Then
$T^{v}r\cong (W/S^1\times \bar \taaa)\oplus \bigoplus_{\alpha\in \Delta^+} V_{\kappa^*\alpha}$. 
We thus have 
\begin{equation}\label{eq1}
T^{v}q=(Z\times \bar \taaa)\oplus s^*(V_1\oplus \bigoplus_{\alpha\in \Delta^+} V_{\kappa^*\alpha})\ .
\end{equation}
The first summand is framed, and the second summand has a natural complex structure.
Therefore, $Tq$ has a natural stable complex structure and hence a $Spin^c$-structure.

Let $\nabla^{T^{v}q^s}$ be a stable complex connection on $T^{v}q^s$ which respects the sum-decomposition above and  which restricts to the trivial connection $\nabla^{triv}$ on $W=\partial Z$ determined by the framing.

The bordism formula \ref{jun0839} predicts that
$$e_{\C}(f)=a(\int_{Z/M}  \Td(\hat I(\hat {\tilde \nu})))\ .$$
Here $ \hat {\tilde \nu}$ is a differential $\bKU$-Thom class of the normal bundle
of $q$ which extends the differential $\bS$-Thom class on the boundary. 

We now use the following general fact. Let $\nabla$ be a complex connection on a complex representative of the
stable normal bundle $\tilde N$ of $q$ which is the trivial connection associated to the framing over $W$.
Then by \ref{aug0335} we can find a Thom class
$\hat{\tilde {\nu}}_{\bMU}\in \hbMU^{\dim(\tilde N)}$
with $\Td(\hat{\tilde  {\nu}}_{\bMU})=u(\nabla)$.
We let $\hat{\tilde {\nu}}$ be the image of $\hat{\tilde {\nu}}_{\bMU}$ under the orientation
$\bMU\to \bKU$. Note that with this choice  $\Td(\hat{\tilde  {\nu}})=\Td(\nabla)^{-1}=\Td(\hat I(\hat{\tilde{\nu}}))$, where  $\Td(\nabla)$ is the form representing the characteristic class of $Spin^{c}$-vector bundles considered in \ref{jul3110}. 
We have fixed an isomorphism of complex vector bundles $T^{v}q^{s}\oplus \tilde N\cong Z\times \C^{l}$ which extends the isomorphism given by the framings over $W$.
We therefore have a transgression formula
$$\Td(\nabla^{T^{v}q^{s}})\wedge\Td(\nabla)=1+d\widetilde{\Td}(\nabla^{T^{v}q^{s}}\oplus\nabla, \nabla^{triv})  \ ,$$
where $\tilde \Td(\nabla^{T^{v}q^{s}}\oplus\nabla, \nabla^{triv}) $ vanishes on $W$.
This gives
$$\Td(\nabla)^{-1}=\frac{\Td(\nabla^{T^{v}q^{s}}) }{1+d\widetilde{\Td}(\nabla^{T^{v}q^{s}}\oplus\nabla, \nabla^{triv})  }\ .$$
It follows that
$$\Rham(\Td(\nabla)^{-1})=\Rham(\Td(\nabla^{T^{v}q^{s}}))\in HA^{0}(Z,W)\ .$$ 
In particular we get
$$\Rham(\int_{Z/M}\Td(\hat I(\hat{\tilde{\nu}})))=\Rham(\int_{Z/M}\Td(\nabla^{T^{v}q^{s}}))\ .$$

 Therefore
$$ e_{\C}(f)=a(\int_{Z/M} \Td(\nabla^{T^{v}q^s}))\ .
$$

Let $x:=bc_1 $, where $c_{1}\in H^2(W/S^1;\C)$ be the first real Chern class of the $S^1$-bundle $W\to W/S^1$.
\begin{lem}
We have
\begin{equation}\label{udiqwdqwdwqd}\Rham(\int_{Z/(W/S^1)} \Td(\nabla^{T^{v}q^s}))=\frac{1}{x}\left(1-\frac{x}{ e^{ x}-1}\prod_{\alpha\in \Delta^+,\kappa^*\alpha\not=0}  \frac{\kappa^*\alpha\:\: x}{ e^{ \kappa^*\alpha \:\:x}-1}\right  )\ .\end{equation}
\end{lem}
 \proof
 Let $T^{v}q^{s}$ denote some stabilization of $T^{v}q$.
According to the decomposition (\ref{eq1}) we choose a complex connection
$\tilde \nabla^{Tq^s}$ as follows. 
First we let $$T^{v}q^s:=(Z\times (\bar \taaa\oplus \bar \taaa))\oplus 
 s^*(V_1\oplus \bigoplus_{\alpha\in \Delta^+} V_{\kappa^*\alpha})\ ,$$
where the first summand has the complex structure
$i(t,t^\prime)=(-t^\prime,t)$. We fix a connection $\nabla$ on $V_1$. It induces connections on $V_\beta$ for all weights $\beta$ of $S^{1}$.
We choose a connection $\tilde \nabla_1$ on $q^*V_1$ which is trivialized near $W$ and coincides with $q^{*}\nabla$ near the zero section. This fixes a choice of $\tilde \nabla_\beta$ for every weight $\beta$ since $q^*V_\beta=q^*V_1^{\otimes n_\beta}$.
In particular, with this choice we have
$c_1(\tilde \nabla_\beta)=\beta c_1(\tilde \nabla_1)$.
We let $\nabla^{Tq^s}$ be the  sum of the connections
$\tilde \nabla_{\kappa^*\alpha}$ and the trivial connection on $Z\times (\bar \taaa\oplus \bar \taaa)$.
Then we have
$$\Td(\nabla^{T^{v}q^s})=\Td(\tilde \nabla_1)\prod_{\alpha\in \Delta^+,\kappa^*\alpha\not=0} \Td(\tilde \nabla_{\kappa^*\alpha})\ .$$ 
Using the generating power series 
$$\Td(x)=\frac{x}{e^{x}-1}$$
(this has been calculated in Problem \ref{jul3110})
we get
$$\Td(\nabla^{T^{v}q^s})= \frac{bc_1(\tilde \nabla_1)}{e^{bc_1(\tilde \nabla_1) }-1}\prod_{\alpha\in \Delta^+,\kappa^*\alpha\not=0}  \frac{\kappa^*\alpha\:\: bc_1(\tilde \nabla_1)}{ e^{ \kappa^*\alpha \:\:bc_1(\tilde \nabla_1) }-1}$$
We get, using the calculation \ref{jul3130},
$$\int_{Z/(W/S^1)} \Td(\nabla^{T^{v}q^s})\equiv -\frac{1}{bc_1(  \nabla)}\left( \frac{bc_1(  \nabla)}{ e^{ bc_1( \nabla)}-1}\prod_{\alpha\in \Delta^+,\kappa^*\alpha\not=0}  \frac{\kappa^*\alpha\:\: bc_1(  \nabla)}{ e^{-\kappa^*\alpha \:\:bc_1( \nabla)}-1}-1\right  )$$
modulo exact forms.
 The Lemma follows from the fact that $x$ is represented by the closed form 
$c_1(  \nabla)$. \hB

\subsubsection*{The higher rank case}

 \begin{lem}
If $\rk(G)\ge 2$, then $e_{\C}(f)=0$.
\end{lem}
\proof 
It suffices to calculate in the universal example where $f:W\to M$ is the universal bundle $\pi:EG\to BG$.
We first observe that $EG/S^1\cong BS^1$ and therefore $H^*(EG/S^1;\C)\cong \C[x]$ with $x\in H^2(EG/S^1;\C)$.
If $\rk(G)\ge 2$, then we can choose a decomposition $T\cong \kappa(S^1)\times T_1\times T^\prime$,
where $T^\prime$ is a $\rk(G)-2$-dimensional torus und $T_1\cong S^1$.
We consider the iterated bundle
$$EG/S^1 \to EG/(S^1\times T_1)  \to BG\ .$$

The right-hand side of (\ref{udiqwdqwdwqd}) is a certain formal power series $P(x)\in \C[[x]]$.
We have
$$e_{\C}(\pi)=\int_{(EG/S^1)/BG} P(x)=\int_{(EG/(S^1\times T_1))/BG}\int_{(EG/S^1)/(EG/(S^1\times T_1))} P(x)\  .$$
Note that $EG/(S^1\times T_1)\cong B(S^1\times T_1)$ so that the cohomology of $EG/(S^1\times T_1)$ is concentrated in even degrees. But $P(x)$ is of even degree and the fibre of $EG/S^1 \to EG/(S^1\times T_1) $ is one-dimensional so that 
$\int_{(EG/S^1)/(EG/(S^1\times T_1))} P(x)$ is a class of odd-degree and hence trivial.
\hB

We now assume that $\rk(G)=1$. The non-trivial examples are
$$G\cong S^1\ ,\quad G\cong SU(2)\ ,\quad G\cong SO(3)\ .$$

\subsubsection*{The case $G=S^1$}

We consider the universal case $f:ES^{1}\to BS^{1}$. In this case
$EG/S^1\cong BS^{1}$ so that
Lemma \ref{udiqwdqwdwqd} immediately gives $$\Rham(\int_{Z/BS^1}\Td(\nabla^{T^{v}q^s}))=\frac{1}{x}\left(1- \frac{x}{ e^{ x}-1} \right)\ .$$
Therefore, with $x=b c_{1}$, 
$$e_{\C}(f)=b\left[\frac{1}{x}\left(1-\frac{x}{e^{x}-1}\right)\right]\in H^{-2}(BS^1;\C[b,b^{-1}])/\im(\ch)\ .$$
Here are the first terms of this power series
$$\frac{1}{2}-\frac{1}{12}x+\frac{1}{720}x^3-\frac{1}{30240} x^5+\frac{1}{1209600} x^7-\frac{1}{47900160}x^9+\frac{691}{1307674368000}x^{11}+ \dots\ .$$
 
Already the first two terms are interesting. 

The constant term $\frac{1}{2}$ is the $e_{\C}$-innvariant of $S^1$ as a framed manifold (a bundle over a point). Note that $\pi^s_1=\Z/2\Z$, and  the generator $\eta$ is represented by the framed manifold $S^1$ via the Pontrjagin Thom construction.  This gives $\ord(e_{\C}(\eta))=2$.

We now consider the linear term. To this end we discuss the pull-back of $e_{\C}(f)$ along the canonical generator $S^{2}\to BS^{1}$ of $\pi_{2}(BS^{1})\cong \Z$.
The restriction of the universal bundle along the generator  $S^2\to BS^1$ of $\pi_2(BS^1)\cong \Z$  gives the Hopf bundle $\pi:S^3\to S^2$. Note that on $S^{2}$ we know that $\im(\ch)$ is the image of the integral cohomology in complex cohomology. Therefore we have
$$\bKU\C/\Z^{-2}(S^2)\cong H^{-2}(S^2,\C[b,b^{-1}])/\im(\ch)\cong (\C[x]/(x^2))/(\Z^2[x]/(x^2))\ .$$ The higher $e_\C$-invariant of the Hopf bundle  is
$$e_{\C}(\pi)=[\frac{1}{2}-\frac{x}{12}]\in (\C[x]/(x^2))/(\Z^2[x]/(x^2))\ .$$
This class has order $12$. On the other hand, the Hopf map $\pi:S^3\to S^2$ represents the  stable cohomotopy class $$\kappa\in S^{-1}(S^2)\cong [S^2,S^{-1}]\cong \pi_3^s\cong \Z/24\Z\ .$$ 
It follows that $e_{\C}(\kappa)$ has order $12$.
 
It remains to interpret the higher terms. See the corresponding discussion in the case of $G=SU(2)$.

%
 \subsubsection*{ The case $G=SU(2)$}

In this case we have one root $\alpha$ which is twice the fundamental weight of the maximal torus $S^1\subset SU(2)$. The reason is that the adjoint representation
factors over the two-fold covering
$SU(2)\to SO(3)$, and the adjoint representation of $SO(3)$ is the standard representation on $\R^3$ on which the maximal torus of $SO(3)$ acts by the fundamental representation.
Therefore we have
$$\Rham(\int_{Z/(ESU(2)/S^1)}\Td(\nabla^{T^{v}q^s}))=\frac{1}{x}\left(1-\frac{x}{e^{x}-1}\frac{2x}{e^{2x}-1}\right)\ .$$ We must integrate this further over the bundle
$ESU(2)/S^1\to BSU(2)$ with fibre $SU(2)/S^1\cong S^2$.
In general, of $G$ is a compact Lie group with maximal torus $T\subseteq G$, then the pull-back along $BT\to BG$ indentifies
$H^*(bG;\C)$ with the subring $H^*(BT;\C)^W\subseteq H^*(BT;\C)$, where
$W$ is the Weyl group of the pair $(G,T)$. In the present case
the Weyl group is $\Z/2\Z$ which acts by $x\to -x$ on
$H^*(BT;\C)\cong \C[[x]]$.  The Gysin sequence of the $S^2$-bundle $r:ESU(2)/S^1\to BSU(2)$ in integral cohomology using the absence of odd cohomology
gives 
$$0\to \Z[[x^{2}]]\stackrel{r^*}{\to} \Z[[x]]\stackrel{\int_r}{\to} \Z[[x^{2}]]\to 0\ ,$$
where
$$\Z[[x]]\cong H^*(BS^1;\Z)\ ,\quad  \Z[[x^2]]\cong H^*(BSU(2);\Z)\ .$$
It follows that $\int_r x^{2n+1}=x^{2n}$  (immediate at least up to sign). The same fomula now holds true in complex cohomology.
We conclude that
$$e_{\C}(\pi)=b^{2}\left[\frac{1}{x^2}\left(1-\frac{x}{ e^{x}-1}\frac{2x}{ e^{2x}-1}\right)\right]_{even}
\in H^{-4}(BSU(2);\C[b,b^{-1}])/\im(\ch)\ .$$
The first few terms of this series are
$$-\frac{11}{12}- \frac{1}{240} x^2 +\frac{1}{6048} x^4 -\frac{1}{172800} x^6 +\frac{1}{5322240} x^8 -\frac{691}{118879488000} x^{10} \dots\ .$$
The constant term $-\frac{11}{12}$ is the $e_{\C}$-invariant of the class in $\pi_3^s\cong \Z/24\Z$ represented by the framed manifold $SU(2)$ under the Pontrjagin Thom construction.

In order to discuss the quadratic term we consider the restriction of $e_{\C}(f)$ along the canonical generator $S^{4}\to BSU(2)$ of $\pi_{4}(BSU(2))\cong \Z$.
Note that $\pi_4(BSU(2))\cong \Z$. The generator classifies a $SU(2)$-principal bundle $E\to S^4$ whose second Chern class
$c_2(E)\in H^4(S^4;\Z)$ is a generator. 
The image of the Chern character in $H^{-4}(S^4;\C[b,b^{-1}])$ is the subgroup
$H^{-4}(S^4;\Z[b,b^{-1}])$. Therefore the class $e_{\C}(q)=[-\frac{11}{12}-\frac{1}{240}\ori_{S^4}]$ has order $240$. We can conclude  that $q$ represents a generator of  $S^{-3}(S^4)\cong \pi^s_{7}\cong \Z/240\Z$.

  In order to interpret the higher terms we must understand the image of $\ch$ in greater detail.
  It is clear that
 \begin{equation}\label{jul3141}\bKU\C/\Z^{-4}(BSU(2))\cong \prod_{n\ge 0} \C/\Z\end{equation}
 but it seems to be not so easy to make this isomorphism explicit. 
 One checks that the Chern character of the standard representation of dimension $n+1$ is given by
 $$\ch(V_{n})=\frac{\sinh((n+\frac{1}{2})x)}{\sinh(\frac{1}{2}x)}$$
 for $n=0,1,\dots$ where $b^{2}c_{2}=-x^{2}$.
It remains to use this calculation to produce the isomorphism \eqref{jul3141} explicitly 
and to derive the corresponding basis decomposition of $e_{\C}(f)$.

\subsubsection*{\textbf{The case $SO(3)$}}

In this case we again have one root $\alpha$ which is equal to the fundamental weight of the maximal torus $S^1\subset SO(3)$ since the   adjoint representation of $SO(3)$ is the standard representation on $\R^3$ on which the maximal torus of $SO(3)$ acts by the fundamental representation.
Therefore we have
$$\Rham(\int_{Z/(ESO(3)/S^1)}\Td(\nabla^{T^{v}q^s}))=\frac{1}{x}\left(1-\left(\frac{x}{e^{x}-1} \right)^2\right)\ .$$ We must integrate this further over the bundle
$ESO(3)/S^1\to BSO(3)$ with fibre $SO(3)/S^1\cong S^2$.
As in the case of $SU(2)$ we conclude that  
$$e_{\C}(f)=b^{2}\left[ \frac{1}{x^2}\left( 1-\left(\frac{x}{e^{x}-1}\right)^2\right)\right]_{even}
\in H^{-4}(BSO(3);\C[b,b^{-1}])/\im(\ch)\ .$$
The first few terms in this expansion are
$$-\frac{5}{12}- \frac{1}{240} x^2 +\frac{1}{6048} x^4 -\frac{1}{172800} x^6 +\frac{1}{5322240} x^8 -\frac{691}{118879488000} x^{10} \dots\ .$$
The constant $-\frac{5}{12}$ is again the $e_{\C}$-invariant of the element in $\pi_3^s$ represented by the framed manifold
$SO(3)$ under the Thom Pontrjagin construction. The higher terms coincide with those of $SU(2)$, a fact which can also be seen directly.

 \hB

}
\end{ex}


\begin{ex}{\em 
We consider a proper submersion $f:W\to M$ and set $n:=\dim(M)-\dim(W)$. The choice of an embedding
$\iota$ determines a representative $\cN_{f}$ of the stable normal bundle for $f$. 
Let us fix a complex structure $I_{f}$ and a connection $\nabla^{f}$ on its underlying bundle $N_{f}$.

The zero section $0_{N_{f}}:W\to  N_{f}$ has a canonical representative of the normal bundle $\cN_{0_{N_{f}}}$ with underlying
bundle $N_{f}$. The datum $(0_{N_{f}},\cN(0_{N_{f}}),I_{f})$ is a topological $\bMU$-cycle representing a $\bMU$-Thom class $\nu$ on $N_{f}$ (cf. Example \ref{jul2610}).  By the solution of Problem \ref{aug0720} there is a natural choice $\hat \nu:=\hat \nu(\nabla^{f})$ of a differential Thom class on $N_{f}$ refining $\nu$ with $\Td(\hat \nu)=u(\nabla^{f})$. 
We can write 
$$\hat \nu:=\hcycl(0_{N_{f}},\cN_{0_{N_{f}}},I_{f},\nabla^{f},\eta_{f})$$
for a suitable choice of $\eta_{f}\in \Omega A^{-1}_{prop/W,-\infty}(N_{f})$.
The pair 
$(\iota,\hat \nu)$ represents a differential $\bMU$-orientation for $f$ and hence yields an integration 
$\hat I(\iota,\hat \nu)_{!}$. 

We now define a geometric integration
$$\hat I^{geom}:\Cycle_{\bMU}^{geom,*}(W)\to \Cycle^{geom,*+n}_{\bMU}(M)$$ as follows.
Let $$\tilde g:=(g:B\to W,\cN_{g},I_{\cN_{g}},\nabla^{\cN_{g}},\eta_{g})\in \Cycle^{*}_{\bMU}(W)\ .$$
Then we get a canonical representative of the normal bundle $\cN_{f\circ g}$ whose underlying bundle admits an isomorphism \begin{equation}\label{jul2210}N_{f\circ g}\cong N_{g}\oplus g^{*}N_{f}\end{equation} which is unique up to homotopy. Indeed, the underlying bundle $N_{f\circ g}$ admits a filtration $N_f\subset N_{f\circ g}$ with quotient $N_g$, and the isomorphism is fixed by the choice of a split.
We choose such an isomorphism
and  define $I_{f\circ g}$ and 
$\nabla^{f\circ g}$ using this direct sum decomposition
and set
$$\eta_{f\circ g}:= \int_{W/M}u(\nabla^{f}) \wedge\eta_{g} + \int_{W/M}(\int_{N_{f}/W} \eta_{f}) \wedge R(\hcycl(\tilde g)) \ .$$
Then we define
$$\hat I^{geom}(\tilde g)=(f\circ g,\cN_{f\circ g},I^{g}\oplus g^{*} I^{f},\nabla^{g}\oplus g^{*}\nabla^{f},\eta_{f\circ g})\ .$$
Observe that
$\hcycl(\hat I^{geom}(\tilde g))$ does not depend the choice of the isomorphism \eqref{jul2210}.

The following is the prototype of a differential index theorem. It compares the integration with respect to a differential orientation with a geometric push-forward construction.
\begin{prob}\label{aug0751}
Show the index theorem:
$$\xymatrix{\Cycle_{\bMU}^{geom,*}(W)\ar[d]^{\hat I^{geom}}\ar[r]^{\hcycl}&\hbMU^{*}(W)\ar[d]^{\hat I(\iota,\hat \nu)_{!}}\\\Cycle_{\bMU}^{geom,*+n}(M)\ar[r]^{\hcycl}&\hbMU^{*+n}(M)}\ .
$$
\end{prob}
\proof
\textcolor{red}{I am interested to see the details.}
\hB

}\end{ex}

 \begin{ex}\label{aug0701}
 {\em
We consider a proper submersion $f:W\to M$ with a fiberwise Riemmannian metric   and a horizontal distribution. We get an induced connection $\nabla^{T^{v}f}$ on the vertical bundle $T^{v}f:=\ker(df)$ by  \cite[Ch. 9]{MR1215720}.
We assume that $T^{v}f$ has a $Spin^{c}$-structure and  that we have chosen a $Spin^{c}$-extension $\tilde \nabla$ of 
$\nabla^{T^{v}f}$. Let $S(T^{v}f)$ be the fibrewise spinor bundle. It acquires the structure of a family of Dirac bundles. In this way we get a geometric family $\cW$ over $M$. 
Let $V\to W$ be a hermitean vector bundle with metric connection $\nabla^{V}$. The geometric family determined by the twisted spinor bundle $S(T^{v}f)\otimes V$ will be denoted by $\cW\otimes V$.
Its associated family of Dirac operators will be denoted by  $D(\cW\otimes V)$. We refer to \cite{MR2191484e} for details on geometric families.

The $Spin^{c}$-structure on $T^{v}f$ determines an equivalence class of a $\bKU$-orientation $[\iota,\nu]$ of $f$. For a complex vector bundle with connection $(V,\nabla^{V})$ on $W$
we define the topological index by
$$\ind^{top}(V,\nabla^{V}):=I[\iota,\nu]_{!}(\cycl(V))\in \bKU^{n}(M)\ ,$$
where $n=\dim(M)-\dim(W)$.
The analytical index is the index of the family of twisted Dirac operators 
$$\ind^{an}(V,\nabla^{V}):=\ind(D(\cW\otimes V))$$
The Atiyah-Singer index theorem for families
states that
$$\ind^{top}(V,\nabla^{V})=\ind^{an}(V,\nabla^{V})\ .$$
Note that both indices do not depend on the connections.
The goal of the following is to state a differential refinement of this index theorem which takes the geometric structures into account.

The $Spin^{c}$-connection $\tilde \nabla$ determines a form
$$\Td(\tilde \nabla)\in \Omega \C[b,b^{-1}]_{cl}^{0}(W)^{*}\ .$$
 By \ref{aug0721} we can find a natural Thom class $\hat \nu(T^{v}f):=\hat \nu(\tilde \nabla)$ on
 $T^{v}f$ which is unique up to homotopy and satisfies $\Td(\hat \nu(T^{v}f))=\Td(\tilde \nabla)^{-1}$.
 We now choose an embedding $\iota$ and get an associated representative of the stable normal bundle $\cN$ with underlying bundle $N$. In particular we have a trivialization of $T^{v}f\oplus N\cong W\times \R^{l}$.

Let  $k:=l+n
$ be the dimension of $N$.
Using the two-out-of-three principle \ref{aug0670} and the canonical Thom class $\hat \nu_{l}$ on the trivial bundle we get a well-defined homotopy class of differential Thom classes
 $\hat \nu\in \hbKU^{k}_{prop/W}(N)$ such that
 $\hat \nu_{l}=\pr_{T^{v}f}\hat \nu(T^{v}f)\cup \pr_{N}^{*} \hat \nu$.
 In particular we have $\Td(\hat \nu)=\Td(\tilde \nabla)$.
The pair $(\iota,\hat \nu)$ induces a well-defined integration
$$\hat I_{!}:=\hat I(\iota,\hat \nu)_{!}:\hbKU^{*}(W)\to \hbKU^{*+n}(M)$$
with curvature map $$R(\iota, \hat \nu)(\alpha)=\int_{W/M} \Td(\tilde \nabla)\wedge \alpha\ .$$ 
 This differential integration only depends on the $Spin^{c}$-structure on $T^{v}f$ and the
choice of the $Spin^{c}$-connection $\tilde \nabla$.
We define the differentially refined topological index of a vector bundle with connection $(V,\nabla^{V})$ on $W$ by 
\begin{equation}\label{sep2205}\hind^{top}(V,\nabla^{V}):=\hat I_{!}(\hcycl(V,\nabla^{V}))\in \hbKU^{n}(M)\ .\end{equation}
We assume that $n$ is even and that the family of twisted Dirac operators
$\cD(\cW\otimes V)$ has a kernel bundle $U$. This bundle acquires an induced metric and a metric connection $\nabla^{U}$ (see \cite[Ch. 10]{MR1215720} for details).
We consider $$\hind^{an}(V,\nabla^{V}):=b^{-\frac{n}{2}}\hcycl(U,\nabla^{U})\in \hbKU^{n}(M)$$
as the differentially refined analytical index. Using the methods of \cite{MR2664467} one can drop the assumptions that $n$ is even and  that there is a kernel bundle so that the analytic index is defined in general.

The comparison between  the differentially refined topological and the analytical indices is the contents of  the differential version of the Atiyah-Singer index theorem.
The local index \cite{MR1215720} theorem asserts that
$$\ch(\nabla^{U})=\int_{W/M} \Td(\tilde \nabla)\wedge \ch(\nabla^{V})+d\eta(\cW\otimes V)\ ,$$
where $\eta(\cW\otimes V)\in \Omega \C[b,b^{-1}]^{n-1}(M)$  is the $\eta$-form of \cite{MR966608}. 
\begin{prob}\label{aug0750}
Show the differential index theorem:  In  $\hbKU^{n}(M)$ we have
\begin{equation}\label{aug2010}\hind^{top}(V,\nabla^{V})=\hind^{an}(V,\nabla^{V})+a(\eta(\cW\otimes V)) \ .\end{equation}
\end{prob}
\proof
The differential index theorem follows from the following proposition. 
\begin{prop}\label{sep2201}
Let $\hat I_{!},\hat I_{!}^{\prime}$ be two integrations in differential $K$-theory defined for differentially $K$-oriented (in the sense of \cite{MR2664467}) proper submersions which satisfy:
\begin{enumerate}
\item $\hat I_{!}$ and $\hat I_{!}^{\prime}$  are compatible with curvature.
\item They are compatible with cartesian diagrams.
\item They are compatible with composition.
\item They satisfy the bordism formula.
\item They coincide with the topological integration $I_{!}$ on flat classes.
\end{enumerate}
Then $\hat I_!=\hat I_!^{\prime}$.
\end{prop}
Indeed, the integration defined in \cite{MR2664467} coincides with the right-hand side of \eqref{aug2010}
and satisfies the assumptions of Proposition \ref{sep2201}.

The integration in differential $K$-theory used in the definition of the
topological index \eqref{sep2205} also satisfies these assumptions by \ref{jun1401} and \ref{jun0839}. \hB

\proof(Prop. \ref{sep2201})
We consider the difference $\delta:=\hat I_{!}-\hat I_{!}^{\prime}$, We must show that $\delta=0$.
The first lemma reduces the proof of Proposition \ref{sep2201} to the case where the base of the bundle is a point.
\begin{lem}
If $\delta=0$ in the case where $M$ is a point, then $\delta=0$ in general.
\end{lem}
\proof
Let $f:W\to M$ be a proper submersion with a differential $K$-orientation.
We let $\hat I_{f,!}$ and $\hat I_{f,!}^{\prime}$ denote the associated integrations
and write
$\delta_{f}:=\hat I_{f,!}-\hat I_{f,!}^{\prime}$.
We consider a class
 $x\in \hbKU^{*}(W)$. 
The assertion that the integrations are compatible with the curvature in detail means that
$$R(\hat I_{!,f}(x))=\int_{W/M} \Td(\tilde \nabla) \wedge  R(x)=R(\hat I_{!,f}^{\prime}(x))\ .$$
 This implies that
$R(\delta_{f}(x))=0$ so that $\delta_{f}(x)\in \hbKU^{*+n}_{flat}(M)\cong \bKU \C/\Z^{*+n-1}(M)$.
In particular, $\delta_f
(x)$ is homotopy invariant and therefore only depends on the underlying topological class $I(x)$ of $x$.
   
As a consequence of the Anderson selfduality of $\bKU$ one can detect classes in $\hbKU\C/\Z^{*}(M)$ by evaluating them against $\bKU_{*}(M)$ classes.
 So if  $\langle \delta_{f}(x), u\rangle=0$ for all $\bKU$-homology classes $u\in \hbKU_{*+n-1}(M)$, then
$\delta_{f}(x)=0$. By \cite{MR679698} or \cite{MR1624352} we can realize $u$ geometrically in the form $u=g_{*}(v\cap [B])$, where
$g:B\to M$ is a smooth map from a closed manifold $B$ of dimension $*+n-1$ with a $Spin^{c}$-structure and $v\in \bKU^{0}(B)$.  The projection $p:B\to *$ is then  $\bKU$-oriented
 and we have
 $$\langle \delta_{f}(x), u\rangle=I_{p,!}(g^{*}\delta_{f}(x)\cup v)\in \bKU\C/\Z^{0}(*)\cong \C/\Z\ .$$
 
We choose a differential refinement of the $\bKU$-orientation of $p$ and a lift $\hat v\in \hbKU^{0}(B)$ of $v$.
We let $\tilde f$ and $\tilde g$ be defined by the cartesian diagram
$$\xymatrix{\tilde W\ar[d]^{\tilde f}\ar[r]^{\tilde g}&W\ar[d]^{f}\\
B\ar[r]^{g}&M}$$
where $\tilde f$ has its induced differential $\bKU$-orientation. 
Then we have the sequence of equalities
\begin{eqnarray*}
\langle \delta_{f}(x), u\rangle&=&I_{p,!}(g^{*}\delta_{f}(x)\cup v)\\
&\stackrel{!!}{=}&\hat I_{p,!}( g^{*}\delta_{f}(x) \cup \hat v)\\
&\stackrel{!}{=}&\hat I_{p,!}(g^{*} \hat I_{f,!}(x)\cup \hat v)-\hat I^{\prime}_{p,!}(g^{*} \hat I_{f,!}^{\prime}(x)\cup \hat v)\\&\stackrel{!!!}{=}&\hat I_{p\circ \tilde f,!}(\tilde g^{*}x\cup \tilde f^{*}\hat v)-\hat I^{\prime}_{p\circ \tilde f,!}(\tilde g^{*}x\cup \tilde f^{*}\hat v)\\
&\stackrel{!}{=}&0\ .
\end{eqnarray*}
At the equalities marked by $!$ we use the assumption that $\delta=0$ for maps to a point. The equality marked by $!!$ uses the assumption \ref{sep2201}, 5. At the equality marked by $!!!$ we use
\ref{sep2201}, 2. und 3.  \hB

In order to finish the proof of Proposition \ref{sep2201} we consider the special case where $M=*$.
Since $\bKU\C/\Z^{odd}(*)=0$ and by periodicity of $\bKU$ it remains to consider the case where $\dim(W)=-n$ and $\deg(x)=1-n$.  We now reduce to the case of odd $n$.

\begin{lem}If we
assume that $\delta=0$ in the case that $M$ is a point and  $n$ is odd, then
$\delta=0$ in general.
\end{lem}
\proof Assume that $n$ is even. We consider the map $p:S^{1}\to *$ with some differential $K$-orientation and a class $e\in \hbKU^{1}(S^{1})$ such that $\hat I_{p,!}(e)=1$.
 Since $\dim(S^{1})$ is odd  we have by assumption that $\hat I_{p,!}^{\prime}(e)=1$.
We now consider the diagram
$$\xymatrix{W\times S^{1}\ar[r]^{f\times \id_{S^{1}}}\ar[d]^{\id_{W}\times p}&S^{1}\ar[d]^{p}\\
W\ar[d]^{f}\ar[r]^{f}&\mbox{$*$}\\
\mbox{$*$}&}\ .$$
By a simple calculation again using the compatibility of the integrations with pull-back, composition and
the projection formula we get
$$
\delta_{f}(x)=\delta_{f\circ (\id_{W}\times p)}(x\cup (f\times \id_{S^{1}})^{*}e)\ .$$
The right-hand side vanishes since $\dim(W\times S^{1})$ is odd. \hB 

We have  not yet used that both integrations satisfy the bordism formula.
Let $Z$ with $W\cong \partial Z$ be a $Spin^{c}$-zero-bordim of $W$. 
 We extend the differential $\bKU$-orientation of $W$ across $Z$. Note that this is possible
 since we use the notion of a differential orientation defined in \cite{MR2664467}.
  Furthermore we assume that $x=y_{|W}$ for some class $y\in \hbKU^{1-n}(Z)$.
The assumption \ref{sep2201}, 4. in detail means that
$$\hat I_{f,!}(x)=a(\int_{Z} \Td(\tilde \nabla)\wedge R(y))=\hat I^{\prime}_{f,!}(x)\ .$$
It follows that $\delta_{f}(x)=0$ if $W$ is zero-bordant and the class $x$ extends over the zero bordism. 

By the Pontrjagin-Thom construction the group of bordism classes of $Spin^{c}$-manifolds of dimension $-n$ equipped with $\bKU^{1-n}$-classes is isomorphic to  $\mathbf{MSpin^{c}}_{-n}(\Omega^{\infty+n-1}\bKU)$.  The difference $\delta$ therefore induces a homomorphism
$\tilde \delta:\mathbf{MSpin^{c}}_{-n}(\Omega^{\infty+n-1}\bKU)\to \C/\Z$. In order to finish the proof of Proposition \ref{sep2201} it suffices to show that $\tilde \delta=0$ in the case of odd $n$.

 The Atiyah-Bott-Shapiro orientation 
$\mathbf{MSpin^{c}}\to \bKU$ induces   a homomorphism
$$ABS:  \mathbf{MSpin^{c}}_{-n}(\Omega^{\infty+n-1}\bKU)\to \bKU_{-n}(\Omega^{\infty+n-1}\bKU)\ .$$
 If $\tilde \delta$ factors over $ABS$, then it is induced by  a map
 $\hat \delta:\bKU_{-n}(\Omega^{\infty+n-1}\bKU)\to \C/\Z$, i.e. a class $\hat \delta\in \bKU\C/\Z^{-n}(\Omega^{*+n-1}\bKU)$.  If $n$ is odd, then
 $\bKU\C/\Z^{-n}(\Omega^{\infty+n-1} \bKU)\cong \bKU\C/\Z^{1}(\Z\times  BU)\cong 0$.
 Therefore the following Lemma completes the proof of Proposition \ref{sep2201}.\footnote{This part of the argument was suggested by Stephan Stolz.}

  \begin{lem}
  The homomorphism $\tilde \delta$ factors over $ABS$.
  \end{lem}
  \proof 
  We let $[W,x]\in  \mathbf{MSpin^{c}}_{-n}(\Omega^{\infty+n-1}\bKU)$
  denote the class represented by the closed $-n$-dimensional $Spin^{c}$-manifold
 $W$ and the class $x\in \bKU^{-n+1}(W)$.
 We have already seen that $\tilde \delta([W,x])$ only depends on the bordism class of $(W,x)$.
 Since $\delta$ is linear we have
 $\tilde \delta([W,x+x^{\prime}])=\delta([W,x])+\delta([W,x^{\prime}])$.
 Following \cite{MR679698}, in order to show that $\tilde \delta$ factors over $ABS$, it remains to show that $\tilde \delta$ preserves vector bundle modifications.
  
Let $g:S\to W$ be the sphere bundle of a real $Spin^{c}$-vector bundle of dimension $2k+1$.  The $Spin^{c}$-structure on $W$ induces a $Spin^{c}$-structure on $S$. The map $\tilde \delta$ preserves vector bundle modifications if in this case we have the equality  $\tilde \delta([S, b^{-k}g^{*}x])=\tilde \delta([W,x])$, where $b^{-k}\in \bKU^{2k}$ is a power of the Bott element.
We choose differential orientations of the maps and a differential lift $\hat x$ of $x$.
We assume:
\begin{lem}\label{sep2210}
We have $\delta_{g}(1)=0$ in the case that $g:S\to W$ is a sphere bundle of a $Spin^{c}$-vector bundle.
\end{lem}
Using this Lemma we can calculate
\begin{eqnarray*}
\tilde \delta([S,b^{-k}g^{*}x])&=&\hat I_{f\circ g,!}(b^{-k}g^{*}\hat x)-\hat I^{\prime}_{f\circ g,!}(b^{-k}g^{*}\hat x)\\
&=&(\hat I_{f,!}\circ \hat I_{g,!})(b^{-k}g^{*}\hat x)-(\hat I^{\prime}_{f,!} \circ \hat I^{\prime}_{g,!})(b^{-k}g^{*}\hat x)\\
&=&\hat I_{f,!}(\hat x\cup b^{-k} \hat I_{g,!}(1))-\hat I^{\prime}_{f,!} (\hat x\cup b^{-k}\hat I^{\prime}_{g,!}(1))\\
&=&\delta_{f}(x\cup b^{-k}\hat I_{g,!}(1))- \hat I^{\prime}_{f,!}(x\cup b^{-k}\delta_{g}(1))\\
&=&\delta_{f}(x)\ .
\end{eqnarray*}
The second summand of the last line vanishes by Lemma \ref{sep2210}.
In the first summand we use the identity of topological $\bKU$-classes $b^{-k}I_{g,!}(1)=1$ and that $\delta_{f}$ factors over the underlying topological $\bKU$-class.

It remains to show Lemma \ref{sep2210}. Since $\delta$ is compatible with pull-back
the class $\delta_{g}(1)$ is pulled by from a universal class in $\bKU\C/\Z^{-2k-1} (BSpin^{c}(2k+1))$.
Since $ BSpin^{c}(2k+1)$ is the classifying space of a compact Lie group its $\bKU$ and $\bKU\C/\Z$-theories are concentrated in even degrees by \cite{MR0139181}. 
\hB

 

%
\hB

A differentially refined Atiyah-Singer index theorem was first shown by Freed and Lott in 
\cite{MR2602854}. The analytic indices of Freed-Lott and in \ref{aug0750} essentially 
coincide. But the construction of the topological indices is different.
So \ref{aug0750} is not completely equivalent to \cite{MR2602854}.  
Differential $\bKU$-index theorems were also studied in \cite{MR2711943} and \cite{2009PhDT} (even equivariantly).

If one applies the differential lift of the Chern character $\hat \ch:\hbKU^{*}\to \widehat{H\Q[b,b^{-1}]}$ to  both sides of the index formula \eqref{aug2010}, then one obtains the equality
$$\hat \ch(\hind^{top}(V,\nabla^{V}))=\hat \ch(\hind^{an}(V,\nabla^{V}))+a(\eta(\cW\otimes V))$$
in $\widehat{H\Q[b,b^{-1}]}^{n}(M)$. The validity of this equality is a mixture if a Riemann-Roch and an index theorem which has been verified in \cite{MR2129894}, \cite{MR2664467}.

 }
 \end{ex}

\subsection{Geometrization}
Let $G$ be a Lie group and $R(G)^{+}\subset R(G)$ be the semiring of isomorphism classes of finite-dimensional representations under $\oplus$ and $\otimes$. Let $P\to M$ be a $G$-principal bundle and $\omega$ be a connection on $P$. If $(\rho,H_{\rho})$ is a finite-dimensional  representation of $G$, then we can form the associated vector bundle $P\times_{G,\rho}H_{\rho}=:V_{\rho}$. It acquires a connection $\nabla^{V_{\rho}}$ from $\omega$.
Below we write $V_{\rho}=:B_{p}(\rho,H_{\rho})$ and $(V_{\rho},\nabla^{V_{\rho}})=:B_{P,\omega}(\rho,H_{\rho})$.
The associated vector bundle   construction gives  a diagram of  morphisms of semirings
$$\xymatrix{& \Vect^{geom}(M)\ar[dd]\\R(G)^{+}\ar[ur]^{B_{P,\omega}}\ar[dr]^{B_{P}}&\\ &\Vect(M)}\ .$$

The construction $B_{P}$ can be applied to the universal $G$-bundle $EG\to BG$.
Since $\bKU^{0}(BG)$ is a group we get an extension
$$\xymatrix{R(G)^{+}\ar[d]\ar[r]^{B_{EG}}&\Vect(BG)\ar[d]^{\cycl}\\
R(G)\ar[r]^{^{\tilde B_{EG}}}&\bKU^{0}(BG)}$$
of $B_{EG}$ to a map $\tilde B_{EG}$.

Let $E^{*}$ be a generalized cohomology theory. If $X$ is a space, then $E^{*}(X)$
has a natural topology such that the open neighbourhoods of $0$ are the kernels of restrictions
$E^{*}(X)\to E^{*}(Y)$ for all maps $Y\to X$ from finite $CW$-complexes $Y$. 
We apply this to $\bKU^{0}(BG)$.

Let $I:=\ker(\dim R(G)\to \Z)$ be the augmentation ideal and
$R(G)^{\hat{}}_{I}$ be the $I$-adic completion of the ring $R(G)$.
For the following see e.g. \cite{MR0259946}.
\begin{theorem}[Atiyah-Hirzebruch]
Assume that $G$ is compact. Then the map $\tilde B_{EG}$ induces a topological isomorphism
$$\hat B_{BG}:R(G)^{\hat{}}_{I}\stackrel{\sim}{\to} \bKU^{0}(BG)\ .$$
\end{theorem}

\begin{ex}\label{jun0741}{\em
Let $G=U(1)$.
We have $BU(1)\cong \C\P^{\infty}$ and
$\bKU^{*}( \C\P^{\infty})\cong \Z[[x]][b,b^{-1}]$,
where $x\in \bKU^{0}(\C\P^{\infty})$ is the class $x=\cycl(L)-1$ with  $L\to \C\P^{\infty}$  being  the  tautological bundle.
We have an identification
$$R(U(1))\cong \Z[u,u^{-1}]\ ,$$
where $u$ is the defining representation $U(1)\stackrel{\id}{\to} U(1)$. With this identification 
we have  $\cycl(L)=B_{EU(1)}(u)$. 
\begin{prob}
Verify explicitly that the map $B_{EU(1)}:u\mapsto \cycl(L)$ extends to a topological isomorphism
$$R(U(1))^{\hat{}}_{I}\stackrel{\sim}{\to} \bKU^{0}(\C\P^{\infty})\ .$$
\end{prob}
\proof
It is natural to consider the restrictions along the subspaces $\C\P^{n}\subset \C\P^{\infty}$.
Note that on $\C\P^{n}$ we have $(\cycl(L)-1)^{n+1}=0$. \hB

Consider now the group $G=\Z/2\Z$. We have $B\Z/2\Z\cong \R\P^{\infty}$.
\begin{prob}\label{jun0742}
Verify the Atiyah-Hirzebruch theorem for the group $\Z/2\Z$ by explicit calculation.
\end{prob}
\proof
Show be calculation that $K(\R\P^{\infty})\cong \Z\oplus \Z_{2}$ and 
$R(S^{1})^{\hat{}}_{I}\cong \Z\oplus \Z_{2}$.
Show then that
$B_{E\Z/2\Z}$ is a topological isomorphism.
\hB

}\end{ex}

We now consider the diagram
$$\xymatrix{R(G)^{+}\ar[d]\ar[r]^{B_{P,\omega}}&\Vect^{geom}(M)\ar[d]^{\hcycl}\\
R(G)\ar[r]^{^{\tilde B_{P,\omega}}}&\hbKU^{0}(M)}\ ,$$
where the factorization $\tilde B_{P,\omega}$ of $B_{P,\omega}$  exists by the universal property of the
 ring completion $R(G)$ of $R(G)^{+}$.
 
 \begin{prob}
Assume that $M$ is compact.
Show that $\tilde B_{P,\omega}$ is continuous with respect to the
$I$-adic topology on $R(G)$ and the discrete topology on $\hbKU^{0}(M)$.
We let
$$\hat B_{P,\omega}:R(G)^{\hat{}}_{I}\to \hbKU^{0}(M)$$ be the extension by
$\tilde B_{P,\omega}$ by continuity.
\end{prob}
\proof
It suffices to show that $\tilde B_{P,\omega}$ annihilates
$I^{2n+2}$, where $n:=\dim(M)$. First of all, $R\circ \tilde B_{P,\omega}$ annihilates $I^{n+1}$.
It follows that we have a factorization
 $$\tilde B_{P,\omega}:I^{n+1}\stackrel{\bar B_{P}}{\to} K\C/\Z^{-1}(M)\to \hbKU^{0}(M)\ .$$
 The first map is a morphism of $R(G)$-modules.
 Hence the restriction of $\tilde B_{P,\omega}$ to $I^{2n+2}$ factors over
 $$I^{n+1}\otimes I^{n+1}\stackrel{\bar B_{P}\otimes B_{P}}{\to} K\C/\Z^{-1}(M)\otimes \bKU^{0}(M) \stackrel{\cup}{\to} K\C/\Z^{-1}(M)\to \hbKU^{0}(M)$$
 which vanishes since the first map vanishes. 
 \hB

We interpret the principal bundle $P$ as a map (determined up to homotopy)
$P:M\to BG$. A connection on $P$ gives rise to a lift 
$$\xymatrix{ &\hbKU^{0}(M)\ar[d]^{I}\\\bK^{0}(BG) \ar[r]^{P^{*}}\ar@{.>}[ur]^{\nabla:= \hat B_{P,\omega}\circ\hat B_{BG}^{-1}}&\bKU^{0}(M)}\ .$$

We complete this diagram as follows:
$$\xymatrix{ &\Omega A^{0}_{cl}(M)\\H(A)^{0}(BG)\ar@{-->}[ur]^{\cC_{\nabla}}&\hbKU^{0}(M)\ar[u]^{R}\ar[d]^{I}\\\bK^{0}(BG)\ar[u]^{c}\ar[r]^{P^{*}}\ar@{.>}[ur]^{\nabla}&\bKU^{0}(M)}\ .$$
The dashed arrow exists since $c$ is an equivalence after complexification and completion, and
$\Omega A^{0}_{cl}(M)$ is a complex vector space. 
Recall that $A=\C[b,b^{-1}]$ and note that $H(A)^{*}(M)$ is bigraded by cohomological and internal degree where $\deg(b)=-2$.
In the discussion above we used the total degree. The domain and target of $\cC_{\nabla}$ are still
graded by the internal degree.
\begin{prob}
Show that
$\cC_{\nabla}$ preserves the internal degree.
\end{prob}
\proof This is a consequence of \ref{may2820}. 
\hB 

The condition that $\cC_{\nabla}$ preserves the internal degree is equivalent to the condition that
$\nabla$ rationally preserves all  Adams operations.

We now turn to the general case. 
Let $B$ be some space not necessarily of the form $BG$ and consider a map $P:M\to B$ from a compact manifold $M$.
We equip $\hbKU^{0}(M)$ and $\Omega A^{0}_{cl}(M)$ with the discrete topology.
Let $\nabla$ be a continuous lift
 \begin{equation}\label{may2831}\xymatrix{&\Omega A^{0}_{cl}(M)\\H(A)^{0}(B)\ar@{-->}[ur]^{\cC_{\nabla}}&\hbKU^{0}(M)\ar[u]^{R} \ar[d]^{I}\\\bK^{0}(B)\ar[u]^{c}\ar[r]^{P^{*}}\ar@{.>}[ur]^{\nabla}&\bKU^{0}(M)}\ .\end{equation}
The map $\cC_{\nabla}$ exists again since $c$ becomes an isomorphism after complexification of the domain and completion.
\begin{ddd}
We call $\cC_{\nabla}$ the cohomological character of $\nabla$.
\end{ddd}
This leads to the following generalization of the notion of a connection on a $G$-principal bundle.

\begin{ddd}
Let $B$ be some space.
A geometrization of a map $P:M\to B$ from a compact manifold $M$ is a continuous lift $\nabla$ in the diagram of topological groups \eqref{may2831}
such that the cohomological character $\cC_{\nabla}$ preserves the internal degree.
 \end{ddd}

\begin{ex}
{\em
Let us  consider geometrizations of the identity maps $S^{2n}\to S^{2n}$, $n\ge 0$.

\begin{enumerate}
\item 
We first consider the case that $B=*$. Then $K^{0}(*)\cong \Z$. We get a geometrization
$\nabla$ be setting $\nabla(1):=1\in \hbKU^{0}(M)$.
\item
The next case is  $n=1$.
We have $\bKU^{0}(S^{2})\cong \Z 1\oplus \Z x$, where $x= \cycl(L)-1$ and
$L\to S^{2}\cong \C\P^{1}$ is the tautological bundle.
We define a geometrization by 
$$\nabla(1):=1 \ ,\quad \nabla(x):=\cycl(L,\nabla^{L})-1\ .$$
Here we can take any connection.
\item
We now consider the geometrization of the identity map $S^{4}\to S^{4}$.
We have $\bKU^{0}(S^{4})\cong \Z 1\oplus \Z x$ with $x:=\cycl(V)-2$, where $V\to S^{4}$ is associated to the
Hopf $SU(2)$-bundle $S^{7}\to S^{4}$ and the defining $2$-dimensional representation of $SU(2)$.
Let us choose some connection $\nabla^{V}$.
We define
$$\nabla(1):=1\ , \quad \nabla(x):=\hcycl(V,\nabla^{V})-2\ .$$
We now calculate the cohomological character. Note that
$\ch(1)=1$ and $\ch(x)=b^{2} \ori_{S^{4}}$. This gives
$$C_{\nabla}(1)=1\ , \quad C_{\nabla}(b^{2}\ori_{S^{4}})=\ch(\nabla^{V})-2=b\ch_{2}(\nabla^{V})+b\ch_{4}(\nabla^{V})\ .$$
This cohomological character preserves the internal degree if any only if $\ch_{2}(\nabla^{V})=0$.
We can  achieve this by choosing a $SU(2)$-connection $\nabla^{V}$.
\item
Let us now consider the identity map $S^{6}\to S^{6}$.
We have $\bKU^{0}(S^{6})\cong \Z 1\oplus \Z x$, with $x=\cycl(V)-d$, where $V\to S^{6}$ is an appropriate bundle.
We again choose a connection $\nabla^{V}$ and  set
$\nabla_{0}(1):=1$ and $\nabla_{0}(x):=\hcycl(V,\nabla^{V})-d$.
We have $\ch(1)=1$ and $\ch([V]-d)=b^{3} \ori_{S^{6}}$.
The cohomological character is now given by 
$$C_{\nabla_{0}}(1)=1\  , \quad 
C_{\nabla_{0}}(b^{3}\ori_{S^{6}})=\ch(\nabla^{V})-d=b\ch_{2}(\nabla^{V})+b^{2}\ch_{4}(\nabla^{V})+b^{3}\ch_{6}(\nabla^{V})\ .$$
It is not clear whether we can find a connection with
$\ch_{2}(\nabla^{V})=0$ and $\ch_{4}(\nabla^{V})=0$. 
But we can choose forms $\alpha_{1}\in \Omega^{1}(S^{6},\C)$  and $\alpha_{3}\in \Omega^{3}(S^{6},\C)$ with $d\alpha_{1}=\ch_{2}(\nabla^{V})$ and $d\alpha_{3}=\ch_{4}(\nabla^{V})$.
Then we define
$$\nabla(1):=1\ , \quad \nabla(x)=\nabla_{0}(x)-a(b\alpha_{1}+b^{2}\alpha_{3})\ .$$
Then 
$$C_{\nabla}(1)=1\ , \quad C_{\nabla}(b^{3}\ori_{S^{6}})= b^{3}\ch_{6}(\nabla^{V})$$
preserves degree.
\end{enumerate}
}
\end{ex}

\begin{ex}{\em
We now consider the geometrization of the inclusion $f:\C\P^{1}\to \C\P^{\infty}$.
We have
$\bKU^{0}(\C\P^{\infty})\cong \Z[[x]]$, where
$x=\cycl(L)-1$ (see \ref{jun0741}). 
We choose any connection $\nabla^{f^{*}L}$ on $f^{*}L$ and define
a geometrization by 
$\nabla(x^{n}):=(\hcycl(f^{*}L,\nabla^{f^{*}L})-1)^{n}$.
\begin{prob}Show that
this prescription defines a geometrization.
\end{prob}
\proof
We have $R(\nabla(x^{n}))=0$ for $n\ge 2$.
It follows that $\nabla(x^{n})=0$ for $n\ge 4$.
This implies that $\nabla$ extends by continuity to a map
$$\nabla:\bKU^{0}(\C\P^{\infty})\to \hbKU^{0}(\C\P^{1})\ .$$
We have $\ch(x)=c_{1}b$ and hence
$\ch(x^{n})= c_{1}^{n}b^{n}$
It follows that
$$C_{\nabla}(x^{n})=b^{n}c_{1}(\nabla^{f^{*}L})^{n}\ .$$
Hence the cohomological character preserves degree.
\hB

\begin{prob}
Construct a geometrization of the inclusion
$f:\R\P^{3}\to \R\P^{\infty}$.
\end{prob} 
\proof
We have (see \ref{jun0742})
$$\bKU^{0}(\R\P^{\infty})\cong \Z 1\oplus \Z_{2} x\ .$$
A topological generator of the $\Z_{2}$-summand is
the class $x:=\cycl(L)-1$, where $L$ is the complexification of the tautological real bundle.
Its restriction to $\R\P^{3}$ has a flat connection $\nabla^{f^{*}L}$.
$$\nabla(1):=1\ , \quad \nabla(x):=\hcycl(f^{*}L,\nabla^{f*L})$$
does the job. Verify!
\hB 
}\end{ex}

\begin{ex}\label{jun0801}{\em 
Let us fix a number $p\in \nat$ and consider the space $X$ defined as the homotopy fibre $$X\to BSU(2)\stackrel{pc_{2}}{\to} K(\Z,4)\ .$$
Let $f:M\to X$ be some map from a compact manifold.
\begin{prob}
Construct a geometrization of $f$.
\end{prob}
\proof
 Since $c_{2}p$ is a rational isomorphism the space $X$ is rationally acyclic.
Note that $X$ has finite skeleta.
We can choose a  smooth compact $\dim(M)+1$-connected approximation
$M^{\prime}\to X$ and a factorization 
$$f:M\stackrel{g}{\to} M^{\prime} \stackrel{f^{\prime}}{\to} X\ .$$
Note that $$0=g^{*}:H(A)^{-1}(M^{\prime})\to H(A)^{-1}(M)\ .$$

We can obtain $X$ as an iterated pull-back
$$\xymatrix{X\ar[d]^{h}\ar[r]&K(\Z/p\Z,3)\ar[d]\ar[r]&{*}\ar[d]\\
BSU(2)\ar[r]^{c_{2}}&K(\Z,4)\ar[r]^{p}&K(\Z,4)}\ .$$
This makes clear that the fibre of $h$ is $K(\Z,3)$.
By \cite{MR0231369} the group $\bKU^{0}(K(\Z,3))$ is divisible. 
It follows, from the Serre spectral sequence that $$h^{*}:R(SU(2))^{\hat{}}_{I}\cong \bKU^{0}(BSU(2))\to \bKU^{0}(X)$$
is injective and has dense range (see \cite[Prop. 5.14]{2011arXiv1103.4217B}  for similar arguments).

We know that $R(SU(2))=\Z[x]$ and $I=(x)$, where $x+2$ is the defining representation $SU(2)\hookrightarrow U(2)$. 
Hence $R(SU(2))^{\hat{}}_{I}\cong \Z[[x]]$.

We choose a bundle $V:=f^{\prime *}h^{*}B_{ESU(2)}(x)$  with connection $\nabla^{V}$.
Note that $\ch(h^{*}x)=2$.
Hence we can find a form $\alpha\in \Omega A^{-1}(M^{\prime})$ such that
$d\alpha=\ch(\nabla^{V})-2$. 
Then we define
$$\nabla^{\prime}(x):=\hcycl(V,\nabla^{V})-2-a(\alpha)\in \hbKU^{0}(M^{\prime})$$
and extend this to a ring map
$$\nabla^{\prime}:R(SU(2))\cong \Z[x]\to \hbKU^{0}(M^{\prime})\ .$$
Finally we define 
$$\nabla(y)=g^{*}\nabla^{\prime}(y)\ , \quad  y\in R(SU(2))\ .$$

We check continuity. 
Let $(y_{i})$ be a sequence in $R(SU(2))$ such that $h^{*}y_{i}\to 0$. We must show that $\nabla(y_{i})\to 0$. 
We observe that $R(\nabla^{\prime}(x))=0$ and
$I(h^{*}y_{i})=0$ for $i>>0$. Hence
$\nabla^{\prime}(y_{i})\in H(A)^{-1}(M^{\prime})/\im(\ch)$ for $i>> 0$.
This implies $\nabla(y_{i})=0$ for $i>> 0$.

Since $R(\nabla(x))=0$ the cohomological character is determined by 
$C_{\nabla}(1)=1$.  It preserves the internal degree. \hB 
}
\end{ex}

\begin{prob}\label{jun0803}
If $B$ is a space of the homotopy type of a $CW$-complex with finite skeleta and $f:M\to B$ is a map from a compact manifold,
then $f$ admits a geometrization. Given $k\ge \dim(M)+1$ we can in addition require that
the geometrization is pulled-back from a $k$-connected approximation
$M^{\prime}\to B$. Such a geometrization will be called good.  \end{prob}
\proof
One proceeds as in Example \ref{jun0801}. We use the same notation.
First decompose $\bKU^{0}(M^{\prime})$ into torsion and free part. Define
$\nabla^{\prime}$ on the torsion part using flat classes. Then deal with the
free part by defining $\nabla^{\prime}$ on generators.

Similar arguments can be found in \cite[Prop. 4.4]{2011arXiv1103.4217B}.
\hB

 \begin{prob}\label{jun0804}
Let $f:M\to B$ as in Prop. \ref{jun0803} and $\nabla_{f}$ be a good geometrization of $f$.
Let further $F:Z\to B$ by an extension of $f$ to a compact zero bordism $Z$ of $M$.
Show that $\nabla_{f}$ extends to a geometrization $\nabla_{F}$ of $F$.

\end{prob}

\begin{prob}
Show by example that in \ref{jun0804} we can not drop the assumtion
that $\nabla_{f}$ is good.
\end{prob}

\begin{ex}{\em 
We continue the example \ref{jun0801}. Let $M$ be a closed $2n-1$-dimensional
$Spin^{c}$-manifold and $f:M\to X$. We choose a good geometrization $\nabla$ of $f$.
 Let $k\ge 1$.
We choose any $\hbKU$-orientation refining the $\bKU$-orientation of $M\to *$  given by the 
$Spin^{c}$-structure and define 
$$e_{k}(M):=\hat I_{!} \nabla(x^{k})\in \hbKU^{-2n+1}(*)\cong \C/\Z$$
using the associated integration map.
\begin{prob}
Show that $e_{k}(M)$ only depends on $Spin^{c}$-bordism class of $M\to X$.
\end{prob}
\proof
Let $F:Z\to X$ be a $Spin^{c}$-zero bordism. Then we can extend the differential $\bKU$-orientation and $\nabla$.
We have by the bordism formula \ref{jun0839}
$$\hat I_{!} \nabla(x^{k})=a(\int_{Z}  \Td(\hat{\tilde \nu})\wedge R(\nabla(x^{k})))=0\ ,$$
since $R(\nabla(x^{k}))=0$.
 \hB
 
 We therefore have defined a family of homomorphisms
 $$e_{k}:MSpin^{c}_{2n-1}(X)\to \C/\Z \ .$$

Let $\sigma:S^{7}\to S^{4}$ be the Hopf bundle and $g:S^{4}\to BSU(2)$ represent the generator of $\pi_{4}(BSU(2))\cong \Z$.
The composition $pc_{2}\circ g\circ \sigma$ is homotopically constant. In fact, we get a unique lift
$$\xymatrix{&X\ar[d]^{h}\\S^{7}\ar@{..>}[ur]^{f}\ar[r]^{g\circ \sigma}&BSU(2)}$$  Note that $S^{7}$ has a canonical $Spin^{c}$-structure induced by the stable framing. We get a class
$[f:S^{7}\to X]\in MSpin^{c}_{7}(X)$.

\begin{prob}
Calculate
$$e_{k}([f:S^{7}\to X])\in \C/\Z\ .$$
\end{prob}

}
\end{ex}

\newpage

\newpage
\section{Exercises}
\newpage
\subsection{Sheet Nr 1.} 
\begin{enumerate}
\item Show in detail, that every vector bundle admits a connection.
\item Show, that a connection on a bundle $E\to M$ uniquely extends to a linear map $\Omega(M,E)\to \Omega(M,E)$
which satisfies the Leibniz rule.
\item Develop the details of the construction of the pull-back of a connection and prove the formula for the pull-back of the curvature.
\item Show that an invariant connection descends if and only if its momentum map vanishes.
Show that in this case the descent is uniquely determined.
\item Calculate the moment map of the Levi-Civita connection on a compact Lie-group with left action and biinvariant metric.
\item Determine explicitly  the connection on the Poincar\'e bundle on the torus $T^{n}$ whose curvature is an invariant two-form on the $2n$-dimensional torus $J(T^{n})\times T^{n}$.
 \end{enumerate}

\newpage

\subsection{Sheet Nr 2.} 
 \begin{enumerate}
\item Calculate the curvature of the $SU(2)$-invariant connection on the tautological bundle
$L\to \C\P^{1}$.
\item   Verify the identities 
$$\ch(\nabla^{E\oplus F})=\ch(\nabla^{E})+\ch(\nabla^{F})\ ,\quad \ch(\nabla^{E\otimes F})=\ch(\nabla^{E})\wedge \ch(\nabla^{F})\ ,$$
$$\ch(\nabla^{\Hom(E,F)})= \ch(\nabla^{F})\wedge \ch(\nabla^{E^{*}})$$
and
 $$\ch_{2i}(\nabla^{E^{*}})=(-1)^{i}\ch_{2i}(\nabla^{E})\ .$$
\item Show that one can express $\ch_{2k}(\nabla)$ as a rational polynomial in the forms
$c_{i}(\nabla)$ for $i\le k$. Determines these polynomials explicitly for $k=4$ and $6$.
\item
For $F\in \Aut(E)$ show
$$\tilde \omega(F)=\int_{S^{1}\times M/M} \omega(E(F))$$
(see \ref{apr12100} for notation)
\item Show that $R^{\nabla^{*}}=-(R^{\nabla})^{*}$ and conclude that
$c_i(\nabla)$ and $\ch_2i(\nabla)$ are real forms if $\nabla$ is unitary.
\item
Show
$$\ch_{2n}(F_\id)=\frac{(-1)^{n-1}(n-1)!}{(2\pi i)^{n}(2n-1)!} \Tr (g^{-1}dg)^{2n-1}\ .$$
(see Lemma \ref{apr2601} for notation)
\end{enumerate}

\newpage

\subsection{Sheet Nr 3.} 

\begin{enumerate}
\item Calculate $\cs_{c_{1}^{n}}(\nabla)$ for all flat connections on line-bundles on the lense space
$L^{2n-1}_{p}:=S^{2n-1}/\mu_{p}$, $p\in \nat$ prime.
\item Consider the map 
 $v:BSL(k,\C^{\delta})\to BGL(k,\C)$ and observe that we then have
$0=v_{*}:MSO_{3}(BSL(k,\C^{\delta}))\to MSO_{3}(BGL(k,\C))$.  
Let $V^{\delta}\to BSL(k,\C^{\delta})$ be associated to the standard representation of $SL(k,\C^{\delta})$. 
 
Show that $\cs^{V^{\delta}}_{c_{1}^{2}}:MSO_{3}(BSL(k,\C^{\delta}))\to \C/\Z$ vanishes.
Furthermore show that
$\Im(\cs^{V^{\delta}}_{c_{2}}):MSO_{3}(BSL(k,\C^{\delta}))\to \R$
is non-trivial.  
Hint: Calculate  
$\tilde c_{2}(f^{*}\nabla)-\tilde c_{2}((f^{*}\nabla)^{*})$ for suitable $f:M\to BSL(k,\C)$ for some metric on $f^{*}V$,
compare this with the Kamber-Tondeur class and use Borel's result that the latter generates $H^{3}(BSL(k,\C);\R)$.
\item
Let $M\to S^{1}$ be a $2$-torus bundle. Then we have a sequence
\begin{equation} 0\to \pi_1(T^{2})\to \pi_1(M)\to \pi_1(S^{1})\to 0\ .\end{equation}
We consider generators $A,B\in \pi_1(T^{2})$ and an element $T\in \pi_1(M)$ which maps to a generator of $\pi_1(S^{1})$. We consider a representation
$\rho:\pi_1(M)\to Sp(1)$ with $\rho(T)=J$ (quaternionic notation) and $\rho(A)=\exp(2\pi i \phi)$, $\rho(B)=\exp(2\pi i \psi)$ for some $\phi,\psi\in \R$. Let $(V\to M,\nabla)$ be the associated two-dimensional flat bundle.
 Calculate 
$\tilde c_2(\nabla)$ in terms of the data $\psi,\phi$ and the extension \eqref{may0301}.
Discuss, under which conditions $\rho$ exists.
Conclude that $\Ree(\cs^{V^{\delta}}_{c_{2}})$ from exercise 2. is non-trivial.
\item
Calculate the Chern-Simons invariant of the lense space $\CS(S^{3}/\mu_{k},g)$, where  
$g$ is a Riemannian metric of constant positive curvature.
\item Determine the sets
$$\{\CS(S^{3},g)\:|\: \mbox{$g$ a Riemannian metric}\}\subset \C/\Z$$
and
$$\{\CS_{refined}(S^{3},g)\:|\: \mbox{$g$ a Riemannian metric}\}\subset \C/\Z$$
(see \eqref{may0603}).
\item
Determine 
$$\CS(M,g)-[\int_{Z} p_{1}(g^{TZ})]$$
in the case that $g^{TZ}$  extends $g$, but does not necessarily have a product structure.
 \end{enumerate}

 \newpage

\subsection{Sheet Nr 4.} 

\begin{enumerate}
\item Let $X$ be a countable $n$-dimensional $CW$-complex. Show that there exists a smooth manifold $M$ and a homotopy equivalence $M\to X$ (see \cite[Sec. 2]{MR2608479} for an argument).
\item Let $f:X\to Y$ be an $n+1$-connected map and $M$ be an $n$-dimensional manifold.
Show that  $f_{*}:[M,X]\to [M,Y]$ is a bijection.
\item Let $M$ be a connected manifold with fundamental group $\pi_{1}(M)=\Z/p\Z$.
For $\rho\in J(M)$ let $P_{\rho}:=P_{|\{\rho\}\times M}$ be the restriction of the Poincar\'e
bundle $P\to J(M)\times M$ (see Ex. \ref{apr1204}).  
Show that
$$J(M)\ni \rho \mapsto c_{1}^{\Z}(P_{\rho})\in H^{2}(M;\Z)$$ 
induces an isomorphism of groups between $J(M)$ and the torsion subgroup of $H^{2}(M;\Z)$.
\item
For $r\in \nat$ we  let $M_{r}\to S^{4}$ be the $SU(2)$-principal bundle classified by 
$c^{\Z}_2(M_{r})=r[S^{4}]\in H^{4}(S^{4};\Z)$. Calculate $H^{4}(M_{r};\Z)$ and its element $p^{\Z}_{1}(TM_{r})$ explicitly.
\item Show that $0\not=c_{1}(L\otimes \C)$ in $H^{2}(\R\P^{\infty};\Z)$, where $L\to \R\P^{\infty}$ is the tautological bundle.  
\item Find an explicit formula for the numbers
$$\langle p_{1}(T\C\P^{n})^{2} \cup p_{3}(T\C\P^{n}),[\C\P^{10}]\rangle\in \Z\ , \quad  n\ge 10\ .$$
\item For a discrete abelian group $A$ we have the following constructions of the cohomology of a manifold $M$ with coefficients in $A$:
$$H^{n}(M;A)=\left\{\begin{array}{cc}H^{n}(M,\underline{A})& \mbox{sheaf cohomology}\\
\mbox{}[M,K(A,n)]&\mbox{homotopy theoretic definition}\\
H^{n}(\Hom(C_{*}(M),A))&\mbox{singular cohomology}
\end{array}\right. $$
Find equivalences between these definitions which are natural in $M$ and $A$.
 \end{enumerate}

\newpage

\subsection{Sheet Nr 5.} 

 \begin{enumerate}
\item For a lower bounded complex of sheaves $\cF$
construct the hypercohomology spectral sequence with
$$H^{p}(M,H^{q}(\cF))=E_2^{p,q}\ ,\quad  (E_r,d_r)\Rightarrow  H^{p+q}(M,\cF)\ .$$
\item  We consider the Hopf fibration $\pi:S^{3}\to \C\P^{1}$ and the complex of sheaves
$$\C\P^{1}\supseteq U\mapsto \Omega_\C(\pi^{-1}(U))\in \Ch^+\ .$$
Make the hypercohomology spectral sequence explicite.
\item Make the hypercohomology spectral sequence for 
$$\cE(n):=\Cone(\Z\to \sigma^{<n}\Omega_\C)$$
explicite and relate it with the structure of Deligne cohomology. 
\item
Show that
$$x\mapsto \{M\ni m\mapsto \ev(m^{*}x)\in \C/\Z\}$$
defines a natural isomorphism
$$f:\hat H^{1}_{Del}(M;\Z)\stackrel{\sim}{\to} C^{\infty}(M;\C/\Z)\ .$$
Show further that $f^{-1}(x)df(x)=R(x)$.
\item
Show that $$x\mapsto \{L(M)\ni \gamma\mapsto \ev(\gamma^{*}x)\in \C/\Z\}$$
defines a natural injective  homomorphism 
$$f:\hat H^{2}_{Del}(M;\Z)\stackrel{\sim}{\to} C^{\infty}(L(M);\C/\Z)\ .$$
Show that it is not surjective, in general, and calculate $df$.
Here $L(M):=C^{\infty}(S^{1},M)$ is the free loop space and smooth maps
out of  $L(M)$ are understood in the diffeological sense.
\item
Show that the canonical generator $e\in \hat H_{Del}^{1}(S^{1};\Z)$ is primitive.
Let further $G$ be a compact connected and simply connected simple Lie group. Show that there is a unique (up to sign) bi-invariant class $e\in \hat H^{3}_{Del}(G;\Z)$ such that
$I(e)\in H^{3}(G;\Z)$ is a generator. Show that this class is primitive.
\end{enumerate}

\newpage
\subsection{Sheet Nr 6.}

\begin{enumerate}
\item Consider the  complex tangent bundle $T\bH\C^{n}$ of the complex hyperbolic space $\bH\C^{n}$ with the Levi-Civita connection $\nabla$ and a point $x\in \bH\C^{n}$. For $r>0$  let $S_r\subset \bH\C^{n}$ be the distance sphere  at $x$.
Calculate
$$\ev_{S_r}(\hat c_n(\nabla))\in \C/\Z$$ as a function of $r$.
\item Let $E\to M$ be a vector bundle, $\nabla$ a connection on $E$ and $h$ be a hermitean metric. 
 Show that $$\overline{\hat  c_{n}(\nabla)}=\hat c_n(\nabla)\ .$$
\item Let $(E,\nabla)$ be a flat complex vector bundle. 
Then $\Im \hat s_{2n}(\nabla)\in H^{2n-1}(M;\R)$.
Find the relation
between $\Im \hat s_n(\nabla)$ and the Kamber-Tondeur class
$\check \ch_{2n}(\nabla)$.
\item
Show that
a \v{C}ech $1$-cocycle  for the complex $$\cK(1):\underline{\C}^{*}\stackrel{d\log}{\to} \Omega^{1}_\C$$
can naturally be identified with the glueing data for a 
line bundle with connection. Use this to construct the isomorphism 
$$\hat c_{1}:\Line_{\nabla}(M)\stackrel{\sim}{\to} H^{1}(M;\cK(1))$$
explicitly  on the level of \v{C}ech cohomology. Verify that this is compatible with the previous constructioin of $\hat c_1$ and the identification $H^{1}(M;\cK(1))\cong \hat H_{Del}^{2}(M;\Z)$ described in the course 
\item
Let $V\to B$ a vector bundle with structure group $GL(k,\C^{\delta})$ over some space $B$ and $\omega$ be an integral characteristic form. For a map $f:M\to B$ from a closed oriented $n-1$-manifold
define
$$\cs_{\omega}^{V}(f:M\to B):=\ev(\hat \omega (\nabla^{f^{*}V}))\in \C/\Z\ .$$  
Show that this defines a homomorphism $\cs_\omega^{V}:\MSO_{n-1}(B)\to \C/\Z$ which extends the provious definition made under the assumption that $V$ was trivializable as a bundle with structure group $GL(k,\C)$. Show further, that $\cs^{V}_\omega$ factorizes over the orientation transformation $\MSO_{n-1}(B)\to H_{n-1}(M;\Z)$, $[f:M\to B]\mapsto f_*[M]$.
\item
Fix a $p$th root of unity $\xi$ and consider the flat bundle
 $V\to B\Z/p\Z$ with holonomy $[k]\mapsto  \xi^{k}$.
 Let $L^{2n-1}_p:=S^{2n-1}/(\Z/p\Z)$ and $f:L^{2n-1}_p
\to B\Z/p\Z$ be the map which induces the canonical identification of fundamental groups. Calculate $\cs^{V}_{c_1^{n}}(f)\in \C/\Z$.

\end{enumerate}

\newpage
\subsection{Sheet Nr 7.} 

\begin{enumerate}
\item For a closed oriented Riemannian manifold $(M,g)$ of dimension $n-1$ define, using the complxification of the Levi-Cevita connection,
$$\CS(g):=\ev((-1)^{n}\hat c_{2n}(\nabla^{TM\otimes \C}))\in \C/\Z\ .$$
Show that this extends the previous definition of the Chern-Simons invariant  dropping the assumption that $M$ is zero-bordant.
\item
Show that $\CS(g)$ is an invariant of the conformal class of $g$ which vanishes if $M$
totally umbilically bounds a locally conformally flat manifold.  
\item Let $\Z/p\Z$ act on $S^{4n-1}\subset \C^{2n}$ diagonally by
$[1]\mapsto (\zeta^{q_{1}},\dots,\zeta^{q_{2n}})$, where $\zeta=\exp(2\pi i p^{-1})$
is a primitive root of unity and the numbers $g.g.T(q_{1},\dots,q_{2n},p)=1$ and set
$L^{4n-1}_{p}(q_{1},\dots,q_{2n}):=S^{4n-1}/(\Z/p\Z)$ with the induced round metric $g$.
Calculate
$$\CS(L^{4n-1}_{p}(q_{1},\dots,q_{2n}),g)\in \C/\Z\ .$$
\item 
Calculate
$ \hat c_1(\nabla \otimes \C) \in \hat H_{Del}^{2}(\R\P^{n};\Z) $ for the canonical bundle $(V,\nabla)$ of $\R\P^{n}$ in algebro-topological terms.  
\item Show that an oriented  real vector bundle $ V $   admits a spin structure if (and only if ?) $\hat c_1(\nabla \otimes \C )=0$ for some and hence every connection.
\item
Show that the integral characteristic form $\frac{p_1}{2}$ for real spin bundles has a unique
differential refinement. Show that this refinement is additive.\end{enumerate}

\newpage

\subsection{Sheet Nr 8.} 

\begin{enumerate}
\item Let $G$ be a Lie group with finitely many connected components and $P\to M$ be a $G$ principal bundle with connection $\omega$.  
Let $\tilde \phi\in \tilde I(G)$. Show that there exists a  natural (under pull-back of principal bundles  with connection) form $\theta(\omega)\in \Omega^{2n-1}(P;\C)$ such that
$d\theta(\omega)=\phi(R^{\omega})$.
Show that $$\pi^{*}\hat \cw(\tilde \phi)(\omega)=a(\theta(\omega))\ .$$
\item We consider $G:=SU(2)$ and $\tilde \phi \in I(SU(2))$ such that
$\phi_\Z=c_2$. 
Give an explicit  formula for
$\theta(\omega)$ in Exercise 1.
\item Let $(V,h ,\nabla)$ be a flat euclidean vector bundle of
dimension $2n$ and $\pi:S(V)\to M$ be its sphere bundle.
Let $\beta$ be the normalized  fibrewise volume form on $S(V)$, extended to $S(V)$ using the connection.  We define
$f(\nabla)\in H^{2n-1}(M;\C/\Z)$ by the condition
$$f(\nabla)(x)=[\int_{\tilde x} \beta]\in \C/\Z$$
for   $\tilde x\in H_{2n-1}(S(V);\Z)$ such that $\pi_*(\tilde x)=x\in H_{2n-1}(M; \Z)$. Show that $f(\nabla)$ is well-defined and that
$f(\nabla)=\hat \chi(\nabla)$.
\item Complete the arguments in the construction of the $\cup$-product.
\item On $S^{1}\times S^{1}$ show that
$\pr_1^{*}\hat e\cup \pr_2
^{*}\hat e\cup \mu^{*}\hat e=0$, where $\mu$ is the group structure.
\item Consider the smooth $\C/\Z$-valued functions on $\C^{*}\setminus \{1\}$  given by 
$f:=[\frac{1}{2\pi i} \ln(z)]$ and $g:=[\frac{1}{2\pi i} \ln(1-z)]$
as elements in $\hat H^{1}_{Del}(\C^{*}\setminus\{1\};\Z)$.
Show that $f\cup g=0$.

\end{enumerate}

\newpage

\subsection{Sheet Nr  9.}

\begin{enumerate}
\item
We consider a closed oriented $n-1$-dimensional manifold $M$ and
  define the
$$\langle\dots,\dots\rangle\colon \hat H^{p}_{Del}(M;\Z)\otimes \hat H^{n-p}_{Del}(M;\Z)\to \C/\Z$$
by $x\otimes y\mapsto \ev(x\cup y)$.
Show that this pairing is non-degenerate.
\item Complete the details of the construction of the differential cohomology diagram for
$\hat H_{CS}^{*}(M;\Z)$.
\item Show that there is a unique natural integration for $H(\dots;\Z)$ along proper oriented submersions which is compatible with Mayer-Vietoris sequences associated to decompositions of the base.
\item Let $(P,\nabla^{P})$ be the Poincar\'e bundle over $J(S^{1})\times S^{1}$.
 Calculate the integrals of $\hat c_{1}(\nabla^{P})$ over the two projections.
\item Let $p:W\to B$ be an oriented bundle of compact manifolds with boundary $\partial W\to B$. Let $q:\partial W\to B$ denote the restriction of $p$ to the boundary.
We consider   $x\in \hat H_{Del}^{n}(W;\Z)$.
  Show that $$q_{!}(x)=a(\int_{W/B} R(x))\ .$$
\item Let $\underline{\Z}\in \Sm(\Nerve(\Ch))$ and consider its interpretation 
$  \underline{\Z}_\infty  \in\Sm(\Nerve(\Ch)[W^{-1}])$.
Show that $\underline{\Z}$ satisfies descent, but  $  \underline{\Z}_\infty $ does not.

\end{enumerate}

\newpage
\subsection{Sheet Nr  10.}

\begin{enumerate}
\item Show that the tensor structure with $\Nerve(\sSet)[W^{-1}])$ on $\Nerve(\Ch)[W^{-1}]$ is given by $A\otimes X\cong A\otimes C_*(X)$, where $C_*(X)$ is the chain complex associated to $C$.
\item We define a functor
$\bar \bs:\Sm(\bC)\to \Sm(\bC)$ as the composition
of
$\bs :\Sm(\bC)\to \Sm(\Fun(\Nerve(\Delta^{op}),\bC))$
and
$\colim_{\Nerve(\Delta^{op})}: \Sm(\Fun(\Nerve(\Delta^{op}),\bC))\to \Sm(\bC)$,
where $\bs$ is  precomposition by
$\Mf\to \Fun(\Delta,\Mf), \, \quad M\mapsto M\times \Delta^{\bullet}$.

We consider the sheaf $\Omega_{cl}^{n}\in \Sh_{\Ab}(\Mf)$ of closed $n$-forms. 
 Calculate the homotopy and cohomology groups  groups of the evaluations on $M$ of the following presheaves and sheaves derived from $\Omega_{cl}^{n}$ in terms of the de Rham cohomology on $M$.
\begin{enumerate}
\item We let $C(\Omega_{cl}^{n})\in \Sm(\Nerve(\Ch)[W^{-1}])$ be the  presheaf of chain complexes given by  $\Omega_{cl}^{n}$ located in degree $0$.
Calculate $H^{*}(C(\Omega_{cl}^{n}))(M))$.
\item Calculate $H^{*}(L(C(\Omega_{cl}^{n}))(M))$, where $L:\Sm(\Nerve(\Ch)[W^{-1}])\to \Sm^{desc}(\Nerve(\Ch)[W^{-1}])$
  is the sheafification.
\item  
Similarly let $C^{\le 0}(\Omega_{cl}^{n})\in \Sm(\Ch^{\le 0}[W^{-1}])$ be the negatively graded chain complex represented by $\Omega^{n}_{cl}$ and $L:\Sm(\Nerve(\Ch^{\le 0})[W^{-1}])\to \Sm^{desc}(\Nerve(\Ch^{\le 0})[W^{-1}])$. 
Calculate $H^{*}(L(C^{\le 0}(\Omega_{cl}^{n}))(M))$ 
\item Calculate the groups
$H^{*}(\bar \bs C(\Omega^{n}_{cl})(M))$, $H^{*}(\bar \bs C^{\le 0}(\Omega^{n}_{cl})(M))$, 
$H^{*}( L(\bar \bs C(\Omega^{n}_{cl}))(M))$ and  $H^{*}(L (\bar \bs C^{\le 0}(\Omega^{n}_{cl}))(M))$, 
 where $\bs$ is as in \ref{jun1901}.
 \item We consider $\Omega_{cl}^{n}$ as a smoth constant simplicial abelian group
$S(\Omega_{cl}^{n})\in \Sm(\Nerve(\sAb)[W^{-1}])$.
Calculate $\pi_{*}(S(\Omega_{cl}^{n})(M))$.
\item  
Calculate $\pi_{*}(\bar \bs S(\Omega_{cl}^{n})(M))$.
 \end{enumerate}
\item Let $(H\Z,\C,c)$ be the canonical datum of $H\Z$ and $\cD(n)$ the $n$th Deligne complex.
 Show that $\Diff^{n}(H\Z,\C,c)\cong H(\cD(n))$
\item Let $(E,A,c)$ be a datum. Provide a natural morphism of monoids
$$\Omega \map(E,H(A))\to \ed(E,A,c)$$ and calculate the action of
$\pi_{0}(\Omega \map(E,H(A)))\cong H(A)^{-1}(E)$ on $\hat E^{*}$.
 \item Let $(E,A,c)$ be a datum. Show that there are pull-back diagrams
 $$ \xymatrix{\Diff^{n}(E,A,c)\ar[r]\ar[d]&L(H(\Omega A_{cl}^{n}))\ar[d]\\
 \Funk(E)\ar[r]&\Funk(H(A))}\ ,\quad \xymatrix{\Diff^{n}(E,A,c)\ar[r]\ar[d]&L(H(\Omega A_{cl}^{n-1}/\im(d)))\ar[d]\\
 \Funk(E)\ar[r]&\Funk(H(A))}\ .$$

\end{enumerate}

\newpage
\subsection{Sheet Nr  11.}

\begin{enumerate}
\item We consider a complex vector bundle $V\to M$. Show, in the framewok of   $\bKU$-theory, that $\Rham(\ch(\nabla^{V}))=c(\cycl(V))$.
\item Show in detail the existence and uniqueness of the lifts
$\hat c_{n}:\hbKU^{0}\to \widehat{H\Z}^{2n}$ of the Chern classes.
\item Show that there is a unique $SU(2)$-invariant 
class $u\in \hbKU^{1}(SU(2))$ which refines a generator of $\widetilde{\bKU}^{1}(SU(2))\cong \Z$ and restricts to zero at the identity element.
\item We consider the canonical datum $(\bMU,A,c)$. Let $u\in H(A)^{0}(BU)$ correspond to $c\in H(A)^{0}(\bMU)$ under the Thom isomorphism.
Show that for a $\bMU$-cycle $z:=(f:W\to M,\cN,I)$ of degree $n$ we have
$f_!u(N)=c(\cycl(z))\in H(A)^{n}(M)$.
\item Write (see previous exercise) $u=\sum_{n=0}^{\infty} [\C\P^{n}] \Phi_n$ with $\Phi_n\in \C[c_1, \dots,c_n]^{2n}$. Calculate $\Phi_n$ for small $n$ (e.g. $n\in \{0,1,2,3\}$) explicitly.
\item We consider $S^{2}$ with the standard metric. A map $f:S^{1}\to S^{2}$ gives rise to a 
$\bMU$-cycle $Z(f)$ of degree $1$ in a canonical way. Let $x\in S^{2}$ be some point. Calculate
$\harm(Z(f))_{|x}\in \widehat{\bMU}^{1}(*)\cong \C/\Z$.

\end{enumerate}

\newpage

\subsection{Sheet Nr  12.}

\begin{enumerate}
\item Show that the cycle map $\hcycl:i\Vect^{geom} \to \hbKU$
is multiplicative.
\item Show that $\bMU$ is formal over $\C$.
\item Show that $\Diff(\bS,\C,c)$ fits into
a pull-back
$$\xymatrix{\Diff(\bS,\C,c)\ar[d]\ar[r]&\Diff(H\Z,\C,c)\ar[d]\\
\Funk(\bS)\ar[r]&\Funk(H\Z)}\ .$$
\item Let $(E,A,c)$ be a multiplicative data.
Work out the definition of a module data $(F,B,d)$.
Construct
$\Diff(F,B,d)$ as a $\Diff(E,A,c)$-module.
\item Let $\eta\in \pi_{1}(\bS)\cong \Z/2\Z$ and $\sigma\in \pi_{3}(\bS)\cong \Z/24\Z$ be the generators given by the Hopf maps $S^{3}\to S^{2}$ and $S^{7}\to S^{4}$. The have canonical lifts
$\hat \eta\in \hat \bS^{-1}(*)$ and $\hat \sigma\in \hat \bS^{-3}(*)$. 
Let $\hat e:\hat \bS^{*}\to \widehat{\bKU}^{*}$ be induced by the unit. 
Calculate $\hat e(\hat\eta)\in  \widehat{\bKU}^{-1}\cong \C/\Z$ and
$\hat e(\hat\sigma)\in  \widehat{\bKU}^{-3}\cong \C/\Z$.
\item Show that $\mathbf{KO}$ is formal over $\C$.
\end{enumerate}

\newpage

\subsection{Sheet Nr  13.}
\begin{enumerate}
\item Let $W$ be a compact  oriented $n$-manifold with boundary $\partial W$. Calculate
$\widehat{H\Z}^{n}_{\partial W}(W)$, $\widehat{H\Z}_{\partial W}^{n+1}(W)$ and $\widehat{H\Z}_{c}^{n+1}(W\setminus \partial W)$.
\item Let $\mathbf{MSpin^{c}}\to \bKU$ be the ABS-orientation. The universal bundle
$V\to BSpin^{c}(n)$ is thus $\bKU$-oriented with Thom class $\nu$.
Let $0_{V}:BSpin^{c}(n)\to V$ be the zero section. 
Describe $H\C[b,b^{-1}]^{*}(BSpin^{c}(n))$ and calculate $\ch(0_{V}^{*}\nu)$.
\item Calculate $\ch(\rho^{k}(V))$, where $\rho^{k}(V)$ is the cannibalistic class
for the Adams operation $\Psi^{k}$ and $V$ is as in $1.$.
\item Show that there exists a unique differential refinement 
$\hat \rho^{k}:\Vect_{\R}^{Spin^{c}, geom}\to \hbKU[\frac{1}{k}]^{0}$.
\item Calculate $\hat \rho^{k}(L,\nabla)$, where $(L,\nabla)$ is the tautological bundle on $\C\P^{n}$ with its invariant connection.
\item 
Let $\hat \nu\in \hat \bS^{n}_{c}(\R^{n})$ be a differential Thom class. Let $x\in \hat E^{*}(M)$.
 Show that $$\desusp(\pr_{\R^{n}}^{*}\hat \nu\cup \pr_{M}^{*}x)=x\ .$$
\end{enumerate}
\newpage

\subsection{Sheet Nr. 14.}\label{sheet14}

These problems are probably more than just exercises. 
I am interested to see solutions.

\begin{enumerate}
 \item Descide \ref{aug0760}, \ref{aug0754}, and \ref{aug0761}.
\item Descide \ref{aug0753}. 
\item Prove the bordism formula \ref{jun0839}.
\item Solve \ref{aug0870}.
\item Prove the differential $\bMU$-index theorem \ref{aug0751}.
\end{enumerate}

\newpage

\bibliographystyle{alpha}
\bibliography{vorl}
\end{document}